\documentclass[english]{elsarticle}
\usepackage[T1]{fontenc}
\usepackage[latin9]{inputenc}
\usepackage{geometry}
\usepackage{color}
\usepackage{array}
\usepackage{textcomp}
\usepackage{multirow}
\usepackage{amsmath}
\usepackage{amssymb}
\usepackage{graphicx}
\usepackage{esint}
\PassOptionsToPackage{normalem}{ulem}
\usepackage{ulem}

\makeatletter

\providecommand{\tabularnewline}{\\}

\journal{Computers and Fluids}
\renewcommand\citet[1]{\cite{#1}}
\usepackage[mathlines]{lineno}
\newcommand*\patchAmsMathEnvironmentForLineno[1]{%
  \expandafter\let\csname old#1\expandafter\endcsname\csname #1\endcsname
  \expandafter\let\csname oldend#1\expandafter\endcsname\csname end#1\endcsname
  \renewenvironment{#1}%
    {\linenomath\csname old#1\endcsname}%
    {\csname oldend#1\endcsname\endlinenomath}}%
\newcommand*\patchBothAmsMathEnvironmentForLineno[1]{%
  \patchAmsMathEnvironmentForLineno{#1}%
  \patchAmsMathEnvironmentForLineno{#1*}}%
\AtBeginDocument{%
\patchBothAmsMathEnvironmentForLineno{equation}%
\patchBothAmsMathEnvironmentForLineno{align}%
\patchBothAmsMathEnvironmentForLineno{flalign}%
\patchBothAmsMathEnvironmentForLineno{alignat}%
\patchBothAmsMathEnvironmentForLineno{gather}%
}

\@ifundefined{showcaptionsetup}{}{%
 \PassOptionsToPackage{caption=false}{subfig}}
\usepackage{subfig}
\makeatother

\usepackage{babel}
\begin{document}

\title{A priori subcell limiting based on compact nonuniform nonlinear weighted
schemes of high-order CPR method for hyperbolic conservation laws}

\author[adress2,adress3]{Huajun Zhu}

\author[adress1]{Huayong Liu}

\author[adress2]{Zhen-Guo Yan\corref{author1}}

\ead{yanzhg@mail.ustc.edu.cn, zgyan@skla.cardc.cn}

\author[adress2]{Guoquan Shi}

\author[adress3,adress4]{\\Xiaogang Deng}

\cortext[author1]{Corresponding author.}

\address[adress2]{State Key Laboratory of Aerodynamics, Mianyang, Sichuan 621000, PR China}

\address[adress1]{Tianfu Engineering-Oriented Numerical Simulation \& Software Innovation
Center, Sichuan University, Chengdu, Sichuan 610000, PR China}

\address[adress3]{College of Aerospace Science and Engineering, National University
of Defense Technology, Changsha, Hunan 410073, PR China}

\address[adress4]{Chinese Academy of Military Science, Beijing 100071, PR China}
\begin{abstract}
This paper develops a shock capturing approach for high-order correction
procedure via reconstruction (CPR) method with Legendre-Gauss solution
points. Shock regions are treated by novel compact nonuniform nonlinear
weighted (CNNW) schemes, which have the same solution points as the
CPR method. CNNW schemes are constructed by discretizing flux derivatives
based on Riemann fluxes at flux points in one cell and using nonuniform
nonlinear weighted (NNW) interpolations to obtain the left and right
values at flux points. Then, a priori subcell p-adaptive CNNW limiting
of the CPR method is proposed for hyperbolic conservation laws. Firstly,
a troubled cell indicator is used to detect shock regions and to quantify
solution smoothness. Secondly, according to the magnitude of the indicator,
CNNW schemes with varying accuracy orders are chosen adaptively for
the troubled cells. The spectral property and discrete conservation
laws are mathematically analyzed. Various numerical experiments show
that the CPR method with subcell CNNW limiting has superiority in
satisfying discrete conservation laws and in good balance between
resolution and shock capturing robustness.\end{abstract}
\begin{keyword}
correction procedure via reconstruction (CPR)\sep shock capturing\sep
compact nonlinear nonuniform weighted (CNNW) schemes\sep subcell
limiting\sep discrete conservation law
\end{keyword}
\maketitle

\section{Introduction}

High-order methods have been widely used in large eddy simulations
(LES) and direct numerical simulations (DNS) of turbulent flows, computational
aeroacoustics (CAA) and shock-induced separation flows \citet{Deng2012,Wang2013,Huynh2014,Wang2017}.
Among high-order methods, high-order finite element (FE) methods are
compact, highly parallelizable, efficient for high-performance computing
and applicable to complex unstructured meshes, such as discontinuous
Galerkin (DG) method \citet{Cockburn1989,Cockburn1998} and correction
procedure via reconstruction method (CPR) \citet{Huynh2007,Wang2009,Huynh2014,Wang2017}.
For conservation laws, since solution may contain discontinuities
even if the initial conditions are smooth, numerical methods need
to be designed carefully to capture discontinuities effectively without
generating obvious oscillations. It is well known that high-order
FE schemes can produce spurious oscillations called the Gibb\textquoteright{}s
instability near discontinuities and may lead to crash of the code
\citet{Wang2013,Huynh2014,Zhong2013,Vilar2019}.

There exist different strategies to deal with spurious oscillations
of FE method. The first strategy is to add artificial viscosity to
the original equations to change properties of PDE and smear out oscillations
near discontinuities \citet{Persson2006,Discacciati2020,Yu2020,Feng2020}.
The second strategy is to limit high-order polynomial approximated
solution distribution in a cell near discontinuity while keep the
rest procedure of the FE method unchanged. Some limiters belong to
this strategy, such as Hermite WENO limiter \citet{Qiu2004,Balsara2007,Zhu2009},
a simple WENO limiter \citet{Zhong2013,Zhu2013,Du2015}, p-weighted
limiter \citet{Li2020} and MLP limiter \citet{Park2016}. The third
strategy is to develop a hybrid method based on different accuracy
orders or different kinds of schemes. The hp-adaption method is a
hybrid method based on different accuracy orders. The method reduces
the degree of the polynomials in shock regions and refine the grid
to guarantee the resolution \citet{Baumann1999,Burbeau2001}. Recently,
DG method based on a subcell limiting is developed \citet{Dumbser2014,Dumbser2016},
which is a hybrid method based on different kinds of schemes. This
hybrid method subdivides the DG cell in shock region into subcells
and adopts shock capturing schemes on the subcells. The third strategy
shows good robustness in shock-capturing since they avoid using high-order
finite element methods to capture shocks directly, but utilize schemes
with better shock capturing abilities instead. According to different
subcell splitting, there are two kinds of subcell limiting approaches.

The first kind of subcell limiting approach is based on equally distributed
subcells. In 2014, Dumbser et al. proposed a posteriori subcell limitig
for DG method for the simple Cartesian case, which refines the troubled
cells into equally spaced subcells and use a high-order ADER-WENO
finite volume (FV) scheme to recompute the discrete solution \citet{Dumbser2014}.
The method has the ability to resolve discontinuities at a sub-grid
scale and has been extended to general unstructured triangular and
tetrahedral meshes \citet{Dumbser2016}, to moving unstructured meshes
\citet{Boscheri2017} and to shallow water equations \citet{Ioriatti2019}.
A posteriori correction of DG schemes at the subcell scale was introduced
by Vilar \citet{Vilar2019}. Although this subcell limiting approach
is effective for shock capturing, it is a bit complicate for code
design since it is a posteriori approach. In addition, since subcells
are equally spaced and the solution points of the two methods are
not coincide with each other, the approach needs data transformation
or projection between DG cells and FV subcells, which adds extra computational
costs. 

The second kind of subcell limiting approach is based on nonuniformly
distributed subcells. In 2017, Sonntag and Munz took a proiri strategy
and used an inherent refinement of the DG elements into several nonuniformly
distributed FV subcells with a lower order approximation without changing
the degrees of freedom (DoFs) \citet{Sonntag2017}. First-order FV
subcell schemes and second-order FV TVD subcell schemes were considered.
Each subcell is associated with one degree of freedom within the DG
grid cell. In 2021, Krais et al. combined DG spectral element method
with a subcell FV schemes (of first-order or second-order ) to capture
shocks in their FLEXI framework \citet{Krais2021}. In 2021, Hennemann
et al. extended the subcell idea and proposed a subcell low order
FV type discretization based on the nodal Legendre-Gauss-Lobatto (LGL)
values for high-order entropy stable discontinuous Galerkin spectral
element method (DGSEM) \citet{Hennemann2021}. This approach uses
nonuniformly spaced solution points of DG schemes and shows superiority
in data exchange and discrete conservation law. However, high-order
subcell schemes based on nonuniform solution points have not been
developed. 

In this paper, we will focus on the second kind of subcell limiting
approach and introduce high-order finite difference (FD) schemes on
nonuniform solution points into subcell limiting. There are some high-order
FD schemes which have good properties in capturing discontinuities,
for instance, weighted essentially non-oscillatory schemes (WENO)
\citet{Jiang1996,Shu1999,Shu2003,Shu2009} and weighted compact nonlinear
schemes (WCNS) \citet{Deng2000}. However these schemes have some
difficulties in being applied to subcell limiting for FE schemes.
Firstly, it needs data transformations between FE solution points
and FD solution points, since high-order FD schemes usually are constructed
based on uniformly spaced solution points while FE schemes are usually
on nonuniformly spaced solution points. Secondly, it has difficulty
in satisfying discrete conservation laws since the discrete conservation
laws for the FD and FE schemes are different \citet{Cheng2014,zhu2017,guo2020}.
To our knowledge, there are no high-order FD shock capturing schemes
designed on the nonuniformly spaced solution points of FE method,
such as Legendre-Gauss (LG) points or LGL points. 

To address the above issues, compact nonuniform nonlinear weighted
schemes (abbr. CNNW) are constructed based on LG solution points and
a priori subcell p-adaptive CNNW limiting of CPR method is proposed
for hyperbolic conservation laws in this paper. The main contributions
of the work go as follows:

1. Compact nonuniform nonlinear weighted (CNNW) schemes are constructed
based on Gauss-Legendre solution points. Nonuniform nonlinear weighted
interpolations are taken to introduce nonlinear mechanism and flux
derivatives are discretized based on Riemann fluxes at flux points
in one cell. Both high-order and low-order CNNW schemes are proved
to be satisfying the same discrete conservation laws as the CPR method.
The spectral properties of CNNW are analyzed and compared with high-order
WCNS and high-order CPR. 

2. A priori subcell p-adaptive CNNW limiting of high-order CPR method
is proposed for hyperbolic conservation laws. Firstly, an indicator
considering modal decay of the polynomial representation based on
an extended stencil is used to detect troubled cells. Secondly, the
troubles cells are divided into nonuniformly spaced subcells being
solved by CNNW schemes. To ensure the hybrid scheme being robust and
accurate, CNNW schemes with varying accuracy orders ($p$-adaptive
CNNW) are chosen adaptively according to the magnitude of troubled
cell indicator to accomplish transition from smooth region to discontinuous
region. In troubled cells, $p$-adaptive CNNW is applied by locally
increasing and decreasing accuracy orders of interpolation operators
or difference operators of CNNW. 

3. Various numerical experiments for linear wave equations and Euler
equations are conducted to show the good properties of the proposed
CNNW scheme and the CPR scheme with subcell CNNW limiting in high
resolution, good robustness in shock capturing and satisfying discrete
conservation law. In addition, the CPR with subcell p-adaptive CNNW
limiting has higher resolution than that with subcell second-order
CNNW limiting. Results of the proposed schemes are also compared with
DG schemes with other limiters.

This paper is organized as follows. In Section 2, high-order CPR methods
are recalled. In Section 3, novel compact nonuniform nonlinear weighted
schemes based on nonuniformly spaced solution points are developed.
In Section 4, a priori subcell CNNW limiting of CPR method is proposed.
In Section 5, spectral properties and discrete conservation laws are
analyzed. In Section 6, numerical investigation about the proposed
CNNW scheme and the CPR scheme with subcell CNNW limiting is conducted
to illustrate the effectiveness of the schemes. Finally, concluding
remarks are given in the Section 7.

\section{Review of high-order CPR }

Correction procedure via reconstruction (CPR) method was originally
proposed by Huynh as flux reconstruction (FR) for structured grids
\citet{Huynh2007} and then was generalized to unstructured grids
by Wang et al.\citet{Wang2009}. Here we give a brief review of the
CPR method. For more details we refer to papers\citet{Huynh2007,Huynh2009,Wang2009}.

Consider two-dimensional conservation law in physical space

\begin{eqnarray}
\frac{\partial\mathbf{U}}{\partial t}+\frac{\partial\mathbf{F}(\mathbf{U})}{\partial x}+\frac{\partial\mathbf{G}(\mathbf{U})}{\partial y} & = & \mathbf{0},\label{eq:conservation-law-physical}
\end{eqnarray}
where $\mathbf{U}$ is the conservative variable vector, and $\mathbf{F}$
is the inviscid flux vector. After transformation into the computational
space, conservation law (\ref{eq:conservation-law-physical}) becomes 

\begin{equation}
\frac{\partial\mathbf{\widehat{\mathbf{U}}}}{\partial t}+\frac{\partial\widehat{\mathbf{F}}}{\partial\xi}+\frac{\partial\widehat{\mathbf{G}}}{\partial\eta}=\mathbf{0},\label{eq:2D-CL}
\end{equation}
where $\mathbf{\widehat{\mathbf{U}}}=J\mathbf{U}$, $\widehat{\mathbf{F}}=\mathbf{F}\widehat{\xi}_{x}+\mathbf{G}\widehat{\xi}_{y}$,
$\widehat{\mathbf{G}}=\mathbf{F}\widehat{\eta}_{x}+\mathbf{G}\widehat{\eta}_{y}$.
Here grid metrics are 
\begin{equation}
\left\{ \begin{alignedat}{3}\widehat{\xi}_{x} & = & J\xi_{x} & = & y_{\eta},\\
\widehat{\xi}_{y} & = & J\xi_{y} & = & -x_{\eta},
\end{alignedat}
\right.\quad\left\{ \begin{alignedat}{3}\widehat{\eta}_{x} & = & J\eta_{x} & = & -y_{\xi},\\
\widehat{\eta}_{y} & = & J\eta_{y} & = & x_{\xi},
\end{alignedat}
\right.\label{eq:grid derivatives}
\end{equation}
and Jacobian is 
\begin{eqnarray}
J & = & \left|\frac{\partial(x,y)}{\partial(\xi,\eta)}\right|=x_{\xi}y_{\eta}-x_{\eta}y_{\xi}.\label{eq:Jacobian-2D-1}
\end{eqnarray}

In the CPR method, the solution inside one element is approximated
by polynomials, for example the following degree $K$ Lagrange interpolation
polynomial

\begin{eqnarray}
\mathbf{U}_{i,j}^{h}(\xi,\eta) & = & \sum_{l=1}^{K+1}\sum_{m=1}^{K+1}\mathbf{U}_{i,j,l,m}L_{l}(\xi)L_{m}(\eta),
\end{eqnarray}
where $\mathbf{U}_{i,j,l,m}$ are the state variables at the solution
point $(l,m)$ of the $(i,j)$ cell , $L_{l}(\xi)$ and $L_{m}(\eta)$
are the 1D Lagrange polynomials in the $\xi$ and $\eta$ directions.
Then, Lagrange Polynomial (LP) approach is applied to approximate
the second term and the third term in (\ref{eq:2D-CL}),

\begin{equation}
\widehat{\mathbf{F}}_{i,j}(\xi,\eta)=\sum_{l=1}^{K+2}\sum_{m=1}^{K+1}\widehat{\mathbf{F}}_{i,j,l,m}L_{l}(\xi)L_{m}(\eta),\quad\widehat{\mathbf{G}}_{i,j}(\xi,\eta)=\sum_{l=1}^{K+1}\sum_{m=1}^{K+2}\widehat{\mathbf{F}}_{i,j,l,m}L_{l}(\xi)L_{m}(\eta).
\end{equation}
Then, the nodal values of the state variable $\mathbf{U}$ at the
solution points are updated by the following equations

\begin{eqnarray}
\frac{\partial\widehat{\mathbf{U}}_{i,j,l,m}}{\partial t}+\frac{\partial\widehat{\mathbf{F}}_{i,j}(\xi_{l},\eta_{m})}{\partial\xi}+\frac{\partial\widehat{\mathbf{G}}_{i,j}(\xi_{l},\eta_{m})}{\partial\eta}+\delta_{i,j}(\xi_{l},\eta_{m}) & = & 0,\quad1\leq l,m\leq K+1\label{eq:CPR-2D}
\end{eqnarray}
where 

\begin{eqnarray*}
\delta_{i,j}(\xi_{l},\eta_{m}) & = & \left[\overline{F}_{i,j}(-1,\eta_{m})-\widehat{F}_{i,j}(-1,\eta_{m})\right]g'{}_{L}(\xi_{l})+\left[\overline{F}_{i,j}(1,\eta_{m})-\widehat{F}_{i,j}(1,\eta_{m})\right]g'{}_{R}(\xi_{l})\\
 &  & +\left[\overline{G}_{i,j}(\xi_{l},-1)-\widehat{G}_{i,j}(\xi_{l},-1)\right]g'{}_{L}(\eta_{m})+\left[\overline{G}_{i,j}(\xi_{l},1)-\widehat{G}_{i,j}(\xi_{l},1)\right]g'{}_{R}(\eta_{m}).
\end{eqnarray*}
Here $\delta_{i,j}(\xi,\eta)$ is a correction flux polynomial, $g_{L}(\xi)$
and $g_{R}(\xi)$ are both the degree $K+1$ polynomials called correction
functions. $\overline{F}$ and $\overline{G}$ are the common fluxes.
Riemann solvers can be used to compute common fluxes, such as Lax-Friedrichs,
Roe, Osher, AUSM, HLL, and their modifications. We refer to papers
\citet{Qu2020,zhu2016} and references therein. 

In this paper, the CPR method takes Legendre-Gauss (LG) points as
solution points and Radau polynomials as correction function, which
is equivalent to a specific DG method. For the equivalence of the
CPR and DG, we refer to details in \citet{Huynh2007,Mengaldo2016}.
Thus, correction functions are $g_{L}=R_{R,K+1},\,\, g_{R}=R_{L,K+1}$.
Here $R_{R,K+1}$ and $R_{L,K+1}$ are the right Radau polynomials
$R_{R,K+1}=\frac{\left(-1\right)^{K+1}}{2}\left(P_{K+1}-P_{K}\right)$
and the left Radau polynomial $R_{L,K+1}=\frac{1}{2}\left(P_{K+1}+P_{K}\right)$,
correspondingly. $P_{K}$ is the Legendre polynomial of order $K$.
For the case $K=4$, we have 
\[
g'_{L}(\xi)=-\frac{1}{16}(315\xi^{4}-140\xi^{3}-210\xi^{2}+60\xi+15),\quad g'{}_{R}(\xi)=\frac{1}{16}(315\xi^{4}+140\xi^{3}-210\xi^{2}-60\xi+15).
\]
In addition, Legendre-Gauss-Lobatto (LGL) points are taken as flux
points, as shown in Fig. \ref{fig:SP-FP of CPR-stru}. For quadrilateral
cells, the operations are in fact one-dimensional. Thus, for two-dimensional
case, each element has $K+1$ solution points and $K+2$ flux points
in each direction. In this paper, we mainly consider 5th-order CPR
(CPR5) with $K=4$. In this case, we have solution points $sp_{1}=-\frac{1}{3}\sqrt{\frac{1}{7}\left(35+2\sqrt{70}\right)}$,
$sp_{2}=-\frac{1}{3}\sqrt{\frac{1}{7}\left(35-2\sqrt{70}\right)}$,
$sp_{3}=0$, $sp_{4}=-sp2$, $sp_{5}=-sp1$, and flux points $fp_{1}=-1$,
$fp_{2}=-\sqrt{\frac{1}{3}+\frac{2}{3\sqrt{7}}}$, $fp_{3}=-\sqrt{\frac{1}{21}\left(7-2\sqrt{7}\right)}$,
$fp_{4}=-fp_{3}$, $fp_{5}=-fp_{2}$, $fp_{6}=-fp_{1}$. The solution
points and flux points are shown in Fig. \ref{fig:Distribution-of-solution}(a). 

\begin{center}
\begin{figure}
\centering{}\includegraphics[scale=0.4]{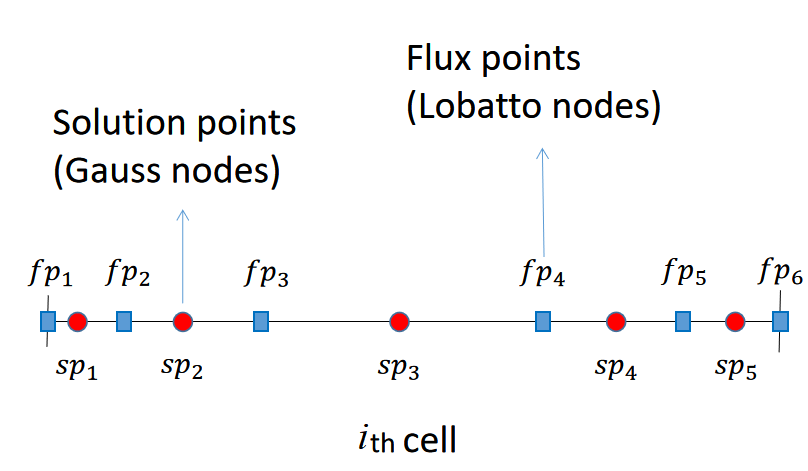}\caption{Solution points and flux points in 1D ($K=4$). \label{fig:SP-FP of CPR-stru}}
\end{figure}

\par\end{center}

\begin{center}
\begin{figure}
\begin{centering}
{\footnotesize }\subfloat[CPR5 cells]{\begin{centering}
{\footnotesize \includegraphics[width=0.4\textwidth]{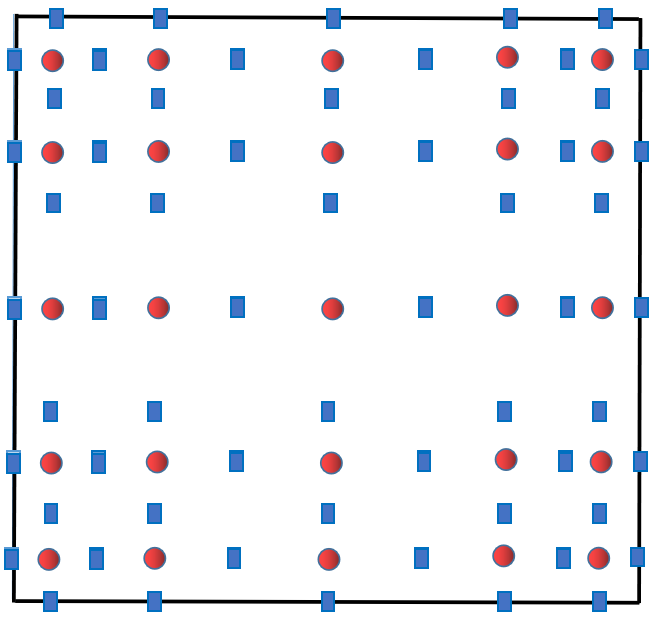}}
\par\end{centering}{\footnotesize \par}

{\footnotesize }}\subfloat[C5NNW5 subcells]{\begin{centering}
\includegraphics[width=0.394\textwidth]{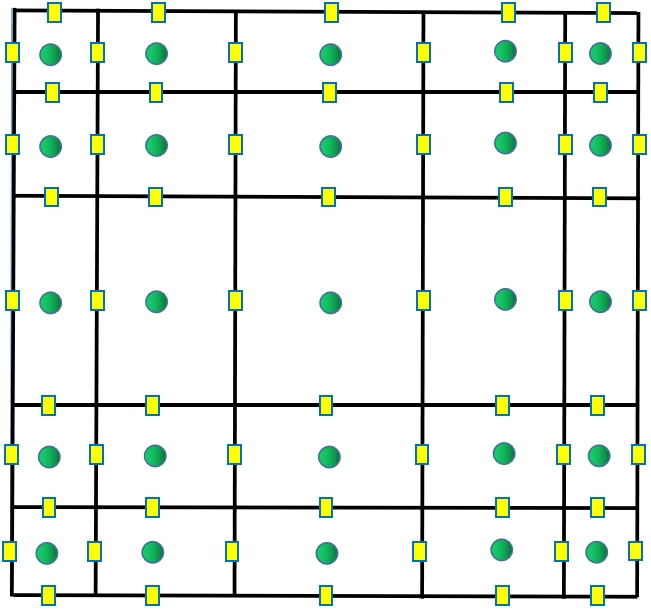}
\par\end{centering}

}
\par\end{centering}

\caption{Solution points and flux points for CPR5 cells and C5NNW5 subcells
($K=4$), where red circle nodes and square blue nodes in Figure (a)
denote solution points and flux points for CPR5, while green circle
nodes and yellow square nodes in Figure (b) denote solution points
and flux points for C5NNW5. \label{fig:Distribution-of-solution}}
\end{figure}

\par\end{center}

\section{Novel subcell schemes }

This section is devoted to develop novel shock capturing schemes based
on Legendre-Gauss (LG) solution points, which will be applied in subcell
limiting for the high-order CPR method. To ensure the subcell limiting
being accurate and robust, both high-order and low-order shock capturing
schemes are constructed.

\subsection{High-order CNNW schemes }

For capturing shock effectively, nonlinear interpolations have been
used in MUSCL \citet{Leer1974,Leer1979}, WCNS \citet{Deng2000} and
WENO \citet{Shu2009} schemes to prevent interpolation across discontinuities.
Inspired by nonlinear interpolation in WCNS and compact differencing
in CPR, we develop new shock capturing schemes by combining nonlinear
interpolation with compact flux differencing. The new shock capturing
schemes are constructed through discretizing flux derivative by a
compact difference operator based on Riemann fluxes within one cell
to make scheme compact, and obtaining the left and right variable
values used in Riemann fluxes by nonuniform nonlinear weighted (NNW)
interpolation from solution points to flux points on computational
space. In the following, the new schemes are called compact nonuniform
nonlinear weighted schemes (CNNW). 

For a one-dimensional element, solution points and flux points of
CNNW schemes are located in staggered form, which means each solution
point locates between two flux points. In order to combine with a
CPR with LG solution points, we takes $\left(K+1\right)$ LG solution
points and $\left(K+2\right)$ Legendre-Gauss-Lobatto (LGL) points
as flux points, as shown in Fig. \ref{fig:Distribution-of-solution}.
In the following, a fifth-order shock capturing scheme is constructed
by using a $5$th-order NNW interpolation and a $5$th-order compact
flux difference operator.

\subsubsection{The fifth-order compact nonuniform nonlinear weighted scheme (C5NNW5)}

Nonuniform nonlinear weighted (NNW) interpolation is taken to obtain
the left and right flow-field variable values used in Riemann fluxes
at flux points. NNW interpolation takes a stencil of several adjacent
solution points, as shown in Fig. \ref{fig:Stencil-for-interpolation}.
A $5$th-order NNW interpolation in one-dimensional case for obtaining
the right values at the first flux point of the $i$th cell is given
in Appendix A. The NNW interpolation procedure is similar as those
in WENO and WCNS. It is worth noticing that there are five LG solution
points for each cell and the distances between two adjacent solution
points are different. 

\begin{center}
\begin{figure}
\begin{centering}
\includegraphics[scale=0.3]{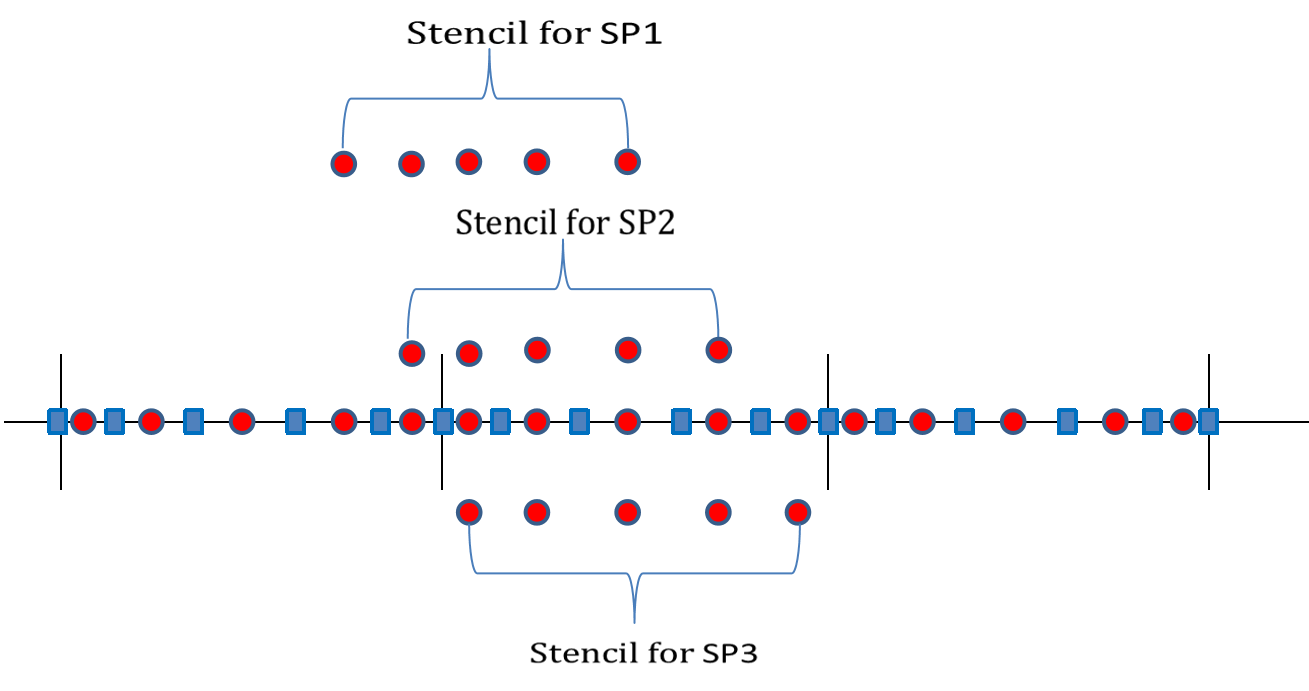}
\par\end{centering}

\caption{Stencil for high-order interpolation from solution points to flux
points\label{fig:Stencil-for-interpolation}}
\end{figure}

\par\end{center}

Smoothness indicators used in nonlinear interpolation need to be carefully
calculated in the case of nonuniformly spaced solution points. Suppose
grid transformation from physical coordinates to computational coordinates
is a linear transformation. Then, for the $l$th subcell $[\xi_{fp_{l}},\xi_{fp_{l+1}}]$
we have $x_{\xi}=\frac{\Delta x_{l}}{\Delta\xi_{l}}$ with $\Delta x_{l}=x_{fp_{\left(l+1\right)}}-x_{fp_{l}}$
and $\Delta\xi_{l}=\xi_{fp_{(l+1)}}-\xi_{fp_{l}}$. A smoothness indicator
used in WENO schemes {[}\citet{Jiang1996,Shi2002}{]} for the $m$th
small stencil $S_{i,sp_{l}}^{(m)}$ of the $l$th solution point $sp_{l}$
has the form

\begin{eqnarray}
IS_{m,sp_{l}} & =\sum_{\alpha=1}^{2} & \int_{x_{fp_{l}}}^{x_{fp_{\left(l+1\right)}}}\left(\Delta x_{l}\right)^{2\alpha-1}\left(\frac{\partial^{\alpha}p_{sp_{l}}^{(m)}}{\partial x^{\alpha}}\right)^{2}dx,\label{eq:WENO indicator}
\end{eqnarray}
where $p_{sp_{l}}^{(m)}(\xi)$ is a Lagrange interpolation polynomial
of degree 2 in the small stencil $S_{i,sp_{l}}^{(m)}$. Then, $\frac{\partial p_{sp_{l}}^{(m)}}{\partial\xi}$
is a linear polynomial and $\frac{\partial^{2}p_{sp_{l}}^{(m)}}{\partial\xi^{2}}$
is a constant. According to $\frac{\partial p_{sp_{l}}^{(m)}}{\partial x}=\frac{\partial p_{sp_{l}}^{(m)}}{\partial\xi}\cdot\frac{1}{x_{\xi}}$,
the smoothness indicator in (\ref{eq:WENO indicator}) becomes
\begin{eqnarray*}
IS_{m,sp_{l}} & = & \left\{ \int_{x_{fp_{l}}}^{x_{fp_{\left(l+1\right)}}}\Delta x_{l}\left(\frac{\partial p_{sp_{l}}^{(m)}}{\partial x}\right)^{2}dx\right\} +\int_{x_{fp_{l}}}^{x_{fp_{(l+1)}}}\left(\Delta x_{l}\right)^{3}\left(\frac{\partial^{2}p_{sp_{l}}^{(m)}}{\partial x^{2}}\right)^{2}dx\\
 & = & \left\{ \frac{1}{6}\left[\left(\left.\frac{\partial p_{sp_{l}}^{(m)}}{\partial\xi}\right|_{\xi_{fp_{l}}}\right)^{2}+4\left(\left.\frac{\partial p_{sp_{l}}^{(m)}}{\partial\xi}\right|_{\frac{\xi_{fp_{l}}+\xi_{fp_{(l+1)}}}{2}}\right)^{2}+\left(\left.\frac{\partial p_{sp_{l}}^{(m)}}{\partial\xi}\right|_{\xi_{fp_{(l+1)}}}\right)^{2}\right]\left(\Delta\xi_{l}\right)^{2}\right\} \\
 &  & +\left[\left(\left.\frac{\partial^{2}p_{sp_{l}}^{(m)}}{\partial\xi^{2}}\right|_{sp_{l}}\right)^{2}\left(\Delta\xi_{l}\right)^{4}\right].
\end{eqnarray*}
In this paper, to reduce computation we approximate the term$\int_{x_{fp_{l}}}^{x_{fp_{\left(l+1\right)}}}\Delta x_{l}\left(\frac{\partial p_{sp_{l}}^{(m)}}{\partial x}\right)^{2}dx$
by $\left(\left.\frac{\partial p_{sp_{l}}^{(m)}}{\partial\xi}\right|_{sp_{l}}\right)^{2}\left(\Delta\xi_{l}\right)^{2}$
and obtain the following new simple smoothness indicator:
\begin{eqnarray}
IS_{m,sp_{l}}^{new} & = & \left(\left.\frac{\partial p_{sp_{l}}^{(m)}}{\partial\xi}\right|_{sp_{l}}\right)^{2}\left(\Delta\xi_{l}\right)^{2}+\left(\left.\frac{\partial^{2}p_{sp_{l}}^{(m)}}{\partial\xi^{2}}\right|_{sp_{l}}\right)^{2}\left(\Delta\xi_{l}\right)^{4}.\label{eq:new-smoothness-indicator}
\end{eqnarray}

After the left and right values at six flux points are obtained by
5th-order NNW interpolation, a $5$th-order compact flux difference
operator is used to discretize the flux derivative. Lagrange polynomial
based on flux points is 
\begin{equation}
\widetilde{F}(\xi)=\sum_{l=1}^{6}\overline{F}_{i,fp_{l}}L_{l}(\xi),\label{eq:flux-polynomial}
\end{equation}
where $\overline{F}_{i,fp_{l}}$ is Riemann flux at LGL flux points
$\overline{F}_{i,fp_{l}}(u_{i,fp_{l}}^{L},u_{i,fp_{l}}^{R})$. Here
$u_{i,fp_{l}}^{L}$ and $u_{i,fp_{l}}^{R}$ are obtained by the $5$th-order
NNW interpolation. Then, the 5th-order compact flux difference operator
is obtained by calculating the first-order derivative of the Lagrange
polynomial (\ref{eq:flux-polynomial}) at solution points, 

\[
\frac{\partial\widetilde{F}}{\partial\xi}|_{i,sp_{m}}=\sum_{l=1}^{6}a_{m,l}\overline{F}_{i,fp_{l}}.
\]

The difference between high-order CNNW and high-orde CPR is that CNNW
uses nonlinear interpolation based on solution points of the cell
and its neighbor cells, uses Riemann fluxes for each flux points,
and does not use correction function.

\subsection{Low-order CNNW schemes}

\subsubsection{C2NNW5}

A low-order shock capturing scheme with high resolution is constructed
by taking 5th-order NNW interpolation proposed in subsection 3.1.1
and Appendix A with following 2nd-order finite difference operator
(C2NNW5) 

\begin{equation}
\frac{\partial\widetilde{F}}{\partial\xi}|_{i,sp_{l}}=\frac{\overline{F}_{i,fp_{(l+1)}}-\overline{F}_{i,fp_{l}}}{\Delta\xi_{l}},\quad l=1,2,\cdots,K+1\label{eq:FD2}
\end{equation}
where $sp_{l}$ are LG solution points, $fp_{l}$ are flux points
and $\Delta\xi_{l}=\xi_{fp_{(l+1)}}-\xi_{fp_{l}}$.

\subsubsection{C2NNW2}

A low-order shock capturing scheme with good robustness is constructed
by taking a 2nd-order nonlinear weighted interpolation (NNW2) with
the 2nd-order finite differential operator in (\ref{eq:FD2}) (C2NNW2).
The NNW2 interpolation based on a stencil of three adjacent nonuniformly
spaced solution points is constructed by using inverse distance weighted
interpolation \citet{Frink1991}  to obtain values at flux points
and then using Birth limiter \citet{Birth1989} to limit linear reconstruction,
as shown in Fig. \ref{fig:C2NNW2}. Details of the NNW2 interpolation
are given in Appendix E. 

\begin{center}
\begin{figure}
\begin{centering}
\includegraphics[scale=0.4]{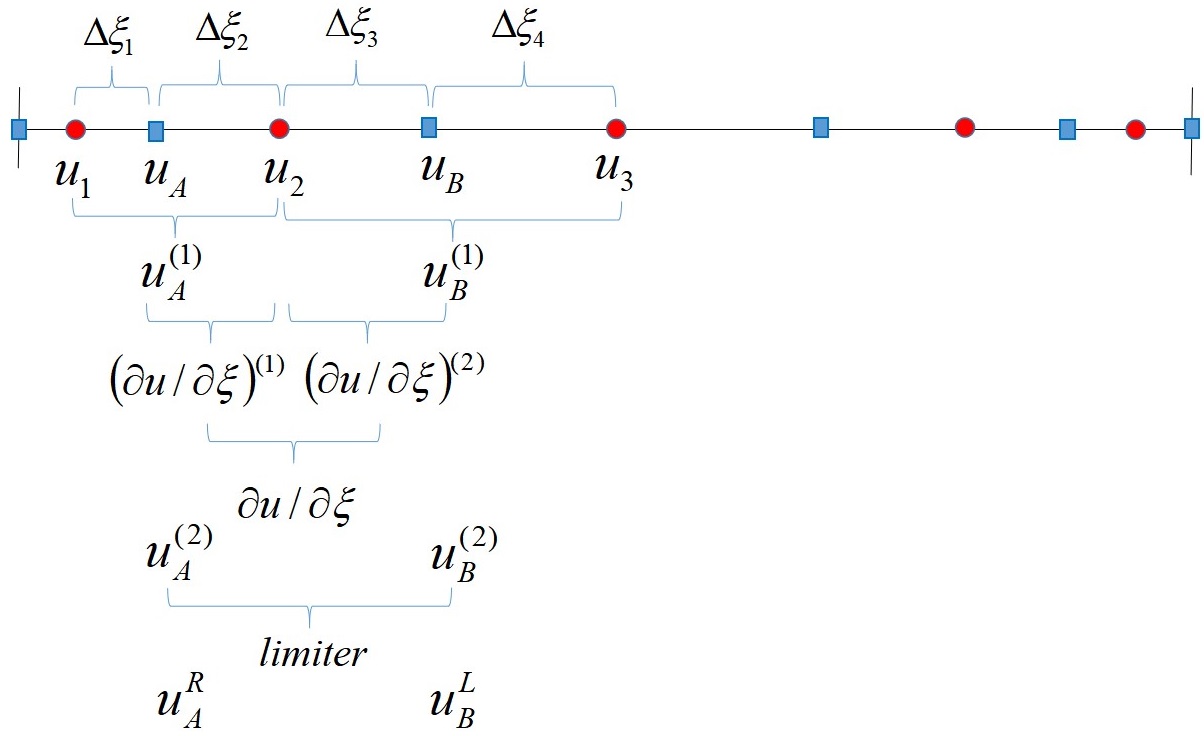}
\par\end{centering}

\caption{Stencil for NNW2 interpolation from solution points to flux points\label{fig:C2NNW2}}
\end{figure}

\par\end{center}

\subsection{Comparison of interpolation methods and difference operators in CNNW
and CPR}

Stencils of interpolation and difference operator for CPR5, C5NNW5,
C2NNW5, C2NNW2 in solving 1D conservation law are shown in Table \ref{tab:different schemes}.
For comparison, the fifth-order weighted compact nonlinear scheme
(WCNS5) in \citet{Deng2000} with hybrid cell-edge-node finite difference
operator \citet{Deng2000} is also shown in the Table \ref{tab:different schemes}. 

\begin{flushleft}
\begin{table}
\begin{centering}
\begin{tabular}{|c|l|c|}
\hline 
\multicolumn{1}{|c|}{{\small Schemes}} & \multicolumn{1}{c|}{{\small Interpolation and FD operator}} & Stencil\tabularnewline
\hline 
\multirow{4}{*}{{\small CPR5}} & {\small Lagrange interpolation on Legendre-Gauss SPs: } & \multirow{2}{*}{{\small \includegraphics[width=0.3\textwidth]{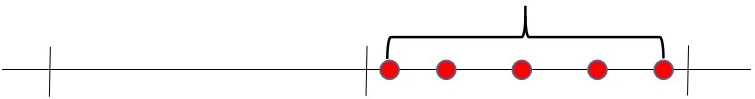}}}\tabularnewline
 & {\small $l_{1}u_{i,sp_{1}}+l_{2}u_{i,sp_{2}}+l_{3}u_{i,sp_{3}}+l_{4}u_{i,sp_{4}}+l_{5}u_{i,sp_{5}}$} & \tabularnewline
\cline{2-3} 
 & {\small Compact FD5 + Correction function: } & \tabularnewline
 & {\small $\frac{\partial\widetilde{F}}{\partial\xi}|_{i,sp_{2}}=\sum_{k=1}^{6}a_{k}\widehat{F}_{i,fp_{k}}$}+$\delta_{i}$ & {\small \includegraphics[width=0.3\textwidth]{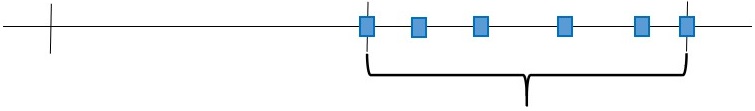}}\tabularnewline
\hline 
\multirow{4}{*}{{\small C5NNW5}} & {\small NNW5 on Legendre-Gauss SPs: } & \multirow{2}{*}{{\small \includegraphics[width=0.3\textwidth]{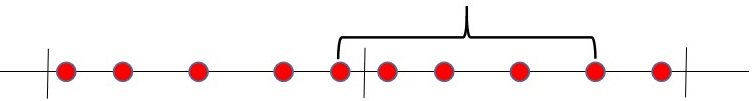}}}\tabularnewline
 & {\small $\omega_{1}p_{sp_{2}}^{(1)}+\omega_{2}p_{sp_{2}}^{(2)}+\omega_{3}p_{sp_{2}}^{(3)}$} & \tabularnewline
\cline{2-3} 
 & {\small Compact FD5: } & \multicolumn{1}{l|}{}\tabularnewline
 & {\small $\frac{\partial\widetilde{F}}{\partial\xi}|_{i,sp_{2}}=\sum_{k=1}^{6}a_{k}\overline{F}_{i,fp_{k}}$} & {\small \includegraphics[width=0.3\textwidth]{fig/bplot-cfd5}}\tabularnewline
\hline 
\multirow{4}{*}{{\small C2NNW5}} & {\small NNW5 on Legendre-Gauss SPs} & \multirow{2}{*}{{\small \includegraphics[width=0.3\textwidth]{fig/bplot-nnw5}}}\tabularnewline
 & {\small $\omega_{1}p_{sp_{2}}^{(1)}+\omega_{2}p_{sp_{2}}^{(2)}+\omega_{3}p_{sp_{2}}^{(3)}$} & \tabularnewline
\cline{2-3} 
 & {\small FD2: } & \multirow{2}{*}{{\small \includegraphics[width=0.3\textwidth]{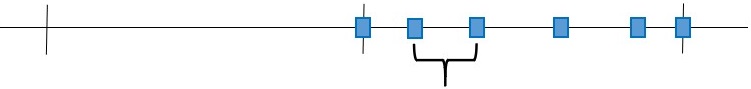}}}\tabularnewline
 & {\small $\frac{\partial\widetilde{F}}{\partial\xi}|_{i,sp_{2}}=\frac{1}{\Delta\xi_{2}}(\overline{F}_{i,fp_{3}}-\overline{F}_{i,fp_{2}})$} & \tabularnewline
\hline 
\multirow{4}{*}{{\small C2NNW2}} & {\small NNW2 on Legendre-Gauss SPs: } & \multirow{2}{*}{{\small \includegraphics[width=0.3\textwidth]{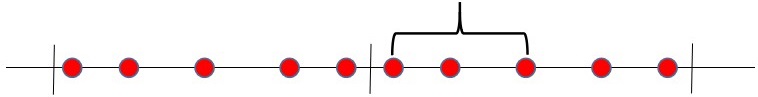}}}\tabularnewline
 & {\small $au_{i,sp_{1}}+bu_{i,sp_{2}}+cu_{i,sp_{3}}$} & \tabularnewline
\cline{2-3} 
 & FD2{\small : } & \multirow{2}{*}{{\small \includegraphics[width=0.3\textwidth]{fig/bplot-fd2}}}\tabularnewline
 & {\small $\frac{\partial\widetilde{F}}{\partial\xi}|_{i,sp_{2}}=\frac{1}{\Delta\xi_{2}}(\overline{F}_{i,fp_{3}}-\overline{F}_{i,fp_{2}})$} & \tabularnewline
\hline 
\multirow{4}{*}{{\small WCNS5}} & {\small WCNS interpolation on uniformly-spaced SPs: } & \multirow{2}{*}{{\small \includegraphics[width=0.3\textwidth]{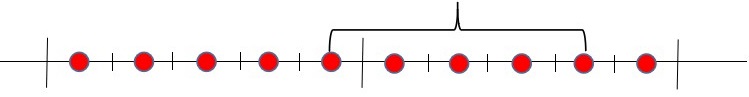}}}\tabularnewline
 & {\small $\hat{\omega}_{1}\hat{p}_{sp2}^{(1)}+\hat{\omega}_{2}\hat{p}_{sp2}^{(2)}+\hat{\omega}_{3}\hat{p}_{sp2}^{(3)}$ } & \tabularnewline
\cline{2-3} 
 & {\small Hybrid FD6: } & \multirow{2}{*}{{\small \includegraphics[width=0.3\textwidth]{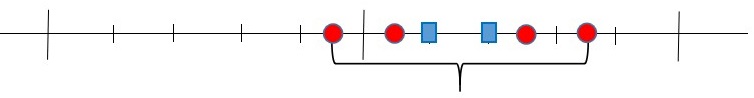}}}\tabularnewline
 & {\small $\begin{array}{ccc}
\frac{\partial\widetilde{F}}{\partial\xi}|_{sp2}= & b_{1}(\overline{F}_{i,fp3}-\overline{F}_{i,fp2})\\
 & +b_{2}(F_{i,sp3}-F_{i,sp1}) & +b_{3}(F_{i,sp4}-F_{i-1,sp5})
\end{array}$} & \tabularnewline
\hline 
\end{tabular}
\par\end{centering}

\noindent \centering{}\caption{Stencil of interpolation and FD operators for different schemes at
the 2nd solution point of a cell.\label{tab:different schemes}}
\end{table}

\par\end{flushleft}

\section{A priori subcell CNNW limiting approach for CPR method}

In this section, a priori subcell limiting approach based on the proposed
CNNW schemes is developed for the fifth-order CPR scheme (CPR5) with
five Legendre-Gauss solution points presented in Subsection 2.1. Firstly,
troubled cell indicators are used to detect troubled cells which may
have discontinuities. Then, the troubled cells are decomposed into
subcells and computed by the CNNW schemes while other cells are computed
by the CPR scheme.

\subsection{Troubled cell indicator }

In order to find troubled cells, we takes the indicator proposed in
\citet{Hennemann2021}, which follow ideas presented by Persson and
Peraire \citet{Persson2006} and consider the rate of the highest
mode to the overall modal energy. Firstly, the representation of the
quantity $\epsilon=\rho p$ with Lagrange interpolation polynomials
of degree $N$ is transformed to a modal representation with Legendre
interpolation polynomials. Secondly, the maximum of proportion of
the highest modes and proportion of the second highest mode to the
total energy of the Legendre interpolation polynomial is calculated
as

\begin{eqnarray}
EI & = & max\left(\frac{m_{N}^{2}}{\sum_{j=0}^{N}m_{j}^{2}},\frac{m_{N-1}^{2}}{\sum_{j=0}^{N-1}m_{j}^{2}}\right),\label{eq:proportion of high modes}
\end{eqnarray}
where $\{m_{j}|j=0,1,\cdots,N\}$ are the modal coefficients. 

In this paper, to consider the jump in cell interfaces, $EI$ is calculated
by a higher degree polynomial based on the ``extended'' stencil
consisting of five solution points in the cell and two end points
at cell interfaces. Thus, for the $5$th-order CPR, the indicator
of the $i$th cell is calculated based on the stencil with seven points
$\{\epsilon_{L},\epsilon_{i,1},\epsilon_{i,2},\epsilon_{i,3},\epsilon_{i,4},\epsilon_{i,5},\epsilon_{R}\}$,
where $\epsilon_{i,1},\epsilon_{i,2},\epsilon_{i,3},\epsilon_{i,4},\epsilon_{i,5}$
are the quantity $\epsilon=\rho p$ at five solution points, $\epsilon_{L}=aver(\epsilon_{i-1,5},\epsilon_{i,1})$
and $\epsilon_{R}=aver(\epsilon_{i,5},\epsilon_{i+1,1})$ are Roe
average values at cell interfaces. Here $aver(\epsilon_{1},\epsilon_{2})=aver(\rho_{1},\rho_{2})\cdot aver(p_{1},p_{2})$
and $aver(\cdot,\cdot)$ is the Roe average function. 

We take a threshold value
\begin{equation}
T(N)=a\cdot10^{-b(N+1)^{1/4}}.\label{eq:TN}
\end{equation}
The parameter $b$ is predetermined as $1.8$, which is the same as
those in \citet{Hennemann2021}. 

It is worth noticing that for CPR5 with $K=4$ we take $N=K+2=6$
in (\ref{eq:proportion of high modes}) and (\ref{eq:TN}). In this
paper, we set the threshold value in the MDHE indicator (\ref{eq:TN})
to be 

\begin{eqnarray}
c(a) & = & a\cdot10^{-1.8(6+1)^{1/4}}.\label{eq:da}
\end{eqnarray}
If $EI\geq c(a)$, the element is denoted as a troubled cell. Thus,
$c(a)$ control the size of CNNW area.

The indicator based on the rate of the highest mode \citet{Persson2006}
is usually called highest modal decay (MDH) indicator. For simplicity,
we denote the highest modal decay indicator based on the ``extended''
stencil as MDHE indicator in this paper.

\subsection{Subcell limiting based on CNNW}

After troubled cell detection, the troubled cells are decomposed into
subcells and solved by CNNW, as shown in Fig. \ref{fig:subcell-limiting-1D},
while other cells are computed by CPR. Thus, a CPR scheme based on
subcell CNNW limiting (abbr. CPR-CNNW) is a hybrid scheme. Subcell
schemes in troubled cells are chosen according to the magnitude of
the MDHE indicator $EI$ in (\ref{eq:proportion of high modes}) on
the extended stencil of $N=K+2$ nodes. To ensure CPR-CNNW having
high-resolution and having good robustness in shock capturing, a $p$-adaptive
limiting procedure is suggested by using both high-order accurate
shock capturing schemes (C5NNW5) and low-order robust shock capturing
schemes (C2NNW5 and C2NNW2) to accomplish transition from smooth region
to discontinuous region, as shown in Fig. \ref{fig:CPR-subcell-limiting}. 

We define the partition vector $\mathbf{dv}=(S_{1},S_{2},S_{3})$
with three parameters $S_{1}$, $S_{2}$ and $S_{3}$ controlling
region division, as shown in Fig. \ref{fig:CPR-subcell-limiting}.
Then, the hybrid CPR-CNNW scheme (abbr. HCCS) can be expressed as

\begin{equation}
\mathrm{HCCS}=\begin{cases}
\mathrm{CPR5}, & 0\leq EI\leq S_{1},\\
\mathrm{C5NNW5}, & S_{1}<EI\leq S_{2},\\
\mathrm{C2NNW5}, & S_{2}<EI\leq S_{3},\\
\mathrm{C2NNW2}, & S_{3}<EI\leq1.
\end{cases}\label{eq:HS}
\end{equation}
The hybrid CPR-CNNW scheme is also denoted by $\mathrm{HCCS}(d_{1},d_{2},d_{3},d_{4})$
with $d_{1},d_{2},d_{3},d_{4}$ marking status of CPR, C5NNW5, C2NNW5,
C2NNW2 correspondingly. Here $d_{i}=1$ means that the corresponding
scheme is contained by the hybrid scheme, otherwise not included.
The values of $d_{1},d_{2},d_{3},d_{4}$ are determined by the relationship
between the three parameters $S_{1},S_{2},S_{3}$. Thus, by controlling
the vector $\mathbf{dv}$, the CPR-CNNW scheme (\ref{eq:HS}) can
contain some of the four schemes. 

In this paper, we will mainly test a CPR-CNNW scheme with $p$-adaption
($\mathrm{HCCS}(1,1,1,1)$) which contains all of the four schemes
and a CPR-CNNW scheme without $p$-adaption ($\mathrm{HCCS}(1,0,0,1)$)
which only contains CPR and C2NNW2, as shown in Table \ref{tab:HS}.

\begin{center}
\begin{figure}
\begin{centering}
\includegraphics[width=0.7\textwidth]{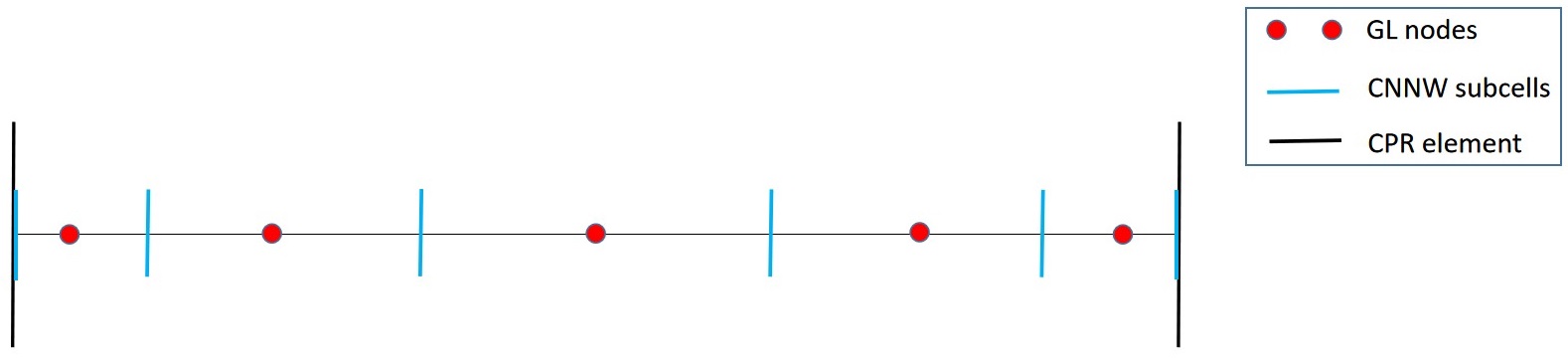}
\par\end{centering}

\caption{A fifth-order CPR element with CNNW subcells \label{fig:subcell-limiting-1D}}
\end{figure}

\par\end{center}

\begin{center}
\begin{figure}
\begin{centering}
\includegraphics[width=0.8\textwidth]{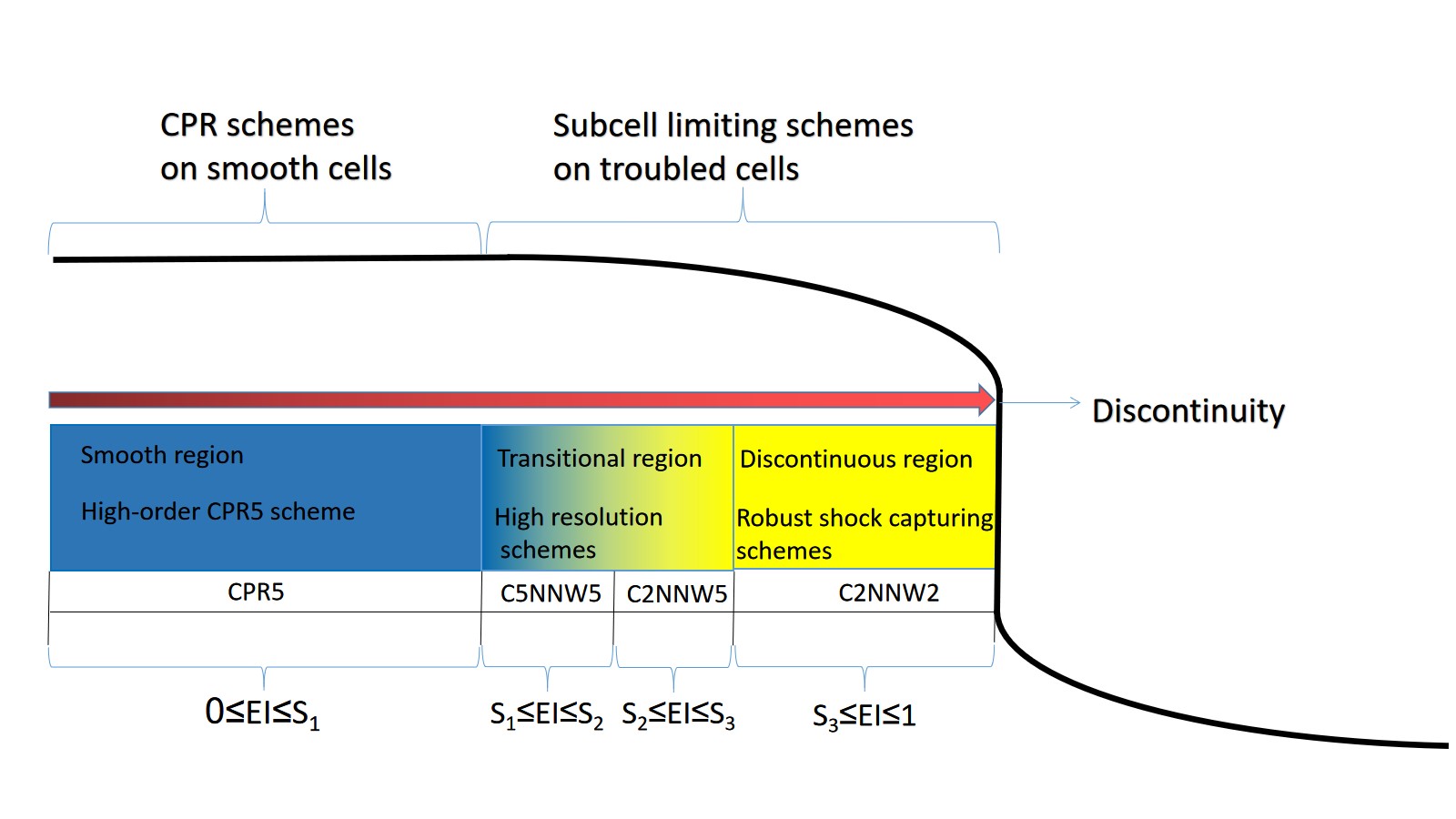}
\par\end{centering}

\caption{Subcell limiting of CPR-CNNW \label{fig:CPR-subcell-limiting}}
\end{figure}

\par\end{center}

\begin{center}
\begin{table}
\begin{centering}
\begin{tabular}{cccc}
\hline 
\multirow{1}{*}{{\small Schemes}} & {\small $\mathrm{HCCS}(d_{1},d_{2},d_{3},d_{4})$} & {\small $dv=(S1,S2,S3)$} & {\small containing schemes }\tabularnewline
\hline 
{\small CPR-CNNW with $p$-adaption} & {\small $\mathrm{HCCS}(1,1,1,1)$} & {\small 0<S1<S2<S3<1} & {\small CPR5,}\textit{\small C5NNW5,C2NNW5,C2NNW2}\tabularnewline
\hline 
{\small CPR-CNNW without $p$-adaption} & {\small $\mathrm{HCCS}(1,0,0,1)$} & {\small 0<S1=S2=S3<1} & {\small CPR5,}\textit{\small C2NNW2}\tabularnewline
\hline 
\end{tabular}
\par\end{centering}

\noindent \centering{}\caption{{\small Two CPR-CNNW} schemes.\label{tab:HS}}
\end{table}

\par\end{center}

\subsection{Interface treatment}

CNNW has the same solution points as CPR, which makes the CPR based
on subcell CNNW limiting approach have some merits. Firstly, there
is no data exchange between solution points of different schemes and
thus the proposed hybrid scheme can take less computations. Secondly,
extra state values required in interpolation method on troubled cells
are taken from neighboring cells directly, and thus there is no need
to add ghost cells for exchanging the state values. 

The only thing needed to do in interface treatment is calculation
of Riemann fluxes at the interface of different schemes. The Riemann
fluxes at the interface of scheme A and scheme B are calculated based
on the left and right values interpolated from scheme A and scheme
B, correspondingly. For example, Riemann fluxes at the interface of
CPR and C5NNW5 are computed based on one side from CPR cell and the
other from C5NNW5.

\section{Theoretical analysis on spectral properties and conservation}

\subsection{Spectral properties of high-order CNNW and CPR}

Finite difference schemes usually obtain the spectrum by Fourier method
while DG-type method which locates several solution points in one
cell usually calculate the spectrum based on local discrete matrices.
To make fair comparisons, the eigenvalues of spatial discretization
matrix of different schemes are calculated by the same method on local
discrete matrices. We analyze the spectrum by local discrete matrices
and prove that all eigenvalues comes from the same function and each
scheme has a unique spectrum curve. 

Suppose the computational domain is decomposed to $M$ cells. The
semi-discretization form of one-dimensional linear advection equation
$u_{t}+u_{x}=0$ with periodic boundary condition can be written as
the first form:

\begin{equation}
\frac{\partial}{\partial t}\mathbf{U}=\frac{1}{\Delta x}E\mathbf{U}\label{eq:Form1}
\end{equation}
where $\mathbf{U}=(u_{1},u_{2},\cdots,u_{(K+1)M})^{T}$, $\Delta x$
is spatial step and $-\frac{1}{\Delta x}E$ is spatial discretization
matrix of first-order derivative. The semi-discretization form can
also be written as the second form:

\begin{eqnarray}
\frac{\partial}{\partial t}\left[\begin{array}{c}
u_{j,1}\\
u_{j,2}\\
\vdots\\
u_{j,K+1}
\end{array}\right] & =\frac{1}{(K+1)\Delta x} & \left(A\left[\begin{array}{c}
u_{j-1,1}\\
u_{j-1,2}\\
\vdots\\
u_{j-1,K+1}
\end{array}\right]+B\left[\begin{array}{c}
u_{j,1}\\
u_{j,2}\\
\vdots\\
u_{j,K+1}
\end{array}\right]+C\left[\begin{array}{c}
u_{j+1,1}\\
u_{j+1,2}\\
\vdots\\
u_{j,K+1}
\end{array}\right]\right)\label{eq:Form2}
\end{eqnarray}
 where $j=1,2,\cdots,M$. Then, the matrix $E$ can be written as 

\begin{eqnarray}
E & = & \frac{1}{K+1}\left[\begin{array}{cccccc}
B & C & 0 & 0 & 0 & A\\
A & B & C & 0 & 0 & 0\\
0 & A & B & C & 0 & 0\\
0 & 0 & A & B & C & 0\\
0 & 0 & 0 & A & B & C\\
C & 0 & 0 & 0 & A & B
\end{array}\right]_{(K+1)M\times(K+1)M},\label{eq:matrixE}
\end{eqnarray}
where $A$, $B$, $C$ are $(K+1)\times(K+1)$ matrix. The matrix
$E$ is a block circulant matrix. 

In the following Theorem 2.1, we prove that all eigenvalues of the
spatial discretization matrix can be obtained by collecting the eigenvalues
of local spatial matrices. In addition, all eigenvalues comes from
the same function and thus all the eigenvalues are on a unique spectrum
curve. 

\label{thm:eigenvalues}\textbf{Theorem 2.1} \textit{The matrix $E$
in spacial discretization matrix with form (\ref{eq:matrixE}) has
following properties:}

\textit{(1)All the eigenvalues of $E$ are given by }

\textit{
\begin{eqnarray*}
\{x|EX=xX,\, X\in\mathbb{C}^{(K+1)M}\} & = & \sum_{m=0}^{M-1}\{x|H_{m}Y_{m}=xY_{m},\, Y_{m}\in\mathbb{C}^{K+1}\},
\end{eqnarray*}
where $H_{m}=H(\phi_{m})$, $\phi_{m}=m\frac{2\pi}{M},\,\, m=0,1,2,\cdots,M-1$
and
\begin{eqnarray}
H(\phi) & = & \left(Ae^{-i\phi}+B+Ce^{i\phi}\right)/(K+1),\,\,\,\,0\leq\phi<2\pi.\label{eq:Hfai}
\end{eqnarray}
In other word, $SH\triangleq Spec\left(E\right)=\left\{ Spec\left(H_{0}\right),Spec\left(H_{1}\right),\cdots,Spec\left(H_{M-1}\right)\right\} $.}

\textit{(2) Suppose $G_{m}=G(\phi_{m})$ with 
\begin{eqnarray}
G(\phi) & = & \left(Ae^{-i\phi(K+1)}+B+Ce^{i\phi(K+1)}\right)/(K+1),\,\,\,\,0\leq\phi<2\pi,\label{eq:Gfai}
\end{eqnarray}
and $Spec\left(G_{m}\right)=\{\lambda^{(l)}(G(\phi_{m}))|l=1,2,\cdots,K+1\}$,
$SG\triangleq\left\{ Spec\left(G_{0}\right),Spec\left(G_{1}\right),\cdots,Spec\left(G_{M-1}\right)\right\} $,
where $\phi_{m}=m\frac{2\pi}{M},\,\, m=0,1,2,\cdots,M-1$. It can
be proved that if $mod(M,K+1)\neq0$ then}

\textit{
\begin{eqnarray*}
SG & = & SH,
\end{eqnarray*}
else
\begin{eqnarray*}
SG=\left\{ Spec\left(H_{0}\right),Spec\left(H_{(K+1)}\right),Spec\left(H_{2(K+1)}\right),\cdots,Spec\left(H_{L(K+1)}\right)\right\}  & \subset & SH,\,\,\,\, SG\neq SH.
\end{eqnarray*}
}

\textit{(3) It can be proved that the $(K+1)$ eigenvalues of $G_{m}$
are }

\textit{
\begin{eqnarray*}
\lambda^{(l)}\left(G_{m}\right) & =\lambda^{(l)}\left(\phi_{m}\right)= & \lambda^{(1)}\left(\phi_{m}-(l-1)\frac{2\pi}{\left(K+1\right)}\right),\,\, l=1,2,\cdots,K+1,
\end{eqnarray*}
 where $\phi_{m}=m\frac{2\pi}{M},\,\, m=0,1,2,\cdots,M-1$. Classify
$SG$ as $(K+1)$th groups, 
\[
SG=\cup_{l=1}^{K+1}Group{}^{(l)}
\]
 with $Group^{(l)}=\left\{ \left.\lambda^{(l)}\left(G_{m}\right)\right|\phi_{m}=m\frac{2\pi}{M},\,\, m=0,1,2,\cdots,M-1\right\} $.
If $mod(M,K+1)=0$, then eigenvalues in each group are the same $Group^{(1)}$=$\cdots$=$Group^{(K+1)}$.}

\textit{(4) $SG$ can be written as 
\begin{eqnarray*}
SG & = & \begin{cases}
\left\{ \left.\lambda^{(1)}(G(\psi_{(K+1)m}))\right|m=0,1,\cdots,M-1\right\} ,\,\, & if\,\,\, mod(M,K+1)=0,\\
\left\{ \left.\lambda^{(1)}(G(\psi_{j}))\right|j=0,1,2\cdots,M(K+1)\right\}  & else,
\end{cases}
\end{eqnarray*}
which means that all eigenvalues comes from the same function. Here
}$\psi_{j}=\frac{2\pi j}{M(K+1)}$.

The proof of the Theorem 2.1 is given in Appendix B. 

A pure upwind flux $\overline{u}=u^{L}$ is used to compute the common
fluxes, where $u^{L}$ is the left value at flux points. An example
is given to explain the properties of the eigenvalues of local discrete
matrices and the unique spectral curve of a scheme in Appendix C.
Comparisons on spectrum of different high-order schemes are given
in Appendix D. 

Dispersion and dissipation characteristics for fifth-order schemes
are shown in Fig. \ref{fig:comparison-of-dispersion-fifth-order},
where real part and imaginary part of eigenvalues computed from local
matrix $G$ in (\ref{eq:Gfai}) with \textit{$\phi_{m}=m\frac{2\pi}{M},\,\, m=0,1,2,\cdots,M-1$}
and $M=40$. We can see that all eigenvalues come from the same function
and the eigenvalue curves can coincide with each other after a shift
of $\frac{2\pi}{K+1}$, which agrees with the property (3) in Theorem
\ref{thm:eigenvalues} under the case $M=40$, $K+1=5$ and $mod(M,K+1)=0$.
This translation phenomenon was also found by Moura in \citet{Moura2015}.
The spectral properties of third-order schemes are given in Appendix
D.

\begin{center}
\begin{figure}
\subfloat[C5NNW5, dispersion]{\includegraphics[width=0.48\textwidth]{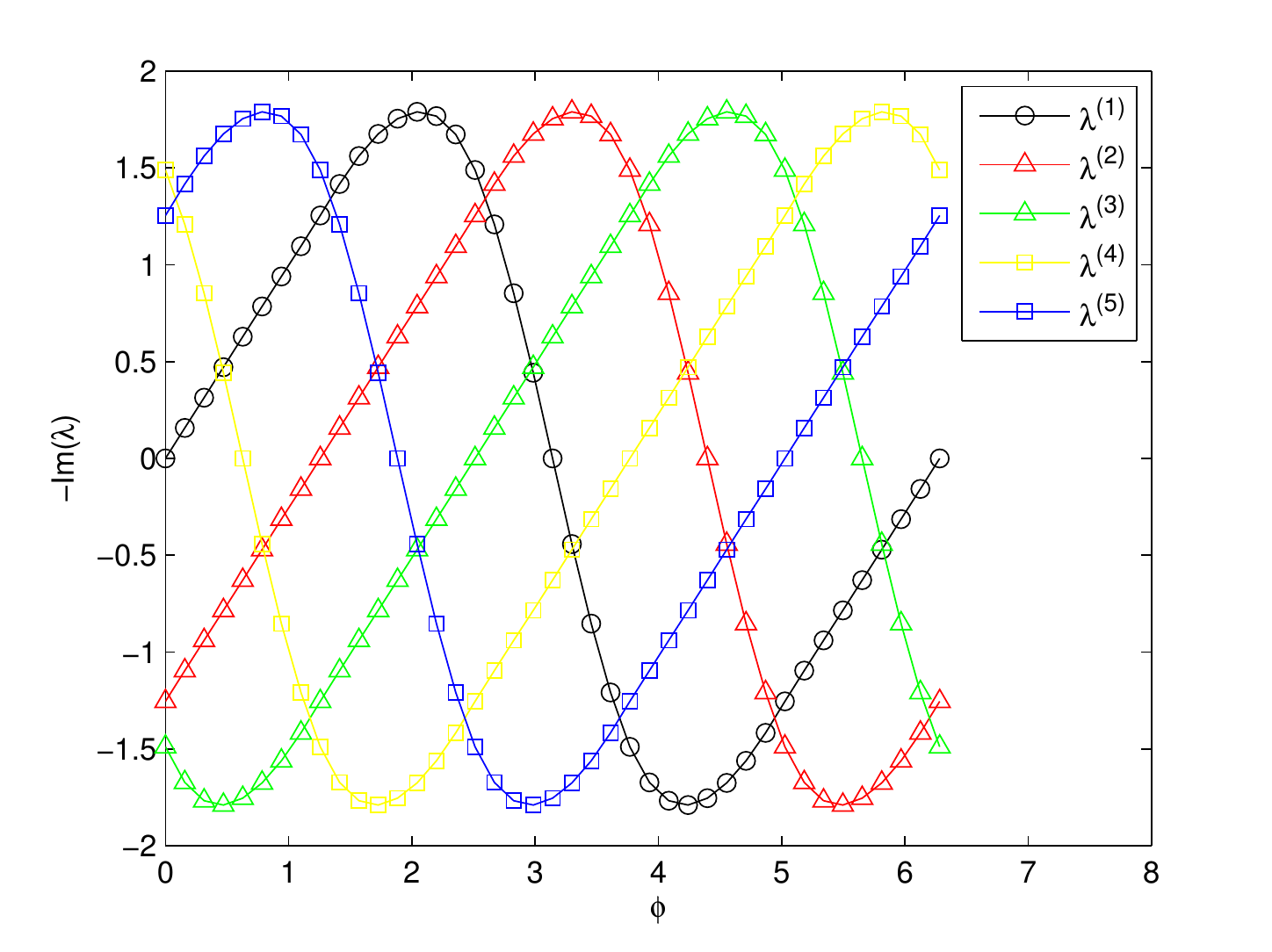}}\subfloat[C5NNW5, dissipation]{\includegraphics[width=0.48\textwidth]{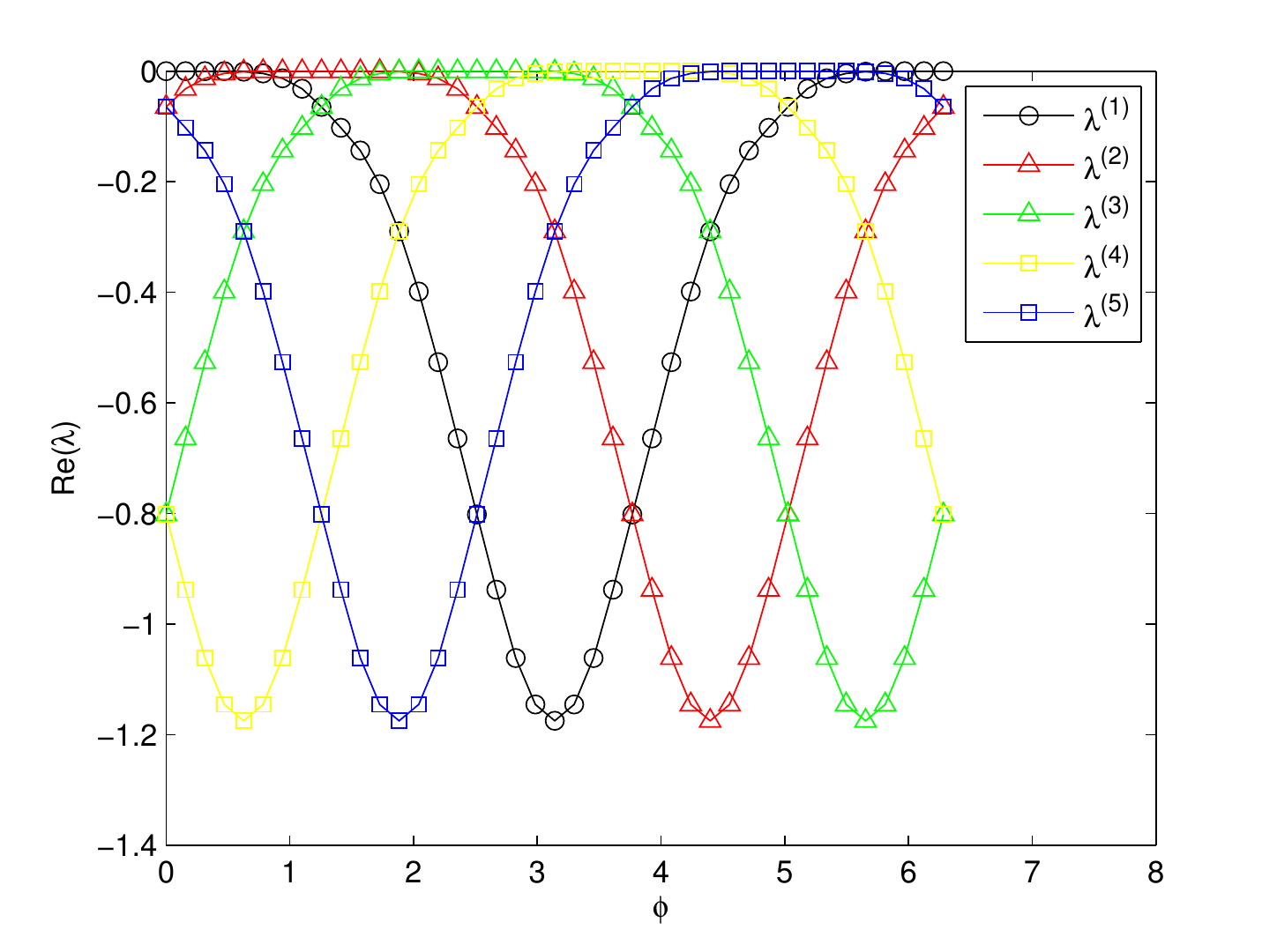}}

\subfloat[WCNS5, dispersion]{\includegraphics[width=0.48\textwidth]{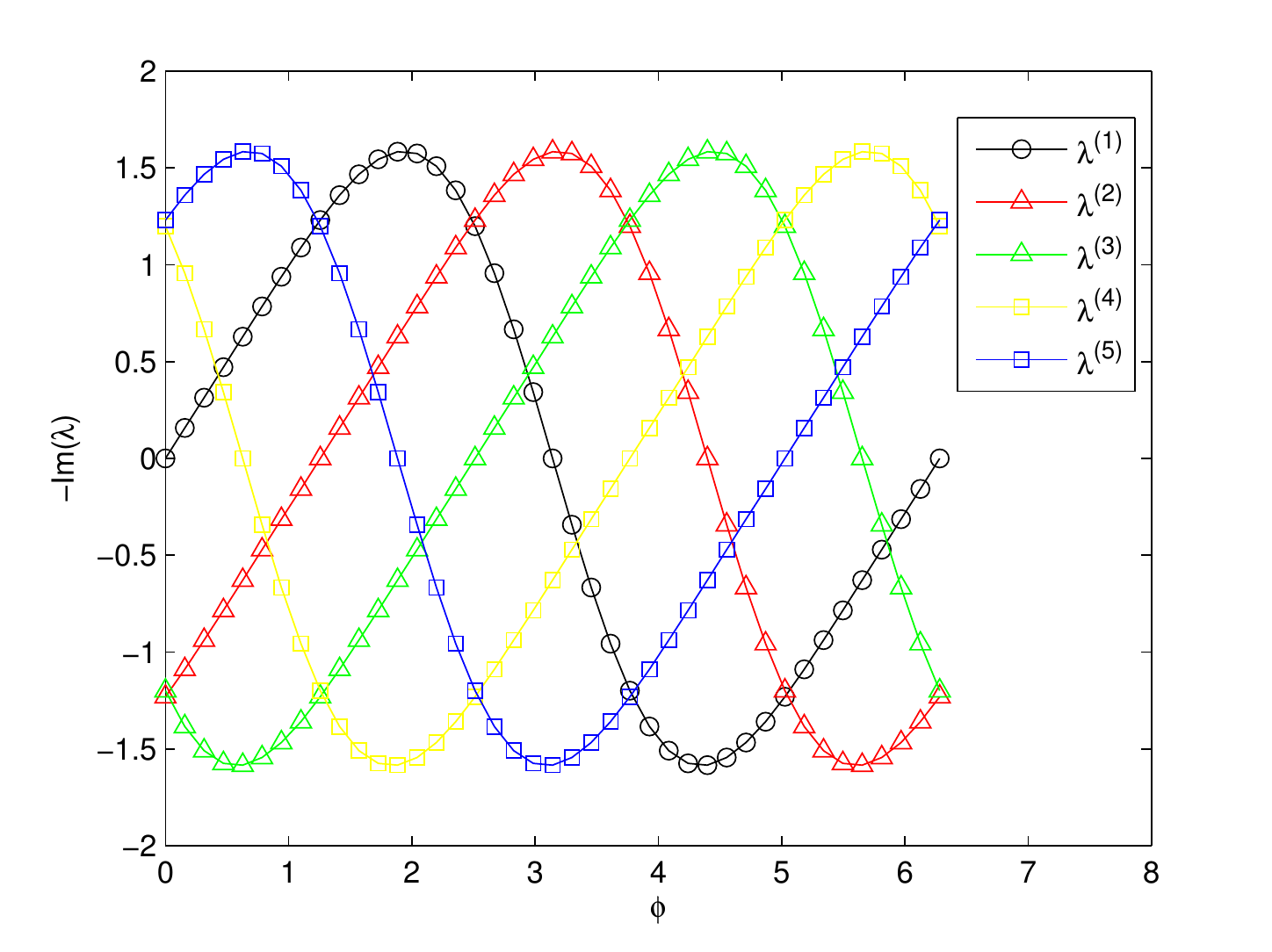}}\subfloat[WCNS5, dissipation]{\includegraphics[width=0.48\textwidth]{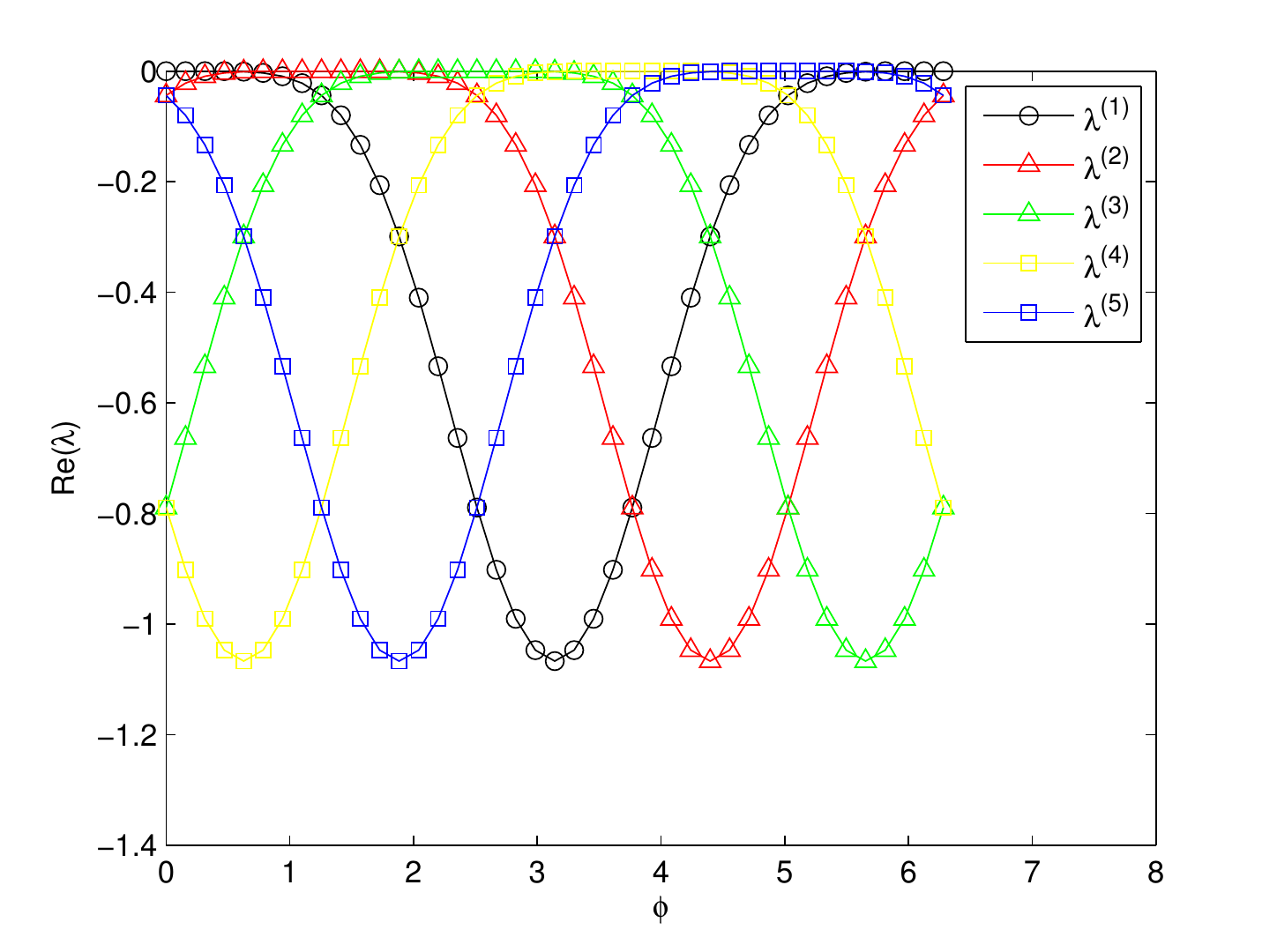}}

\subfloat[CPR5, dispersion]{\includegraphics[width=0.48\textwidth]{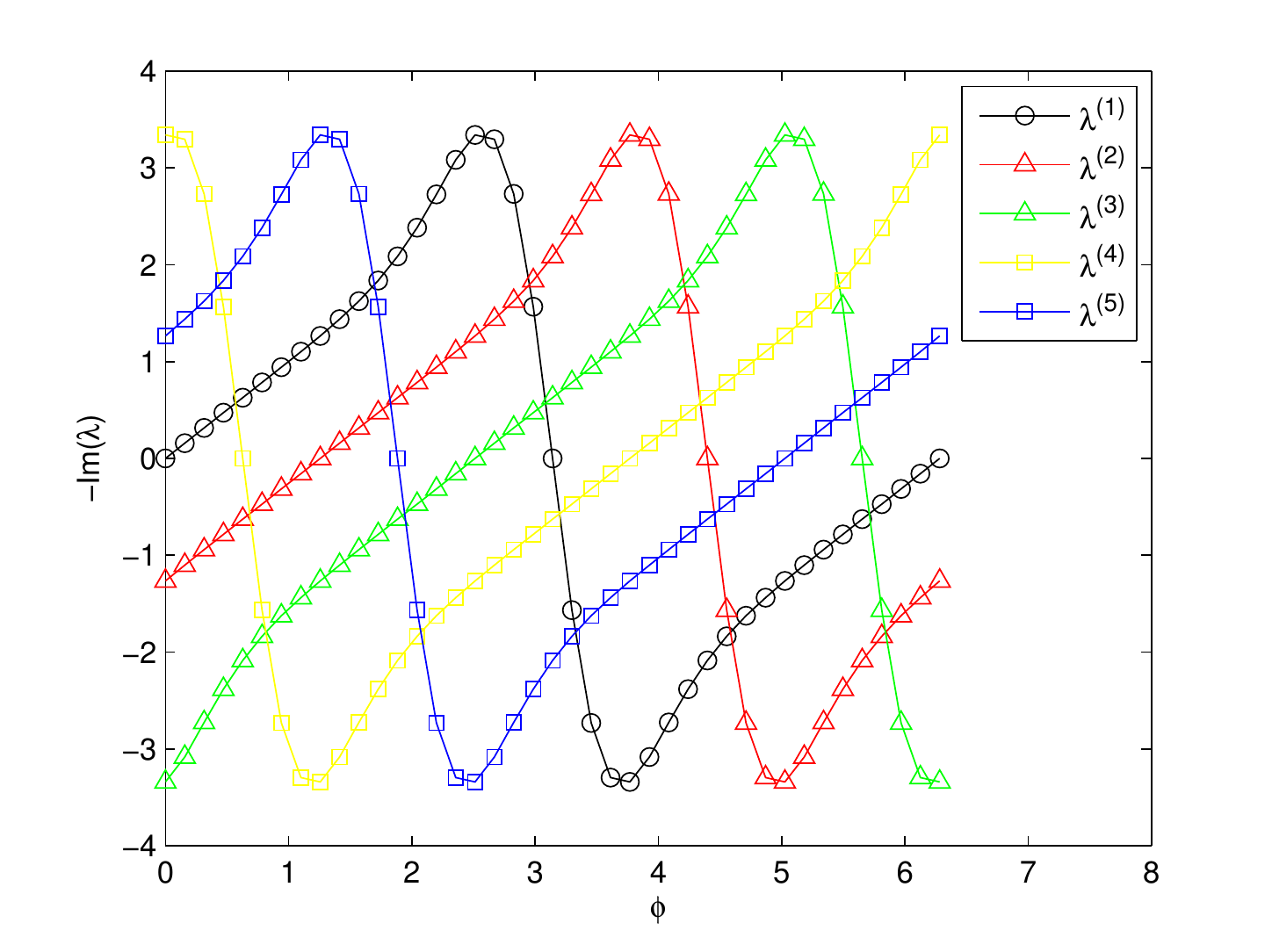}}\subfloat[CPR5, dissipation]{\includegraphics[width=0.48\textwidth]{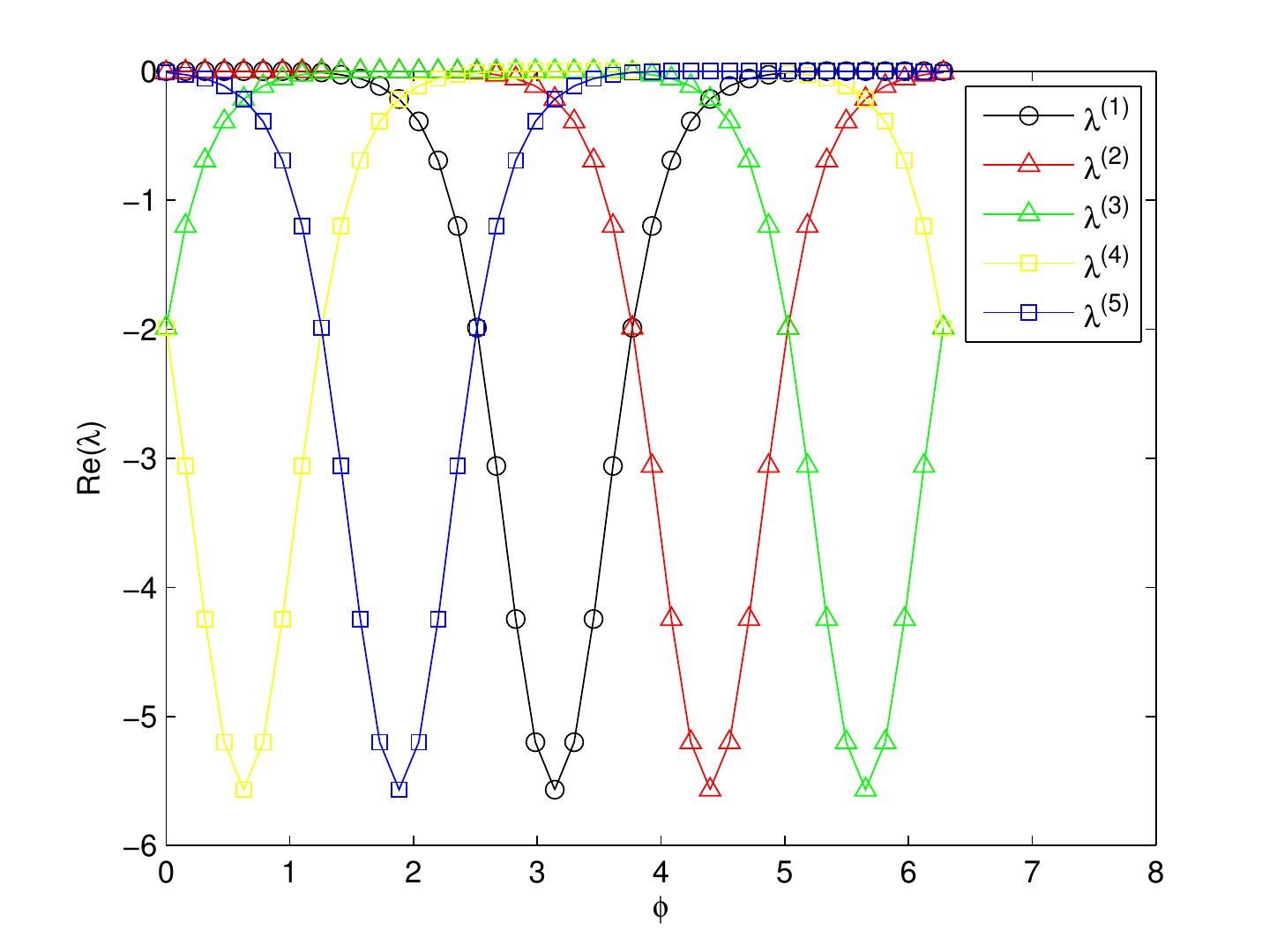}}

\caption{Comparison of dispersion (left) and dissipation (right) eigencurves
for fifth-order schemes, where $\left\{ \lambda_{m}^{(l)}\left|\lambda_{m}^{(l)}=\lambda^{(l)}(G(\phi_{m})),l=1,2,\cdots K+1;m=0,1,2,\cdots,M\right.\right\} $,
$\phi_{m}=m\frac{2\pi}{M}$, $K=4$ and $M=40$.\label{fig:comparison-of-dispersion-fifth-order}}
\end{figure}

\par\end{center}

Comparisons of different schemes show that the spectral property of
the proposed CNNW is closer to WCNS than CPR. Form Fig. \ref{fig:comparison }(a)(c)
we can see that C3NNW3 has smaller dispersion errors than WCNS3. For
fifth-order schemes, C5NNW5 has smaller dispersion errors than WCNS5
at area $\phi_{0}\leq\phi<\pi$ but bigger dispersion errors for $0\leq\phi<\phi_{0}$,
where $\phi_{0}\approx0.3772$. As for dissipation errors, CNNW has
similar errors as WCNS for both of the third-order schemes and fifth-order
schemes, which can been seen from Fig. \ref{fig:comparison }(b)(d). 

\begin{center}
\begin{figure}
\begin{centering}
\subfloat[Dispersion]{\begin{centering}
\includegraphics[width=0.48\textwidth]{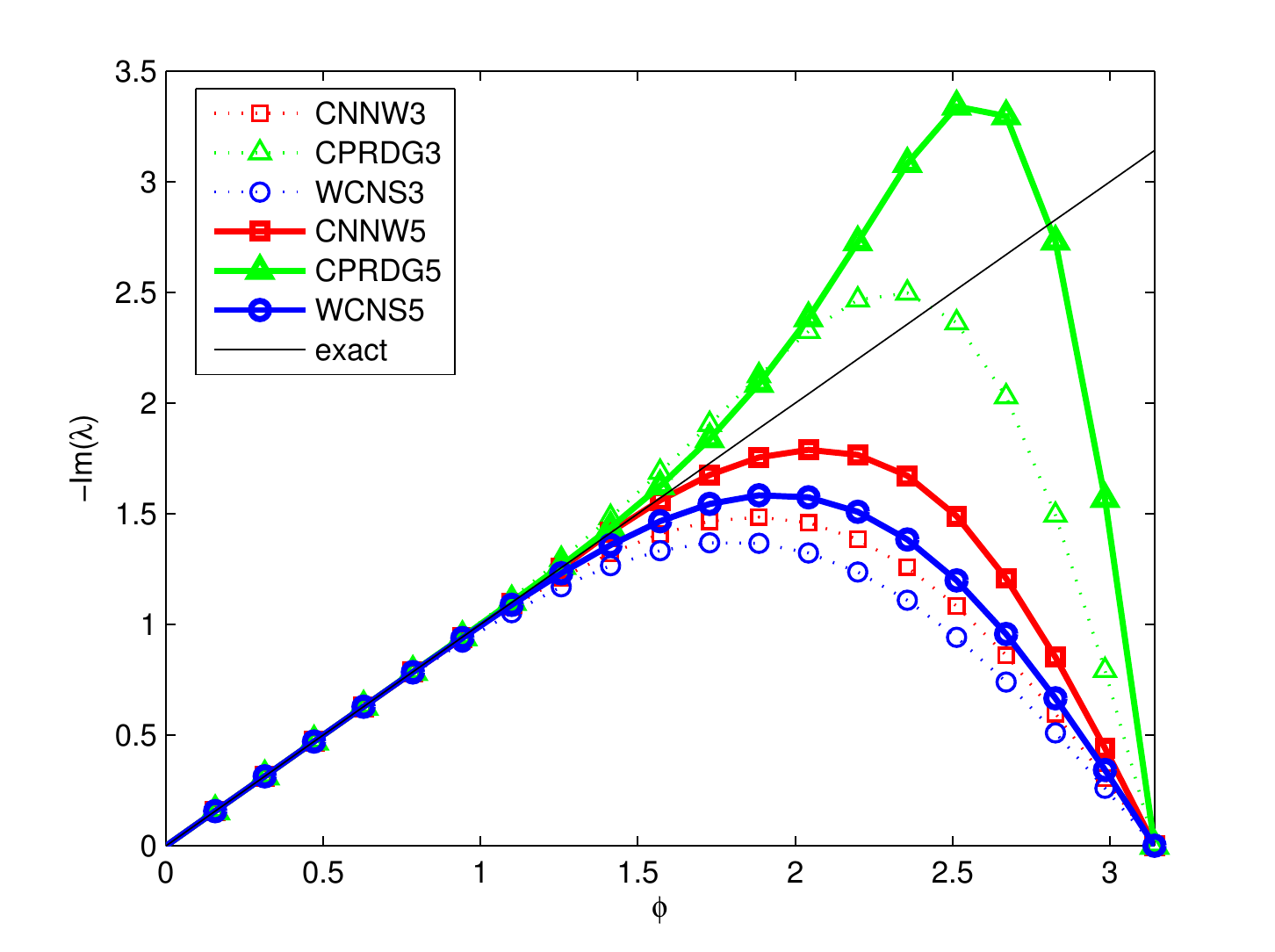}
\par\end{centering}

}\subfloat[Dissipation]{\begin{centering}
\includegraphics[width=0.48\textwidth]{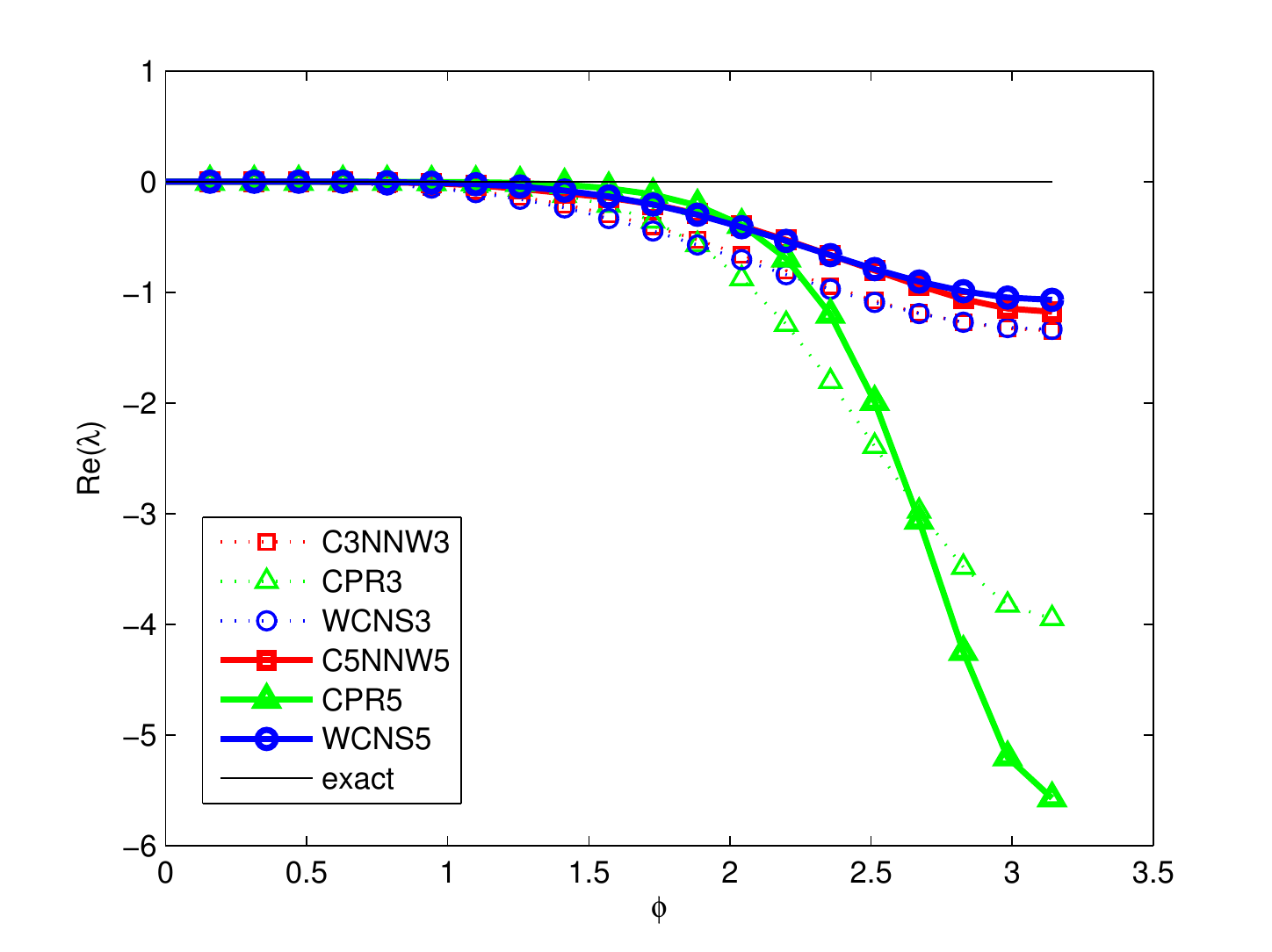} 
\par\end{centering}

}
\par\end{centering}

\begin{centering}
\subfloat[Dispersion errors]{\begin{centering}
\includegraphics[width=0.48\textwidth]{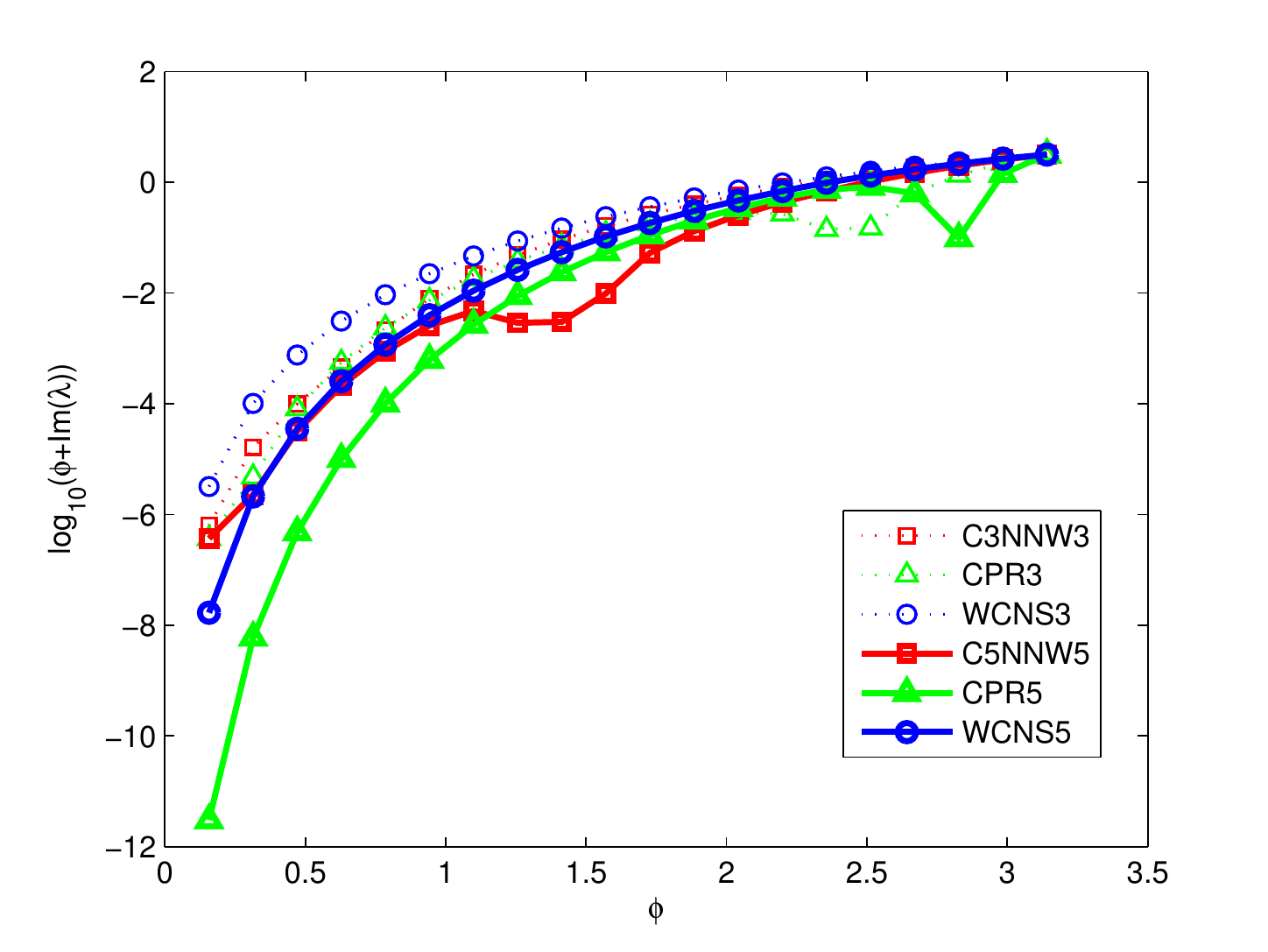}
\par\end{centering}

}\subfloat[Dissipation errors]{\begin{centering}
\includegraphics[width=0.48\textwidth]{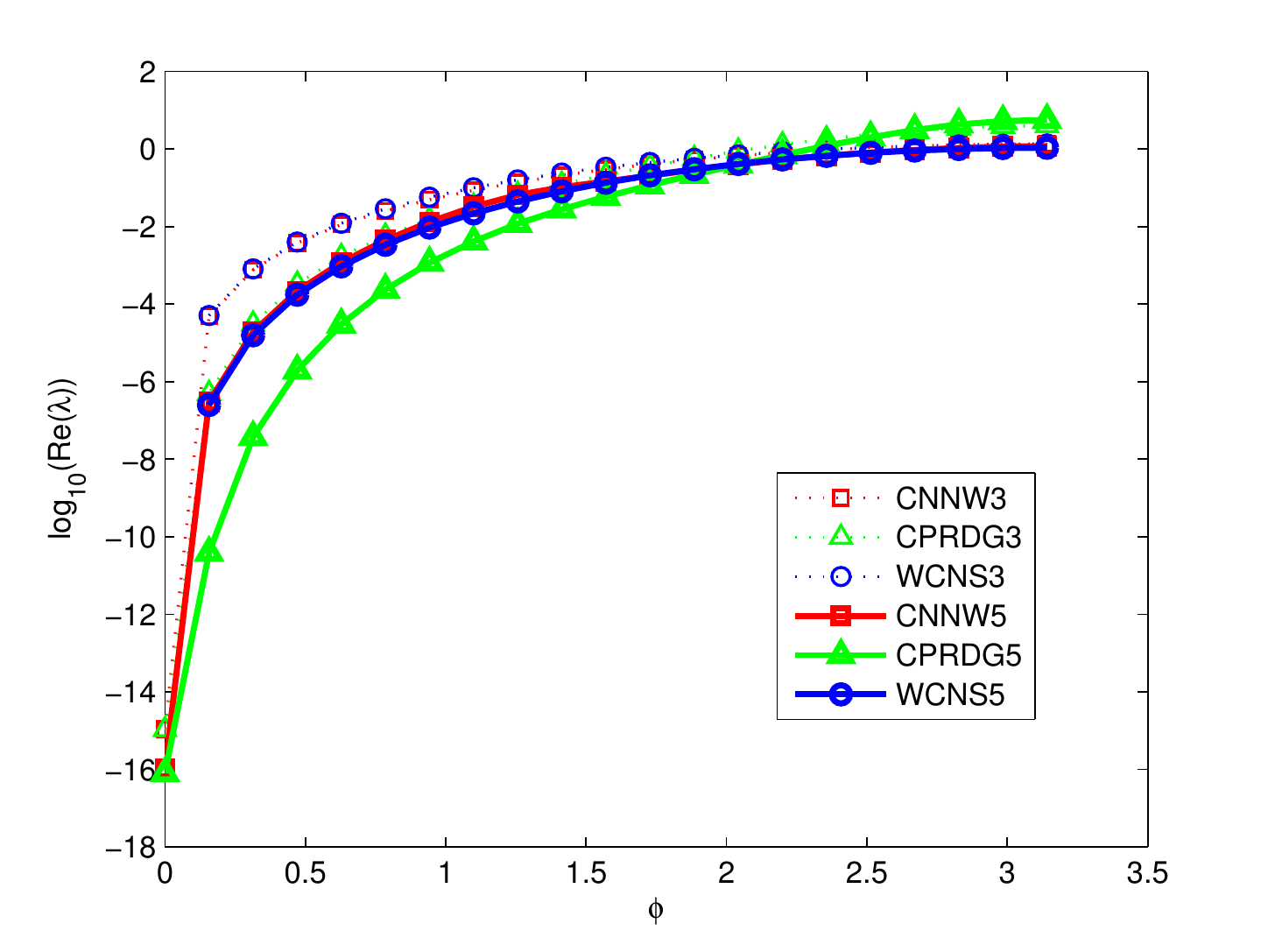} 
\par\end{centering}

}
\par\end{centering}

\caption{comparison of dispersion and dissipation for third-order schemes and
fifth-order schemes\label{fig:comparison }}
\end{figure}

\par\end{center}

\subsection{Discrete conservation laws}

\subsubsection{Discrete conservation laws of high-order CNNW}

In order to satisfy one-dimensional conservation law (CL), the following
integral conservation law in a cell $[x_{i},x_{i+1}]$ should be satisfied,
i.e.,

\begin{eqnarray*}
\int_{x_{i}}^{x_{i+1}}\frac{\partial u^{h}}{\partial t}dx+\left(\widetilde{F}(x_{i+1})-\widetilde{F}(x_{i})\right) & = & 0.
\end{eqnarray*}
For each solution point, CNNW with $\left(K+1\right)$th-order of
accuracy reads

\[
\frac{\partial u_{i,l}}{\partial t}=-\frac{2}{h}\frac{\partial\widetilde{F}}{\partial\xi}|_{i,l}=-\frac{2}{h}\sum_{j=1}^{K+2}a_{j}\overline{F}_{i,j},\quad l=1,2,\cdots,K+1.
\]
Thus, we can obtain

\begin{eqnarray}
\frac{\partial\left(\int_{x_{i}}^{x_{i+1}}u^{h}dx\right)}{\partial t} & = & \frac{h}{2}\frac{\partial\left(\int_{-1}^{1}u^{h}d\xi\right)}{\partial t}=\frac{h}{2}\frac{\partial\left(\sum_{l=1}^{K+1}2w_{l}u_{i,l}\right)}{\partial t}\label{eq:conservation-law}\\
 & = & \frac{h}{2}\sum_{l=1}^{K+1}2w_{l}\left(\frac{\partial u_{i,l}}{\partial t}\right)=\sum_{l=1}^{K+1}2w_{l}\left(-\frac{\partial\widetilde{F}}{\partial\xi}|_{i,l}\right).\nonumber 
\end{eqnarray}
where $w_{l}$ are the weights in Gaussian quadrature formulas. When
$K=4$, we have

\begin{eqnarray}
w_{1}=\frac{322-13\sqrt{70}}{1800},\quad w_{2}=\frac{322+13\sqrt{70}}{1800},\quad & w_{3}=\frac{64}{225},\quad & w_{4}=w_{2},\quad w_{5}=w_{1}.\label{eq:Gauss-weights}
\end{eqnarray}

Since flux polynomial is the Lagrange polynomial based on Riemann
fluxes at $(K+2)$th flux points, 
\[
\widetilde{F}(\xi)=\sum_{l=1}^{K+2}\overline{F}_{i,fp_{l}}L_{l}(\xi),
\]
$\widetilde{F}$ is a $K+1$ degree polynomial. Then, $\frac{\partial\widetilde{F}}{\partial\xi}$
belongs to $P^{K}$. Since the quadrature rule based on $K+1$ points
has at least $K$ algebraic accuracy, the rule is exact for degree
$K$ polynomial $\left(\frac{\partial\widetilde{F}}{\partial\xi}\right)$.
Thus, we have
\begin{eqnarray}
\sum_{l=1}^{K+1}2w_{l}\left(-\frac{\partial\widetilde{F}}{\partial\xi}|_{i,l}\right) & = & -\int_{-1}^{1}\left(\frac{\partial\widetilde{F}}{\partial\xi}\right)d\xi=-\left(\widetilde{F}(1)-\widetilde{F}(-1)\right)\label{eq:integral-exact}\\
 & = & -\left(\overline{F}_{i,fp_{K+2}}-\overline{F}_{i,fp_{1}}\right)=-\left(\overline{F}_{i}(1)-\overline{F}_{i}(-1)\right).\nonumber 
\end{eqnarray}
According to (\ref{eq:conservation-law}) and (\ref{eq:integral-exact}),
we obtain following relation

\begin{eqnarray*}
\frac{\partial\left(\sum_{l=1}^{K+1}w_{l}u_{i,l}\right)}{\partial t} & = & -\left(\overline{F}_{i}(1)-\overline{F}_{i}(-1)\right).
\end{eqnarray*}
Therefore, high-order CNNW satisfies discrete conservation laws.

\textbf{Remark 5.1.} Suppose flux derivatives at solutions points
are calculated by the Lagrange interpolation polynomial using all
flux points in the cell. If a cell has $K+1$ solution points and
flux points are less equal to $K+2$, then CNNW satisfies discrete
conservation law no matter which kind of solution points and flux
points are selected. When $K+1$ Legendre-Gauss solution points are
used, then CNNW satisfies discrete conservation law if flux points
are less equal to $2K+1$.

For two-dimensional case, we have

\begin{eqnarray*}
\frac{\partial u_{i,j,l,m}}{\partial t} & = & -\frac{4}{h^{2}}\left(\frac{\partial\widetilde{F}}{\partial\xi}|_{i,j,l,m}+\frac{\partial\widetilde{G}}{\partial\eta}|_{i,j,l,m}\right)\\
 & = & -\frac{4}{h^{2}}\left(\sum_{k=1}^{K+2}c_{l,k}\overline{F}_{i,j,fp_{k},m}+\sum_{k=1}^{K+2}c_{m,k}\overline{G}_{i,j,l,fp_{k}}\right).
\end{eqnarray*}
Then, it can be easily proved that the discrete conservation law is
\begin{eqnarray}
 &  & \frac{\partial\sum_{l=1}^{K+1}\sum_{m=1}^{K+1}w_{l}w_{m}u_{i,j,l,m}}{\partial t}\label{eq:CNNW5-conservation}\\
 & = & -\frac{4}{h^{2}}\left[\sum_{m=1}^{K+1}w_{m}\left(\sum_{l=1}^{K+1}w_{l}\frac{\partial\widetilde{F}}{\partial\xi}|_{i,j,l,m}\right)+\sum_{l=1}^{K+1}w_{l}\left(\sum_{m=1}^{K+1}w_{m}\frac{\partial\widetilde{G}}{\partial\eta}|_{i,j,l,m}\right)\right]\nonumber \\
 & = & -\frac{4}{h^{2}}\left[\sum_{m=1}^{K+1}w_{m}\left(\overline{F}(1,\eta_{m})-\overline{F}(-1,\eta_{m})\right)+\sum_{l=1}^{K+1}w_{l}\left(\overline{G}(\xi_{l},1)-\overline{G}(\xi_{l},-1)\right)\right].\nonumber 
\end{eqnarray}

\subsubsection{Discrete conservation law of C2NNW5 and C2NNW2}

For each solution point, we have

\[
\frac{\partial u_{i,sp_{l}}}{\partial t}=-\frac{2}{h}\frac{\partial\widetilde{F}}{\partial\xi}|_{i,sp_{l}}=-\frac{2}{h}\frac{\overline{F}_{i,fp_{(l+1)}}-\overline{F}_{i,fp_{l}}}{\Delta\xi_{l}}.
\]
Then, we can obtain

\begin{eqnarray*}
\frac{\partial\left(\int_{x_{i}}^{x_{i+1}}u^{h}dx\right)}{\partial t} & = & \frac{h}{2}\frac{\partial\left(\int_{-1}^{1}u^{h}d\xi\right)}{\partial t}=\frac{h}{2}\frac{\partial\left(\sum_{m=1}^{K+1}\int_{\xi_{i,fp_{l}}}^{\xi_{i,fp_{(l+1)}}}u^{h}d\xi\right)}{\partial t}\\
 & = & \frac{h}{2}\frac{\partial\left(\sum_{l=1}^{K+1}u_{i,sp_{l}}\Delta\xi_{l}\right)}{\partial t}=\frac{h}{2}\left(\sum_{l=1}^{K+1}\Delta\xi_{l}\frac{\partial u_{i,sp_{l}}}{\partial t}\right)\\
 & = & \frac{h}{2}\left(\sum_{l=1}^{K+1}\Delta\xi_{l}\left(-\frac{2}{h}\frac{\overline{F}_{i,fp_{(l+1)}}-\overline{F}_{i,fp_{l}}}{\Delta\xi_{l}}\right)\right)=\sum_{l=1}^{K+1}-\left(\overline{F}_{i,fp_{(l+1)}}-\overline{F}_{i,fp_{l}}\right)\\
 & = & -\left(\overline{F}_{i,fp_{(K+2)}}-\overline{F}_{i,fp_{1}}\right)=-\left(\overline{F}_{i}(1)-\overline{F}_{i}(-1)\right).
\end{eqnarray*}
Therefore, C2NNW5 and C2NNW2 satisfy following discrete conservation
law

\[
\frac{\partial\left(\sum_{l=1}^{K+1}u_{i,sp_{l}}\Delta\xi_{l}\right)}{\partial t}=-\left(\overline{F}_{i}(1)-\overline{F}_{i}(-1)\right).
\]
For two-dimensional case, it can be easily proved that the discrete
conservation law is
\begin{eqnarray}
 &  & \frac{\partial\sum_{l=1}^{K+1}\sum_{m=1}^{K+1}\Delta\xi_{l}\Delta\eta_{m}u_{i,j,l,m}}{\partial t}\nonumber \\
 & = & -\frac{4}{h^{2}}\left[\sum_{m=1}^{K+1}\Delta\eta_{m}\left(\sum_{l=1}^{K+1}\Delta\xi_{l}\frac{\partial\widetilde{F}}{\partial\xi}|_{i,j,l,m}\right)+\sum_{l=1}^{K+1}\Delta\xi_{l}\left(\sum_{m=1}^{K+1}\Delta\eta_{m}\frac{\partial\widetilde{G}}{\partial\eta}|_{i,j,l,m}\right)\right]\nonumber \\
 & = & -\frac{4}{h^{2}}\left[\sum_{m=1}^{K+1}\Delta\eta_{m}\left(\overline{F}(1,\eta_{m})-\overline{F}(-1,\eta_{m})\right)+\sum_{l=1}^{K+1}\Delta\xi_{l}\left(\overline{G}(\xi_{l},1)-\overline{G}(\xi_{l},-1)\right)\right].\label{eq:C2-conservation}
\end{eqnarray}
Therefore, the discrete conservation law holds for both C2NNW5 and
C2NNW2.

\subsubsection{Discrete conservation law of the CPR-CNNW scheme}

CPR and C5NNW5 has discrete conservation law in (\ref{eq:CNNW5-conservation}),
while the discrete conservation law for C2NNW5 and C2NNW2 has the
form in (\ref{eq:C2-conservation}). To make the CPR-CNNW scheme satisfying
discrete conservation law, it is necessary to have the same form of
discrete conservation laws for different schemes. Thus, we choose
flux points $\xi_{fp_{m}}=-1+\sum_{l=1}^{m-1}w_{l}$ for C2NNW5 and
C2NNW2. Here $w_{l}$ are the Legendre-Gauss weights in (\ref{eq:Gauss-weights}).
Then $\Delta\xi_{l}=w_{l}$ and $\Delta\eta_{m}=w_{m}$, the discrete
conservation law in (\ref{eq:C2-conservation}) of C2NNW5 and C2NNW2
becomes the same as the form (\ref{eq:CNNW5-conservation}). Then,
the C2NNW5 and C2NNW2 have the same form of discrete conservation
law with that of C5NNW5, which means that the Riemann fluxes at cell
interfaces will be eliminated during total summation. Therefore, the
hybrid schemes also satisfy discrete conservation law. The splitting
of a CPR element into subcells based on Gauss weights is similar as
that done by Sonntag et al. \citet{Sonntag2017} for DG method with
Gauss solution points and by Hennemann et al. \citet{Hennemann2021}
for DG method with Legendre-Gauss-Lobatto solution points.

\section{Numerical investigation}

In this section, seven test cases are adopted for both one-dimensional
and two-dimensional hyperbolic conservation laws to illustrate the
following properties of the CPR-CNNW schemes: (1) high resolution;
(2) good shock capturing robustness; (3) satisfying discrete conservation.
The test cases adopted and the properties studied in these test cases
have been summarized in Tab. \ref{tab:Summary-of-test}.

\begin{table}
\begin{centering}
\emph{}%
\begin{tabular}{|l|l|l|l|}
\hline 
\textit{Sections} & \textit{Equations} & \textit{Test cases} & \textit{Purpose}\tabularnewline
\hline 
\hline 
{\footnotesize \ref{sub:1D-linear-advection} } & {\footnotesize 1D advection} & {\footnotesize Advection of smooth wave} & {\footnotesize Order of accuracy}\tabularnewline
\hline 
{\footnotesize \ref{sub:1D-Shock-tube}} & {\footnotesize 1D Euler} & {\footnotesize Shock tube problems} & {\footnotesize Shock capturing (Numerical oscillations)}\tabularnewline
\hline 
\multirow{3}{*}{{\footnotesize \ref{sub:Shu-Osher-problem}}} & \multirow{3}{*}{{\footnotesize 1D Euler}} & \multirow{3}{*}{{\footnotesize Shu-Osher problem}} & {\footnotesize Resolution, Scheme adaption strategy}\tabularnewline
\cline{4-4} 
 &  &  & {\footnotesize Influence of indicating on resolution}\tabularnewline
 &  &  & {\footnotesize Shock capturing (Comparison with others)}\tabularnewline
\hline 
\multirow{2}{*}{{\footnotesize \ref{sub:2D-Euler-vortex}}} & \multirow{2}{*}{{\footnotesize 2D Euler}} & \multirow{2}{*}{{\footnotesize Isentropic vortex problem}} & \multirow{1}{*}{{\footnotesize Order of accuracy, influence of indicating on accuracy}}\tabularnewline
 &  &  & {\footnotesize Time evolution of errors, Discrete conservation,}\tabularnewline
\hline 
\multirow{2}{*}{{\footnotesize \ref{sub:2D-Riemann-problem}}} & \multirow{2}{*}{{\footnotesize 2D Euler}} & \multirow{2}{*}{{\footnotesize Riemann problem}} & \multicolumn{1}{l|}{{\footnotesize Resolution, Influence of indicating on resolution}}\tabularnewline
 &  &  & {\footnotesize Shock capturing (Quantitative study)}\tabularnewline
\hline 
{\footnotesize \ref{sub:Double-Mach-reflection}} & {\footnotesize 2D Euler} & {\footnotesize Double Mach reflection} & {\footnotesize Resolution, Shock capturing (Strong shocks)}\tabularnewline
\hline 
{\footnotesize \ref{sub:Shock-vortex-interaction}} & {\footnotesize 2D Euler} & {\footnotesize Shock vortex interaction} & {\footnotesize Resolution, Shock capturing}\tabularnewline
\hline 
\end{tabular}
\par\end{centering}

\caption{Summary of test cases\label{tab:Summary-of-test}}
\end{table}

Unless otherwise specified, $a=0.5$ is used and then $c_{0}=c(0.5)=0.5\cdot10^{-1.8(6+1)^{1/4}}\approx5.90\times10^{-4}$
is applied to control the size of CNNW area. Partition vector $\mathbf{dv}=(c_{0},0.05,0.1)$
and $\mathbf{dv}=(c_{0},c_{0},c_{0})$ in (\ref{eq:HS}) will be applied
to obtain HCCS(1,1,1,1) and HCCS(1,0,0,1) respectively. The Z-weights
in Formula (\ref{eq:Zweight}) in Appendix A with $\varepsilon=10^{-10}$
is used for computing nonlinear weights in the NNW procedure. By default,
in shock capturing test cases, characteristic variables are applied
in the NNW interpolations. The Lax-Friedrichs flux is adopted to compute
common fluxes. The explicit third-order TVD Runge-Kutta scheme is
used for time integration.

\subsection{1D linear advection equation with a smooth solution\label{sub:1D-linear-advection}}

An accuracy test is conducted based on the one-dimensional (1D) linear
advection equation

\begin{eqnarray*}
u_{t}+u_{x} & = & 0.
\end{eqnarray*}
The problem is solved in a spatial domain of $[-3,3]$ with initial
condition $u(x,0)=sin\left(\frac{\pi x}{3}\right)$ till time $T=3$.
The $L_{2}$ and $L_{\infty}$ errors, as well as the numerical order
of accuracy, are summarized in Tabs. \ref{tab:high-order-linear}
and \ref{tab:high-order-nonlinear} for the fifth-order schemes and
in Tab. \ref{tab:second-order-schemes} for the second-order schemes.
As expected, C5NNW5 reaches fifth-order of accuracy. In addition,
C5NNW5 has smaller numerical errors than the fifth-order CPR-g2 scheme
with Legendre-Gauss-Lobatto solution points and $g_{2}$ correction
function (see references \citet{Huynh2007,zhu2017} for details),
while C5NNW5 has bigger errors than WCNS5 and CPR5. For the second-order
schemes, C2NNW5 has smaller errors than C2NNW2 for both linear schemes
and nonlinear schemes.

\begin{center}
\begin{table}
\begin{centering}
\begin{tabular}{cccccccccc}
\hline 
\multirow{2}{*}{{\footnotesize $Norm$}} & \multirow{2}{*}{{\footnotesize DoFs}} & \multicolumn{2}{c}{{\footnotesize C5NNW5}} & \multicolumn{2}{c}{{\footnotesize WCNS5}} & \multicolumn{2}{c}{{\footnotesize CPR-g2}} & \multicolumn{2}{c}{{\footnotesize CPR5}}\tabularnewline
\cline{3-10} 
 &  & \textit{\footnotesize error} & \textit{\footnotesize order} & \textit{\footnotesize error} & \textit{\footnotesize order} & \textit{\footnotesize error} & \textit{\footnotesize order} & \textit{\footnotesize error} & \textit{\footnotesize order}\tabularnewline
\multirow{5}{*}{{\footnotesize $L_{\infty}$}} & {\footnotesize $15$} & {\footnotesize 9.04E-04} & {\footnotesize -} & {\footnotesize 6.55E-04} & {\footnotesize -} & {\footnotesize 5.28E-03} & {\footnotesize -} & \textbf{\textcolor{black}{\footnotesize 5.72E-04 }} & \textcolor{black}{\footnotesize -}\tabularnewline
 & {\footnotesize $30$} & {\footnotesize 2.74E-05} & {\footnotesize 5.04} & {\footnotesize 2.09E-05} & {\footnotesize 4.97} & {\footnotesize 1.77E-04} & {\footnotesize 4.90} & \textbf{\textcolor{black}{\footnotesize 1.28E-05 }} & \textcolor{black}{\footnotesize 5.48}\tabularnewline
 & {\footnotesize $60$} & {\footnotesize 8.10E-07} & {\footnotesize 5.08} & {\footnotesize 6.58E-07} & {\footnotesize 4.99} & {\footnotesize 5.93E-06} & {\footnotesize 4.90} & \textbf{\textcolor{black}{\footnotesize 4.38E-07 }} & \textcolor{black}{\footnotesize 4.87}\tabularnewline
 & {\footnotesize $120$} & {\footnotesize 2.54E-08} & {\footnotesize 5.00} & {\footnotesize 2.06E-08} & {\footnotesize 5.00} & {\footnotesize 1.84E-07} & {\footnotesize 5.01} & \textbf{\textcolor{black}{\footnotesize 1.41E-08 }} & \textcolor{black}{\footnotesize 4.95}\tabularnewline
 & \uline{\footnotesize $240$} & {\footnotesize 7.97E-10} & {\footnotesize 4.99} & {\footnotesize 6.44E-10} & {\footnotesize 5.00} & {\footnotesize 5.72E-09} & {\footnotesize 5.01} & \textbf{\textcolor{black}{\footnotesize 4.55E-10 }} & \textcolor{black}{\footnotesize 4.96}\tabularnewline
\hline 
\multirow{5}{*}{{\footnotesize $L_{2}$}} & {\footnotesize $15$} & {\footnotesize 5.90E-04} & {\footnotesize -} & {\footnotesize 4.64E-04} & {\footnotesize -} & {\footnotesize 1.94E-03} & {\footnotesize -} & \textbf{\textcolor{black}{\footnotesize 3.18E-04 }} & \textcolor{black}{\footnotesize -}\tabularnewline
 & {\footnotesize $30$} & {\footnotesize 1.80E-05} & {\footnotesize 5.03} & {\footnotesize 1.48E-05} & {\footnotesize 4.97} & {\footnotesize 6.43E-05} & {\footnotesize 4.91} & \textbf{\textcolor{black}{\footnotesize 9.48E-06 }} & \textcolor{black}{\footnotesize 5.07}\tabularnewline
 & {\footnotesize $60$} & {\footnotesize 5.63E-07} & {\footnotesize 5.00} & {\footnotesize 4.65E-07} & {\footnotesize 4.99} & {\footnotesize 2.03E-06} & {\footnotesize 4.99} & \textbf{\textcolor{black}{\footnotesize 2.98E-07 }} & \textcolor{black}{\footnotesize 4.99}\tabularnewline
 & {\footnotesize $120$} & {\footnotesize 1.77E-08} & {\footnotesize 4.99} & {\footnotesize 1.46E-08} & {\footnotesize 4.99} & {\footnotesize 6.29E-08} & {\footnotesize 5.01} & \textbf{\textcolor{black}{\footnotesize 9.23E-09 }} & \textcolor{black}{\footnotesize 5.01}\tabularnewline
 & \uline{\footnotesize $240$} & {\footnotesize 5.54E-10} & {\footnotesize 5.00} & {\footnotesize 4.55E-10} & {\footnotesize 5.00} & {\footnotesize 1.95E-09} & {\footnotesize 5.01} & \textbf{\textcolor{black}{\footnotesize 2.85E-10 }} & \textcolor{black}{\footnotesize 5.02}\tabularnewline
\hline 
\end{tabular}
\par\end{centering}

\noindent \centering{}\caption{Accuracy test results of the high-order linear schemes in solving
1D linear equation. Here ``linear'' indicates that the optimal linear
weights are adopted for the NNW interpolation procedure. \label{tab:high-order-linear}}
\end{table}

\par\end{center}

\begin{center}
\begin{table}
\begin{centering}
\begin{tabular}{cccccc}
\hline 
\multirow{2}{*}{{\footnotesize $Norm$}} & \multirow{2}{*}{{\footnotesize DoFs}} & \multicolumn{2}{c}{{\footnotesize C5NNW5}} & \multicolumn{2}{c}{{\footnotesize WCNS5}}\tabularnewline
\cline{3-6} 
 &  & \textit{\footnotesize error} & \textit{\footnotesize order} & \textit{\footnotesize error} & \textit{\footnotesize order}\tabularnewline
\multirow{5}{*}{{\footnotesize $L_{\infty}$}} & {\footnotesize $15$} & {\footnotesize 3.94E-02} & {\footnotesize -} & \textbf{\textcolor{black}{\footnotesize 5.67E-03}} & \textbf{\textcolor{black}{\footnotesize -}}\tabularnewline
 & {\footnotesize $30$} & {\footnotesize 4.68E-04} & {\footnotesize 6.40} & \textbf{\textcolor{black}{\footnotesize 2.05E-04}} & \textcolor{black}{\footnotesize 4.79}\tabularnewline
 & {\footnotesize $60$} & {\footnotesize 9.53E-06} & {\footnotesize 5.62} & \textbf{\textcolor{black}{\footnotesize 6.59E-06}} & \textcolor{black}{\footnotesize 4.96}\tabularnewline
 & {\footnotesize $120$} & {\footnotesize 2.39E-07} & {\footnotesize 5.32} & \textbf{\textcolor{black}{\footnotesize 2.04E-07}} & \textcolor{black}{\footnotesize 5.01}\tabularnewline
 & \uline{\footnotesize $240$} & {\footnotesize 7.58E-09} & {\footnotesize 4.98} & \textbf{\textcolor{black}{\footnotesize 6.03E-09}} & \textcolor{black}{\footnotesize 5.08}\tabularnewline
\hline 
\multirow{5}{*}{{\footnotesize $L_{2}$}} & {\footnotesize $15$} & {\footnotesize 1.98E-02} & {\footnotesize -} & \textbf{\textcolor{black}{\footnotesize 3.73E-03}} & \textcolor{black}{\footnotesize -}\tabularnewline
 & {\footnotesize $30$} & {\footnotesize 2.83E-04} & {\footnotesize 6.13} & \textbf{\textcolor{black}{\footnotesize 1.23E-04}} & \textcolor{black}{\footnotesize 4.92}\tabularnewline
 & {\footnotesize $60$} & {\footnotesize 5.16E-06} & {\footnotesize 5.78} & \textbf{\textcolor{black}{\footnotesize 3.74E-06}} & \textcolor{black}{\footnotesize 5.04}\tabularnewline
 & {\footnotesize $120$} & {\footnotesize 1.52E-07} & {\footnotesize 5.09} & \textbf{\textcolor{black}{\footnotesize 1.14E-07}} & \textcolor{black}{\footnotesize 5.04}\tabularnewline
 & \uline{\footnotesize $240$} & {\footnotesize 4.76E-09} & {\footnotesize 5.00} & \textbf{\textcolor{black}{\footnotesize 3.53E-09}} & \textcolor{black}{\footnotesize 5.01}\tabularnewline
\hline 
\end{tabular}
\par\end{centering}

\noindent \centering{}\caption{Accuracy test results of the high-order nonlinear schemes in solving
1D linear equation.\label{tab:high-order-nonlinear}}
\end{table}

\par\end{center}

\begin{center}
\begin{table}
\begin{centering}
\begin{tabular}{cccccccccc}
\cline{3-10} 
\multirow{3}{*}{Norm} & \multirow{3}{*}{N} & \multicolumn{4}{c}{linear schemes} & \multicolumn{4}{c}{nonlinear schemes}\tabularnewline
\cline{3-10} 
 &  & \multicolumn{2}{c}{{\footnotesize C2NNW5}} & \multicolumn{2}{c}{{\footnotesize C2NNW2}} & \multicolumn{2}{c}{{\footnotesize C2NNW5}} & \multicolumn{2}{c}{{\footnotesize C2NNW2}}\tabularnewline
\cline{3-10} 
 &  & \textit{\footnotesize error} & \textit{\footnotesize order} & \textit{\footnotesize error} & \textit{\footnotesize order} & \textit{\footnotesize error} & \textit{\footnotesize order} & \textit{\footnotesize error} & \textit{\footnotesize order}\tabularnewline
\multirow{5}{*}{{\footnotesize $L_{\infty}$}} & {\footnotesize $15$} & {\footnotesize 3.88E-02 } & {\footnotesize -} & {\footnotesize 7.15E-02} & {\footnotesize -} & {\footnotesize 4.25E-02 } & {\footnotesize -} & {\footnotesize 1.38E-01} & {\footnotesize -}\tabularnewline
 & {\footnotesize $30$} & {\footnotesize 8.61E-03 } & {\footnotesize 2.17} & {\footnotesize 1.61E-02 } & {\footnotesize 2.15} & {\footnotesize 9.29E-03} & {\footnotesize 2.19} & {\footnotesize 4.33E-02} & {\footnotesize 1.68}\tabularnewline
 & {\footnotesize $60$} & {\footnotesize 2.17E-03 } & {\footnotesize 1.99} & {\footnotesize 3.95E-03 } & {\footnotesize 2.02} & {\footnotesize 2.22E-03} & {\footnotesize 2.07} & {\footnotesize 1.49E-02} & {\footnotesize 1.53}\tabularnewline
 & {\footnotesize $120$} & {\footnotesize 5.35E-04} & {\footnotesize 2.02} & {\footnotesize 9.57E-04 } & {\footnotesize 2.05} & {\footnotesize 5.35E-04} & {\footnotesize 2.05} & {\footnotesize 6.17E-03} & {\footnotesize 1.28}\tabularnewline
 & \uline{\footnotesize $240$} & {\footnotesize 1.33E-04} & {\footnotesize 2.01} & {\footnotesize 2.38E-04 } & {\footnotesize 2.01} & {\footnotesize 1.33E-04} & {\footnotesize 2.01} & {\footnotesize 2.57E-03} & {\footnotesize 1.26}\tabularnewline
\hline 
\multirow{5}{*}{{\footnotesize $L_{2}$}} & {\footnotesize $15$} & {\footnotesize 2.55E-02} & {\footnotesize -} & {\footnotesize 4.73E-02 } & {\footnotesize -} & {\footnotesize 2.85E-02} & {\footnotesize -} & {\footnotesize 6.20E-02} & {\footnotesize -}\tabularnewline
 & {\footnotesize $30$} & {\footnotesize 6.00E-03} & {\footnotesize 2.09} & {\footnotesize 1.10E-02 } & {\footnotesize 2.10} & {\footnotesize 6.42E-03} & {\footnotesize 2.15} & {\footnotesize 1.97E-02} & {\footnotesize 1.65}\tabularnewline
 & {\footnotesize $60$} & {\footnotesize 1.47E-03} & {\footnotesize 2.03} & {\footnotesize 2.67E-03 } & {\footnotesize 2.04} & {\footnotesize 1.55E-03} & {\footnotesize 2.05} & {\footnotesize 6.30E-03} & {\footnotesize 1.65}\tabularnewline
 & {\footnotesize $120$} & {\footnotesize 3.66E-04} & {\footnotesize 2.01} & {\footnotesize 6.61E-04 } & {\footnotesize 2.01} & {\footnotesize 3.82E-04} & {\footnotesize 2.02} & {\footnotesize 1.85E-03} & {\footnotesize 1.77}\tabularnewline
 & \uline{\footnotesize $240$} & {\footnotesize 9.13E-05} & {\footnotesize 2.00} & {\footnotesize 1.65E-04 } & {\footnotesize 2.00} & {\footnotesize 9.31E-05} & {\footnotesize 2.04} & {\footnotesize 5.67E-04} & {\footnotesize 1.71}\tabularnewline
\hline 
\end{tabular}
\par\end{centering}

\noindent \centering{}\caption{Accuracy test results of the second-order schemes in solving 1D linear
equation.\label{tab:second-order-schemes}}
\end{table}

\par\end{center}

\subsection{1D Shock tube problems\label{sub:1D-Shock-tube}}

The shock tube problems are solved to test the shock-capturing capability
of the CNNW and CPR-CNNW. The Sod problem with initial conditions
\[
\left(\rho,u,p\right)=\begin{cases}
\left(1,0,1\right), & 0\leq x<0.5,\\
(0.125,0,0.1), & 0.5\leq x<1,
\end{cases}
\]
is solved till $t=0.2$ with $40$ cells and $DoFs=200$. In addition,
the Lax problem with initial conditions 
\[
\left(\rho,u,p\right)=\begin{cases}
\left(0.445,0.698,3.528\right), & 0\leq x<0.5,\\
(0.5,0,0.571), & 0.5\leq x<1,
\end{cases}
\]
is solved till $t=0.1$ with $100$ cells and $DoFs=500$. 

Firstly, the CNNW schemes are used to solve these problems. The computed
density distributions are shown in Fig. \ref{fig:single-schemes}.
In the Sod problem, there is no obvious numerical oscillation for
simulations with C5NNW5, C2NNW5 and C2NNW2. Moreover, the C5NNW5 acquires
a similar result as that of the WCNS5. For the second-order schemes,
C2NNW5 results in sharper discontinuities than the C2NNW2. For the
Lax problem with stronger discontinuities, the results of the C5NNW5,
C2NNW5 and C2NNW2 have small oscillations near the location $x=0.74$
while there is no obvious oscillations for the result of the WCNS. 

Secondly, CPR-CNNW schemes are applied to solve these problems. The
results are shown in Fig. \ref{fig:Hybrid-schemes-HS}. We can see
that both HCCS(1,1,1,1) and HCCS(1,0,0,1) can capture shock robustly.
In addition, there are only 1-2 troubled cells near each discontinuity.
Symbols are also given to highlight the values at all solution points.
For both HCCS(1,1,1,1) and HCCS(1,0,0,1), the shock thickness is close
to the element length and 6 solution points can nearly recover the
discontinuity, which indicates that the present limiting method has
a subcell resolution.

\begin{center}
\begin{figure}
\begin{centering}
\subfloat[Sod problem, $DOFs=200$ ]{\begin{centering}
\includegraphics[width=0.45\textwidth]{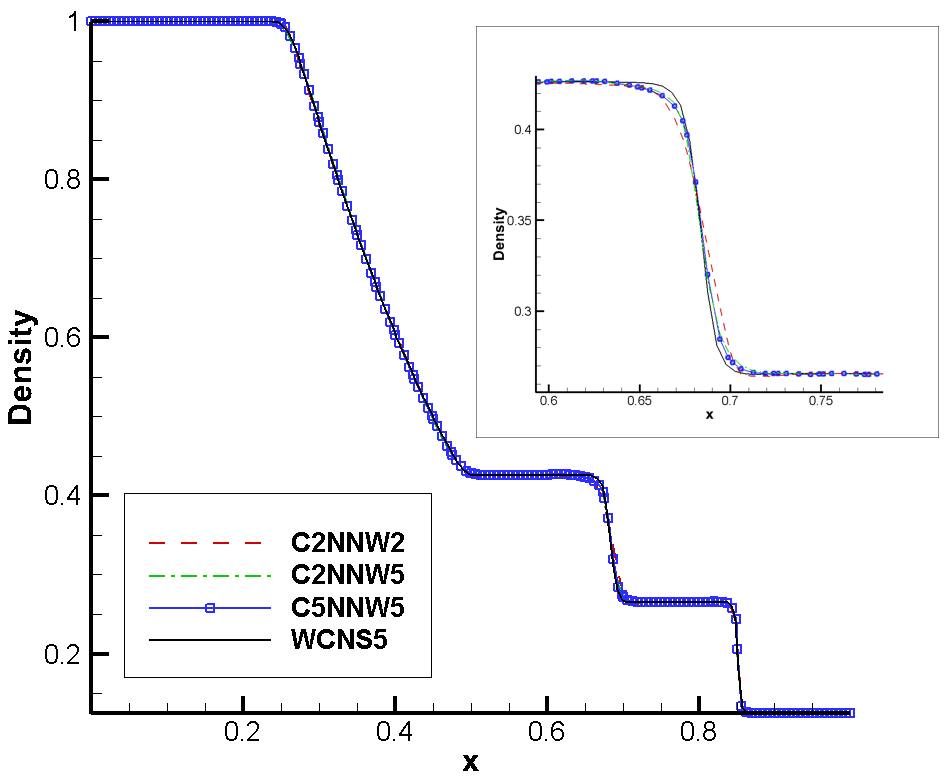}
\par\end{centering}

}\subfloat[Lax problem, $DOFs=500$]{\begin{centering}
\includegraphics[width=0.45\textwidth]{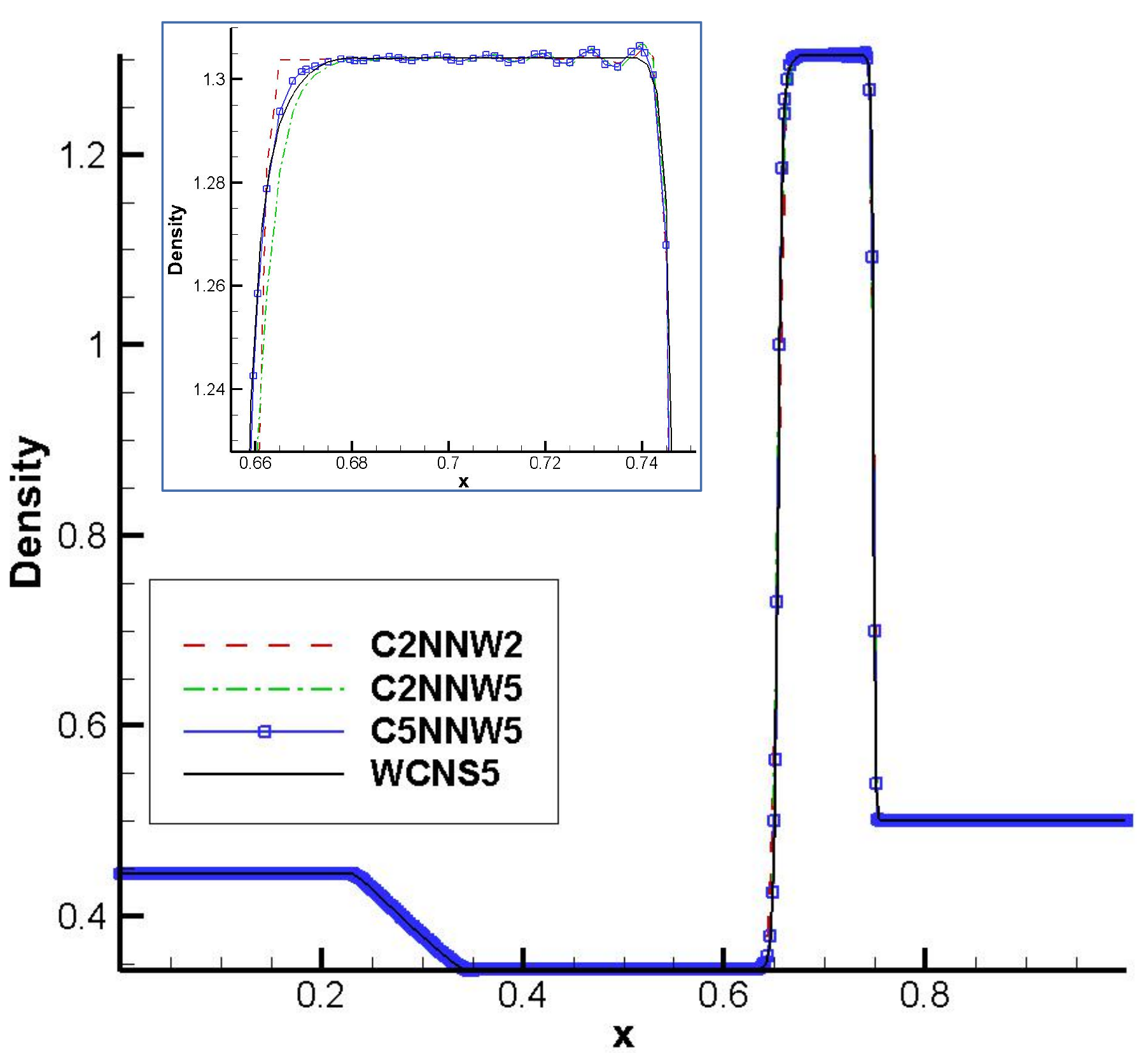}
\par\end{centering}

}
\par\end{centering}

\centering{}\caption{Density distribution of single schemes in solving 1D shock tube problems.\label{fig:single-schemes}}
\end{figure}

\par\end{center}

\begin{center}
\begin{figure}
\begin{centering}
\subfloat[Sod problem, HCCS(1,1,1,1), $DOFs=200$]{\begin{centering}
\includegraphics[width=0.47\textwidth]{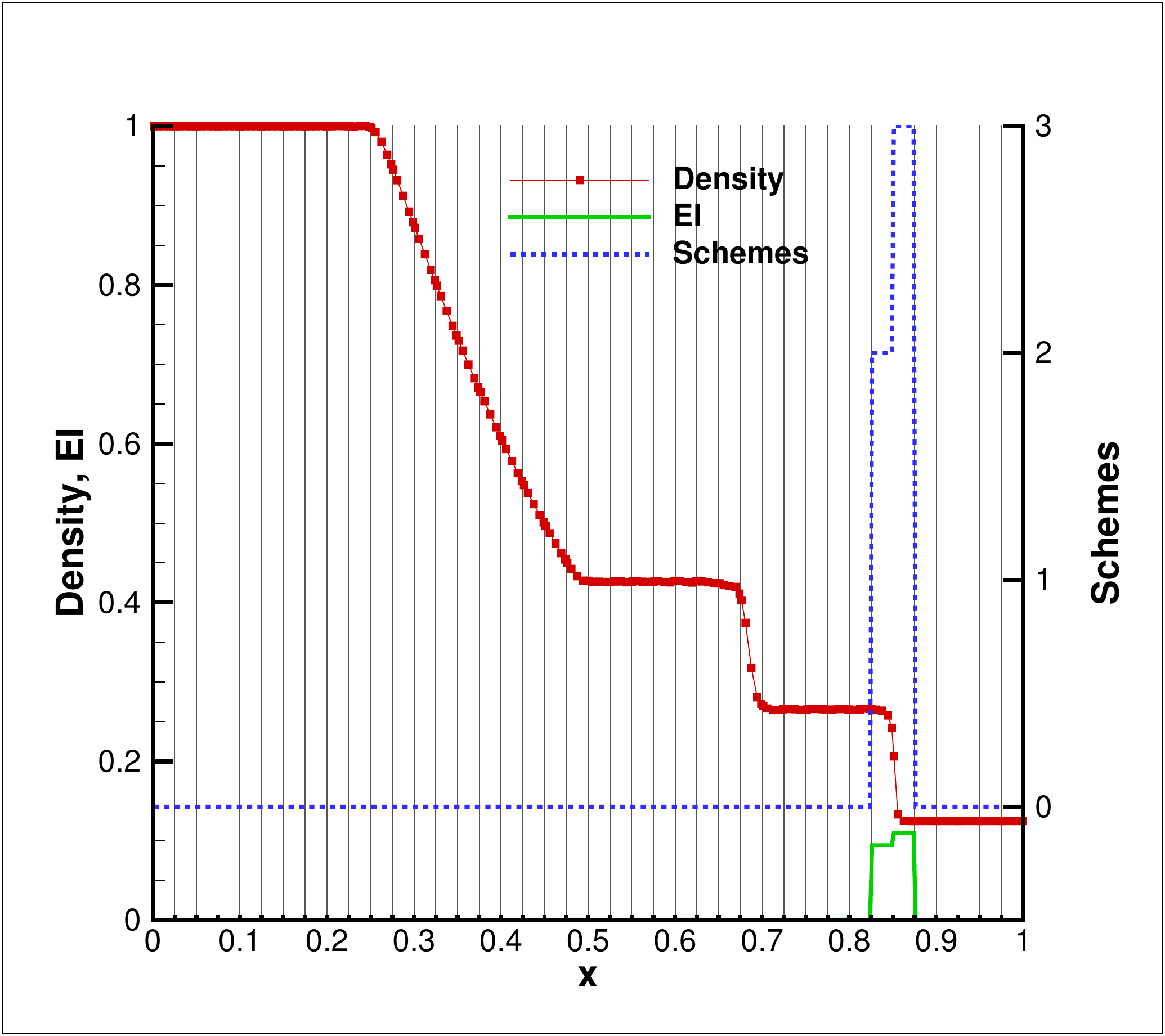}
\par\end{centering}

}\subfloat[Sod problem, HCCS(1,0,0,1), $DOFs=200$ ]{\begin{centering}
\includegraphics[width=0.47\textwidth]{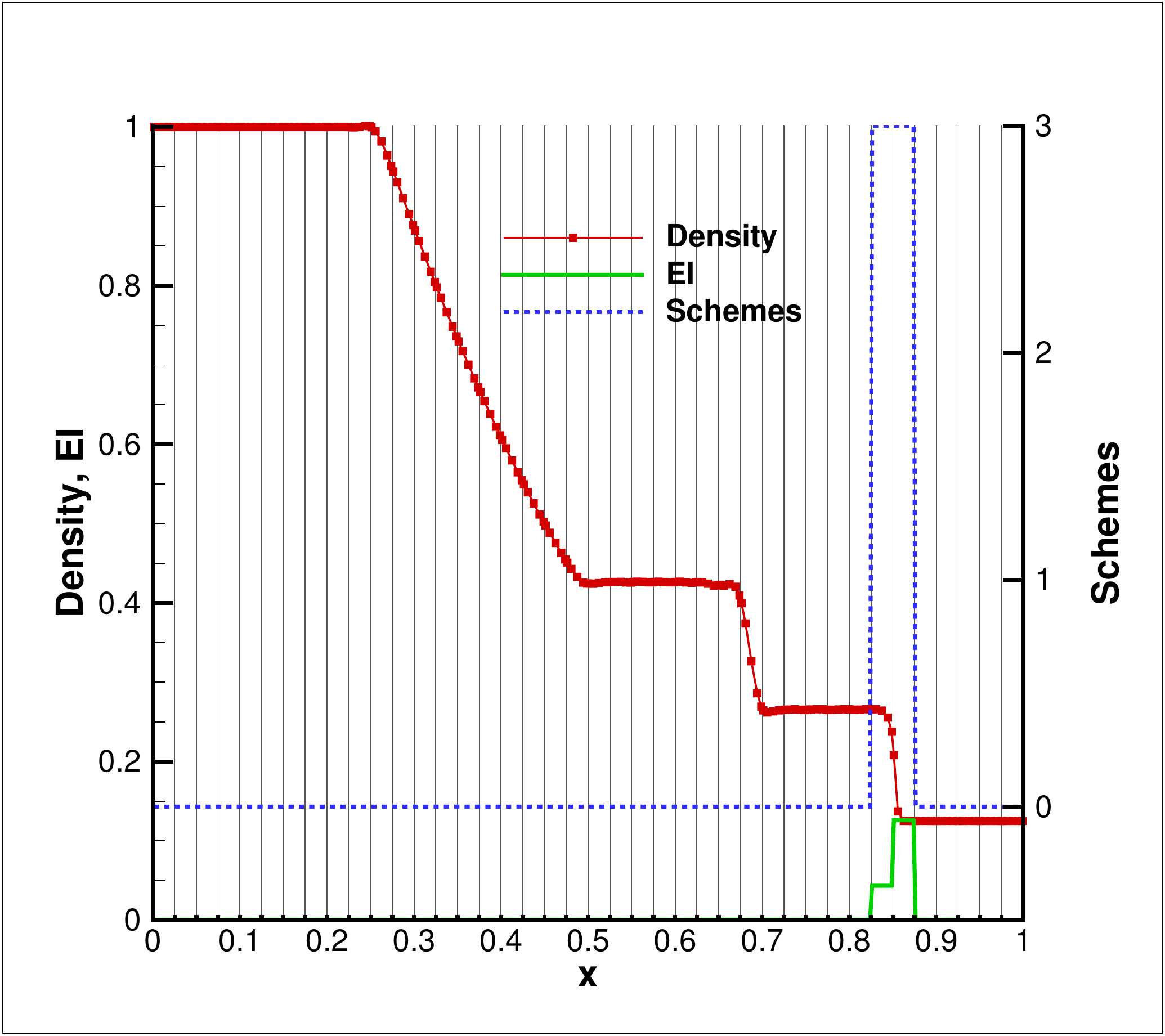}
\par\end{centering}

}
\par\end{centering}

\begin{centering}
\subfloat[Lax problem, HCCS(1,1,1,1), $DOFs=500$]{\begin{centering}
\includegraphics[width=0.47\textwidth]{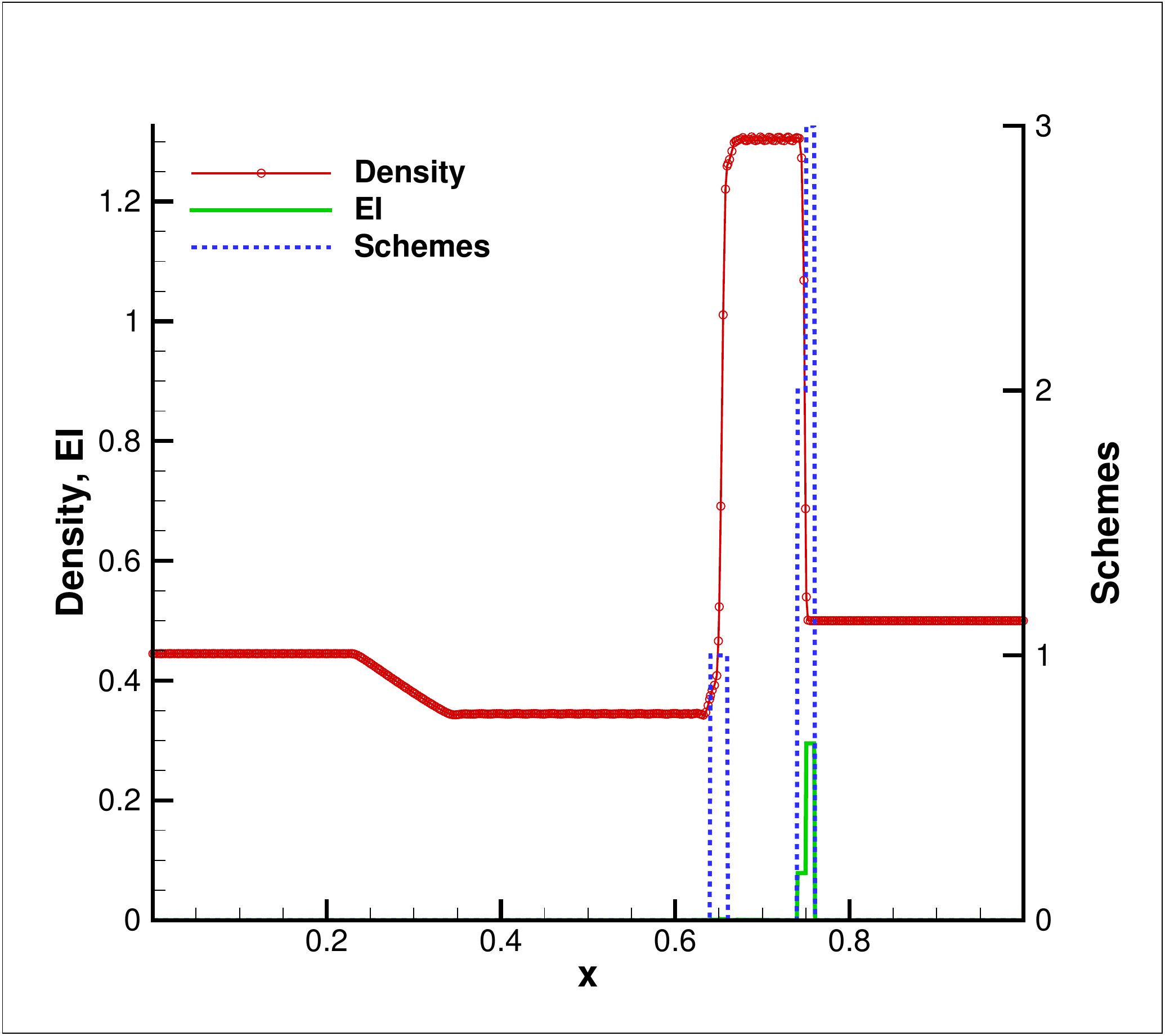} 
\par\end{centering}

}\subfloat[Lax problem, HCCS(1,0,0,1), $DOFs=500$]{\begin{centering}
\includegraphics[width=0.47\textwidth]{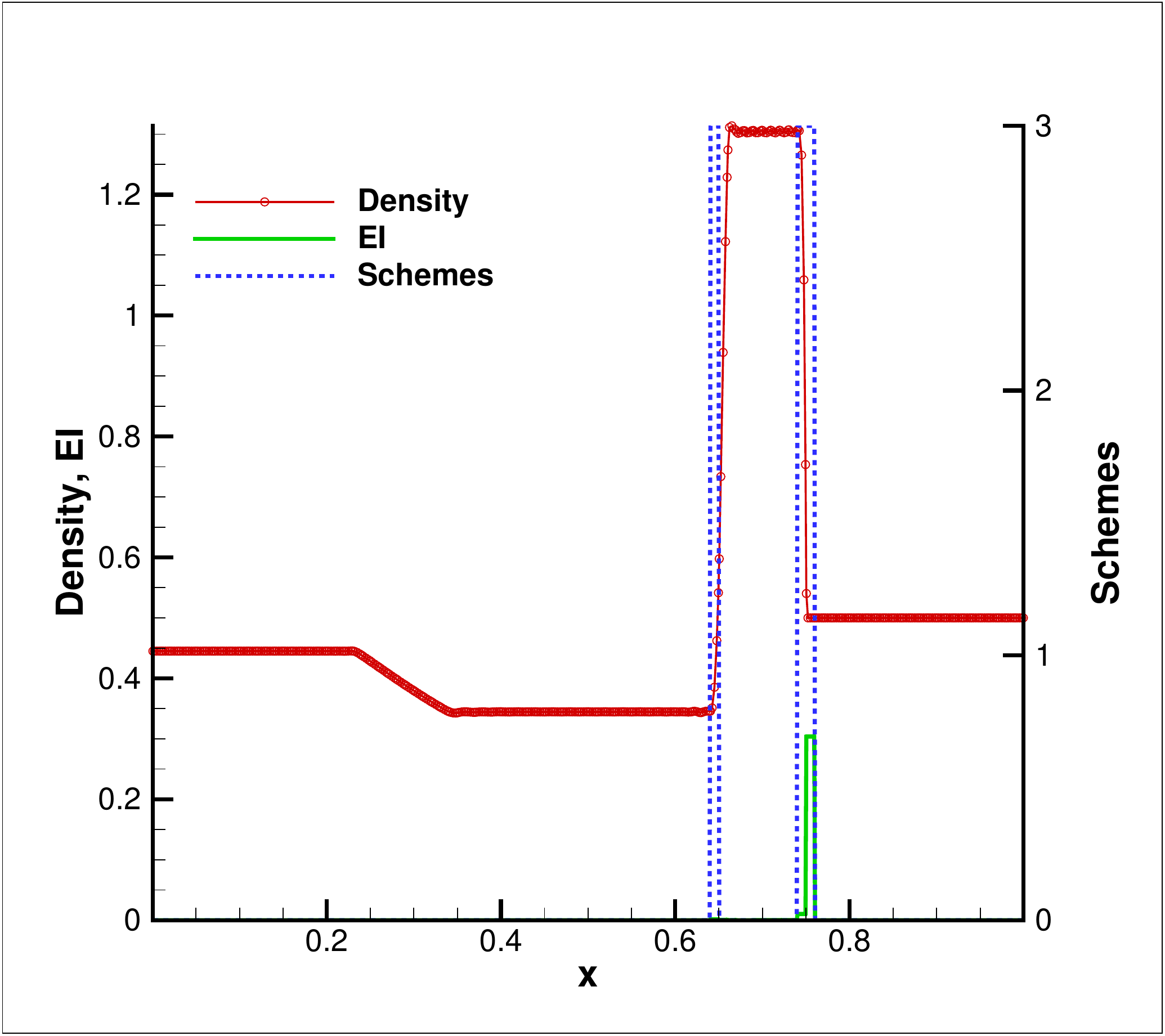}
\par\end{centering}

}
\par\end{centering}

\centering{}\caption{Hybrid schemes HCCS(1,1,1,1) and HCCS(1,0,0,1) in solving 1D shock
tube problems. Here schemes 0, 1, 2, 3 represent CPR, C5NNW5, C2NNW5
and C2NNW2, correspondingly. \label{fig:Hybrid-schemes-HS}}
\end{figure}

\par\end{center}

\subsection{Shu-Osher problem\label{sub:Shu-Osher-problem}}

The Shu-Osher problem with initial conditions 
\[
\left(\rho,u,p\right)=\begin{cases}
\left(3.857143,2.629369,10.333333\right), & -5\leq x<-4,\\
(1.0+0.2sin(5x),0,1.0), & -4\leq x<5,
\end{cases}
\]
is solved till $t=1.8$ with $DoFs=400$. 

Firstly, the CNNW schemes are used to solve the problem and are compared
with WCNS. The computed density distributions are shown in Fig. \ref{fig:Shu-Osher-single}.
C5NNW5, C2NNW5 and C2NNW2 can capture shocks without obvious oscillations.
In addition, C5NNW5 has the highest resolution, which is similar as
WCNS5. On the contrary, the C2NNW2 has lowest resolution. 

Secondly, CPR-CNNW schemes are tested in this problem. HCCS(1,1,1,1)
with $\mathbf{dv}=(c(a),0.05,0.1)$ and HCCS(1,0,0,1) with $\mathbf{dv}=(c(a),c(a),c(a))$
under different $a$ are applied to solve this problem. We consider
three cases $a=0.5$, $a=0.05$ and $a=0.005$. From Fig. \ref{fig:Hybrid-schemes},
we can see that both HCCS(1,1,1,1) and HCCS(1,0,0,1) can capture shocks
very well without obvious oscillations. The numbers of troubled cells
of HCCS(1,1,1,1) at $t=1.8$ are 2, 8, 19 for $a=0.5$, $a=0.05$
and $a=0.005$, correspondingly and the numbers of troubled cells
of HCCS(1,0,0,1) are 3, 7, 19, correspondingly. Comparing Fig. \ref{fig:Hybrid-schemes}(c)
and (d), we can see that HCCS(1,1,1,1) has much higher resolution
than HCCS(1,0,0,1). HCCS(1,1,1,1) can still obtain similar resolution
during increasing the number of troubled cells (by increasing $a$)
while the resolution of HCCS(1,0,0,1) decrease dramatically. Thus,
indicating wrongly will affect the resolution of CPR with second-order
CNNW limiting (HCCS(1,0,0,1)) since second-order CNNW scheme may be
used in smooth region. However it has less influences on the resolution
of CPR with p-adaptive CNNW limiting (HCCS(1,1,1,1)) since just few
troubled cells in shocks are computed by second-order schemes and
most troubled cells will be computed by C5NNW5 which also have high-order
of accuracy in smooth region. This results indicate that it is better
to include the high-order scheme C5NNW5 or the high resolution scheme
C2NNW5 in subcell limiting to keep high resolution. 

Thirdly, the CPR-CNNW schemes are compared with other four shock capturing
FE schemes. The following four results are compared: (1) the result
computed by DG-$P4$ with artificial viscosity, the MDH model and
DoFs=600 in \citet{Discacciati2020}; (2) the result of DG-$P3$ with
simple WENO limiter \citet{Zhong2013,Zhu2013,Du2015} under DoFs=400
from Fig. 1(a) in \citet{Li2020}; (3) the result of DG-$P5$ with
$p$-weighted limiter under DoFs=600 in Fig. 1(d) in \citet{Li2020};
(4) the result of DG-$P4$ with subcell shock limiting based on LGL
nodes under DoFs=640 in Fig. 5(c) in \citet{Hennemann2021}. as Comparisons
are shown in Fig. \ref{fig:comparison-schemes}. We can see that HCCS(1,1,1,1)
can obtain similar results in capturing the high-frequent waves with
fewer DoFs (DoFs=400). In addition, both the simple WENO limiter and
the $p$-weighted limiter lead to obvious deviations in wave phase
in Fig. \ref{fig:comparison-schemes}(b) and Fig. \ref{fig:comparison-schemes}(c).
The simple WENO limiter has obvious oscillations near $x=-2$ as shown
in Fig. \ref{fig:comparison-schemes}(b) while the $p$-weighted limiter
has overshoots near $x=2.45$ as shown in Fig. \ref{fig:comparison-schemes}(d).
Meanwhile the artificial viscosity approach has some oscillations
near shock, as shown in Fig. \ref{fig:comparison-schemes}(b). From
the figure we can also see that the subcell shock capturing approach
based on LGL points of \citet{Hennemann2021} results in obvious oscillations
and has overshoots near $x=-1.45$. On the contrary, HCCS(1,1,1,1)
and HCCS(1,0,0,1) obtain results with higher resolution without obvious
numerical oscillations. These results indicate that the CNNW subcell
limiting approach has better performances than simple WENO limiter
proposed in \citet{Zhong2013}, $p$-weighted limiter proposed in
\citet{Li2020}, artificial viscosity approach in \citet{Discacciati2020}
and subcell shock capturing on LGL nodes proposed in \citet{Hennemann2021}
in both resolution and shock capturing.

\begin{center}
\begin{figure}
\begin{centering}
\subfloat[{Density on the area $[-5.00,5.00]$,}]{\begin{centering}
\includegraphics[width=0.48\textwidth]{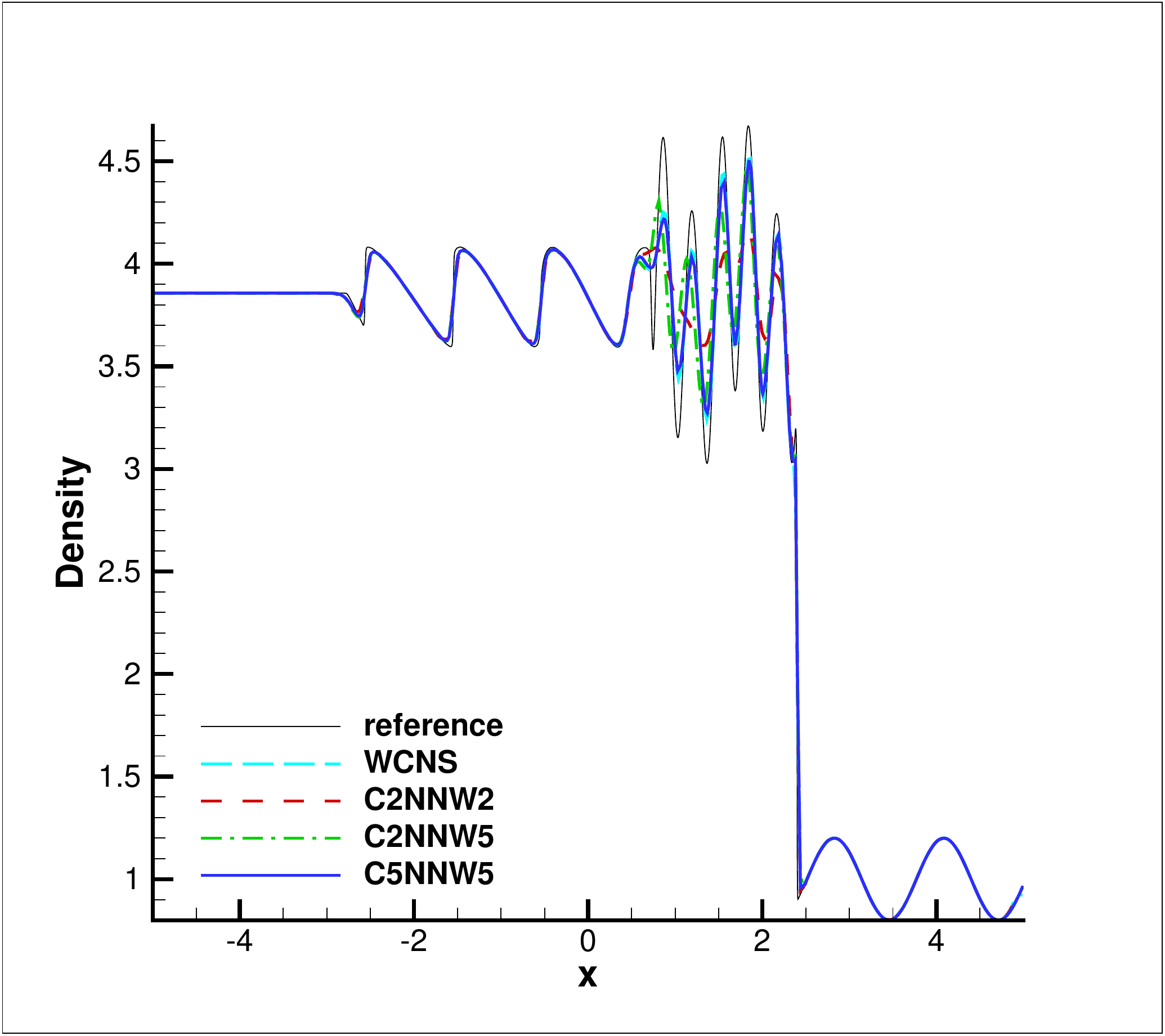}
\par\end{centering}

}\subfloat[{Density on the area $[0.5,2.4]$,}]{\begin{centering}
\includegraphics[width=0.48\textwidth]{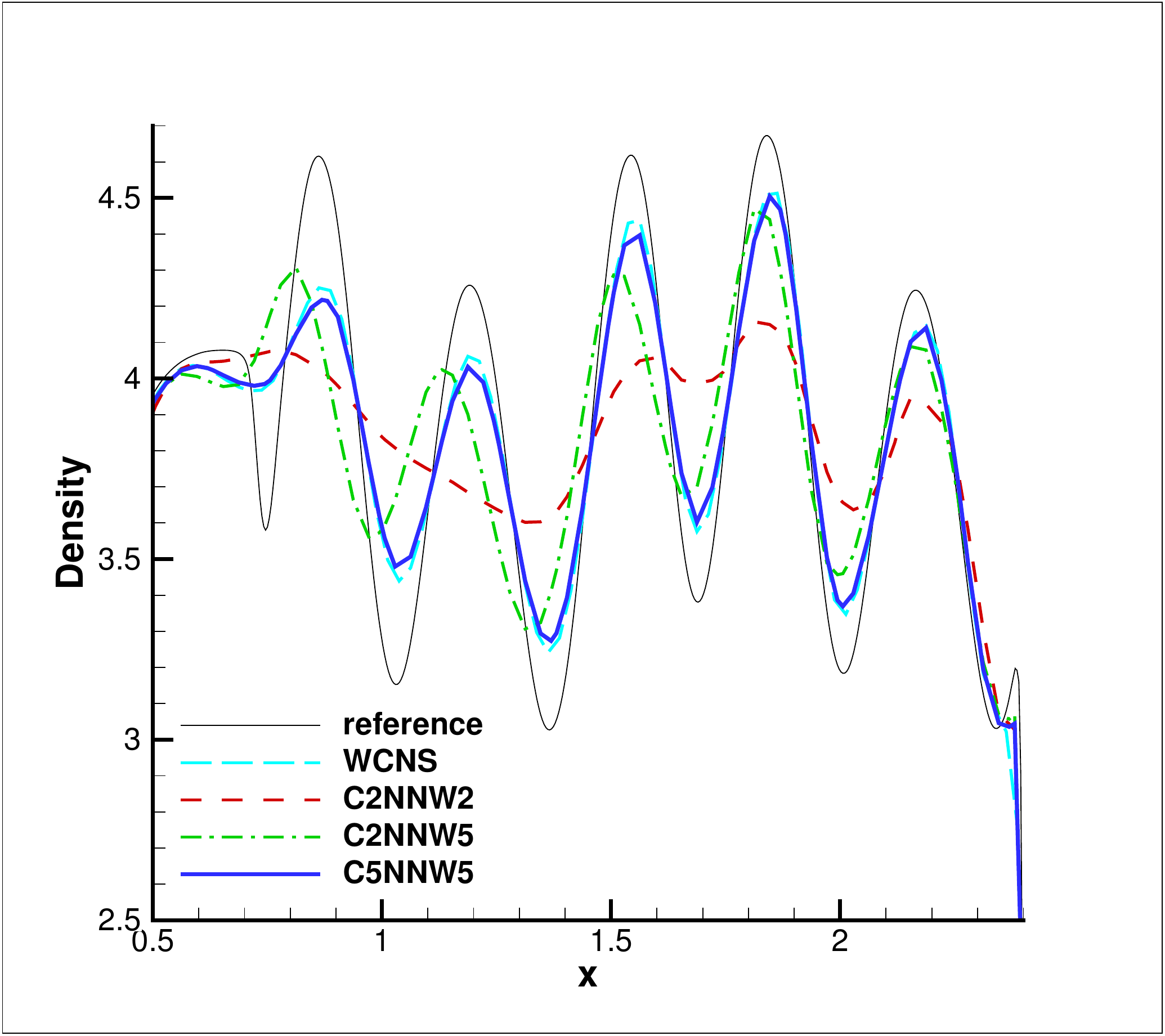} 
\par\end{centering}

}
\par\end{centering}

\caption{Density distributions of single schemes in solving Shu-Osher problem
($T=1.8$, $DOFs=400$).\label{fig:Shu-Osher-single}}
\end{figure}

\par\end{center}

\begin{center}
\begin{figure}
\begin{centering}
\subfloat[Distribution of density and schemes of HCCS(1,1,1,1)]{\begin{centering}
\includegraphics[width=0.48\textwidth]{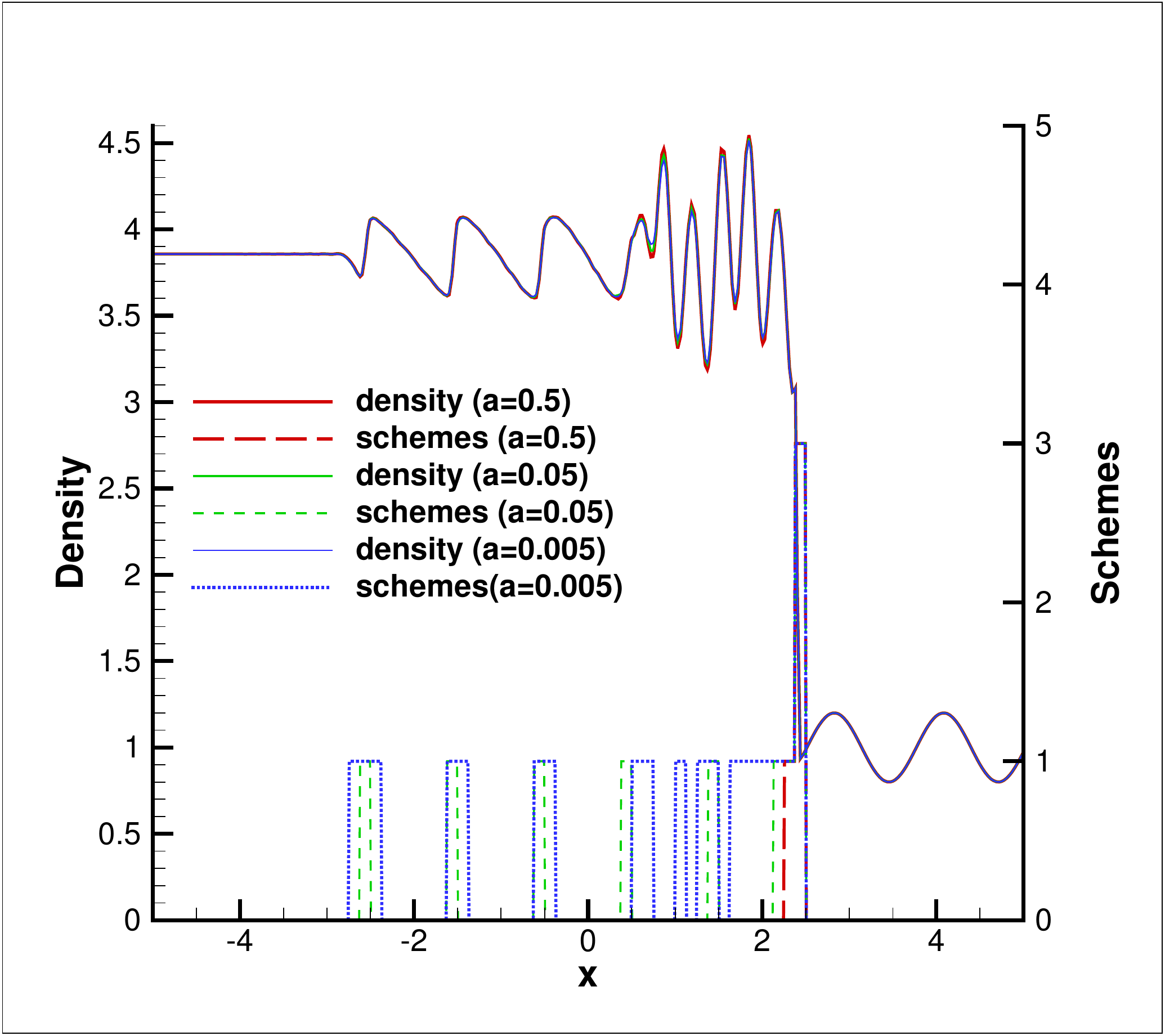}
\par\end{centering}

}\subfloat[Distribution of density and schemes of HCCS(1,0,0,1) ]{\begin{centering}
\includegraphics[width=0.48\textwidth]{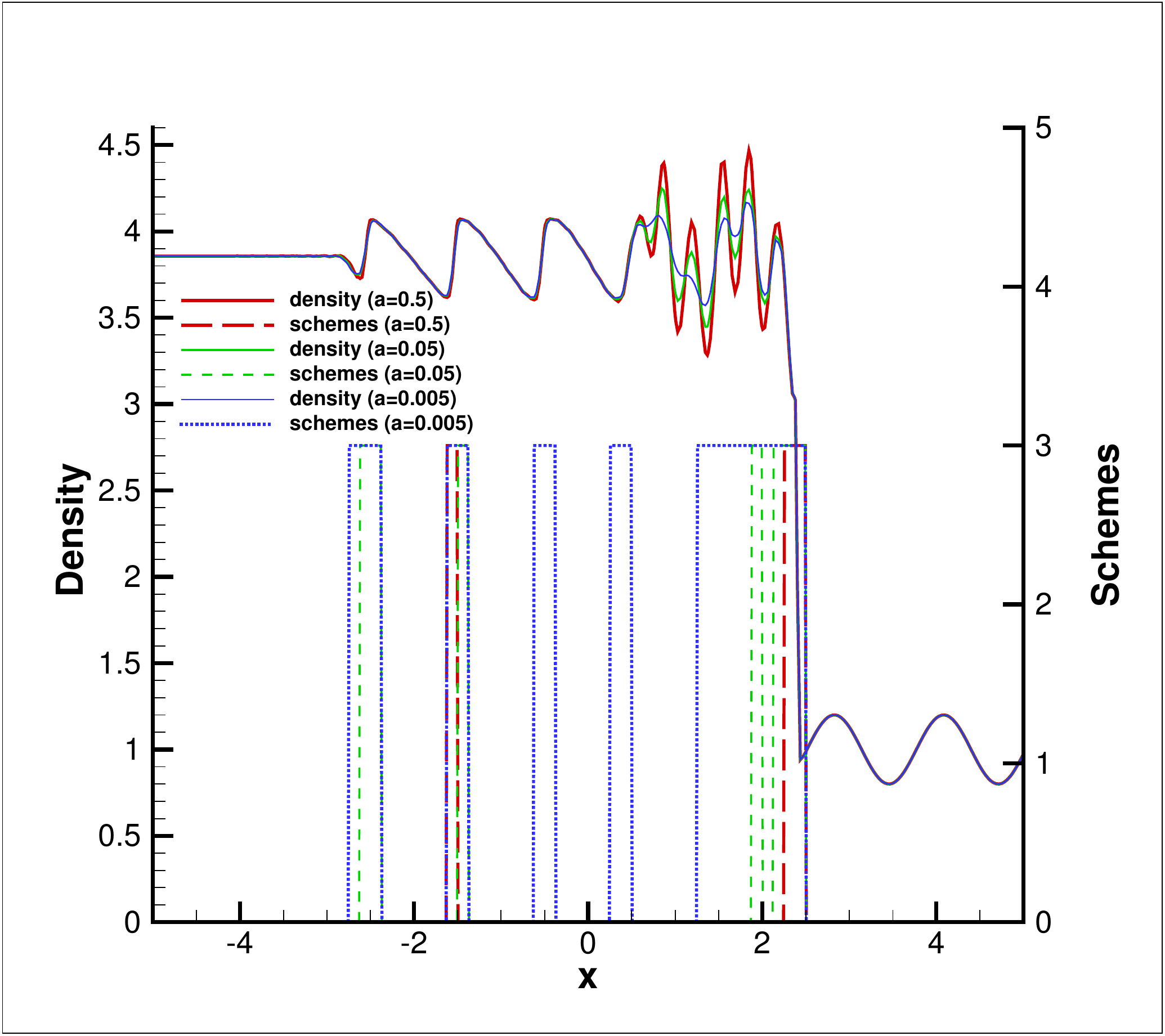}
\par\end{centering}

}
\par\end{centering}

\centering{}\subfloat[{Density on the area {[}\textminus{}0.4,2.40{]} of HCCS(1,1,1,1)}]{\begin{centering}
\includegraphics[width=0.48\textwidth]{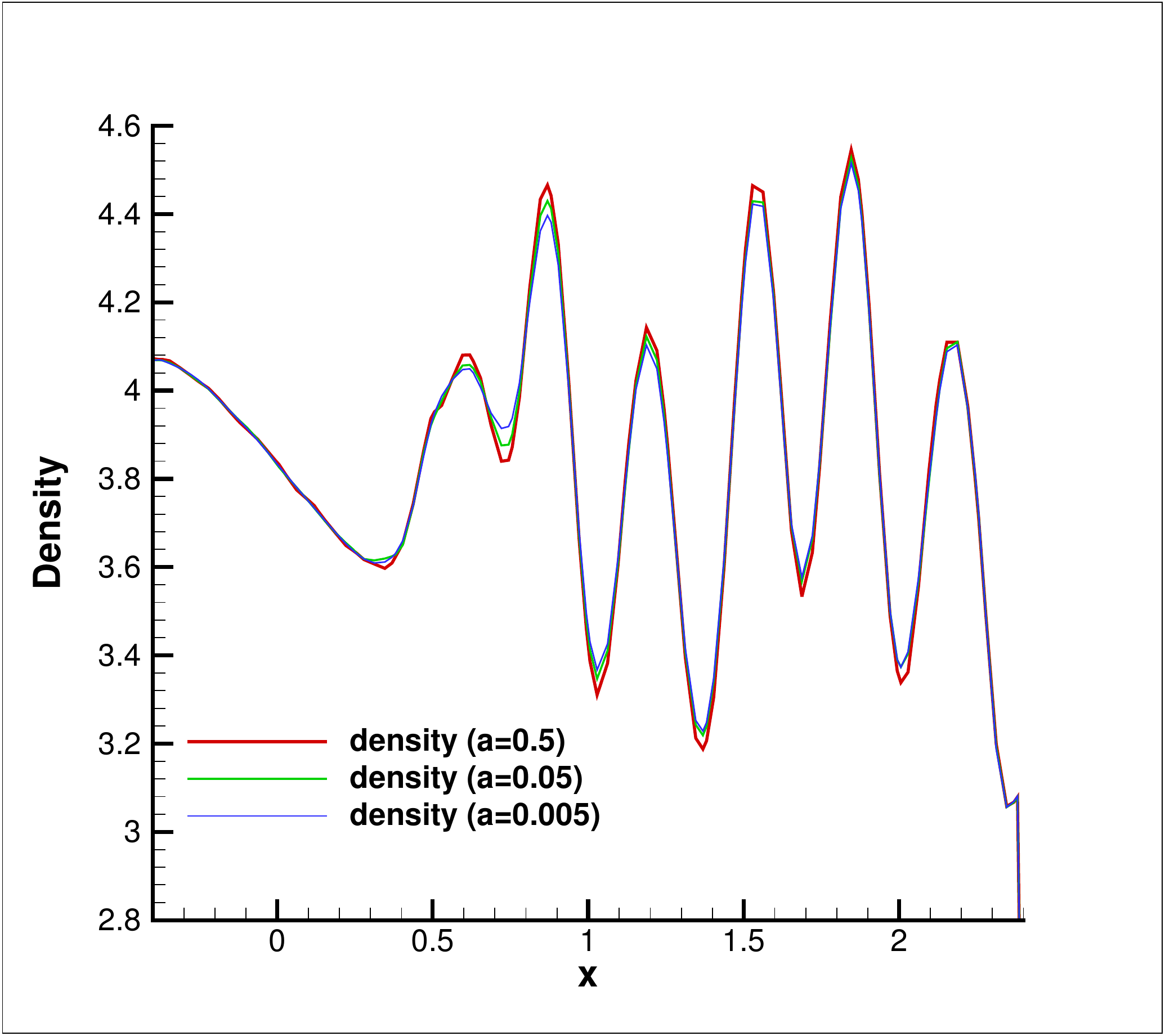}
\par\end{centering}

}\subfloat[{Density on the area {[}\textminus{}0.4,2.40{]} of HCCS(1,0,0,1) }]{\begin{centering}
\includegraphics[width=0.48\textwidth]{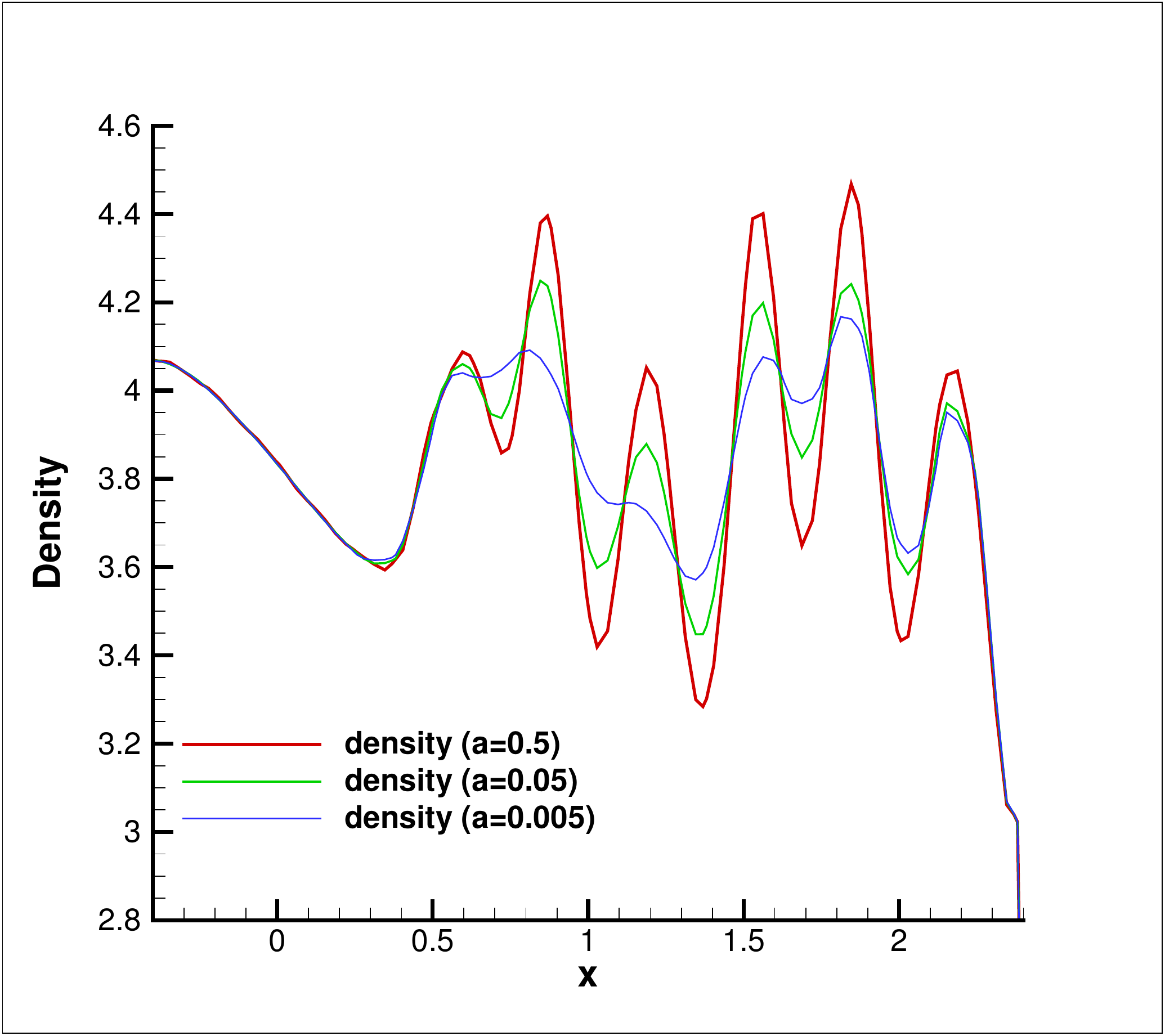}
\par\end{centering}

}\caption{Distribution of density and schemes for HCCS(1,1,1,1) and HCCS(1,0,0,1)
in solving Shu-Osher problem. Here schemes 0, 1, 2, 3 represent CPR,
C5NNW5, C2NNW5 and C2NNW2, correspondingly. ($T=1.8$, $DOFs=400$)
\label{fig:Hybrid-schemes}}
\end{figure}

\par\end{center}

\begin{center}
\begin{figure}
\begin{centering}
\subfloat[{Density on the area $[-5.00,5.00]$,}]{\begin{centering}
\includegraphics[width=0.48\textwidth]{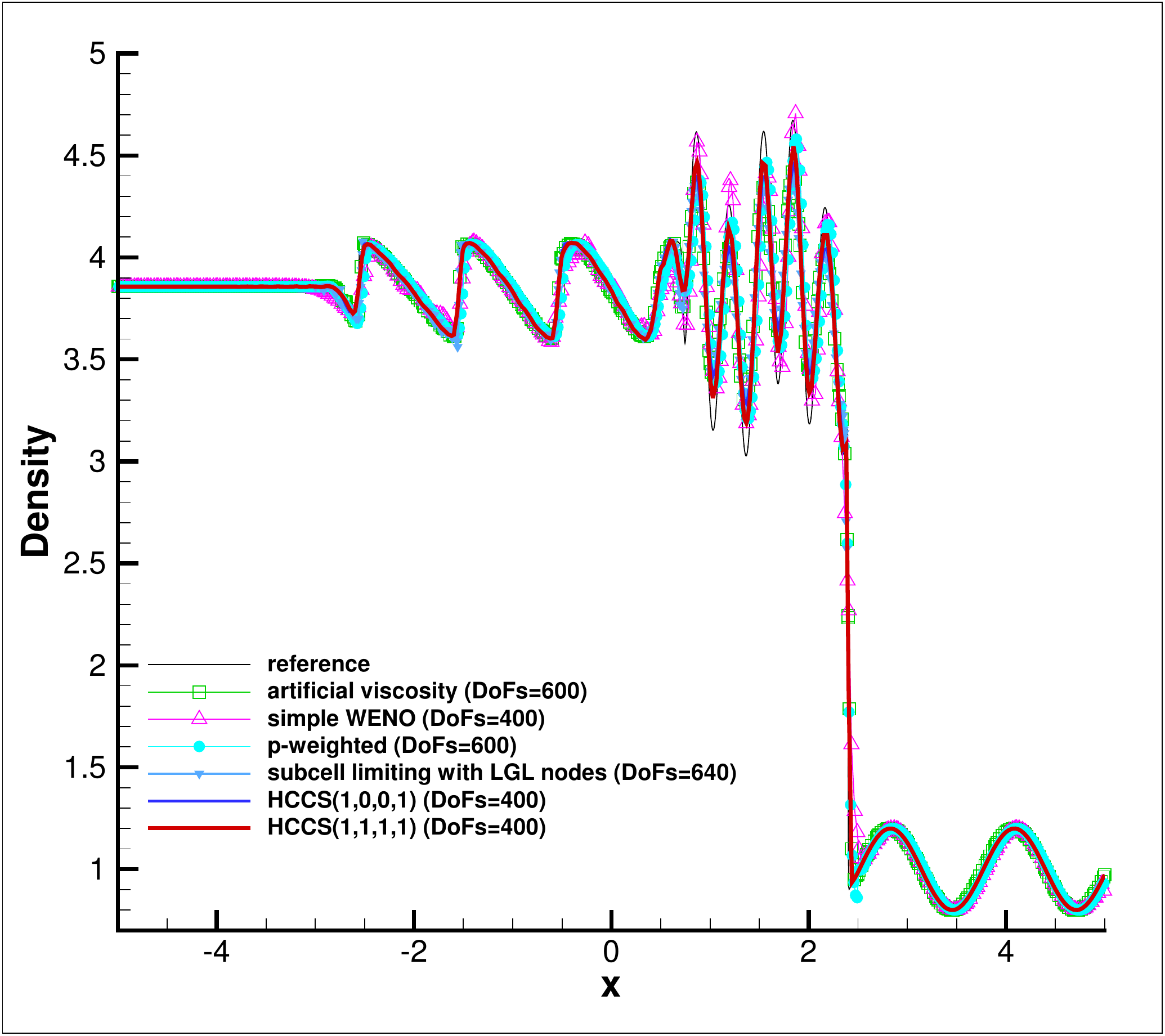}
\par\end{centering}

}\subfloat[{Density on the area $[-3.00,-1.40]$}]{\begin{centering}
\includegraphics[width=0.48\textwidth]{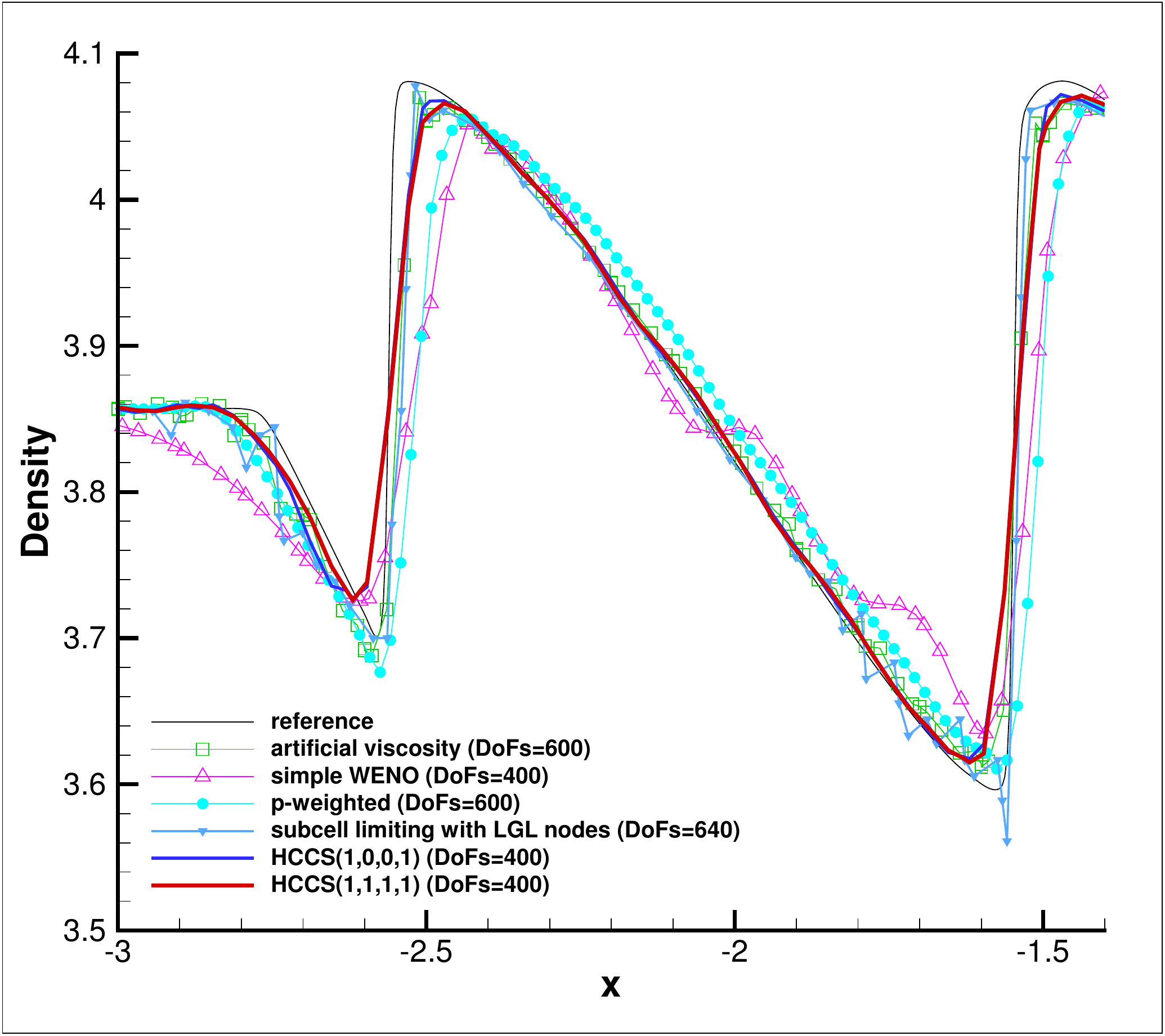}
\par\end{centering}

}
\par\end{centering}

\begin{centering}
\subfloat[{Density on the area $[0.50,2.40]$}]{\begin{centering}
\includegraphics[width=0.48\textwidth]{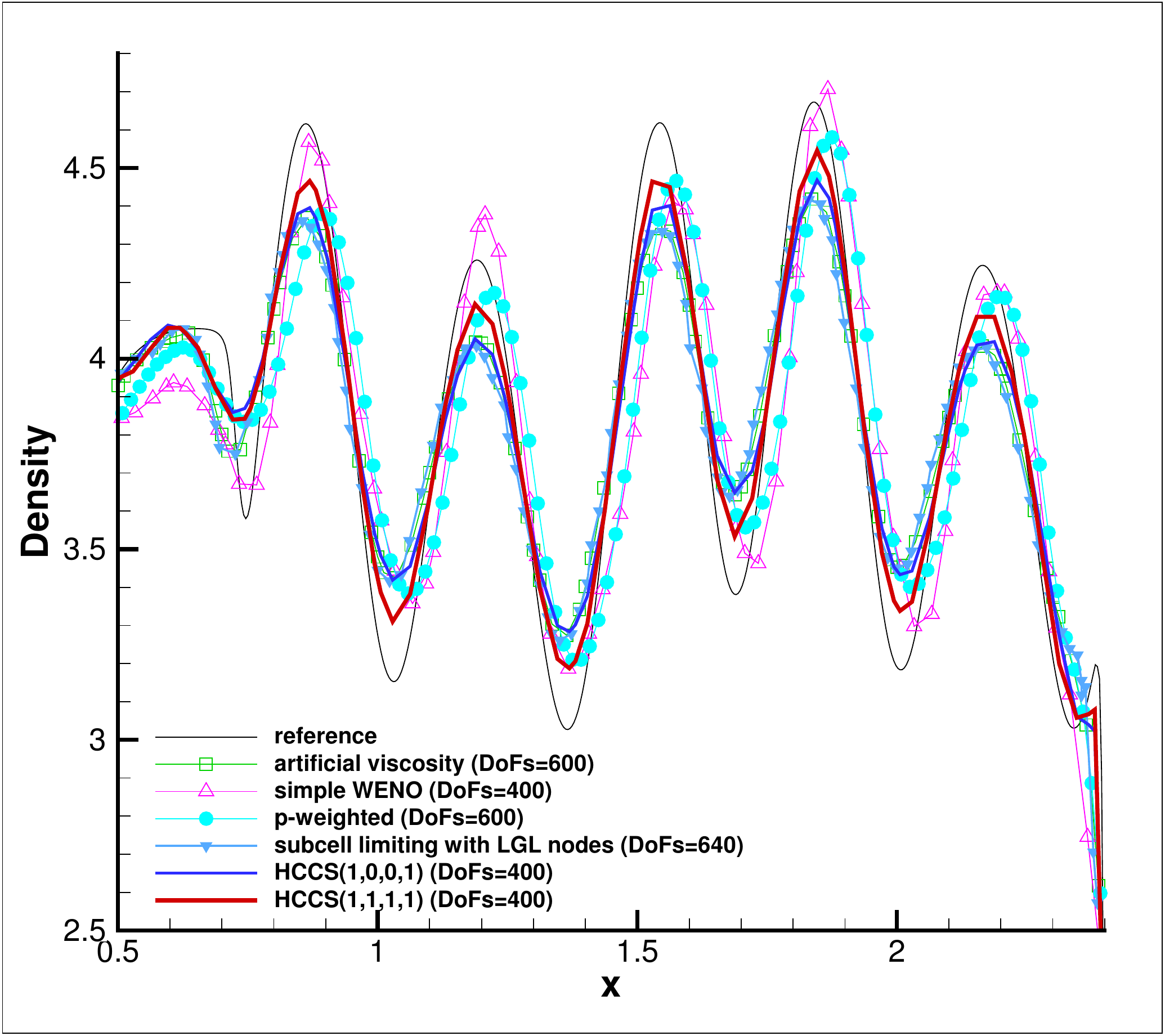}
\par\end{centering}

}\subfloat[{Density on area $[2.29,3.00]$ }]{\begin{centering}
\includegraphics[width=0.48\textwidth]{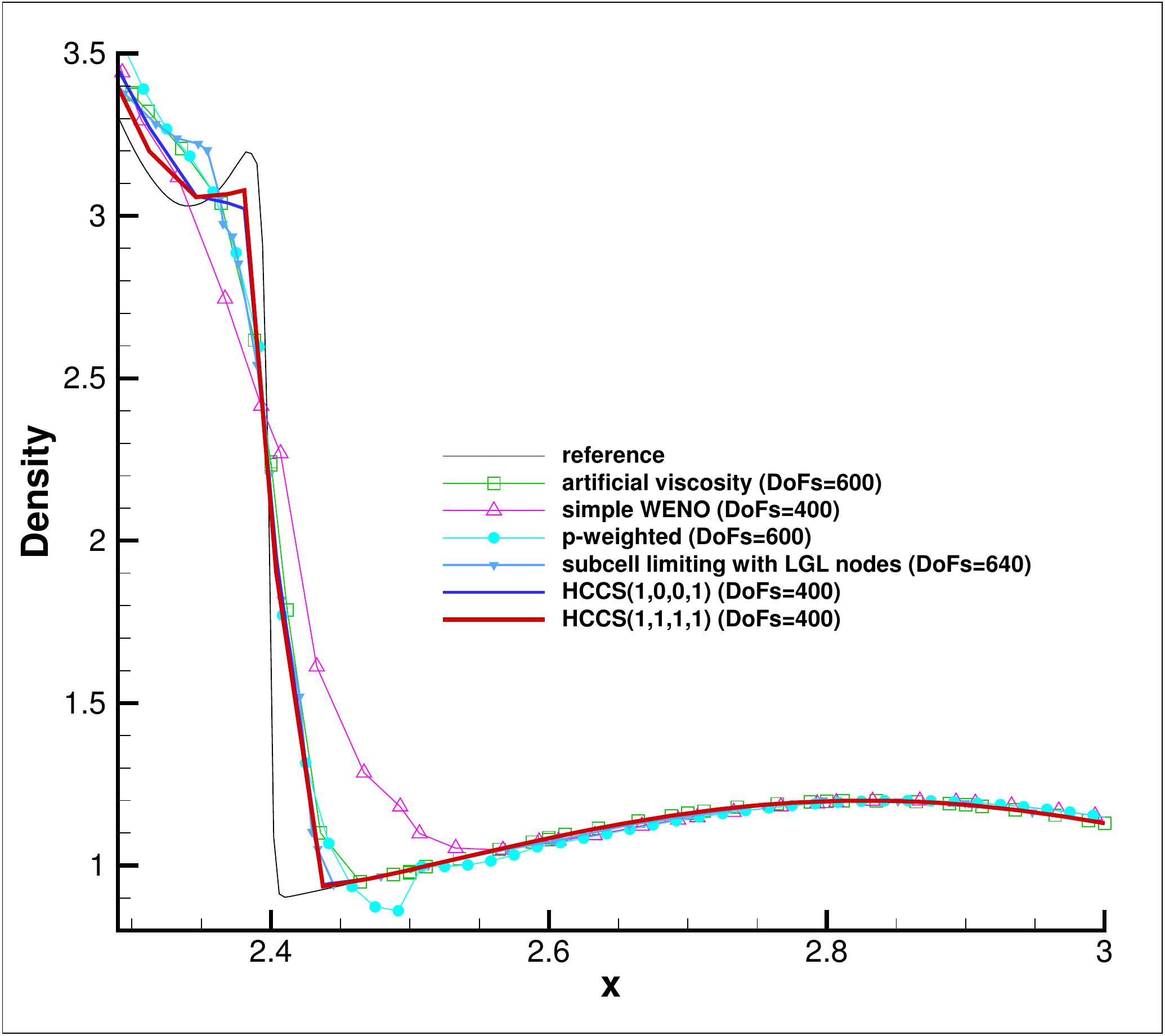}
\par\end{centering}

}
\par\end{centering}

\centering{}\caption{Comparison of hybrid schemes HCCS(1,1,1,1) and HCCS(1,0,0,1) with
artificial viscosity approach in \citet{Discacciati2020}, simple
WENO limiter in \citet{Li2020}, $p$-weighted limiter in \citet{Li2020}
and subcell shock capturing based on LGL nodes in \citet{Hennemann2021}
for solving 1D Euler equations.\label{fig:comparison-schemes}}
\end{figure}

\par\end{center}

\subsection{2D Euler vortex problem\label{sub:2D-Euler-vortex}}

In this subsection, the isentropic vortex problem in reference \citet{Shu1999}
is solved to test order of accuracy, time evolution of numerical errors
and discrete conservation properties of the proposed schemes. The
initial condition is a mean flow with $\{\rho,u,v,p\}=\{1,1,1,1\}$.
An isotropic vortex is then added to the mean flow with perturbations
in $u,v$ and $T=p/\rho$ but a constant entropy $S=p/\rho^{\gamma}$:
\begin{eqnarray*}
\left(\Delta u,\Delta v\right) & = & \frac{\varepsilon}{2\pi}e^{0.5(1-r^{2})}(-y,x),\\
\Delta T & = & -\frac{(\gamma-1)\varepsilon^{2}}{8\gamma\pi^{2}}e^{1-r^{2}},
\end{eqnarray*}
with $r=\sqrt{x^{2}+y^{2}}$ and the vortex strength $\varepsilon=5$.
In the numerical simulations, the computational domain is taken to
be $[-10,10]\times[-10,10]$ with periodic boundary conditions. Here
errors of density are adopted in calculating numerical errors.

Firstly, the fifth-order linear schemes, CPR5, C5NNW5 with linear
weights and WCNS5 with linear weights, are used to solve this problem
till $T=dt$ and $T=2$ with a time step of $dt=0.0001$. Tabs. \ref{tab:2d-vortex-problem-dt}
and \ref{tab:2d-vortex-problem-T2} present numerical errors of density
and orders of accuracy for different schemes. At $T=dt$, C5NNW5 and
CPR5 have fourth-order of accuracy while WCNS5 has fifth-order of
accuracy and the smallest numerical errors. At $T=2$, C5NNW5 and
WCNS5 can obtain fifth-order of accuracy which is higher than that
of CPR5. However, the CPR5 has the smallest numerical errors. To further
study the time evolution of errors, the time evolution of numerical
errors and orders of accuracy are presented in Fig. \ref{fig:Time-evolution}
for the three schemes. We can see that the numerical errors of the
CPR5 grows slower in time than the WCNS5 and C5NNW5, which illustrates
the reason that WCNS5 has smaller error than CPR5 at $T=dt$ while
has larger numerical errors than CPR5 at $T=2$. 

Secondly, the fifth-order nonlinear schemes C5NNW5 and WCNS5 are taken
to solve this problem till $T=2$ by taking JS weights with $\epsilon=10^{-6}$.
From Table \ref{tab:2d-vortex-nonlinear} we can see that schemes
using nonlinear weights on characteristic variables have larger numerical
errors than that on primary variables. In addition, both of C5NNW5
and WCNS5 can still acquire fifth-order of accuracy and C5NNW5 has
a bit larger numerical errors than WCNS5. 

Thirdly, the influence of indicating on accuracy of CPR-CNNW schemes
is tested. A special threshold value $c_{1}=c(0.0005)$ is considered,
which makes some cells in smooth regions being indicated wrongly as
troubled cells. HCCS(1,1,1,1) with $\mathbf{dv}=(c_{1},0.05,0.1)$
and HCCS(1,0,0,1) with $\mathbf{dv}=(c_{1},c_{1},c_{1})$ are applied
to solve this problem by using characteristic{\small{} variables} in
NNW interpolations. The numerical errors and CPU times are given in
Table \ref{tab:accuracy-CPUtime}. We can see that HCCS(1,1,1,1) still
have nearly fifth-order of accuracy while HCCS(1,0,0,1) has second-order
of accuracy, which illustrates that indicating wrongly has less influences
on accuracy of HCCS(1,1,1,1) than that of HCCS(1,0,0,1). In addition,
HCCS(1,1,1,1) and HCCS(1,0,0,1) have similar CPU time, which is nearly
111\% of pure CPR5 schemes for the grid with $\sqrt{DOFs}=800$.

\begin{center}
\begin{table}
\begin{centering}
\begin{tabular}{cccccccc}
\hline 
\multirow{2}{*}{{\small $Norm$}} & \multirow{2}{*}{{\small $\sqrt{DOFs}$}} & \multicolumn{2}{c}{{\small C5NNW5}} & \multicolumn{2}{c}{{\small CPR5}} & \multicolumn{2}{c}{{\small WCNS5}}\tabularnewline
\cline{3-8} 
 &  & \textit{\small error} & \textit{\small order} & \textit{\small error} & \textit{\small order} & \textit{\small error} & \textit{\small order}\tabularnewline
\multirow{4}{*}{{\small $L_{\infty}$}} & {\small $100$} & {\small 1.38E-07} & {\small -} & {\small 7.33E-07} & {\small -} & \textbf{\textcolor{black}{\small 6.52E-08}} & \textcolor{black}{\small -}\tabularnewline
 & {\small $200$} & {\small 1.03E-08} & {\small 3.74} & {\small 9.24E-08} & {\small 2.99} & \textbf{\textcolor{black}{\small 2.46E-09}} & \textcolor{black}{\small 4.73}\tabularnewline
 & {\small $400$} & {\small 6.69E-10} & {\small 3.94} & {\small 7.36E-09 } & {\small 3.65} & \textbf{\textcolor{black}{\small 8.09E-11}} & \textcolor{black}{\small 4.93}\tabularnewline
 & \uline{\small $800$} & {\small 4.11E-11} & {\small 4.02} & {\small 4.71E-10} & {\small 3.97} & \textbf{\textcolor{black}{\small 2.57E-12}} & \textcolor{black}{\small 4.98}\tabularnewline
\hline 
\multirow{4}{*}{{\small $L_{2}$}} & {\small $100$} & {\small 7.15E-09} & {\small -} & {\small 4.57E-08} & {\small -} & \textbf{\textcolor{black}{\small 4.12E-09}} & \textcolor{black}{\small -}\tabularnewline
 & {\small $200$} & {\small 3.91E-10} & {\small 4.19} & {\small 3.65E-09} & {\small 3.65} & \textbf{\textcolor{black}{\small 1.47E-10}} & \textcolor{black}{\small 4.81}\tabularnewline
 & {\small $400$} & {\small 2.31E-11} & {\small 4.08} & {\small 2.50E-10} & {\small 3.87} & \textbf{\textcolor{black}{\small 4.76E-12}} & \textcolor{black}{\small 4.95}\tabularnewline
 & \uline{\small $800$} & {\small 1.43E-12} & {\small 4.01} & {\small 1.59E-11} & {\small 3.97} & \textbf{\textcolor{black}{\small 1.50E-13}} & \textcolor{black}{\small 4.99}\tabularnewline
\hline 
\end{tabular}
\par\end{centering}

\noindent \centering{}\caption{Comparisons of linear schemes based on original physical variables
in solving 2D vortex problem ($T=dt=0.0001$).\label{tab:2d-vortex-problem-dt}}
\end{table}

\par\end{center}

\begin{center}
\begin{table}
\begin{centering}
\begin{tabular}{cccccccc}
\hline 
\multirow{2}{*}{{\small $Norm$}} & \multirow{2}{*}{{\small $\sqrt{DOFs}$}} & \multicolumn{2}{c}{{\small C5NNW5}} & \multicolumn{2}{c}{{\small CPR5}} & \multicolumn{2}{c}{{\small WCNS5}}\tabularnewline
\cline{3-8} 
 &  & \textit{\small error} & \textit{\small order} & \textit{\small error} & \textit{\small order} & \textit{\small error} & \textit{\small order}\tabularnewline
\multirow{4}{*}{{\small $L_{\infty}$}} & {\small $100$} & {\small 1.71E-3} & {\small -} & \textbf{\textcolor{black}{\small 3.00E-4}} & {\small -} & {\small 1.56E-3} & {\small -}\tabularnewline
 & {\small $200$} & {\small 7.08E-5} & {\small 4.59} & \textbf{\textcolor{black}{\small 7.58E-6}} & {\small 5.31} & {\small 5.24E-5} & {\small 4.90}\tabularnewline
 & {\small $400$} & {\small 2.88E-6} & {\small 4.62} & \textbf{\textcolor{black}{\small 5.55E-7}} & {\small 3.77} & {\small 1.64E-6} & {\small 5.00}\tabularnewline
 & \uline{\small $800$} & {\small 9.69E-8} & {\small 4.89} & \textbf{\textcolor{black}{\small 2.70E-8}} & {\small 4.36} & {\small 5.13E-8} & {\small 5.00}\tabularnewline
\hline 
\multirow{4}{*}{{\small $L_{2}$}} & {\small $100$} & {\small 6.83E-5} & {\small -} & \textbf{\textcolor{black}{\small 1.29E-5}} & {\small -} & {\small 5.30E-5} & {\small -}\tabularnewline
 & {\small $200$} & {\small 2.35E-6} & {\small 4.86} & \textbf{\textcolor{black}{\small 4.44E-7}} & {\small 4.86} & {\small 1.89E-6} & {\small 4.81}\tabularnewline
 & {\small $400$} & {\small 8.68E-8} & {\small 4.76} & \textbf{\textcolor{black}{\small 2.57E-8}} & {\small 4.11} & {\small 6.09E-8} & {\small 4.96}\tabularnewline
 & \uline{\small $800$} & {\small 2.62E-9} & {\small 5.05} & \textbf{\textcolor{black}{\small 1.31E-9}} & {\small 4.29} & {\small 1.91E-9} & {\small 4.99}\tabularnewline
\hline 
\end{tabular}
\par\end{centering}

\noindent \centering{}\caption{Comparisons of linear schemes based on original physical variables
in solving 2D vortex problem ($T=2$).\label{tab:2d-vortex-problem-T2}}
\end{table}

\par\end{center}

\begin{center}
\begin{figure}
\begin{centering}
\subfloat[Time evolution of L2 accuracy order based on $DOFs=400$ and $DOFs=800$,
$t_{k}=0.0001,0.2,0.4,0.6,0.8,...,2.0$.]{\begin{centering}
\includegraphics[width=0.45\textwidth]{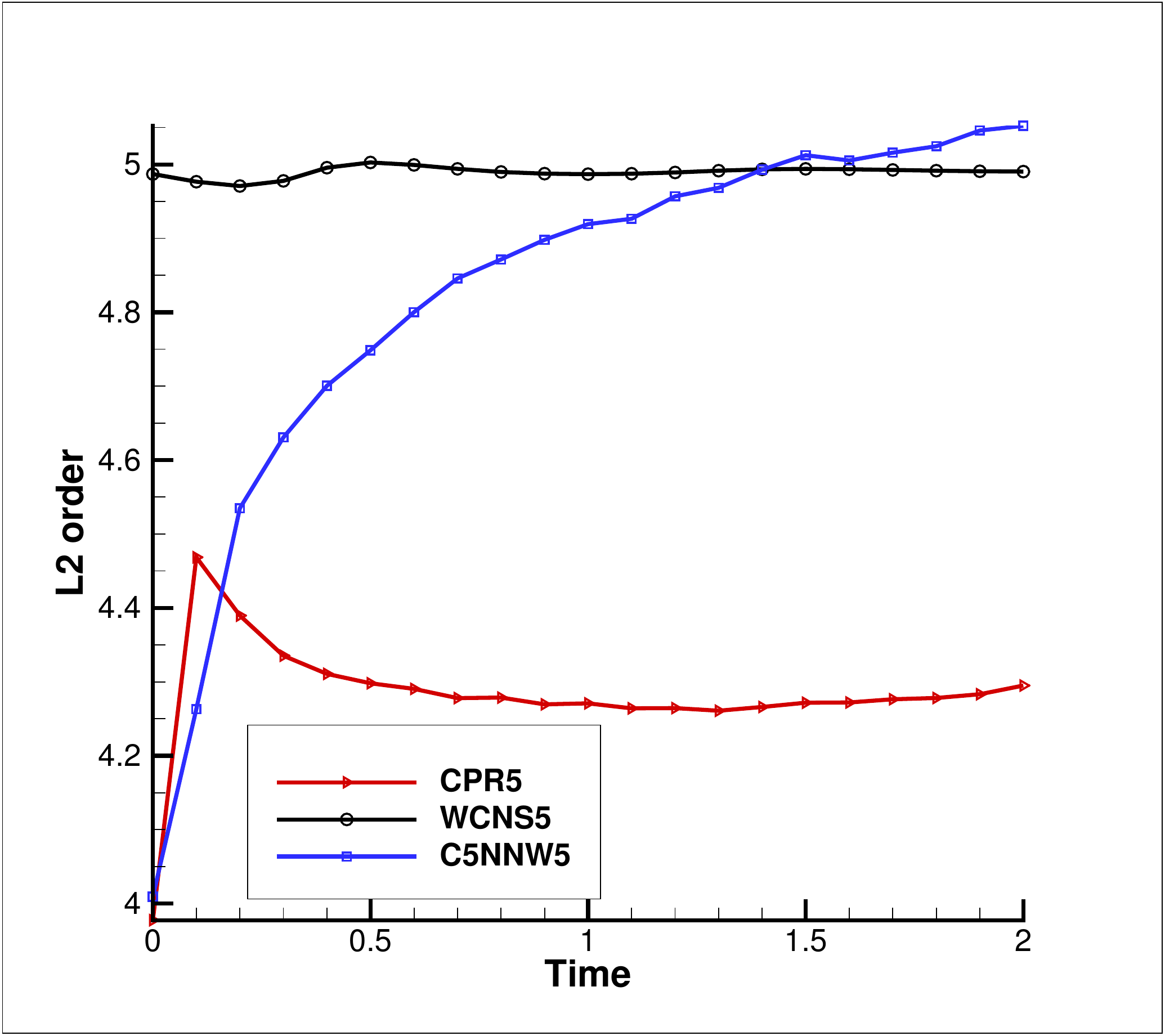}
\par\end{centering}

}\subfloat[Time evolution of L2 numerical error ($DOFs=800$)]{\begin{centering}
\includegraphics[width=0.45\textwidth]{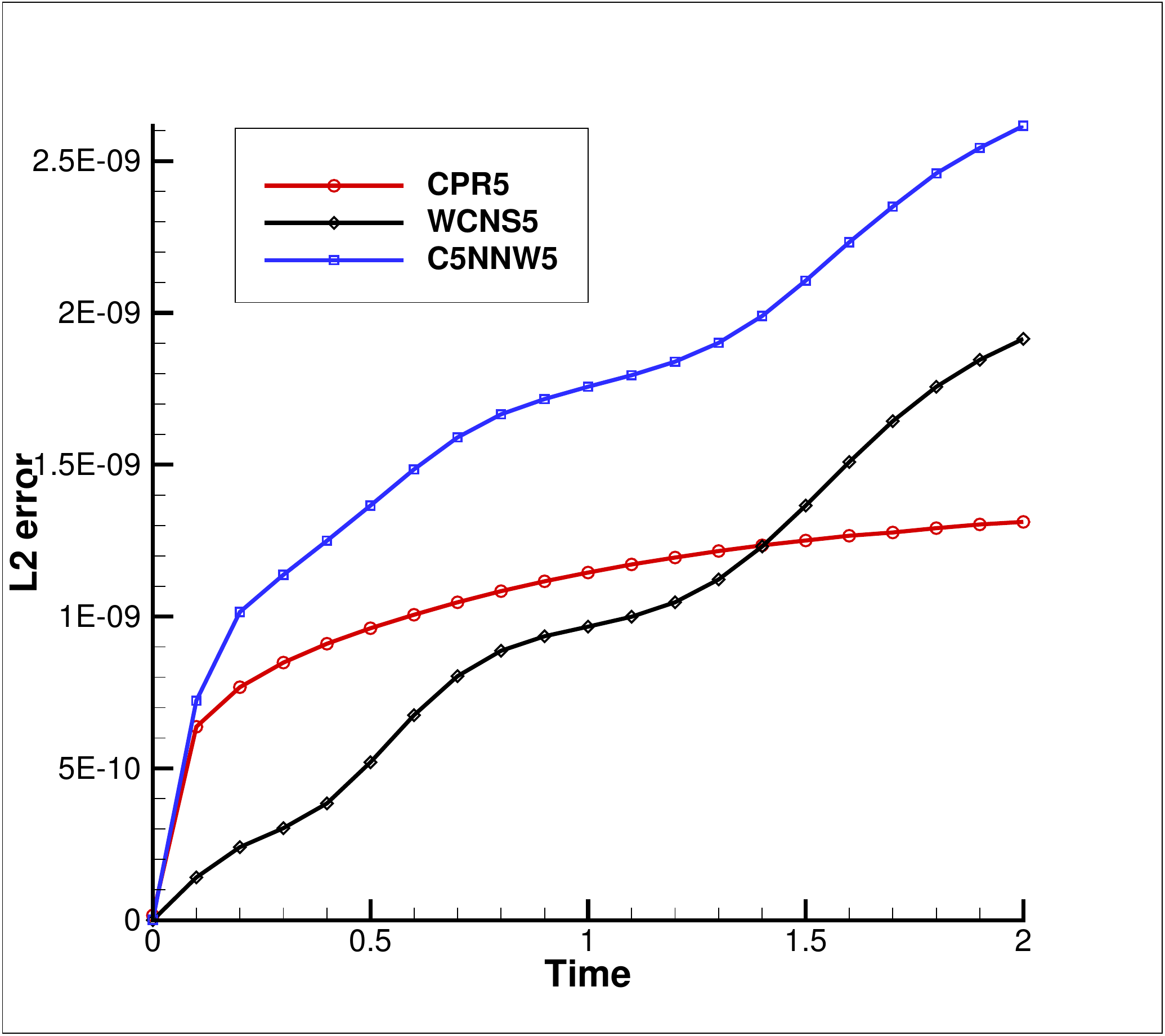}
\par\end{centering}

}
\par\end{centering}

\caption{Time evolution of numerical accuracy order and numerical errors ($T=2$)\label{fig:Time-evolution}}
\end{figure}

\par\end{center}

\begin{center}
\begin{table}
\begin{centering}
\begin{tabular}{cccccccccc}
\cline{3-10} 
\multirow{3}{*}{{\small $Norm$}} & \multirow{3}{*}{{\small $\sqrt{DOFs}$}} & \multicolumn{4}{c}{{\small Primary variables}} & \multicolumn{4}{c}{{\small Characteristic variables}}\tabularnewline
\cline{3-10} 
 &  & \multicolumn{2}{c}{{\small C5NNW5}} & \multicolumn{2}{c}{{\small WCNS5}} & \multicolumn{2}{c}{{\small C5NNW5}} & \multicolumn{2}{c}{{\small WCNS5}}\tabularnewline
\cline{3-10} 
 &  & \textit{\small error} & \textit{\small order} & \textit{\small error} & \textit{\small order} & \textit{\textcolor{black}{\small error}} & \textit{\textcolor{black}{\small order}} & \textit{\small error} & \textit{\small order}\tabularnewline
\multirow{4}{*}{{\small $L_{\infty}$}} & {\small $100$} & {\small 4.41E-03} & {\small -} & \textbf{\textcolor{black}{\small 3.08E-03}} & {\small -} & \textcolor{black}{\small 1.12E-02} & {\small -} & \textbf{\textcolor{black}{\small 7.63E-03}} & {\small -}\tabularnewline
 & {\small $200$} & {\small 1.97E-04} & {\small 4.48} & \textbf{\textcolor{black}{\small 1.55E-04}} & {\small 4.31} & \textcolor{black}{\small 4.43E-04} & {\small 4.66} & \textbf{\textcolor{black}{\small 2.80E-04}} & {\small 4.77}\tabularnewline
 & {\small $400$} & {\small 9.18E-06} & {\small 4.42} & \textbf{\textcolor{black}{\small 5.65E-06}} & {\small 4.78} & \textcolor{black}{\small 8.40E-06} & {\small 5.72} & \textbf{\textcolor{black}{\small 3.98E-06}} & {\small 6.14}\tabularnewline
 & \uline{\small $800$} & {\small 3.06E-07} & {\small 4.91} & \textbf{\textcolor{black}{\small 1.79E-07}} & {\small 4.98} & {\small 3.57E-07} & {\small 4.56} & \textbf{\textcolor{black}{\small 1.46E-07}} & {\small 4.77}\tabularnewline
\hline 
\multirow{4}{*}{{\small $L_{2}$}} & {\small $100$} & {\small 2.17E-04} & {\small -} & \textbf{\textcolor{black}{\small 1.39E-04}} & {\small -} & {\small 3.70E-04} & {\small -} & \textbf{\textcolor{black}{\small 2.44E-04}} & {\small -}\tabularnewline
 & {\small $200$} & {\small 9.27E-06} & {\small 4.55} & \textbf{\textcolor{black}{\small 6.73E-06}} & {\small 4.37} & {\small 1.95E-05} & {\small 4.25} & \textbf{\textcolor{black}{\small 1.15E-05}} & {\small 4.41}\tabularnewline
 & {\small $400$} & {\small 3.54E-07} & {\small 4.71} & \textbf{\textcolor{black}{\small 2.43E-07}} & {\small 4.79} & {\small 4.27E-07} & {\small 5.51} & \textbf{\textcolor{black}{\small 2.54E-07}} & {\small 5.50}\tabularnewline
 & \uline{\small $800$} & {\small 1.11E-08} & {\small 5.00} & \textbf{\textcolor{black}{\small 7.98E-09}} & {\small 4.93} & {\small 1.35E-08} & {\small 4.98} & \textbf{\textcolor{black}{\small 8.95E-09}} & {\small 4.83}\tabularnewline
\hline 
\end{tabular}
\par\end{centering}

\noindent \centering{}\caption{Comparisons of nonlinear schemes based on primary variables and characteristic
variables in solving 2D vortex problem ($T=2$).\label{tab:2d-vortex-nonlinear}}
\end{table}

\par\end{center}

\begin{table}
\begin{centering}
\begin{tabular}{ccccccccc}
\hline 
\multirow{3}{*}{{\small $Norm$}} & \multirow{3}{*}{{\small $\sqrt{DOFs}$}} & \multicolumn{6}{c}{{\small CPR-CNNW schemes}} & \multicolumn{1}{c}{{\small Single scheme}}\tabularnewline
\cline{3-9} 
 &  & \multicolumn{3}{c}{{\small HCCS(1,1,1,1)}} & \multicolumn{3}{c}{{\small HCCS(1,0,0,1)}} & \multicolumn{1}{c}{{\small CPR5}}\tabularnewline
\cline{3-9} 
 &  & \textit{\small error} & \textit{\small order} & {\small CPU (s)} & \textit{\small error} & \textit{\small order} & {\small CPU (s)} & {\small CPU (s)}\tabularnewline
\multirow{3}{*}{{\small $L_{\infty}$}} & {\small $100$} & {\small 6.01E-03} & {\small -} & {\small 64.5312} & {\small 3.3057E-02} & {\small -} & {\small 64.6875} & {\small 52.2656}\tabularnewline
\cline{2-9} 
 & {\small $200$} & {\small 3.41E-04} & {\small 4.13} & {\small 330.5469} & {\small 7.4320E-03} & {\small 2.15} & {\small 330.0938} & {\small 278.1406}\tabularnewline
\cline{2-9} 
 & {\small $400$} & {\small 7.45E-06} & {\small 5.52} & {\small 1421.4531} & {\small 1.5424E-03} & {\small 2.27} & {\small 1412.9844} & {\small 1279.2812}\tabularnewline
\hline 
\end{tabular}
\par\end{centering}

\noindent \centering{}\caption{Comparisons of numerical errors and CPU time of HCCS(1,1,1,1) and
HCCS(1,0,0,1) in solving 2D vortex problem ($T=2$).\label{tab:accuracy-CPUtime}}
\end{table}
Fourthly, discrete conservation properties are tested for both single
schemes and the hybrid CPR-CNNW schemes. The discrete conservation
error is defined as the relative error in the preservation of the
integral of the conservative quantity $\rho$:

\[
\left\langle \rho\right\rangle =\frac{\int_{V}(\rho-\rho_{0})dxdy}{\int_{V}\rho_{0}dxdy}=\frac{\int_{V}(\rho-\rho_{0})Jd\xi d\eta}{\int_{V}\rho_{0}Jd\xi d\eta}=\frac{INT(\rho)-INT(\rho_{0})}{INT(\rho_{0})}
\]
where $\rho$ and $\rho_{0}$ denote the density at the final time
and the initial time, correspondingly. Here $V$ indicates the whole
computational domain and $INT(\rho)=\int_{V}\rho Jd\xi d\eta$, which
is estimated by

\[
\sum_{i,j}\left(\sum_{l}\sum_{m}w_{l}w_{m}\rho_{i,j,l,m}J_{i,j,l,m}\right)\Delta\xi\Delta\eta,
\]
for all schemes, where $w_{l}$ are the weights of the quadrature
rule. In order to test discrete conservation of CPR-CNNW schemes,
a special threshold value $c_{1}=c(0.0005)=0.0005\cdot10^{-1.8(6+1)^{1/4}}\thickapprox5.9\times10^{-7}$
is taken. Then, HCCS(1,1,1,1) with $\mathbf{dv}=(c_{1},5\times10^{-6},10^{-5})$
and HCCS(1,0,0,1) with $\mathbf{dv}=(c_{1},c_{1},c_{1})$ are applied
to maintain that both CPR and CNNW are used in the computations. The
density distribution, error distribution, scheme distribution and
time evolution of the error $\left\langle \rho\right\rangle $ of
HCCS(1,1,1,1) are shown in Fig. \ref{fig:HS(1,1,1,1)-evolution}.
We can see that the obtained flow is smooth and four schemes are adopted
in the computation. In addition, the error $\left\langle \rho\right\rangle $
keeps almost zero. Errors in the conservation $\left\langle \rho\right\rangle $
are also summarized in Tab. \ref{tab:errors of invariant} with $INT(\rho_{0})\approx398.241743560187$
for this problem. All the conservation errors are at the level of
rounding errors, which indicates that all the schemes satisfy discrete
conservation. 

\begin{center}
\begin{figure}
\begin{centering}
\subfloat[Density distribution at $T=20$.]{\begin{centering}
\includegraphics[width=0.45\textwidth]{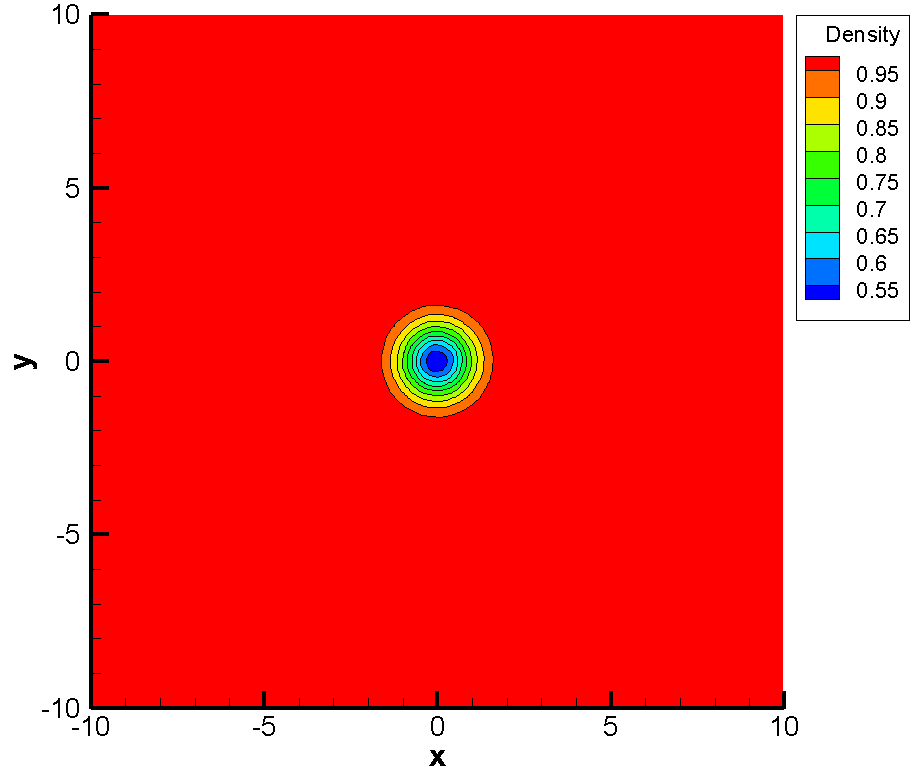}
\par\end{centering}

}\subfloat[Error distribution at $T=20$.]{\begin{centering}
\includegraphics[width=0.45\textwidth]{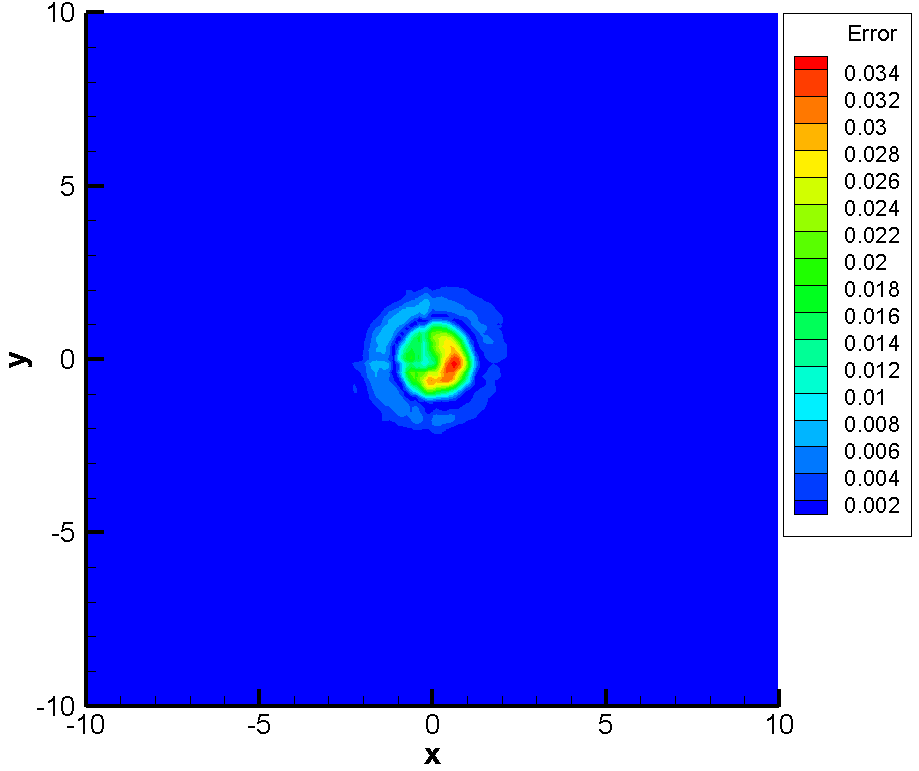}
\par\end{centering}

}
\par\end{centering}

\begin{centering}
\subfloat[Scheme distribution at $T=20$.]{\begin{centering}
\includegraphics[width=0.45\textwidth]{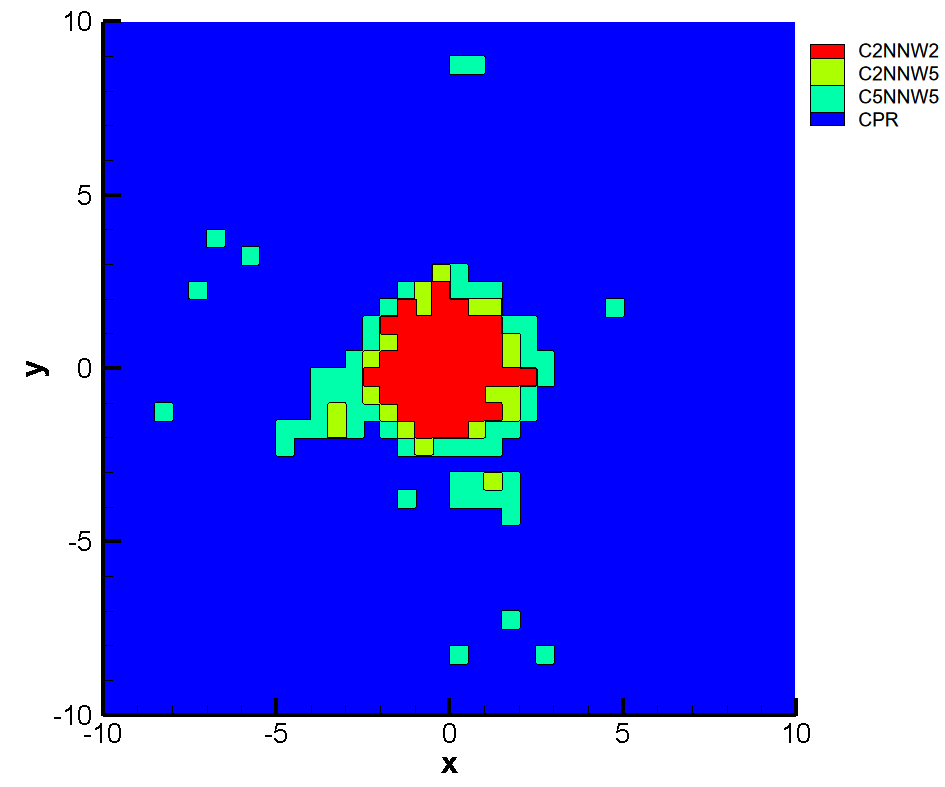}
\par\end{centering}

}\subfloat[Time evolution of the error $\left\langle \rho\right\rangle $.]{\begin{centering}
\includegraphics[width=0.45\textwidth]{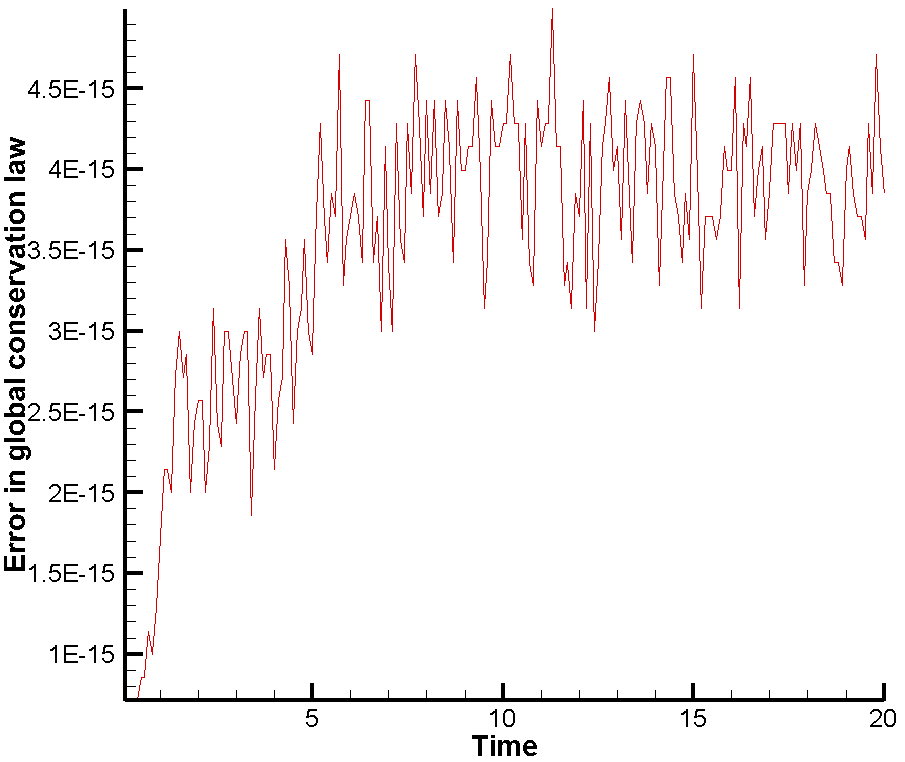}
\par\end{centering}

}
\par\end{centering}

\caption{HCCS(1,1,1,1) with $\mathbf{dv}=(c_{1},0.000005,0.00001)$ for solving
Euler vortex problem ($40\times40$ grid, $DOFs=200\times200$, $dt=0.001$).\label{fig:HS(1,1,1,1)-evolution}}
\end{figure}

\par\end{center}

\begin{center}
\begin{table}
\begin{centering}
\begin{tabular}{ccccccc}
\hline 
{\small Schemes} & {\small CPR} & {\small C5NNW5} & {\small C2NNW5} & {\small C2NNW2} & {\small HCCS(1,1,1,1)} & {\small HCCS(1,0,0,1)}\tabularnewline
\hline 
{\small $\left\langle \rho\right\rangle _{max}$} & 4.42E-15 & 3.57E-15 & 3.71E-15 & 4.42E-15 & 5.00E-15 & 5.00E-15\tabularnewline
\hline 
\end{tabular}
\par\end{centering}

\noindent \centering{}\caption{The maximum of $\left\langle \rho\right\rangle $ on time samples
$t_{i}\in\left\{ 0.1i|i=1,2,\cdots,200\right\} $ for all schemes
($40\times40$ grid, $dt=0.001$, $T=20$)\label{tab:errors of invariant}}
\end{table}

\par\end{center}

\subsection{2D Riemann problem\label{sub:2D-Riemann-problem}}

The CNNW and CPR-CNNW schemes are applied to solve 2D Riemann problem
proposed by Schulz-Rinne \citet{Schulz1993}. The computational domain
$\left[0,1\right]\times[0,1]$ is divided into four quadrants by two
lines $x=0.8$ and $y=0.8$ and the initial constant states on the
four quadrants are 
\begin{eqnarray*}
(\rho,u,\upsilon,p) & = & \begin{cases}
\mathbf{V_{\mathit{\mathrm{1}}}=}(\rho_{1},u_{1},\upsilon_{1},p_{1}), & 0.8\leq x\leq1.0,\,\,0.8\leq y\leq1.0,\\
\mathbf{V}_{2}=(\rho_{2},u_{2},\upsilon_{2},p_{2}), & 0.0\leq x<0.8,\,\,0.8\leq y\leq1.0,\\
\mathbf{V}_{3}=(\rho_{3},u_{3},\upsilon_{3},p_{3}), & 0.0\leq x<0.8,\,\,0.0\leq y<0.8,\\
\mathbf{V}_{4}=(\rho_{4},u_{4},\upsilon_{4},p_{4}), & 0.8<x\leq1.0,\,\,0.0\leq y<0.8.
\end{cases}
\end{eqnarray*}

\begin{figure}
\subfloat[{\footnotesize C5NNW5 }]{\begin{centering}
\includegraphics[width=0.45\textwidth]{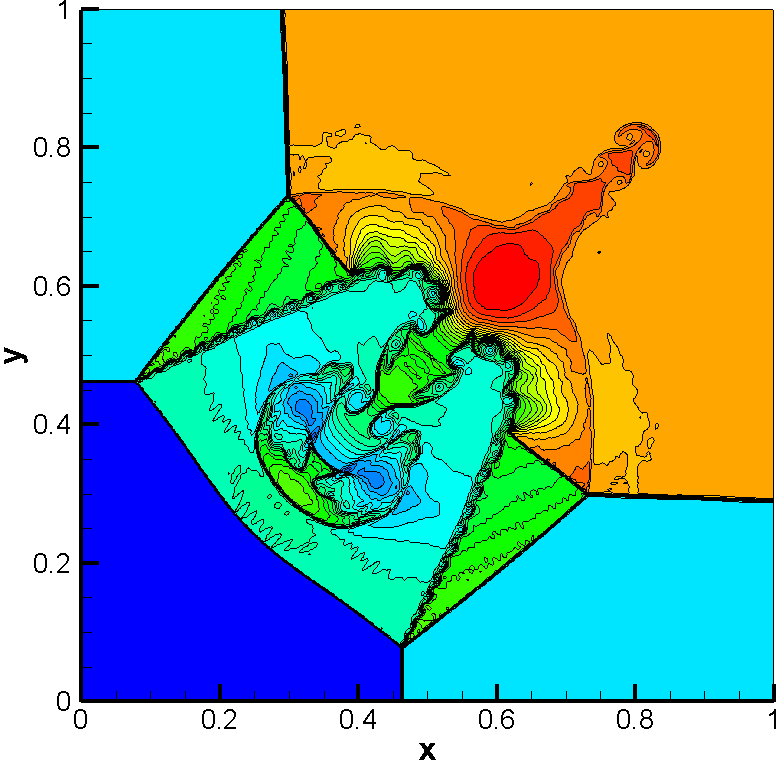}
\par\end{centering}

}\subfloat[{\footnotesize WCNS5}]{\begin{centering}
\includegraphics[width=0.45\textwidth]{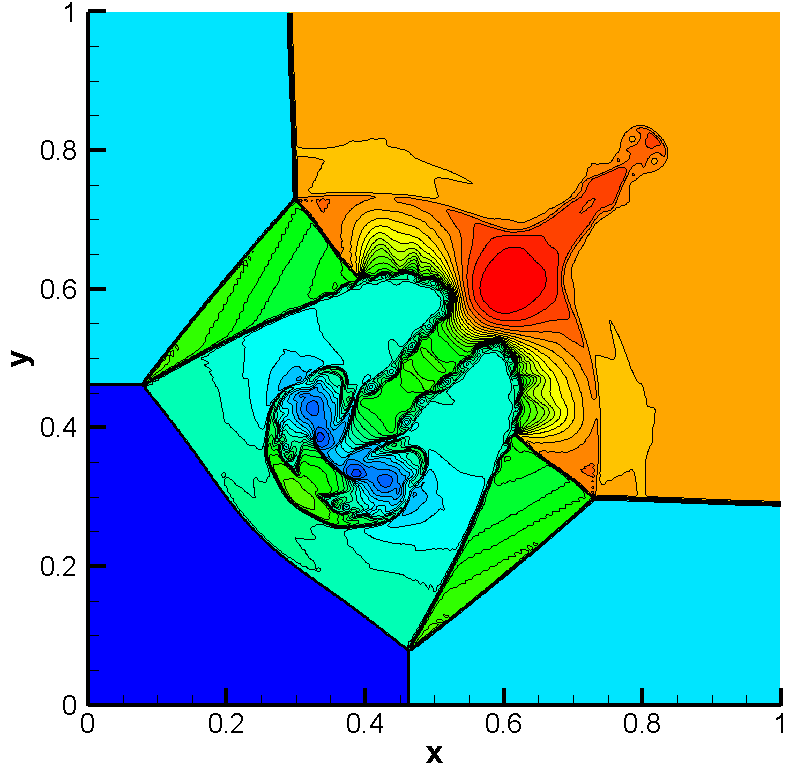}
\par\end{centering}

}

\subfloat[{\footnotesize C2NNW5}]{\begin{centering}
\includegraphics[width=0.45\textwidth]{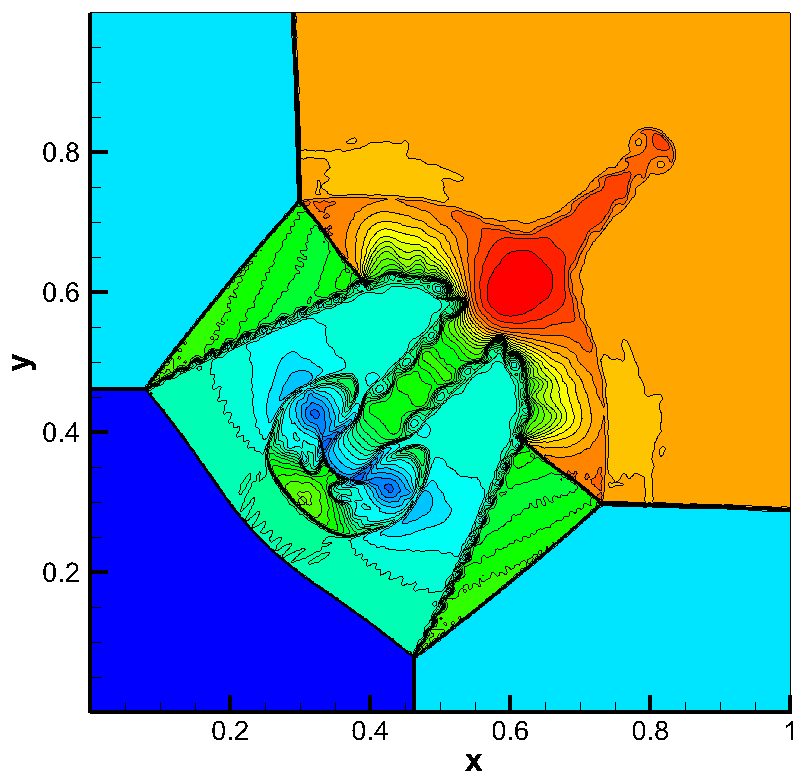}
\par\end{centering}

}\subfloat[{\footnotesize C2NNW2}]{\begin{centering}
\includegraphics[width=0.45\textwidth]{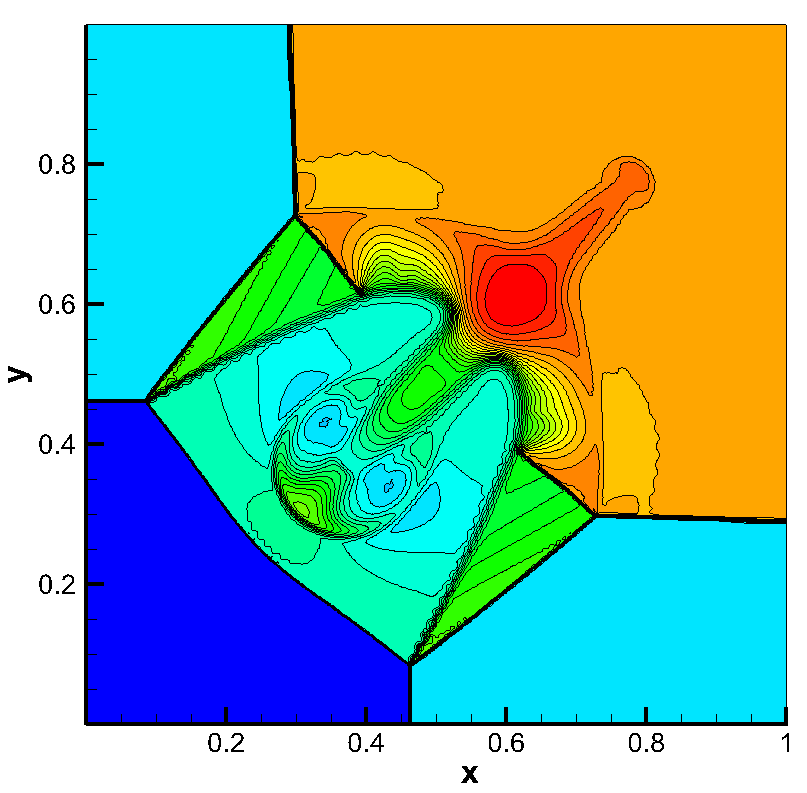}
\par\end{centering}

}

\caption{Density from 0.2 to 1.7 with 31 contours computed by C5NNW5, WCNS5,
C2NNW5, C2NNW2 ($120\times120$ grid, $DOFs=600\times600$, $T=0.8$)\label{fig:CNNW-WCNS}}
\end{figure}

Firstly, we test shock capturing properties for CNNW in solving the
2D Riemann problem with the initial constant states

\begin{equation}
\begin{array}{cccc}
 & \mathbf{V}_{1}=(1.500,0,0,1.500),\,\, &  & \mathbf{V}_{2}=(0.5323,1.206,0,0.3),\\
 & \mathbf{V}_{3}=(0.138,1.206,1.206,0.029),\,\, &  & \mathbf{V}_{4}=(0.5323,0,1.206,0.3),
\end{array}\label{eq:2D-Riemann-ref}
\end{equation}
till $T=0.8$. From Fig. \ref{fig:CNNW-WCNS}, we can see that C5NNW5
has higher resolution than WCNS5 in capturing the vortices along the
shear layers. For second-order schemes, C2NNW5 captures more small
scale features than C2NNW2, which illustrates that the former has
higher resolution than the latter.

\begin{figure}
\begin{centering}
\subfloat[Density contours, $a=0.5$]{\includegraphics[width=0.45\textwidth]{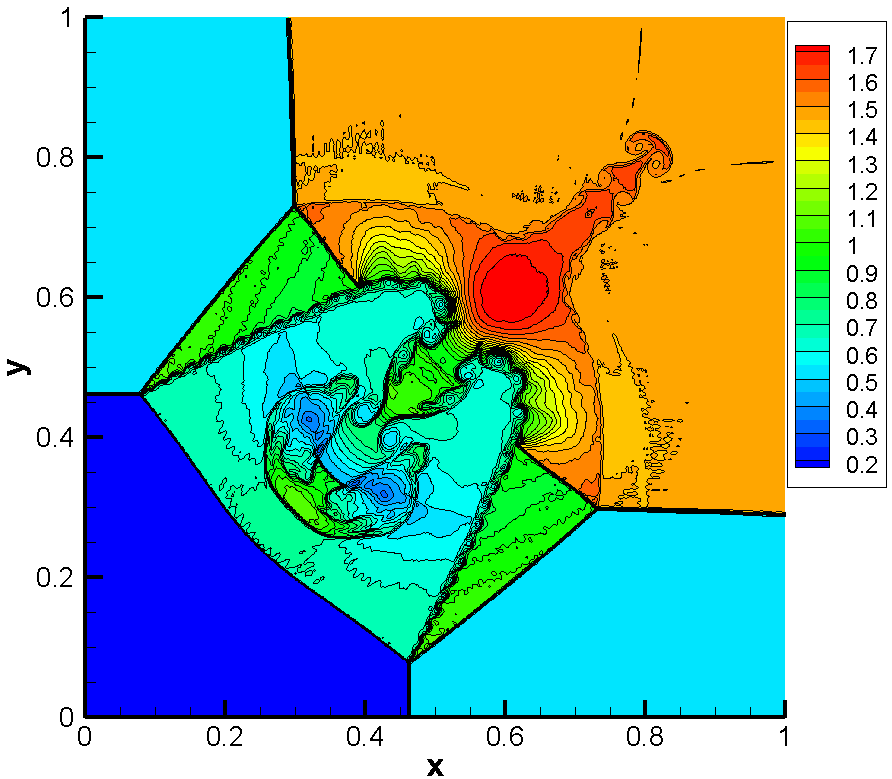}

}\subfloat[Distribution of different cells, $a=0.5$]{\includegraphics[width=0.45\textwidth]{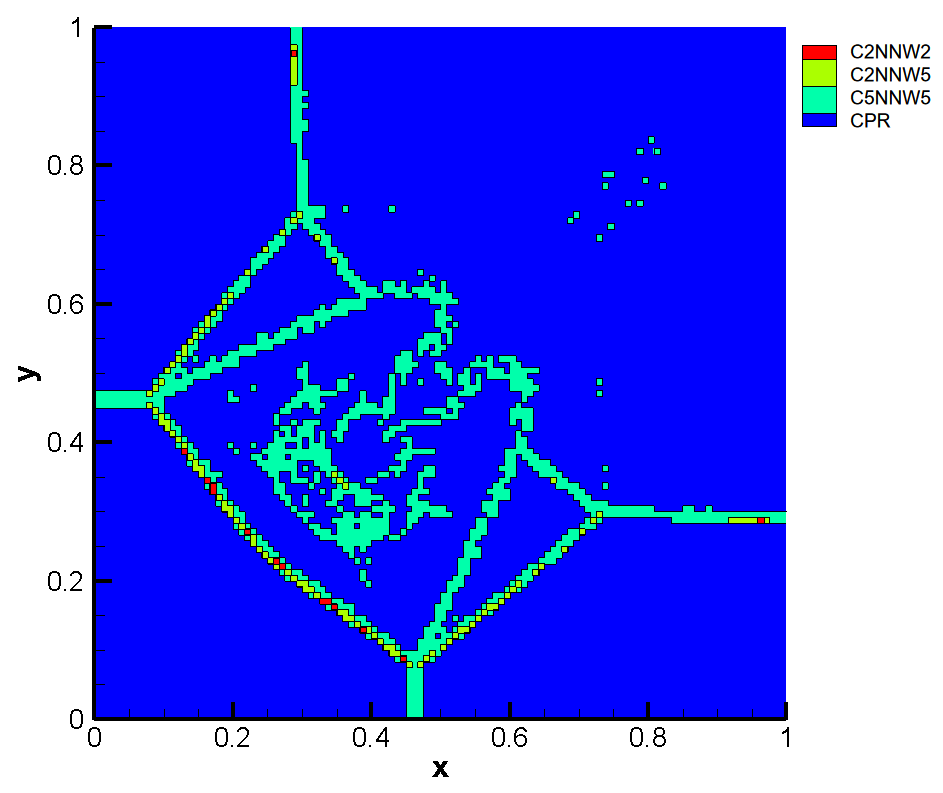}

}
\par\end{centering}

\begin{centering}
\subfloat[Density contours, $a=0.05$]{\includegraphics[width=0.45\textwidth]{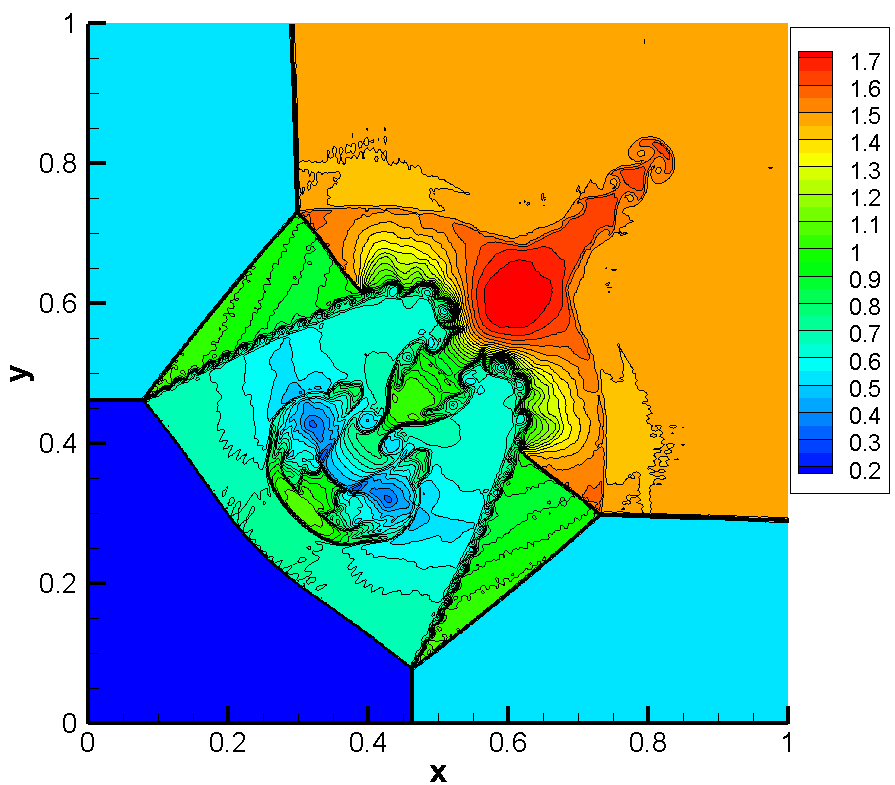}

}\subfloat[Distribution of different cells, $a=0.05$]{\includegraphics[width=0.45\textwidth]{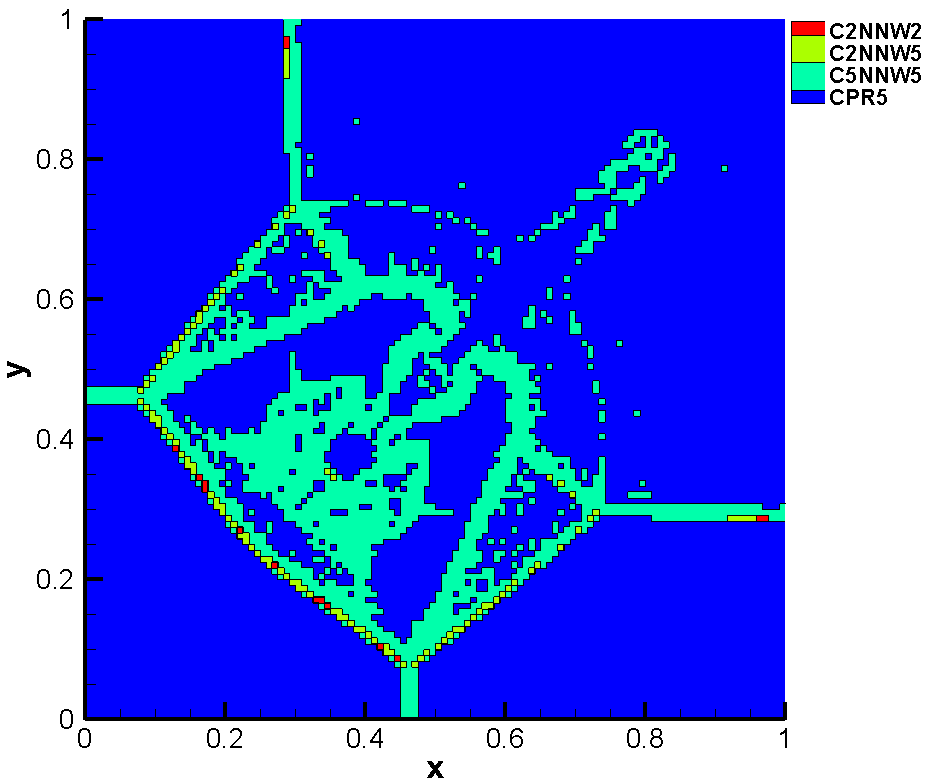}

}
\par\end{centering}

\begin{centering}
\subfloat[Density contours, $a=0.005$]{\includegraphics[width=0.45\textwidth]{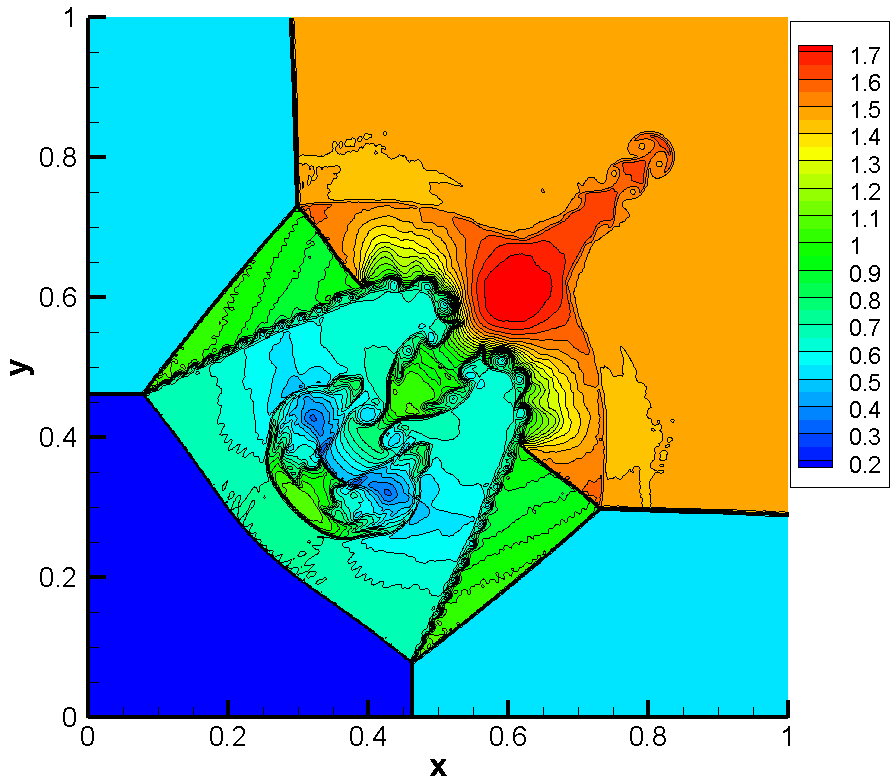}

}\subfloat[Distribution of different cells, $a=0.005$]{\includegraphics[width=0.45\textwidth]{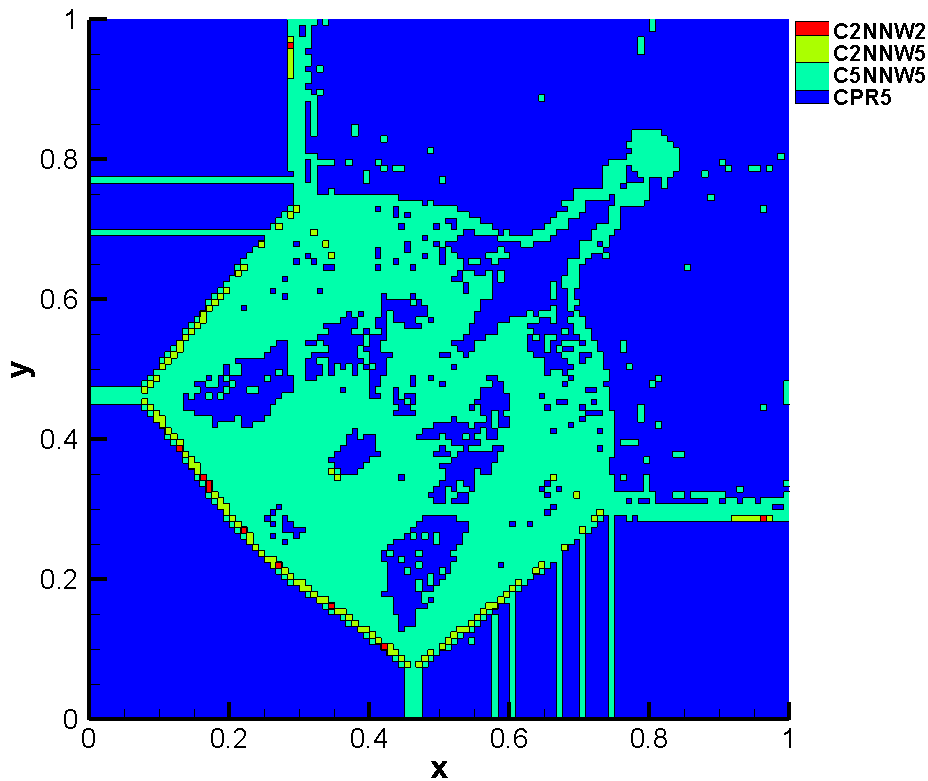}

}
\par\end{centering}

\caption{2D Riemann problem solved by HCCS(1,1,1,1) with $\mathbf{dv}=(c(a),0.3,0.6)$.
Density from 0.2 to 1.7 with 31 contours and distribution of different
cells. ($120\times120$ grid, $DOFs=600\times600$, $T=0.8$)\label{fig:HS1111}}
\end{figure}

\begin{figure}
\begin{centering}
\subfloat[Density contours, $a=0.5$]{\includegraphics[width=0.45\textwidth]{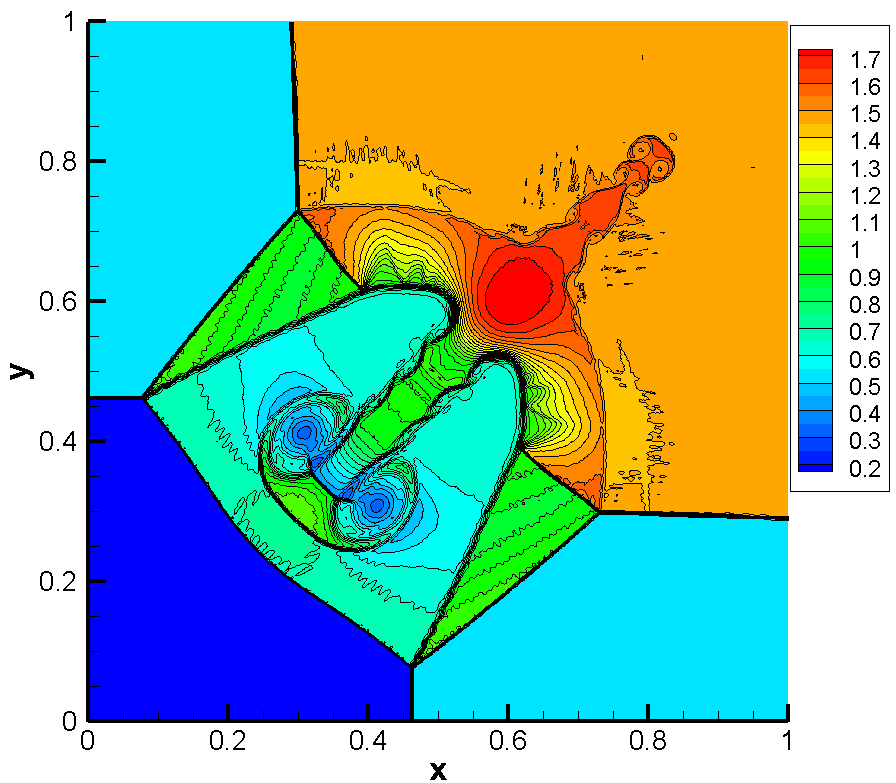}

}\subfloat[Distribution of different cells, $a=0.5$]{\includegraphics[width=0.45\textwidth]{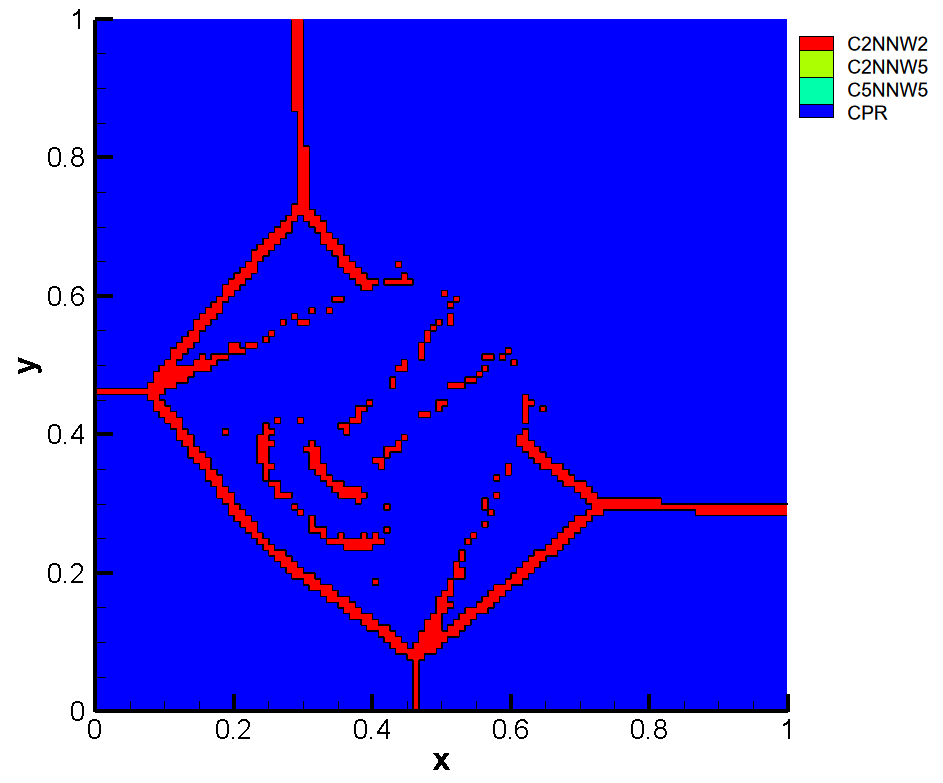}

}
\par\end{centering}

\begin{centering}
\subfloat[Density contours, $a=0.05$]{\includegraphics[width=0.45\textwidth]{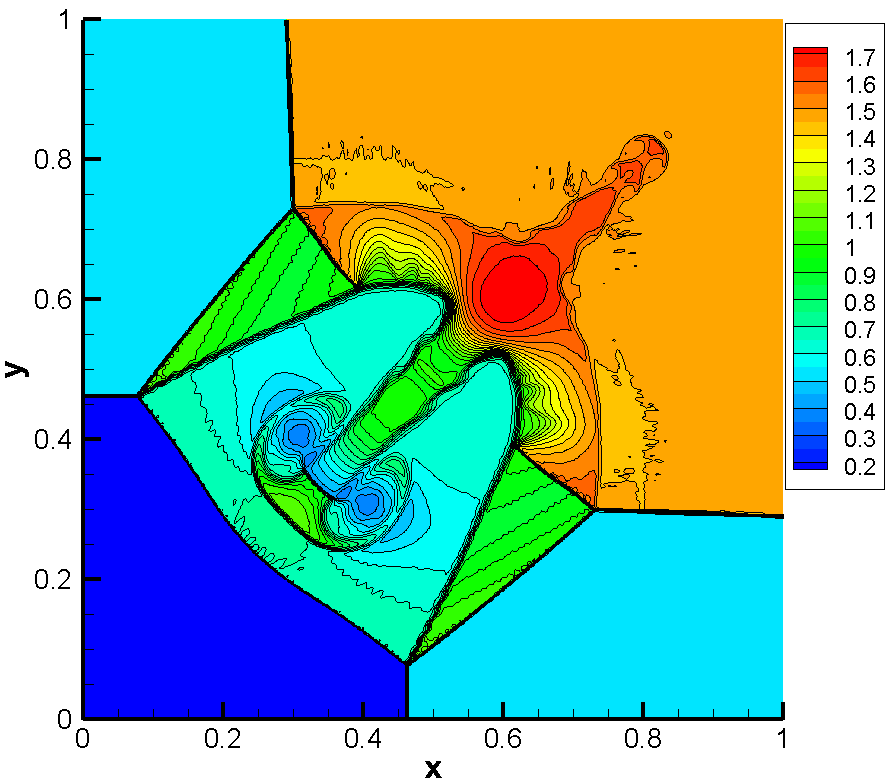}

}\subfloat[Distribution of different cells, $a=0.05$]{\includegraphics[width=0.45\textwidth]{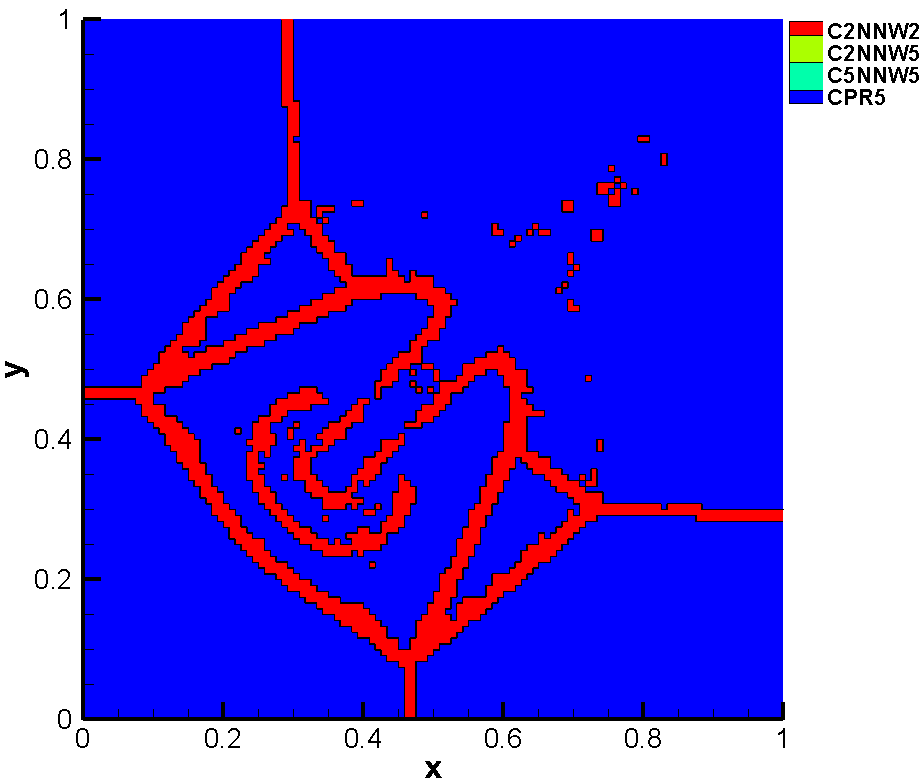}

}
\par\end{centering}

\begin{centering}
\subfloat[Density contours, $a=0.005$]{\includegraphics[width=0.45\textwidth]{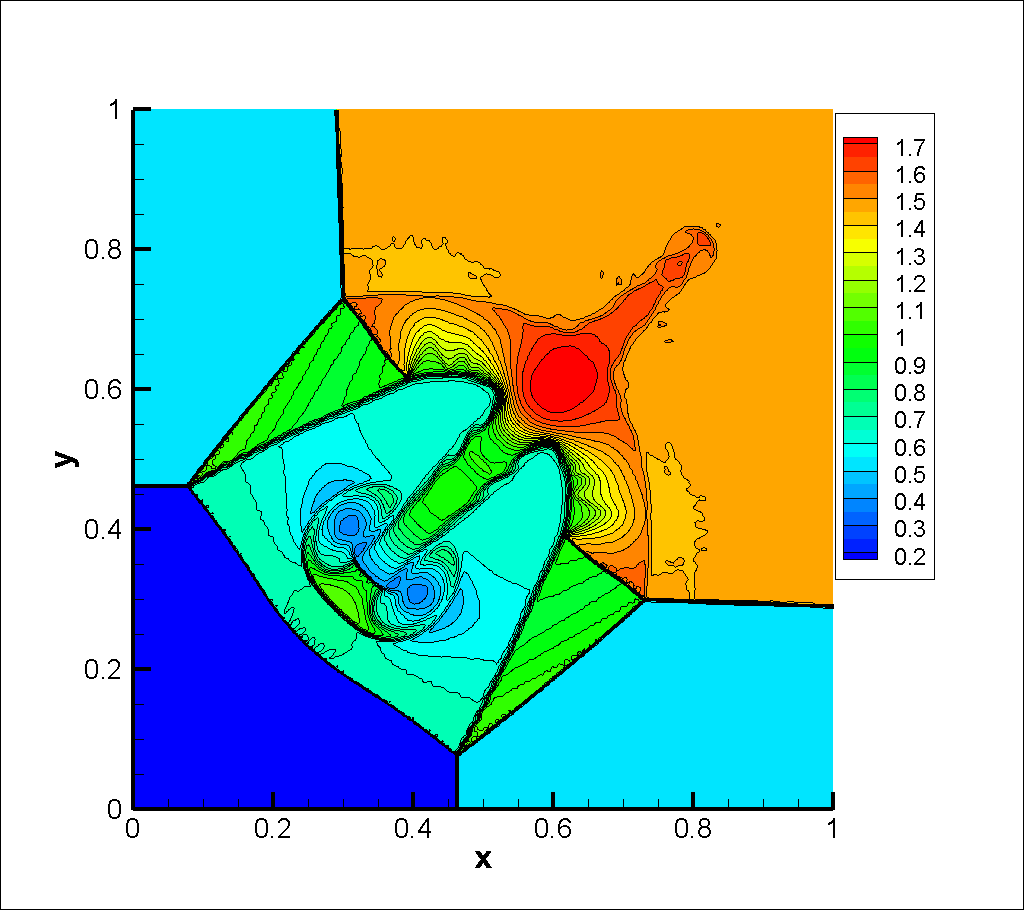}

}\subfloat[Distribution of different cells, $a=0.005$]{\includegraphics[width=0.45\textwidth]{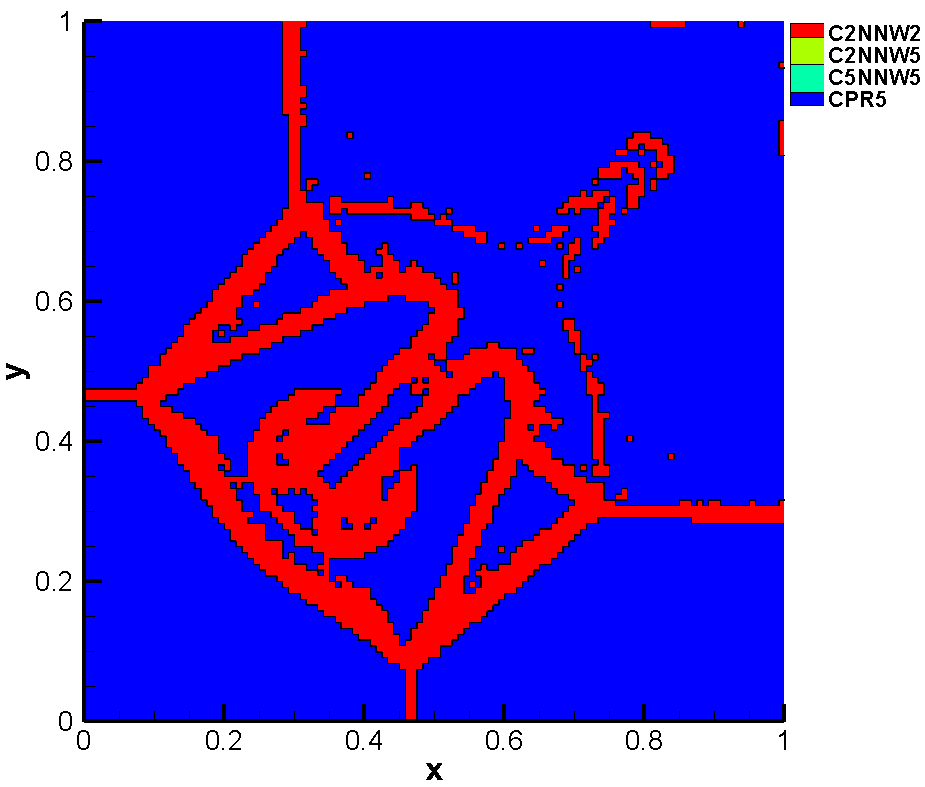}

}
\par\end{centering}

\caption{2D Riemann problem solved by HCCS(1,0,0,1) with $\mathbf{dv}=(c(a),c(a),c(a))$.
Density from 0.2 to 1.7 with 31 contours and distribution of different
cells ($120\times120$ grid, $DOFs=600\times600$, $T=0.8$)\label{fig:HS1001}}
\end{figure}

Secondly, HCCS(1,1,1,1) with $\mathbf{dv}=(c(a),0.3,0.6)$ and HCCS(1,0,0,1)
with $\mathbf{dv}=(c(a),c(a),c(a))$ are adopted to solve this problem.
The results in Fig. \ref{fig:HS1111} and Fig. \ref{fig:HS1001} show
that both of HCCS(1,1,1,1) and HCCS(1,0,0,1) can capture shocks effectively
and the former which contain C5NNW5 and C2NNW5 captures more flow
structures than the latter. In addition, HCCS(1,1,1,1) can keep high
resolution when increasing the number of troubled cells as shown in
Fig. \ref{fig:HS1111}. HCCS(1,0,0,1) can not acquire small-scale
features since when the shear wave or vortex structures are indicated
wrongly, they are limited by the dissipative C2NNW2, which smears
out the small scale vortex structures.

At last, we change shock strength of the Riemann problem to test shock
capturing ability of CNNW and CPR-CNNW. The initial constant states
are

\begin{eqnarray*}
\begin{cases}
\rho_{3} & =0.138,\\
u_{3} & =M\sqrt{\frac{\gamma p_{3}}{\rho_{3}}}/\sqrt{2},\\
\upsilon_{3} & =u_{3},\\
p_{3} & =0.029,
\end{cases} &  & \begin{cases}
s_{23} & =\frac{3-\gamma}{4}\upsilon_{3}-\sqrt{\frac{\gamma p_{3}}{\rho_{3}}+\left[\frac{\left(\gamma+1\right)}{4}\upsilon_{3}\right]^{2}},\\
\rho_{2} & =-\frac{\rho_{3}\upsilon_{3}}{s_{23}}+\rho_{3},\\
u_{2} & =u_{3},\\
\upsilon_{2} & =0,\\
p_{2} & =(\rho_{3}\upsilon_{3}^{2}+p_{3})-s_{23}\rho_{3}\upsilon_{3},
\end{cases}
\end{eqnarray*}
\begin{eqnarray*}
\begin{cases}
\rho_{4} & =\rho_{2},\\
u_{4} & =0,\\
\upsilon_{4} & =u_{2},\\
p_{4} & =p_{2},
\end{cases} &  & \begin{cases}
s_{12} & =\frac{3-\gamma}{4}u_{2}-\sqrt{\frac{\gamma p_{2}}{\rho_{2}}+\left[\frac{\left(\gamma+1\right)}{4}u_{2}\right]^{2}},\\
\rho_{1} & =-\frac{\rho_{2}u_{2}}{s_{12}}+\rho_{2},\\
u_{1} & =0,\\
\upsilon_{1} & =0,\\
p_{1} & =(\rho_{2}u_{2}^{2}+p_{2})-s_{12}\rho_{2}u_{2}.
\end{cases}
\end{eqnarray*}
The shock Mach $M=\sqrt{2}u_{3}/\sqrt{\gamma p_{3}/\rho_{3}}$ here
is an adjustable parameter to control the shock strength and the states
$\mathbf{V}_{2}$,$\mathbf{V}_{4}$ and $\mathbf{V}_{1}$ are determined
based on the Rankine-Hugoniot condition. Please note that the initial
conditions of \ref{eq:2D-Riemann-ref} is a special case of $M=3.145$
with possible rounding errors. In our tests, $M$ is changed gradually
from $3.145$ to $10^{6}$ to test whether a numerical scheme can
capture the shocks without blowing up. The largest $M$ for different
schemes without blowing up are summarized in Table \ref{tab:2d-Riemann-Ma}. 

{\small In the C5NNW5, nonlinear mechanism is introduced in the interpolation
process (from solution points to flux points). However, the }compact
flux difference operator (C5) still relies on an uniform polynomial{\small{}
assumption, which is known to result in Gibbs phenomenon near discontinuities.
The resulting C5NNW5 }can calculate problems with $M<6$. Even if
we introduce constant interpolation (NNW1), which was shown to be
very dissipative in the communities of finite difference and finite
volume studies, such as \citet{Leer1974,Buffard2010,zhu2016}, the
resulting C5NNW1 still blows up at about $M=8.2$. This indicates
that it maybe impossible to construct a high-order FE schemes capable
of simulating very strong shocks by introducing nonlinear mechanism
to only parts of the operators that are based on uniform polynomial
assumption, which is fundamental in high-order FE schemes like CPR
and DG. In another word, a possible reason for the difficulty in constructing
robust shock capturing high-order FE methods using the artificial
viscosity \citet{Persson2006,Yu2020} or limiting solution distribution
\citet{Zhu2009,Zhong2013,Zhu2013,Du2015,Park2016,Li2020} is that
the nonlinear mechanism only breaks the uniform polynomial assumption
in parts of the high-order FE methods. However, the remaining operators
based on uniform polynomial still generates oscillations and easily
leads to blow up.

The C2NNW5 scheme, which breaks the uniform polynomial distribution
assumption for both interpolation and flux difference operators using
NNW5 and C2 respectively, successfully captures shocks with $M=10^{6}$
and a much more robust shock capturing ability is obtained. We can
see that second-order schemes C2NNW5 and C2NNW2 can compute stronger
shocks than high-order schemes C5NNW5. Based on this observation,
it is better for hybrid schemes to contain low-order robust shock
capturing schemes to capture strong shocks. 

In our subcell limiting strategy, we can divide the high-order CPR
cells into subcells computed by second-order schemes C2NNW5 and C2NNW2
for strong shocks, and therefore break the uniform polynomial distribution
assumption for every operators. From Table \ref{tab:2d-Riemann-Ma},
we can see that the HCCS(1,1,1,1) with $\mathbf{dv}=(c_{0},0.005,0.01)$
and HCCS(1,0,0,1) with $\mathbf{dv}=(c_{0},c_{0},c_{0})$ can calculate
this problem with $M=10^{6}$, which illustrate that the proposed
CPR-CNNW have good properties in capturing strong shocks.

\begin{flushleft}
\begin{table}
\begin{centering}
\begin{tabular}{c|ccc|c||cc||c|c||c||cc||c}
\hline 
\multirow{2}{*}{{\small Schemes}} & \multicolumn{3}{c|}{{\small Low-order schemes}} & \multicolumn{4}{c|}{{\small Fifth-order schemes}} & \multicolumn{5}{c}{{\small CPR-CNNW}}\tabularnewline
\cline{2-13} 
 & {\small C5NNW1} & {\small C2NNW5} & {\small C2NNW2} & \multicolumn{2}{c}{{\small C5NNW5}} & \multicolumn{2}{c|}{{\small WCNS5}} & \multicolumn{3}{c}{{\small HCCS(1,1,1,1)}} & \multicolumn{2}{c}{{\small HCCS(1,0,0,1)}}\tabularnewline
\hline 
{\small $M_{max}$} & 8.2 & $10^{6}$ & $10^{6}$ & \multicolumn{2}{c}{{\small 6.0}} & \multicolumn{2}{c|}{{\small 183}} & \multicolumn{3}{c}{$10^{6}$} & \multicolumn{2}{c}{$10^{6}$}\tabularnewline
\hline 
\end{tabular}
\par\end{centering}

\noindent \centering{}\caption{Comparisons of nonlinear schemes in solving 2D Riemann problems with
different Ma.\label{tab:2d-Riemann-Ma}}
\end{table}

\par\end{flushleft}

\subsection{Double Mach reflection\label{sub:Double-Mach-reflection}}

\begin{figure}
\begin{centering}
\subfloat[HCCS(1,1,1,1), density from 1.5 to 21.7 with 31 contours.]{\begin{centering}
\includegraphics[width=0.45\textwidth]{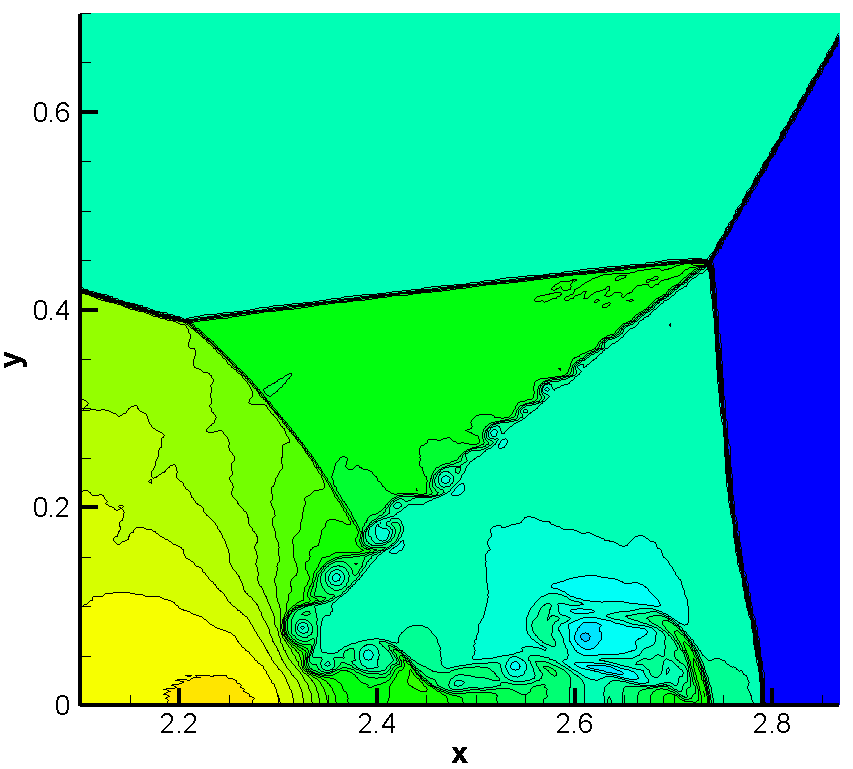}
\par\end{centering}

}\subfloat[HCCS(1,1,1,1), distribution of different cells]{\begin{centering}
\includegraphics[width=0.54\textwidth]{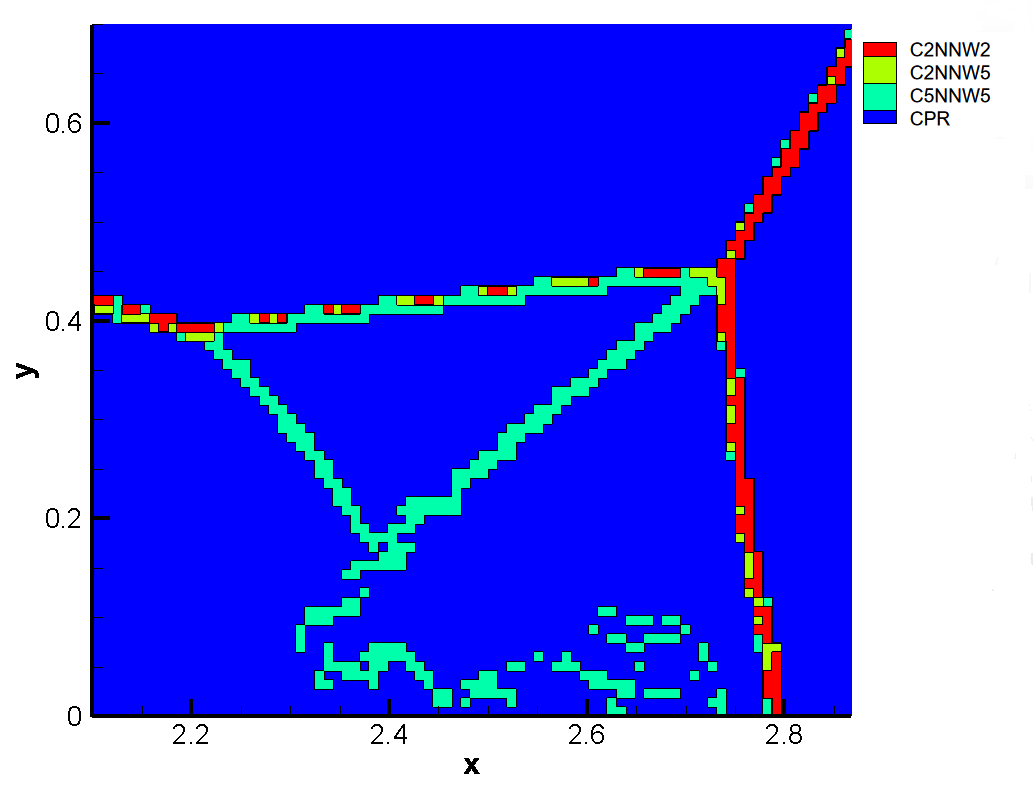}
\par\end{centering}

}
\par\end{centering}

\begin{centering}
\subfloat[HCCS(1,0,0,1), density from 1.5 to 21.7 with 31 contours.]{\begin{centering}
\includegraphics[width=0.45\textwidth]{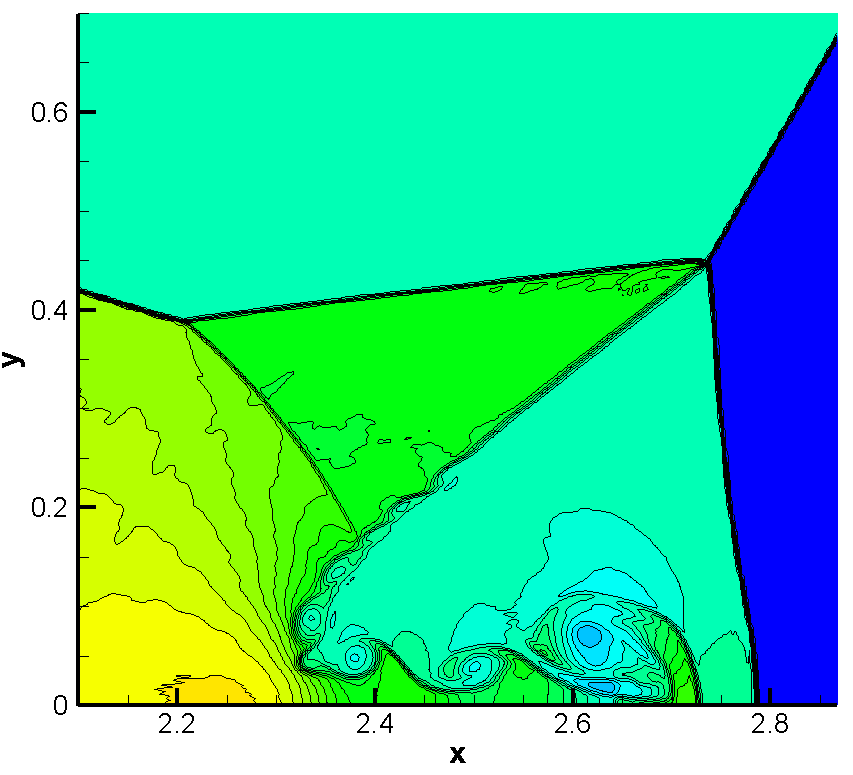}
\par\end{centering}

}\subfloat[HCCS(1,0,0,1), distribution of different cells]{\begin{centering}
\includegraphics[width=0.53\textwidth]{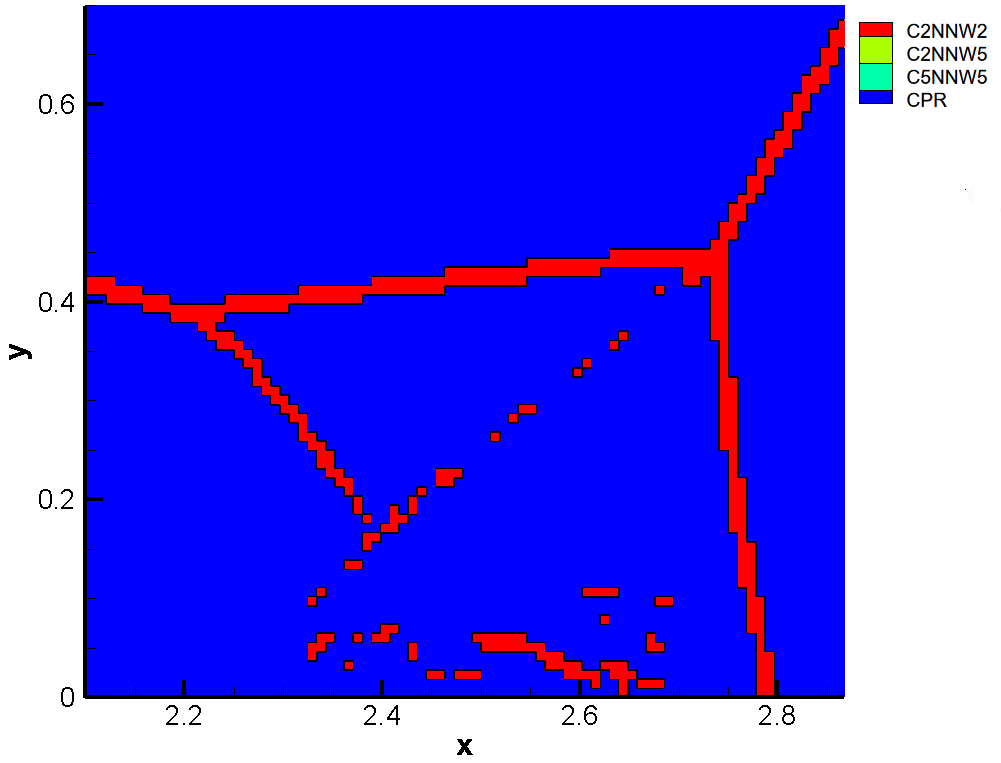}
\par\end{centering}

}
\par\end{centering}

\caption{Hybrid scheme HCCS(1,1,1,1) with $\mathbf{dv}=(c_{0},0.05,0.1)$ and
HCCS(1,0,0,1) with $\mathbf{dv}=(c_{0},c_{0},c_{0})$ in solving Double
Mach reflection problem ($h=1/108$, $DOFs=2160\times540$, $T=0.2$)\label{fig:Hybrid-schemes-tworeflection-1}}
\end{figure}

Double Mach reflection problem described in \citet{Woodward1984}
is a popular test case to test strong shock capturing ability of high-resolution
schemes. The problem is solved by CPR-CNNW schemes on a grid with
$h=1/108$. As shown in Fig. \ref{fig:Hybrid-schemes-tworeflection-1},
CPR-CNNW schemes can capture strong shock robustly. In addition, HCCS(1,1,1,1)
can capture small-scale structures better than HCCS(1,0,0,1). Comparing
Fig. \ref{fig:Hybrid-schemes-tworeflection-1}(b) and (d), we can
see HCCS(1,1,1,1) computes the troubled cells with strong shocks using
C2NNW2 while using C2NNW5 and C5NNW5 for weak discontinuities and
cells next to discontinuities. However, HCCS(1,0,0,1) compute all
these cells by the C2NNW2. These results illustrate that HCCS(1,1,1,1)
containing C5NNW5 and C2NNW5 can acquire higher resolution than HCCS(1,0,0,1). 

\subsection{Shock-vortex interaction\label{sub:Shock-vortex-interaction}}

As a final test, shock-vortex interaction is solved by CPR-CNNW schemes
to test their good balance in shock capturing and high resolution.
This test case is originally proposed by \citet{Rault2003} and is
one of the benchmark problems of the 5th International Workshop on
High Order CFD Methods (HiOCFD5). It has been adopted as a benchmark
for high-order numerical schemes \citet{Dumbser2007,Dumbser2014,You2018},
since it involves a complex flow pattern with both smooth features
and discontinuous waves. The computational domain is $[0,2]\times[0,1]$.
The initial conditions are given by a stationary normal shock wave
placed at $x=0.5$ with shock Mach number $M_{s}$ and by a vortex
placed at $(x_{c},y_{c})=(0.5,0.5)$ with the strength $M_{v}$. We
choose $M_{s}=1.5$ and $M_{v}=0.9$, which is the same as that in
\citet{Rault2003}. A computational grid with $h=1/140$ and $DoFs=1400\times700=980000$
is used. 

\begin{figure}
\begin{centering}
\subfloat[Distribution of density and different cells]{\begin{centering}
\includegraphics[width=0.85\textwidth]{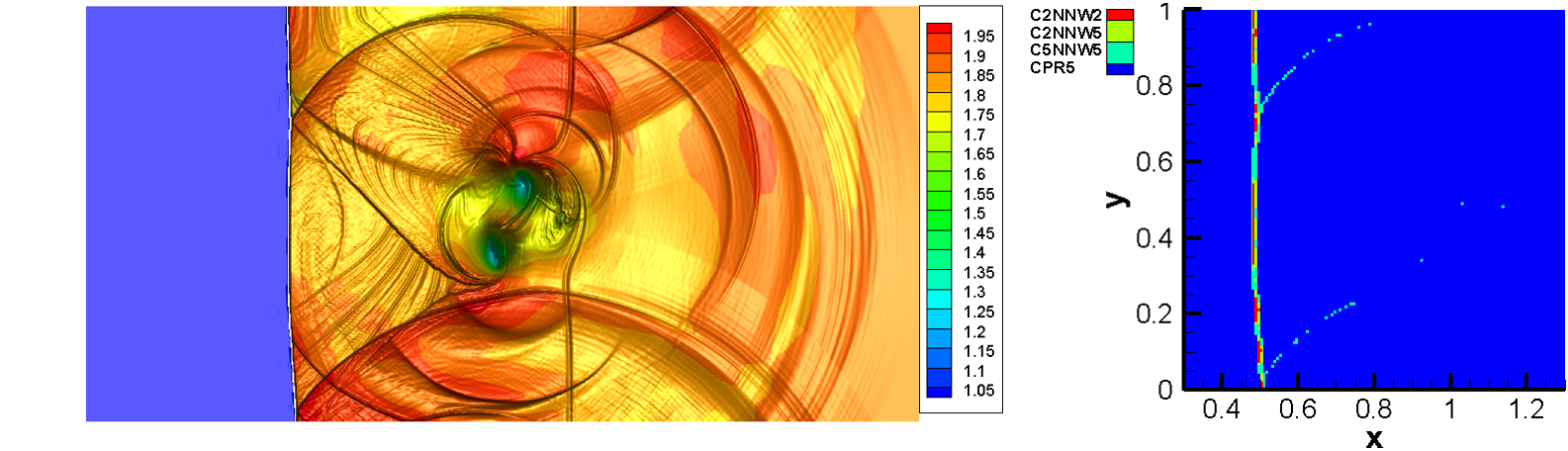}
\par\end{centering}

}
\par\end{centering}

\begin{centering}
\subfloat[Schlieren View of density and its close-up]{\begin{centering}
\includegraphics[width=0.85\textwidth]{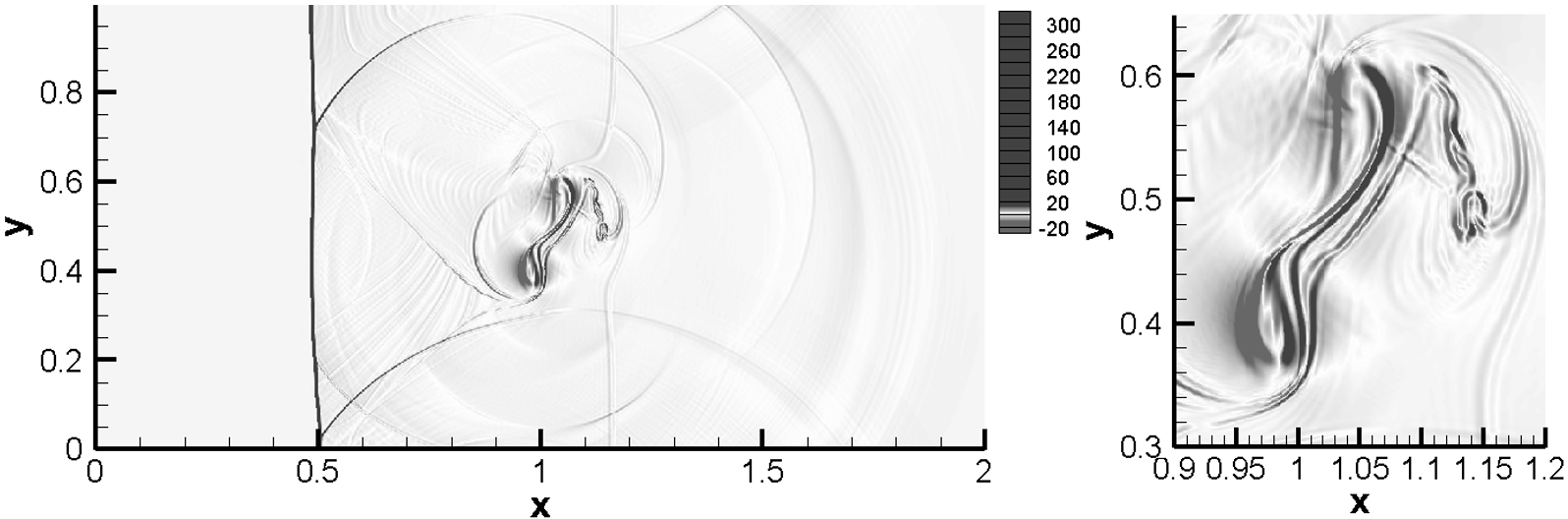}
\par\end{centering}

}
\par\end{centering}

\caption{Shock-vortex interaction problem solved by HCCS(1,1,1,1) with $\mathbf{dv}=(c_{0},0.05,0.1)$
($h=1/140$, $DoFs=1400\times700$, $T=0.7$)\label{fig:HS1111-cell}}
\end{figure}

\begin{figure}
\begin{centering}
\subfloat[Distribution of density and different cells]{\begin{centering}
\includegraphics[width=0.85\textwidth]{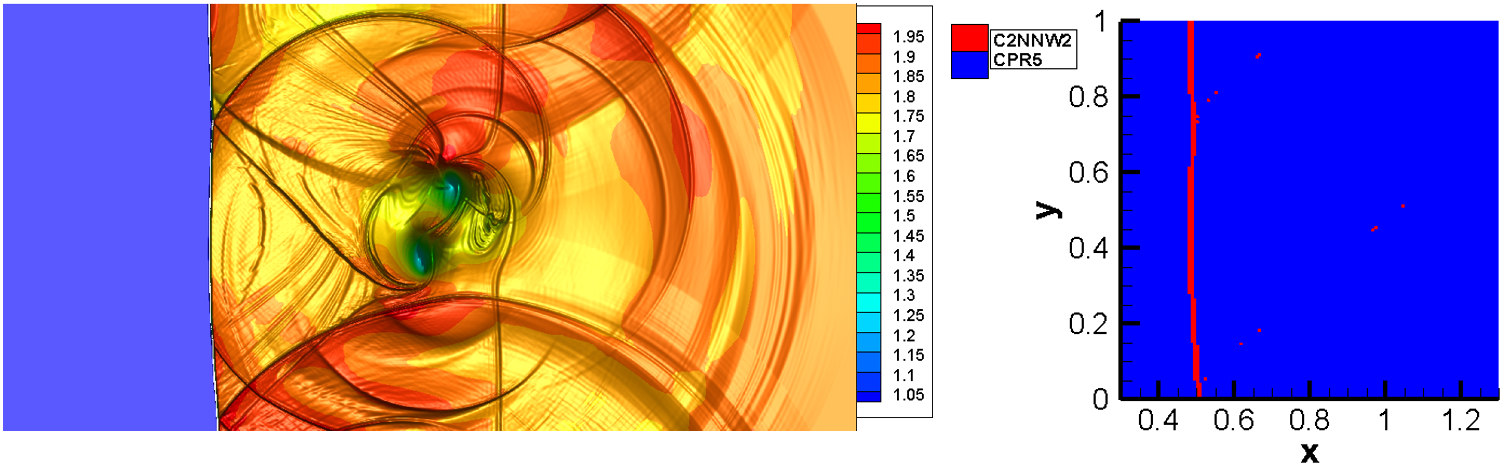}
\par\end{centering}

}
\par\end{centering}

\begin{centering}
\subfloat[Schlieren View of density and its close-up]{\begin{centering}
\includegraphics[width=0.85\textwidth]{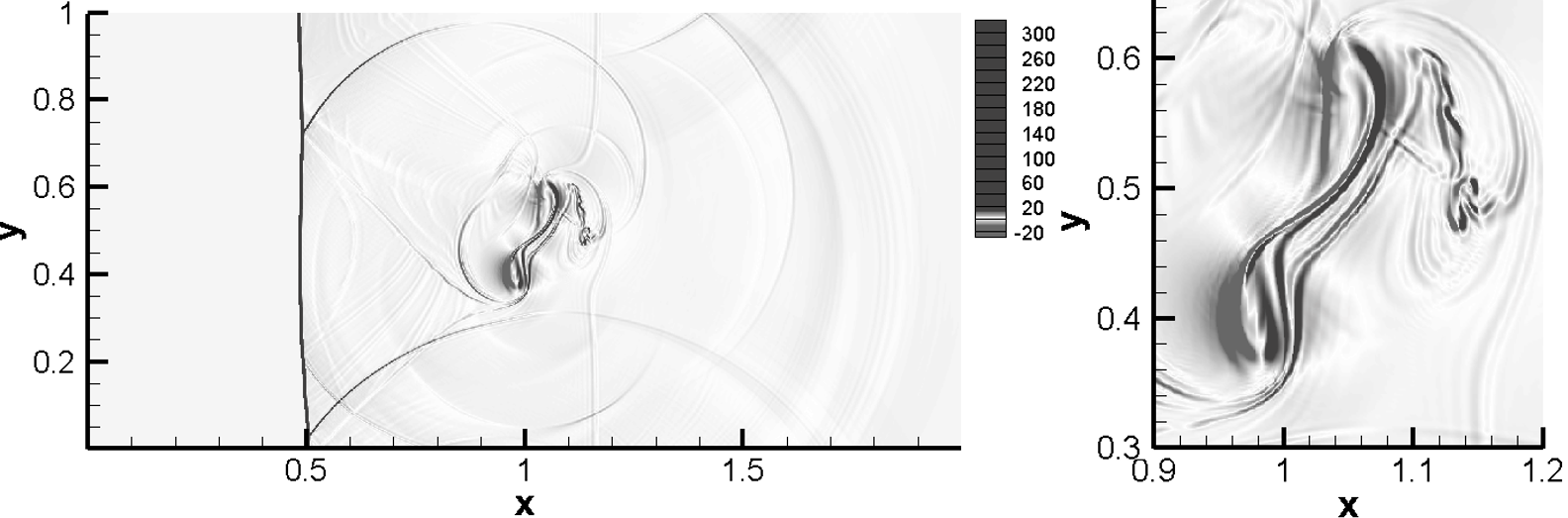}
\par\end{centering}

}
\par\end{centering}

\caption{Shock-vortex interaction problem solved by HCCS(1,0,0,1) with $\mathbf{dv}=(c_{0},c_{0},c_{0})$
($h=1/140$, $DoFs=1400\times700$, $T=0.7$)\label{fig:HS1001-cell}}
\end{figure}

\begin{figure}
\begin{centering}
\subfloat[Density along Line 1 ($y=0.4$) ]{\begin{centering}
\includegraphics[width=0.4\textwidth]{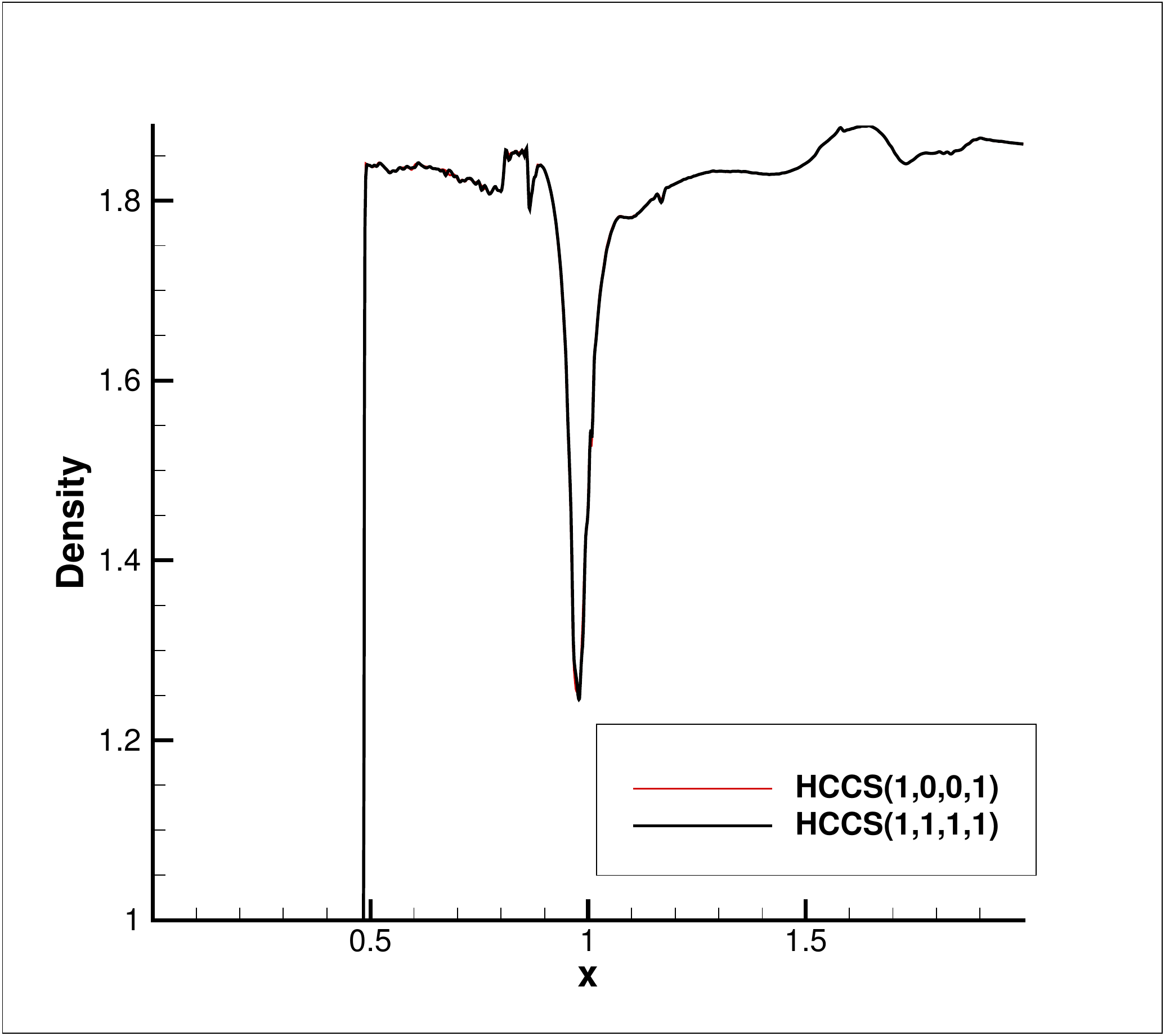}
\par\end{centering}

}\subfloat[{Density along Line 1 ($y=0.4$) near shock on the area $x\in[0.47,0.51]$. }]{\begin{centering}
\includegraphics[width=0.4\textwidth]{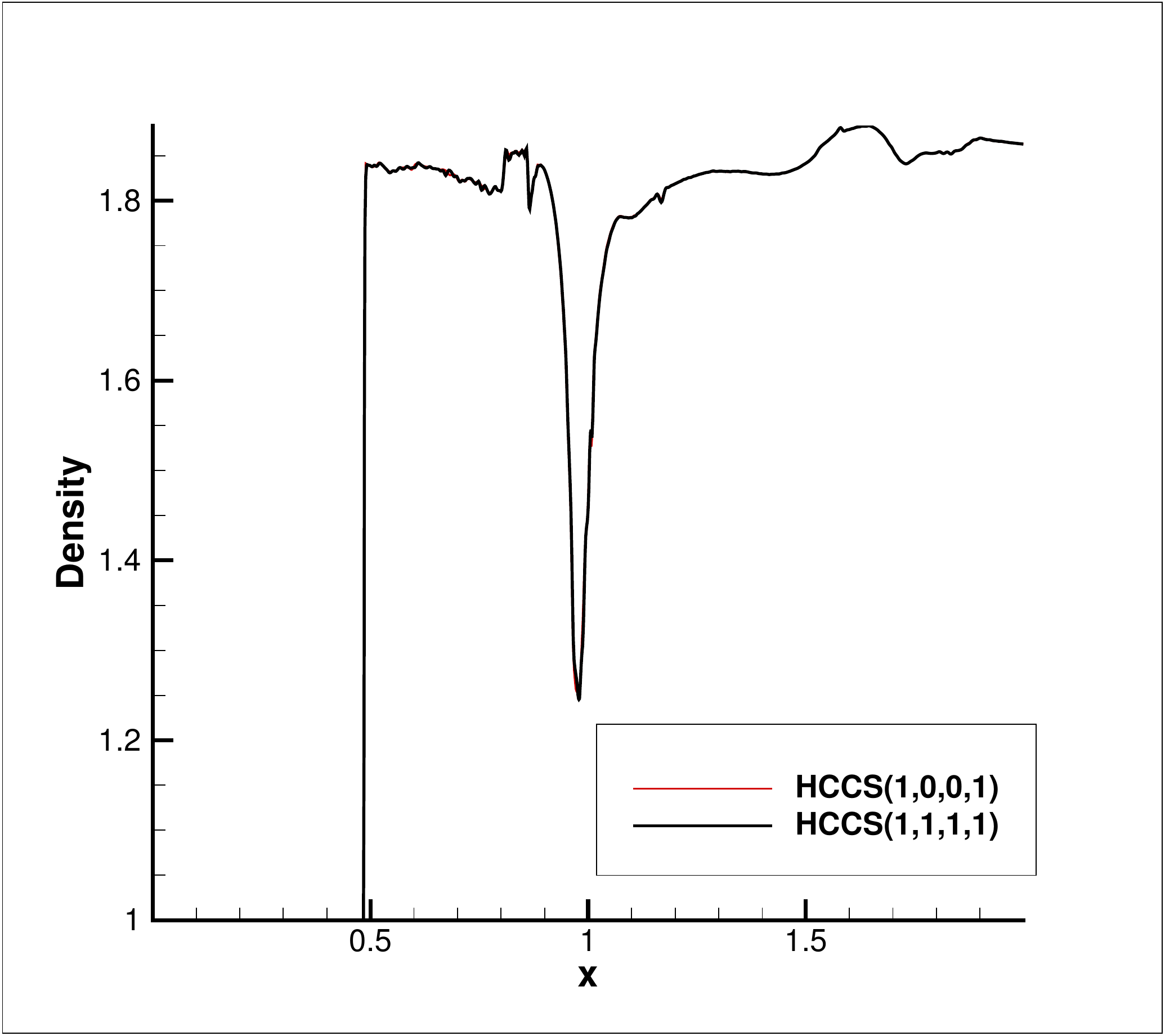}
\par\end{centering}

}
\par\end{centering}

\begin{centering}
\subfloat[Density along Line 2 ($y=0.7$)]{\begin{centering}
\includegraphics[width=0.4\textwidth]{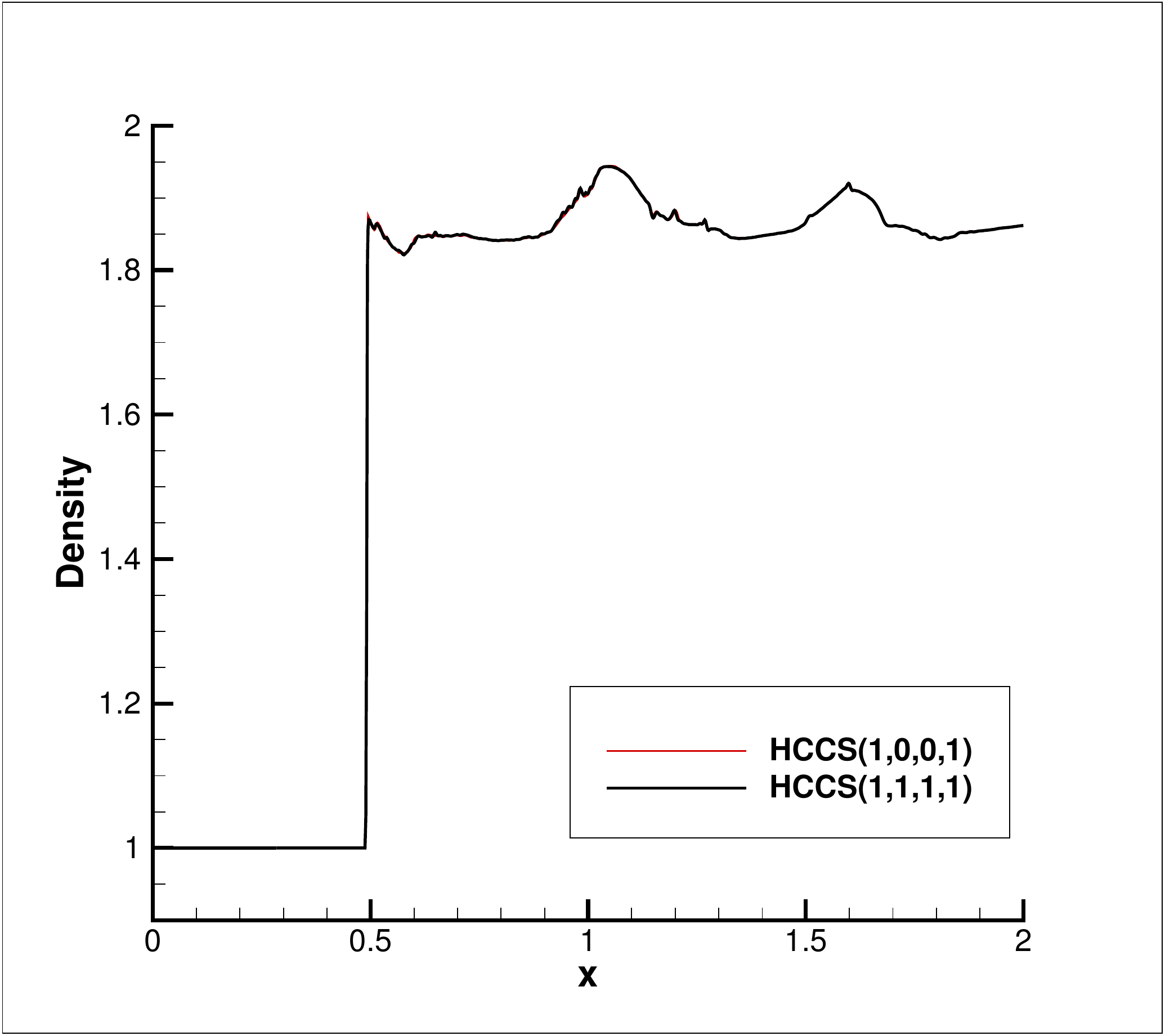}
\par\end{centering}

}\subfloat[Density along Line 3 ($x=0.52$)]{\begin{centering}
\includegraphics[width=0.4\textwidth]{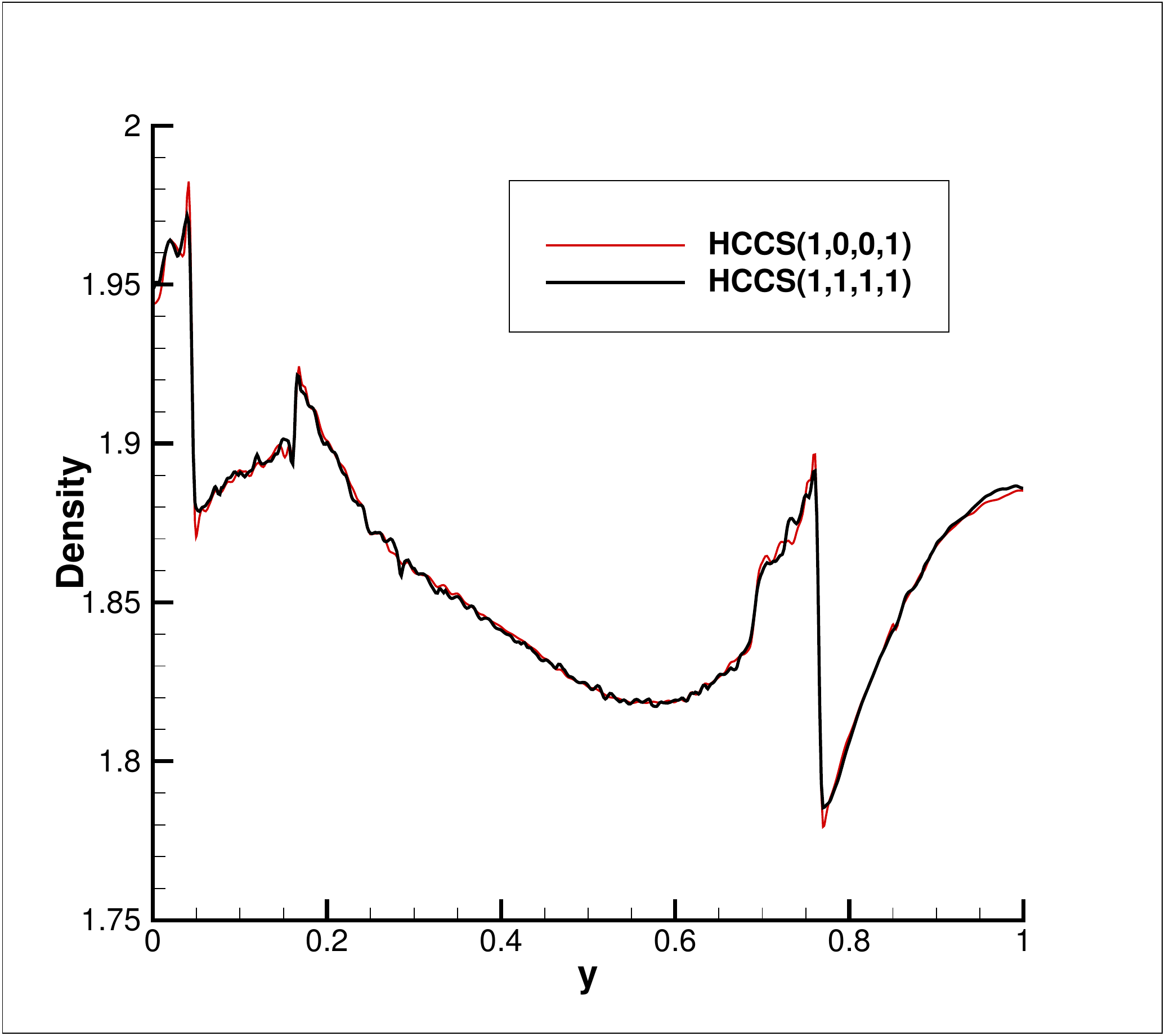}
\par\end{centering}

}
\par\end{centering}

\begin{centering}
\subfloat[Density along Line 4 ($x=1.05$)]{\begin{centering}
\includegraphics[width=0.4\textwidth]{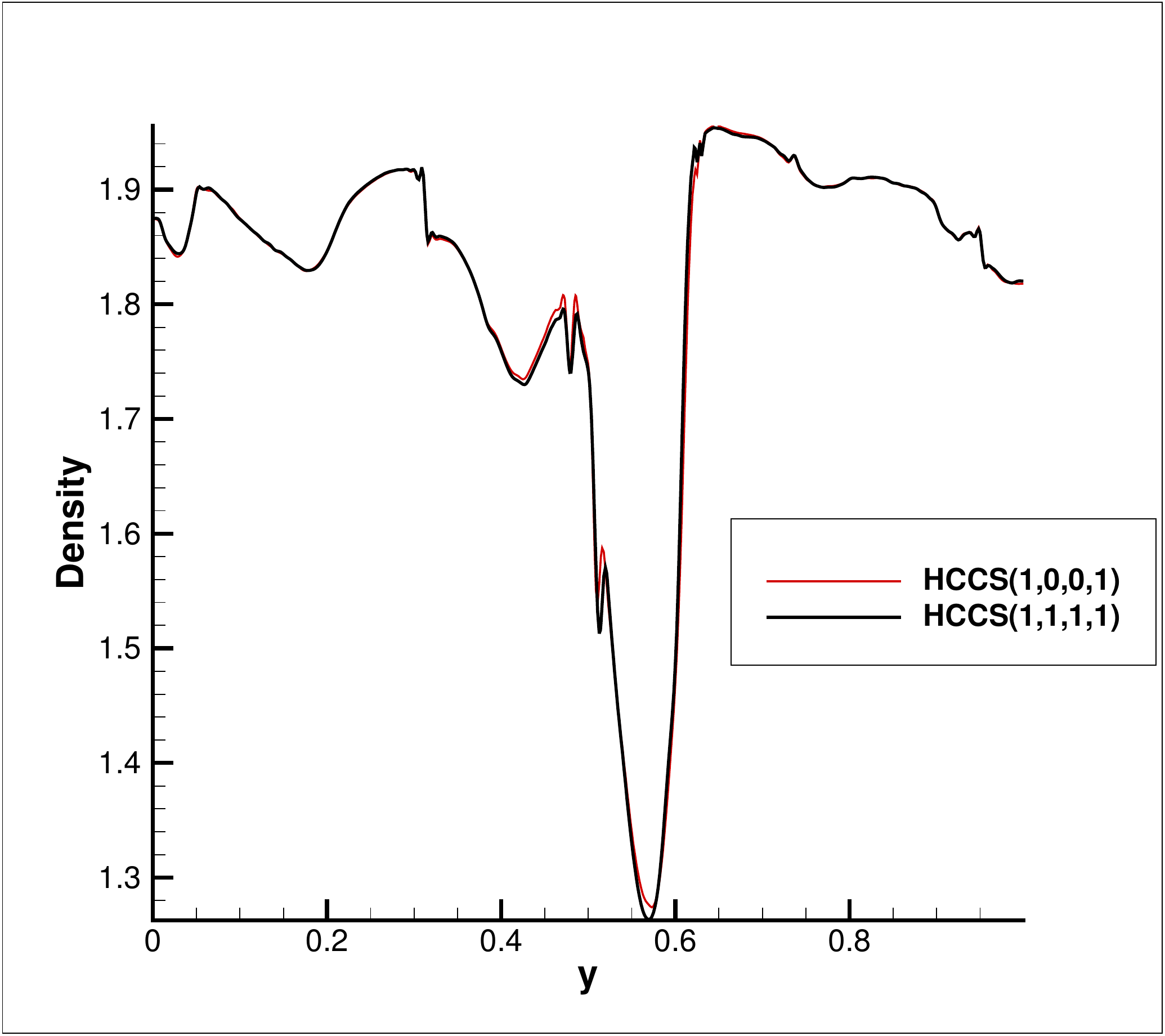}
\par\end{centering}

}\subfloat[Density along Line 5 ($y=1.65$)]{\begin{centering}
\includegraphics[width=0.4\textwidth]{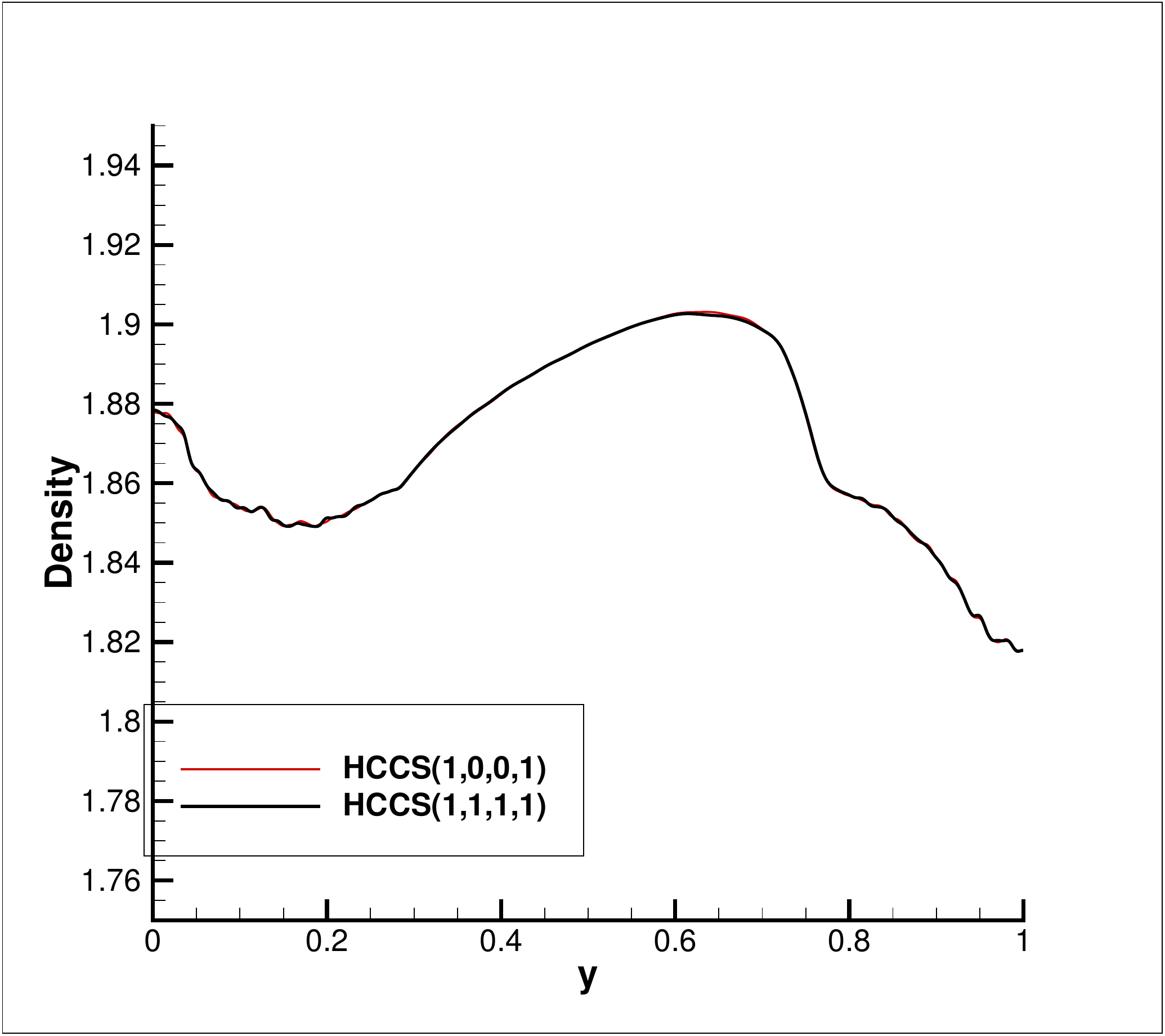}
\par\end{centering}

}
\par\end{centering}

\caption{Results along five lines of HCCS(1,1,1,1) with $\mathbf{dv}=(c_{0},0.05,0.1)$
and HCCS(1,0,0,1) with $\mathbf{dv}=(c_{0},c_{0},c_{0})$ where $a=0.5$
in solving shock-vortex interaction ($h=1/140$, $DoFs=980000$, $T=0.7$).
\label{fig:density-lines}}
\end{figure}

HCCS(1,1,1,1) and HCCS(1,0,0,1) are applied to solve this problem.
The distribution of the density, of different cells and of the density
derivatives at time $T=0.7$ are shown in Fig. \ref{fig:HS1111-cell}
for HCCS(1,1,1,1) and in Fig. \ref{fig:HS1001-cell} for HCCS(1,0,0,1).
The results show that the hybrid CPR-CNNW scheme can capture shock
wave as well as smooth vortex features. Our results are comparable
qualitatively with those computed by ADER-DG-P5 with a posteriori
ADER-WENO3 subcell limiter \citet{Dumbser2014} and with the numerical
solution provided in \citet{Rault2003}. In addition, density along
five lines are also shown in Fig. \ref{fig:density-lines}. From Fig.
\ref{fig:density-lines}(a) (b) (c), we can see that there are no
obvious oscillations near the stationary shock wave at $x=0.5$. Vortical
structures are also well simulated, as shown in Fig. \ref{fig:density-lines}(a)
(e). Large scale flow structures are well resolved although there
are some small shock-driven oscillations, as shown in \ref{fig:density-lines}(c)
(d) (f). Density along Line 4 in Fig. \ref{fig:density-lines}(e)
also shows that HCCS(1,1,1,1) can acquire better resolution in computing
the vortical structure. 

\begin{figure}
\begin{centering}
\subfloat[Distribution of different cells, HCCS(1,1,1,1)]{\begin{centering}
\includegraphics[width=0.6\textwidth]{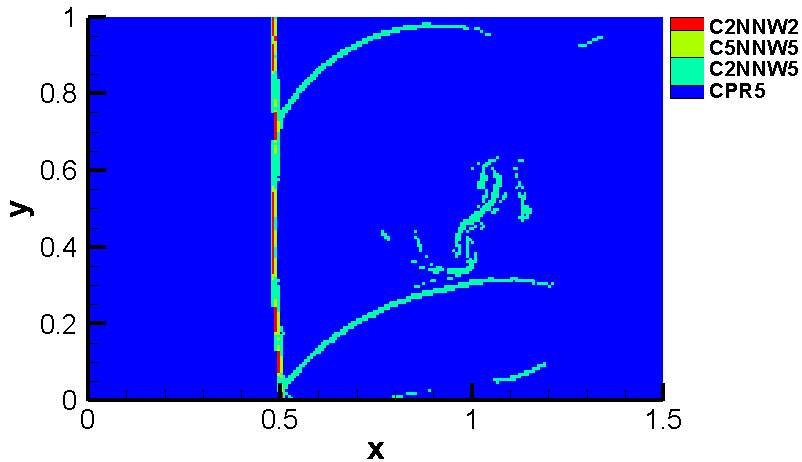}
\par\end{centering}

}\subfloat[Schlieren View, HCCS(1,1,1,1)]{\begin{centering}
\includegraphics[width=0.36\textwidth]{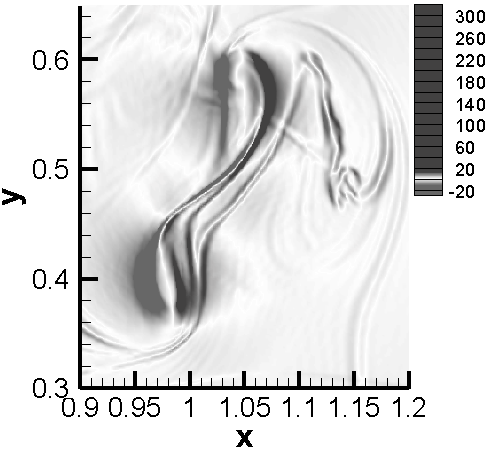}
\par\end{centering}

}
\par\end{centering}

\begin{centering}
\subfloat[Distribution of different cells, HCCS(1,0,0,1)]{\begin{centering}
\includegraphics[width=0.6\textwidth]{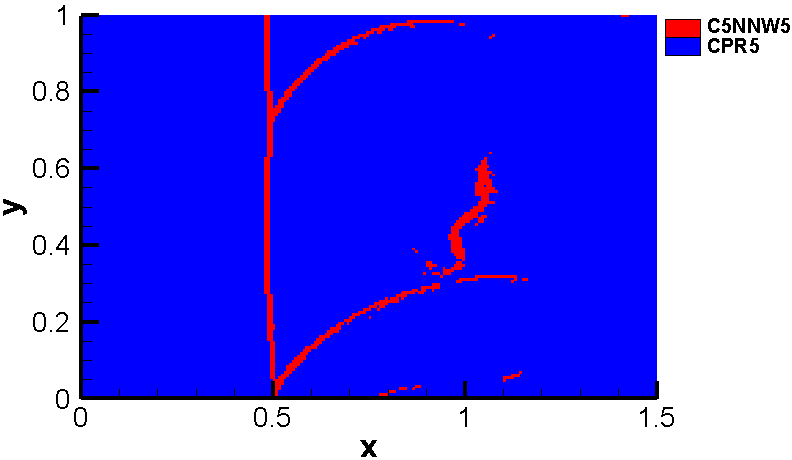}
\par\end{centering}

}\subfloat[Schlieren View, HCCS(1,0,0,1)]{\begin{centering}
\includegraphics[width=0.36\textwidth]{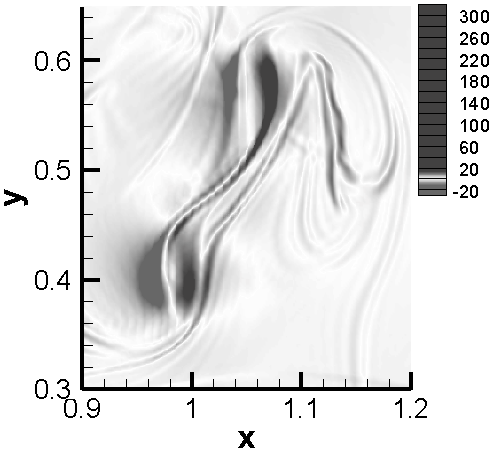}
\par\end{centering}

}
\par\end{centering}

\begin{centering}
\subfloat[Comparisons on Line 1.]{\begin{centering}
\includegraphics[width=0.48\textwidth]{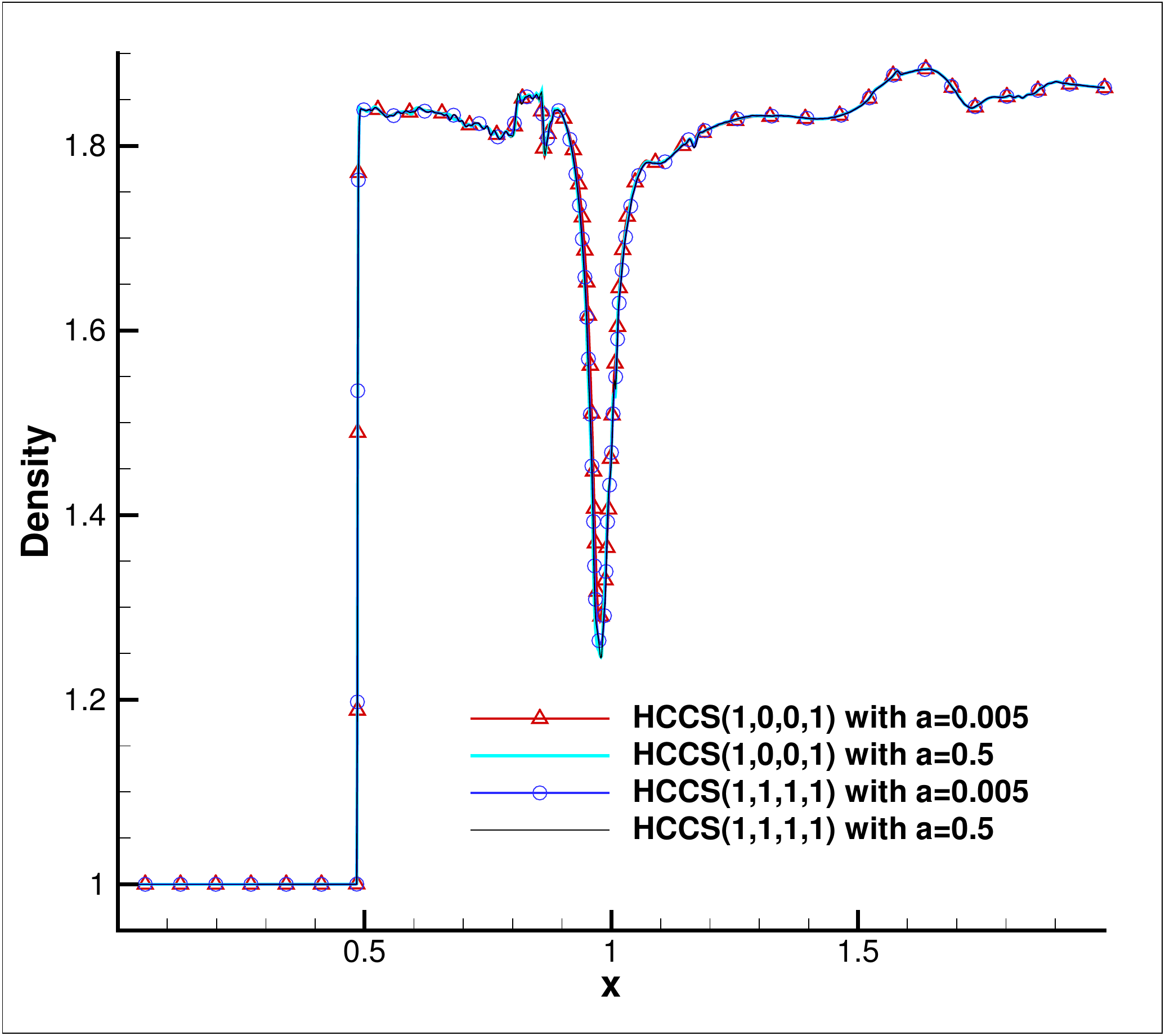}
\par\end{centering}

}\subfloat[Comparisons on Line 4]{\begin{centering}
\includegraphics[width=0.48\textwidth]{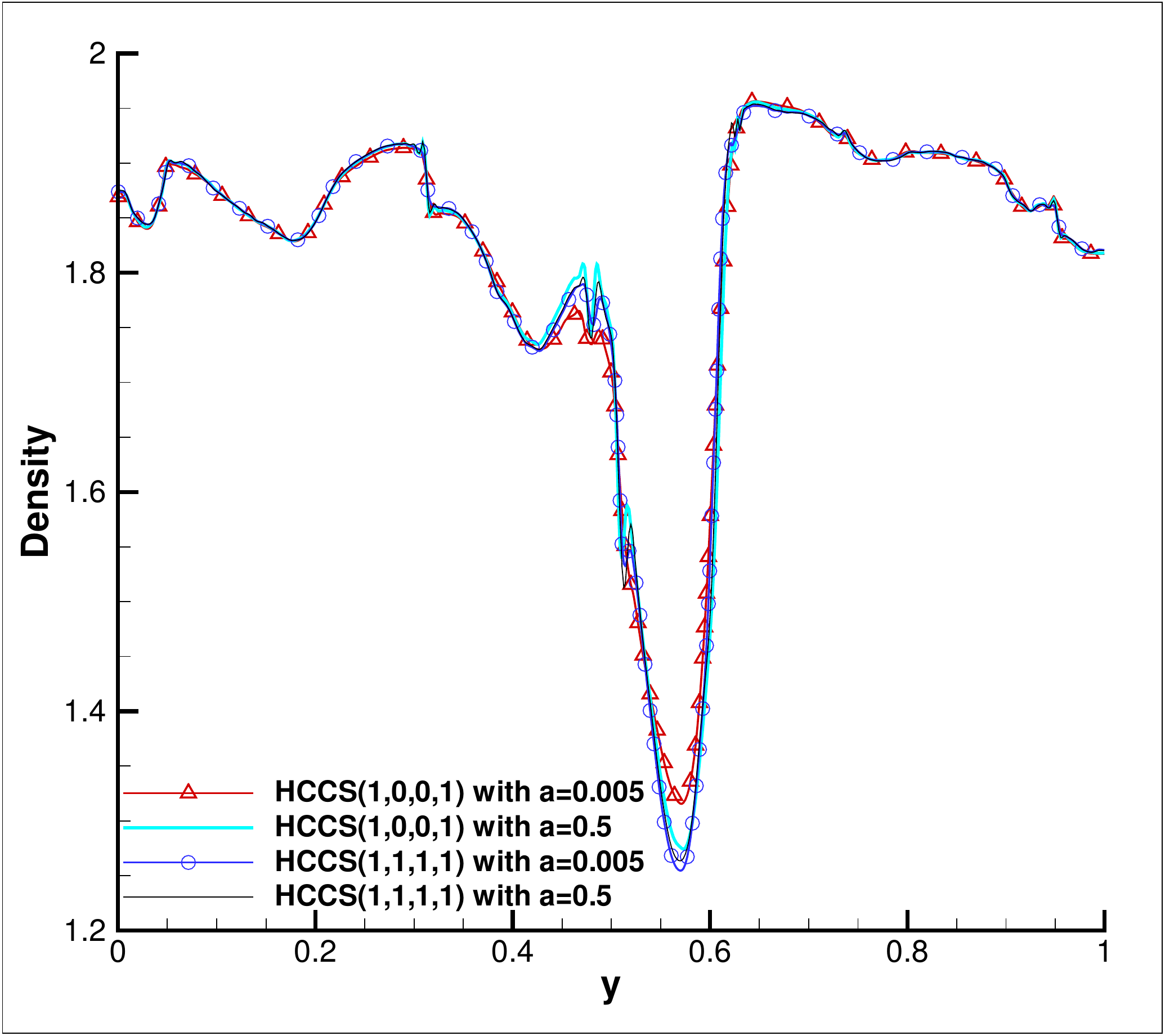}
\par\end{centering}

}
\par\end{centering}

\caption{Shock-vortex interaction problem solved by HCCS(1,1,1,1) with $\mathbf{dv}=(c(a),0.05,0.1)$
and HCCS(1,0,0,1) with $\mathbf{dv}=(c(a),c(a),c(a))$ where $a=0.005$.
($h=1/140$, $DoFs=1400\times700$, $T=0.7$)\label{fig:comparison-a0p005}}
\end{figure}

To view the influence of troubled cell detecting, the results of HCCS(1,1,1,1)
and HCCS(1,0,0,1) with parameter $a=0.005$ are shown in Fig. \ref{fig:comparison-a0p005}.
Under this parameter the vortex area is detected as troubled cells,
as shown in Fig. \ref{fig:comparison-a0p005}(a)(c). We can see that
the HCCS(1,1,1,1) can compute the vortical structure better than HCCS(1,0,0,1),
as shown in Fig. \ref{fig:comparison-a0p005}(b) (d). Moreover, density
along Line 1 and Line 4 are compared in Fig. \ref{fig:comparison-a0p005}(e)
(f). We can see that HCCS(1,1,1,1) can keep high resolution when increasing
the number of troubled cells (by increasing $a$) while the resolution
of HCCS(1,0,0,1) decreases significantly. 

\section{Concluding Remarks}

In this paper, both high-order and second-order shock capturing schemes
based on nonuniform nonlinear weighted interpolation are proposed
and these schemes are applied in subcell limiting for the high-order
CPR method. Eigenvalues of the spatial discretization matrix of the
proposed high-order CNNW are proved to be a collection of eigenvalues
of local matrices. All eigenvalues are computed and compared with
CPR and WCNS. The results show that the proposed high-order CNNW schemes
are stable and have similar spectral properties as those of WCNS.
Then, a priori subcell CNNW limiting approach is developed for the
fifth-order CPR, resulting in a special hybrid scheme and is called
CPR-CNNW. To ensure robustness and accuracy of the hybrid scheme,
CNNW schemes with varying accuracy orders are chosen adaptively according
to the magnitude of troubled cell indicator. The proposed CNNW and
CPR-CNNW schemes are applied to solve linear advection equation and
Euler equations. Numerical investigations show that the proposed C5NNW5
have fifth-order of accuracy. C2NNW5 and C2NNW2 has second-order of
accuracy and the former has higher resolution than the later one.
In addition, C2NNW5 and C2NNW2 are more robust in shock capturing
than fifth-order scheme C5NNW5. It is shown that the CPR with subcell
$p$-adaptive CNNW limiting has higher resolution than subcell second-order
C2NNW2 limiting, which illustrate that high-order interpolation in
subcells can improve resolution. CPR with subcell $p$-adaptive CNNW
limiting has good balance in high resolution and shock capturing.
The scheme combines the advantages of high-order CPR schemes in smooth
regions, with the robustness and accuracy of $p$-adaptive CNNW for
shock capturing. Both analytical and numerical results show that the
CPR-CNNW schemes satisfy discrete conservation law. The proposed hybrid
CPR-CNNW has some merits in less data exchanging for physical variables,
in satisfying discrete conservation laws, and in good balance between
high resolution and good shock capturing robustness. The proposed
subcell limiting approach could be generalized to CPR on unstructured
meshes by making some changes in interpolation procedure. The limiting
approach in theory can be applied to DG or other kind of FE method.

\section*{Acknowledgments}

This study was supported by National Numerical Wind-tunnel project,
National Natural Science Foundation of China (Grant Nos.12172375,
11902344, 11572342), the foundation of State Key Laboratory of Aerodynamics
(Grant No. SKLA2019010101).

\section*{Appendix A. NNW5 interpolation}

Here we give fifth-order nonuniform nonlinear weighted interpolation
in one-dimensional case for obtaining the right values at the first
flux point of the $i$th cell $u_{i,fp_{1}}^{R}$based on the stencil
of the first solution point $sp_{1}$, where $fp_{1}$ and $sp_{1}$
denote the first flux point and the first solution point at the $i$th
cell, correspondingly. The values at other flux points can be obtained
by similar procedure.

Step A: Choose a stencil of five points $S_{i,sp_{1}}=\{u_{i-1,sp_{4}},u_{i-1,sp_{5}},u_{i,sp_{1}},u_{i,sp_{2}},u_{i,sp_{3}}\}$
and divide this stencil into three small stencil $S_{i,sp_{1}}^{(1)}=\{u_{i-1,sp_{4}},u_{i-1,sp_{5}},u_{i,sp_{1}}\}$,
$S_{i,sp_{1}}^{(2)}=\{u_{i-1,sp_{5}},u_{i,sp_{1}},u_{i,sp_{2}}\}$
and $S_{i,sp_{1}}^{(3)}=\{u_{i,sp_{1}},u_{i,sp_{2}},u_{i,sp_{3}}\}$. 

Step B: Construct Lagrange interpolation polynomial $p_{sp_{1}}^{(m)}(\xi)$
in each stencil $S_{i,sp_{1}}^{(m)},m=1,2,3$, we have

\begin{eqnarray}
p_{sp_{1}}^{(1)}(\xi_{fp_{1}}^{R}) & = & c_{11}u_{i-1,sp4}+c_{12}u_{i-1,sp5}+c_{13}u_{i,sp1},\nonumber \\
p_{sp_{1}}^{(2)}(\xi_{fp_{1}}^{R}) & = & c_{21}u_{i-1,sp5}+c_{22}u_{i,sp1}+c_{23}u_{i,sp2},\label{eq:polynomial-small-stencil}\\
p_{sp_{1}}^{(3)}(\xi_{fp_{1}}^{R}) & = & c_{31}u_{i,sp1}+c_{32}u_{i,sp2}+c_{33}u_{i,sp3}.\nonumber 
\end{eqnarray}

Step C: Calculate the linear weights $d_{m}$ for each stencil $S_{i,sp_{1}}^{(m)},m=1,2,3$.
The fifth-order linear interpolation for obtaining $u_{i,fp_{1}}^{R}$can
be obtained by Taylor expansion or Lagrange interpolation polynomial,
\begin{equation}
\begin{alignedat}{2}u_{i,fp_{1}}^{R} & = & a_{11}u_{i-1,sp4}+a_{12}u_{i-1,sp5}+a_{13}u_{i,sp1}+a_{14}u_{i,sp2}+a_{15}u_{i,sp3}\end{alignedat}
\label{eq:fifth-order linear interpolation}
\end{equation}
The linear weight are chosen to make the weighted results has optimal
fifth-order of accuracy, thus

\begin{eqnarray}
 &  & d_{1}(c_{11}u_{i-1,sp4}+c_{12}u_{i-1,sp5}+c_{13}u_{i,sp1})\nonumber \\
 & + & d_{2}(c_{21}u_{i-1,sp5}+c_{22}u_{i,sp1}+c_{23}u_{i,sp2})\label{eq:linear weights}\\
 & + & d_{3}(c_{31}u_{i,sp1}+c_{32}u_{i,sp2}+c_{33}u_{i,sp3}),\nonumber \\
 & = & a_{11}u_{i-1,sp4}+a_{12}u_{i-1,sp5}+a_{13}u_{i,sp1}+a_{14}u_{i,sp2}+a_{15}u_{i,sp3}.\nonumber 
\end{eqnarray}
Then, the linear weights $d_{1}$, $d_{2}$,$d_{3}$ can be determined.
Coefficients and linear weights for fifth-order interpolation in (\ref{eq:polynomial-small-stencil})
and (\ref{eq:linear weights}) are given in Table 11. It is worth
noticing that all of the linear weight are positive.

Step D: Compute smoothness indicator $IS_{m,sp_{1}}$ and nonlinear
weights $\omega_{m}$ to get the NNW interpolation value at the flux
points $\xi_{fp1}$,

\[
u_{i,fp_{1}}^{R}=\sum_{m=1}^{3}\omega_{m}p_{sp_{1}}^{(m)}(\xi_{fp1}),
\]
where $\left\{ \omega_{1},\omega_{2},\omega_{3}\right\} $ are nonlinear
weights. Various types of nonlinear weights have been developed, we
refer to \citet{Jiang1996,Borges2008,Yan2016} and references therein.
We consider two types of nonlinear weights in this paper. The first
one is the JS weights \citet{Jiang1996}, which are defined by
\begin{equation}
\omega_{m}=\frac{\beta_{m}}{\sum_{l=1}^{3}\beta_{l}},\quad\beta_{m}=\frac{d_{m}}{(\varepsilon+IS_{m,sp_{1}})^{2}},\label{eq:JSweight}
\end{equation}
where $\varepsilon=10^{-6}$ is a small number and $IS_{m,sp_{1}}$
is a smoothness indicator. The second one is the Z weights \citet{Borges2008},
which are defined by

\begin{equation}
\omega_{m}=\frac{\beta_{m}}{\sum_{l=1}^{3}\beta_{l}},\quad\beta_{m}=d_{m}\left(1+\frac{\left|IS_{3,sp_{1}}-IS_{1,sp_{1}}\right|^{2}}{(\varepsilon+IS_{m,sp_{1}})^{2}}\right),\label{eq:Zweight}
\end{equation}
where $\varepsilon=10^{-10}$ is a small number and $IS_{m,sp_{1}}$
is a smoothness indicator. 

\begin{flushleft}
\begin{table}
\begin{centering}
\begin{tabular}{c|c|c||cc||c}
\hline 
sp & fp & \multicolumn{4}{c}{$\left\{ d_{l}|l=1,2,3\right\} $, $\left\{ c_{i,j}|i=1,2,3,\; j=1,2,3\right\} $}\tabularnewline
\hline 
\multirow{12}{*}{sp1} & \multirow{6}{*}{$u_{i,fp_{1}}^{R}$} & \multicolumn{2}{c}{{\footnotesize $d_{1}^{1L}$=0.34210708229202832129514593919672}} & \multicolumn{2}{c}{{\footnotesize $c_{11}^{1L}$=-0.043104119062505851129386949103353}}\tabularnewline
 &  & \multicolumn{2}{c}{{\footnotesize $d_{2}^{1L}$=0.6308070429239803035449279980960}} & \multicolumn{2}{c}{{\footnotesize $c_{12}^{1L}$=0.62757338302641182681128438378820}}\tabularnewline
 &  & \multicolumn{2}{c}{{\footnotesize $d_{3}^{1L}$=0.027085874783991375159926062707238}} & \multicolumn{2}{c}{{\footnotesize $c_{13}^{1L}$=0.41553073603609402431810256531515}}\tabularnewline
 &  & \multicolumn{2}{c}{{\footnotesize $c_{21}^{1L}$=$c_{13}^{1L}$}} & \multicolumn{2}{c}{{\footnotesize $c_{31}^{1L}$=1.3850967572035142771188790865197}}\tabularnewline
 &  & \multicolumn{2}{c}{{\footnotesize $c_{22}^{1L}$=$c_{12}^{1L}$}} & \multicolumn{2}{c}{{\footnotesize $c_{32}^{1L}$=-0.47383715518113200994718334103294}}\tabularnewline
 &  & \multicolumn{2}{c}{{\footnotesize $c_{23}^{1L}$=$c_{11}^{1L}$}} & \multicolumn{2}{c}{{\footnotesize $c_{33}^{1L}$=0.088740397977617732828304254513256}}\tabularnewline
\cline{2-6} 
 & \multirow{6}{*}{$u_{i,fp_{2}}^{L}$} & \multicolumn{2}{c}{{\footnotesize $d_{1}^{1R}$=0.12849535271459107836474476379591}} & \multicolumn{2}{c}{{\footnotesize $c_{11}^{1R}$=0.22720313860956266078248678148547}}\tabularnewline
 &  & \multicolumn{2}{c}{{\footnotesize $d_{2}^{1R}$=0.7282735720676514895919532282452}} & \multicolumn{2}{c}{{\footnotesize $c_{12}^{1R}$=-1.4245445419469326647976617868947}}\tabularnewline
 &  & \multicolumn{2}{c}{{\footnotesize $d_{3}^{1R}$=0.14323107521775743204330200795893}} & \multicolumn{2}{c}{{\footnotesize $c_{13}^{1R}$=2.1973414033373700040151750054092}}\tabularnewline
 &  & \multicolumn{2}{c}{{\footnotesize $c_{21}^{1R}$=-0.30686122157984679047945926834897}} & \multicolumn{2}{c}{{\footnotesize $c_{31}^{1R}$=0.52024206914481960337059414500993}}\tabularnewline
 &  & \multicolumn{2}{c}{{\footnotesize $c_{22}^{1R}$=1.0796580829702841296969724868635}} & \multicolumn{2}{c}{{\footnotesize $c_{32}^{1R}$=0.54529095781405304293936341037295}}\tabularnewline
 &  & \multicolumn{2}{c}{{\footnotesize $c_{23}^{1R}$=0.22720313860956266078248678148547}} & \multicolumn{2}{c}{{\footnotesize $c_{33}^{1R}$=-0.065533026958872646309957555382876}}\tabularnewline
\hline 
\multirow{12}{*}{sp2} & \multirow{6}{*}{$u_{i,fp_{2}}^{R}$} & \multicolumn{2}{c}{{\footnotesize $d_{1}^{2L}$=0.5585589126271906359535502604210}} & \multicolumn{2}{c}{{\footnotesize $c_{11}^{2L}$=-0.30686122157984679047945926834897}}\tabularnewline
 &  & \multicolumn{2}{c}{{\footnotesize $d_{2}^{2L}$=0.421765474422970721147577236989}} & \multicolumn{2}{c}{{\footnotesize $c_{12}^{2L}$=1.0796580829702841296969724868635}}\tabularnewline
 &  & \multicolumn{2}{c}{{\footnotesize $d_{3}^{2L}$=0.019675612949838642898872502590336}} & \multicolumn{2}{c}{{\footnotesize $c_{13}^{2L}$=0.22720313860956266078248678148547}}\tabularnewline
 &  & \multicolumn{2}{c}{{\footnotesize $c_{21}^{2L}$=0.52024206914481960337059414500993}} & \multicolumn{2}{c}{{\footnotesize $c_{31}^{2L}$=1.7197296427724247141318458677543}}\tabularnewline
 &  & \multicolumn{2}{c}{{\footnotesize $c_{22}^{2L}$=0.54529095781405304293936341037295}} & \multicolumn{2}{c}{{\footnotesize $c_{32}^{2L}$=-1.0186627618222429820572260495517}}\tabularnewline
 &  & \multicolumn{2}{c}{{\footnotesize $c_{23}^{2L}$=-0.065533026958872646309957555382876}} & \multicolumn{2}{c}{{\footnotesize $c_{33}^{2L}$=0.29893311904981826792538018179744}}\tabularnewline
\cline{2-6} 
 & \multirow{6}{*}{$u_{i,fp_{3}}^{L}$} & \multicolumn{2}{c}{{\footnotesize $d_{1}^{2R}$=0.13159031124797584071088104105971}} & \multicolumn{2}{c}{{\footnotesize $c_{11}^{2R}$=1.5090040675292501024133653209020}}\tabularnewline
 &  & \multicolumn{2}{c}{{\footnotesize $d_{2}^{2R}$=0.6554862204946998473379356400573}} & \multicolumn{2}{c}{{\footnotesize $c_{12}^{2R}$=-2.9677269746646236620838843129099}}\tabularnewline
 &  & \multicolumn{2}{c}{{\footnotesize $d_{3}^{2R}$=0.21292346825732431195118331888295}} & \multicolumn{2}{c}{{\footnotesize $c_{13}^{2R}$=2.4587229071353735596705189920079}}\tabularnewline
 &  & \multicolumn{2}{c}{{\footnotesize $c_{21}^{2R}$=-0.21677318533425579686895436642487}} & \multicolumn{2}{c}{{\footnotesize $c_{31}^{2R}$=0.40514931162836377750297939141122}}\tabularnewline
 &  & \multicolumn{2}{c}{{\footnotesize $c_{22}^{2R}$=0.89451153321243085774192320462416}} & \multicolumn{2}{c}{{\footnotesize $c_{32}^{2R}$=0.71940940490306471718773980440933}}\tabularnewline
 &  & \multicolumn{2}{c}{{\footnotesize $c_{23}^{2R}$=0.32226165212182493912703116180071}} & \multicolumn{2}{c}{{\footnotesize $c_{33}^{2R}$=-0.12455871653142849469071919582055}}\tabularnewline
\hline 
\multirow{8}{*}{sp3} & \multirow{6}{*}{$u_{i,fp_{3}}^{R}$} & \multicolumn{2}{c}{{\footnotesize $d_{1}^{3L}$=0.37482146743990888513676931020257}} & \multicolumn{2}{c}{{\footnotesize $c_{11}^{3L}$=$c_{21}^{2R}$}}\tabularnewline
 &  & \multicolumn{2}{c}{{\footnotesize $d_{2}^{3L}$=0.5651196564082159925965654892928}} & \multicolumn{2}{c}{{\footnotesize $c_{12}^{3L}$=$c_{22}^{2R}$}}\tabularnewline
 &  & \multicolumn{2}{c}{{\footnotesize $d_{3}^{3L}$=0.06005887615187512226666520050463}} & \multicolumn{2}{c}{{\footnotesize $c_{13}^{3L}$=$c_{23}^{2R}$}}\tabularnewline
 &  & \multicolumn{2}{c}{{\footnotesize $c_{21}^{3L}$=$c_{31}^{2R}$}} & \multicolumn{2}{c}{{\footnotesize $c_{31}^{3L}$=2.0112028910174770092546365556601}}\tabularnewline
 &  & \multicolumn{2}{c}{{\footnotesize $c_{22}^{3L}$=$c_{32}^{2R}$}} & \multicolumn{2}{c}{{\footnotesize $c_{32}^{3L}$=-1.7162961134215585397544348429330}}\tabularnewline
 &  & \multicolumn{2}{c}{{\footnotesize $c_{23}^{3L}$=$c_{33}^{2R}$}} & \multicolumn{2}{c}{{\footnotesize $c_{33}^{3L}$=0.70509322240408153049979828727285}}\tabularnewline
\cline{2-6} 
 & \multirow{2}{*}{$u_{i,fp_{4}}^{L}$} & \multicolumn{2}{c}{{\footnotesize $d_{1}^{3R}$=$d_{3}^{3L}$,$d_{2}^{3R}$=$d_{2}^{3L}$,$d_{3}^{3R}$=$d_{1}^{3L}$}} & \multicolumn{2}{c}{{\footnotesize $c_{11}^{3R}$=$c_{33}^{3L}$,$c_{12}^{3R}$=$c_{32}^{3L}$,$c_{13}^{3R}$=$c_{31}^{3L}$}}\tabularnewline
 &  & \multicolumn{2}{c}{{\footnotesize $c_{21}^{3R}$=$c_{23}^{3L}$,$c_{22}^{3R}$=$c_{22}^{3L}$,$c_{23}^{3R}$=$c_{21}^{3L}$}} & \multicolumn{2}{c}{{\footnotesize $c_{31}^{3R}$=$c_{13}^{3L}$,$c_{32}^{3R}$=$c_{12}^{3L}$,$c_{33}^{3R}$=$c_{11}^{3L}$}}\tabularnewline
\hline 
\multirow{4}{*}{sp4} & \multirow{2}{*}{$u_{i,fp_{4}}^{R}$} & \multicolumn{2}{c}{{\footnotesize $d_{1}^{4L}$=$d_{3}^{2R}$,$d_{2}^{4L}$=$d_{2}^{2R}$,$d_{3}^{4L}$=$d_{1}^{2R}$}} & \multicolumn{2}{c}{{\footnotesize $c_{11}^{4L}$=$c_{33}^{2R}$,$c_{12}^{4L}$=$c_{32}^{2R}$,$c_{13}^{4L}$=$c_{31}^{2R}$}}\tabularnewline
 &  & \multicolumn{2}{c}{{\footnotesize $c_{21}^{4L}$=$c_{23}^{2R}$,$c_{22}^{4L}$=$c_{22}^{2R}$,$c_{23}^{4L}$=$c_{21}^{2R}$}} & \multicolumn{2}{c}{{\footnotesize $c_{31}^{4L}$=$c_{13}^{2R}$,$c_{32}^{4L}$=$c_{12}^{2R}$,$c_{33}^{4L}$=$c_{11}^{2R}$}}\tabularnewline
\cline{2-6} 
 & \multirow{2}{*}{$u_{i,fp_{5}}^{L}$} & \multicolumn{2}{c}{{\footnotesize $d_{1}^{4R}$=$d_{3}^{2L}$,$d_{2}^{4R}$=$d_{2}^{2L}$,$d_{3}^{4R}$=$d_{1}^{2L}$}} & \multicolumn{2}{c}{{\footnotesize $c_{11}^{4R}$=$c_{33}^{2L}$,$c_{12}^{4R}$=$c_{32}^{2L}$,$c_{13}^{4R}$=$c_{31}^{2L}$}}\tabularnewline
 &  & \multicolumn{2}{c}{{\footnotesize $c_{21}^{4R}$=$c_{23}^{2L}$,$c_{22}^{4R}$=$c_{22}^{2L}$,$c_{23}^{4R}$=$c_{21}^{2L}$}} & \multicolumn{2}{c}{{\footnotesize $c_{31}^{4R}$=$c_{13}^{2L}$,$c_{32}^{4R}$=$c_{12}^{2L}$,$c_{33}^{4R}$=$c_{11}^{2L}$}}\tabularnewline
\hline 
\multirow{4}{*}{sp5} & \multirow{2}{*}{$u_{i,fp_{5}}^{R}$} & \multicolumn{2}{c}{{\footnotesize $d_{1}^{5L}$=$d_{3}^{1R}$,$d_{2}^{5L}$=$d_{2}^{1R}$,$d_{3}^{5L}$=$d_{1}^{1R}$}} & \multicolumn{2}{c}{{\footnotesize $c_{11}^{5L}$=$c_{33}^{1R}$,$c_{12}^{5L}$=$c_{32}^{1R}$,$c_{13}^{5L}$=$c_{31}^{1R}$}}\tabularnewline
 &  & \multicolumn{2}{c}{{\footnotesize $c_{21}^{5L}$=$c_{23}^{1R}$,$c_{22}^{5L}$=$c_{22}^{1R}$,$c_{23}^{5L}$=$c_{21}^{1R}$}} & \multicolumn{2}{c}{{\footnotesize $c_{31}^{5L}$=$c_{13}^{1R}$,$c_{32}^{5L}$=$c_{12}^{1R}$,$c_{33}^{5L}$=$c_{11}^{1R}$}}\tabularnewline
\cline{2-6} 
 & \multirow{2}{*}{$u_{i,fp_{6}}^{L}$} & \multicolumn{2}{c}{{\footnotesize $d_{1}^{5R}$=$d_{3}^{1L}$,$d_{2}^{5R}$=$d_{2}^{1L}$,$d_{3}^{5R}$=$d_{1}^{1L}$}} & \multicolumn{2}{c}{{\footnotesize $c_{11}^{5R}$=$c_{33}^{1L}$,$c_{12}^{5R}$=$c_{32}^{1L}$,$c_{13}^{5R}$=$c_{31}^{1L}$}}\tabularnewline
 &  & \multicolumn{2}{c}{{\footnotesize $c_{21}^{5R}$=$c_{23}^{1L}$,$c_{22}^{5R}$=$c_{22}^{1L}$,$c_{23}^{5R}$=$c_{21}^{1L}$}} & \multicolumn{2}{c}{{\footnotesize $c_{31}^{5R}$=$c_{13}^{1L}$,$c_{32}^{5R}$=$c_{12}^{1L}$,$c_{33}^{5R}$=$c_{11}^{1L}$}}\tabularnewline
\hline 
\end{tabular}
\par\end{centering}

\noindent \centering{}\caption{Coefficients and linear weights in NNW5 interpolation}
\end{table}

\par\end{flushleft}

\section*{Appendix B. Proofs}

(1) Suppose $D$ is a block circulant matrix $D=\frac{1}{K+1}Circ(C_{0},C_{1},C_{2},\cdots,C_{M-1})$,
then there exists a Fourier matrix $F_{(K+1)M}^{*}$ \citet{Davis1979}
such that

\begin{equation}
F_{(K+1)M}DF_{(K+1)M}^{*}=diag(H_{0},H_{1},\cdots,H_{M-1}),\label{eq:block circulant matrics}
\end{equation}
where 
\begin{equation}
H_{m}=\frac{1}{K+1}\sum_{k=0}^{M-1}(\omega^{m})^{k}C_{k},\,\, m=0,1,2,\cdots,M-1,\,\,\,\omega=e^{\frac{2\pi}{M}i}\label{eq:Hm}
\end{equation}
and $F_{(K+1)M}^{*}=F_{M}^{*}\otimes I_{(K+1)}$, 
\[
F_{M}^{*}=\frac{1}{\sqrt{M}}\left[\begin{array}{ccccc}
1 & 1 & 1 & \cdots & 1\\
1 & \omega & \omega^{2} & \cdots & \omega^{M-1}\\
1 & \omega^{2} & \omega^{4} & \cdots & \omega^{2(M-1)}\\
\vdots & \vdots & \vdots & \vdots & \vdots\\
1 & \omega^{M-1} & \omega^{2(M-1)} & \cdots & \omega^{(M-1)(M-1)}
\end{array}\right].
\]
We have

\begin{eqnarray*}
det(D-\mu I_{(K+1)M}) & = & \prod_{m=0}^{M-1}det(H_{m}-\mu I_{M}),
\end{eqnarray*}
where $I_{(K+1)M}$ and $I_{M}$ are $(K+1)M$ identify matrix and
$M$ identify matrix, correspondingly. Then, all eigenvalues of $D$
are given by 
\begin{eqnarray}
\{\mu|DX=\mu X,\, X\in\mathbb{C}^{(K+1)M}\} & = & \{\mu|H_{0}Y_{0}=\mu Y_{0},\, Y_{0}\in\mathbb{C}^{K+1}\}\nonumber \\
 &  & \cup\{\mu|H_{1}Y_{1}=\mu Y_{1},\, Y_{1}\in\mathbb{C}^{K+1}\}\nonumber \\
 &  & \cdots\nonumber \\
 &  & \cup\{\mu|H_{M-1}Y_{M-1}=\mu Y_{M-1},\, Y_{M-1}\in\mathbb{C}^{K+1}\}.\label{eq:all eigenvalues}
\end{eqnarray}

Denote $(K+1)M$ eigenvalues of the matrix $E$ be $spec\left(E\right)=\{\mu_{j},j=0,1,2,\cdots,(K+1)M-1\}$
and $(K+1)$ eigenvalues of the matrix $H_{m}$ be $Spec\left(H_{m}\right)=\{\mu^{(l)}(H_{m})|l=1,2,\cdots,K+1\}$,
where $m=0,1,2,\cdots,M-1$. For the matrix $E$ in (\ref{eq:matrixE}),
we have $E=\frac{1}{K+1}circ(B,C,0,\cdots,0,A)$. According to (\ref{eq:Hm}),
it can be easily obtained that

\begin{eqnarray*}
H_{m} & = & \frac{1}{K+1}\sum_{k=0}^{M-1}(\omega^{m})^{k}C_{k}\\
 & = & \frac{1}{K+1}\left(B+\omega^{m}C+(\omega^{m})^{M-1}A\right)\\
 & = & \frac{1}{K+1}\left((\omega^{m})^{-1}A+B+\omega^{m}C\right)\\
 & = & \frac{1}{K+1}\left(e^{-i\phi_{m}}A+B+e^{i\phi_{m}}C\right),
\end{eqnarray*}
where $\phi_{m}=m\frac{2\pi}{M}$. Therefore, we have $Spec\left(E_{(K+1)M\times(K+1)M}\right)=\left\{ Spec\left(H_{0}\right),Spec\left(H_{1}\right),\cdots,Spec\left(H_{M-1}\right)\right\} $.

(2) Consider the case $mod(M,(K+1))=0$ and set $M=L(K+1)$. For a
fixed integer $n_{0}\in[0,K]$ and $n_{0}M\leq m(K+1)<(n_{0}+1)M$,
we have 
\begin{eqnarray*}
G_{m}=G(\phi_{m}) & = & \left(Ae^{-i\phi_{m}(K+1)}+B+Ce^{i\phi_{m}(K+1)}\right)/(K+1)\\
 & = & \left(Ae^{-i\frac{m(K+1)}{M}2\pi}+B+Ce^{i\frac{m(K+1)}{M}2\pi}\right)/(K+1)\\
 & = & \left(Ae^{-i\frac{m(K+1)-n_{0}M}{M}2\pi}+B+Ce^{i\frac{m(K+1)-n_{0}M}{M}2\pi}\right)/(K+1)\\
 & = & \left(Ae^{-i\frac{(m-n_{0}L)(K+1)}{M}2\pi}+B+Ce^{i\frac{(m-n_{0}L)(K+1)}{M}2\pi}\right)/(K+1)\\
 & = & H(\phi_{(m-n_{0}L)(K+1)})\\
 & = & H_{(m-n_{0}L)(K+1)}.
\end{eqnarray*}
It is easy to check that $0\leq m-n_{0}L\leq L$ or $n_{0}L\leq m\leq(n_{0}+1)L$,
and 
\begin{eqnarray}
\left\{ G_{n_{0}L},G_{\left(n_{0}L+1\right)},\cdots,G_{\left(n_{0}L+L\right)}\right\}  & = & \left\{ H_{0},H_{\left(K+1\right)},\cdots,H_{\left(L(K+1)\right)}\right\} .\label{eq:GH_relation1}
\end{eqnarray}
Since the relation (\ref{eq:GH_relation1}) is satisfied for every
integer $n_{0}\in[0,K]$. Thus, we have

\begin{eqnarray*}
SG & = & \left\{ Spec\left(G_{0}\right),Spec\left(G_{1}\right),\cdots,Spec\left(G_{M-1}\right)\right\} \\
 & = & \left\{ Spec\left(H_{0}\right),Spec\left(H_{(K+1)}\right),Spec\left(H_{2(K+1)}\right),\cdots,Spec\left(H_{L(K+1)}\right)\right\} .
\end{eqnarray*}
Therefore,
\begin{eqnarray*}
SG & \subset & SH,
\end{eqnarray*}
and 
\begin{eqnarray*}
SG & \neq & SH.
\end{eqnarray*}
For the case $mod(M,(K+1))\neq0$, we set $mod(M,(K+1))=l_{0}$ and
$M=L(K+1)+l_{0}$ with $l_{0}\in[1,K]$. For a fixed integer $n_{0}\in[0,K]$
and $n_{0}M\leq m(K+1)<(n_{0}+1)M$ , we have 
\begin{eqnarray}
G_{m}=G(\phi_{m}) & = & \left(Ae^{-i\phi_{m}(K+1)}+B+Ce^{i\phi_{m}(K+1)}\right)/(K+1),\,\,\,\nonumber \\
 & = & \left(Ae^{-i\frac{m(K+1)}{M}2\pi}+B+Ce^{i\frac{m(K+1)}{M}2\pi}\right)/(K+1)\nonumber \\
 & = & \left(Ae^{-i\frac{m(K+1)-n_{0}M}{M}2\pi}+B+Ce^{i\frac{m(K+1)-n_{0}M}{M}2\pi}\right)/(K+1)\nonumber \\
 & = & \left(Ae^{-i\frac{(m-n_{0}L)(K+1)-n_{0}l_{0}}{M}2\pi}+B+Ce^{i\frac{(m-n_{0}L)(K+1)-n_{0}l_{0}}{M}2\pi}\right)/(K+1)\nonumber \\
 & = & H(\phi_{(m-n_{0}L)(K+1)-n_{0}l_{0}})\nonumber \\
 & = & H_{(m-n_{0}L)(K+1)-n_{0}l_{0}}.\label{eq:GH_relation2}
\end{eqnarray}
Thus, for $n_{0}=0$, we have $G_{0}=H_{0}$,$G_{1}=H_{K+1},\cdots,G_{L}=H_{(K+1)L}$.
For a fixed integer $n_{0}\in[1,K]$, we have $G_{(n_{0}L+1)}=H_{(K+1)-n_{0}l_{0}}$,$G_{(n_{0}L+2)}=H_{2(K+1)-n_{0}l_{0}},\cdots,G_{n_{0}L+L}=H_{L(K+1)-n_{0}l_{0}}$.
Notice that for the case $l_{0}\neq0$ we have $\left\{ \left.mod(n_{0}l_{0},K+1)\right|n_{0}\in[0,K]\right\} =\left\{ 0,1,2,\cdots,K\right\} $.
Thus, according to the relation (\ref{eq:GH_relation2}) it can be
easily checked that 

\begin{eqnarray*}
\left\{ G_{0},G_{1},\cdots,G_{M}\right\}  & = & \left\{ H_{0},H_{1},\cdots,H_{M}\right\} .
\end{eqnarray*}
Therefore, for the case $mod(M,(K+1))\neq0$, 
\begin{eqnarray*}
SG & = & SH.
\end{eqnarray*}

(3) Since $G_{m}$ can be written as 
\begin{eqnarray*}
G_{m} & = & \left(A\omega^{-m(K+1)}+B+C\omega^{m(K+1)}\right)/(K+1),
\end{eqnarray*}
the eigenvalue of $G_{m}$ is a function of $\omega^{m}$ and $\omega=e^{\frac{2\pi}{M}i}$.
Suppose $\lambda^{(1)}\left(\omega^{m}\right)$ be an eigenvalue of
$G_{m}$, then we have
\[
|\lambda^{(1)}\left(\omega^{m}\right)I_{(K+1)}-G(\omega^{m(K+1)})|=0.
\]
Thus
\[
|\lambda^{(1)}\left(e^{i\phi_{m}}\right)I_{(K+1)}-G(e^{i\phi_{m}(K+1)})|=0,
\]

\[
|\lambda^{(1)}\left(e^{i(\phi_{m}-\frac{(l-1)}{K+1}2\pi)}\right)I_{(K+1)}-G(e^{i(\phi_{m}-\frac{(l-1)}{K+1}2\pi)(K+1)})|=0,\,\,\,\, l=2,3,\cdots,K+1.
\]
Since $e^{il2\pi}=1$, we have 

\[
|\lambda^{(1)}\left(e^{i(\phi_{m}-\frac{(l-1)}{K+1}2\pi)}\right)I_{(K+1)}-G(e^{i\phi_{m}(K+1)})|=0,\,\,\,\, l=2,3,\cdots,K+1.
\]
Thus $\lambda^{(1)}(\phi_{m}-\frac{(l-1)}{K+1}2\pi),\,\,2,3,\cdots,K+1$
are also eigenvalues of $G_{m}$. In addition, $\lambda^{(1)}(\phi_{m}-\frac{(l-1)}{K+1}2\pi),\, l=1,2,\cdots,K+1$
are different from each other, thus the collection of them are the
all eigenvalues of $G_{m}$, and we set that

\begin{eqnarray}
\lambda^{(l)}(\phi_{m}) & = & \lambda^{(1)}(\phi_{m}-\frac{(l-1)}{K+1}2\pi),\, l=1,2,\cdots,K+1.\label{eq:eigen_of_G}
\end{eqnarray}
 If $M=(K+1)L$, we have

\begin{eqnarray*}
\lambda^{(l)}\left(\phi_{m}\right) & = & \lambda^{(1)}\left(\phi_{m}-(l-1)\frac{2\pi}{\left(K+1\right)}\right)\\
 & = & \lambda^{(1)}\left(m\frac{2\pi}{M}-(l-1)\frac{2\pi}{\left(K+1\right)}\right)\\
 & = & \lambda^{(1)}\left(m\frac{2\pi}{(K+1)L}-L(l-1)\frac{2\pi}{\left(K+1\right)L}\right)\\
 & = & \lambda^{(1)}\left((-L(l-1)+m)\frac{2\pi}{M}\right).
\end{eqnarray*}
Since $\lambda^{(1)}(\phi)$ is a periodic function, 

\begin{eqnarray*}
Group^{(l)} & = & \left\{ \left.\lambda^{(1)}\left(s\frac{2\pi}{M}\right)\right|\, s=-L(l-1),-L(l-1)+1,-L(l-1)+2,\cdots,-L(l-1)+(M-1)\right\} \\
 & = & \left\{ \left.\lambda^{(1)}\left(m\frac{2\pi}{M}\right)\right|\, m=0,1,2,\cdots,M-1\right\} =Group{}^{(1)}.
\end{eqnarray*}

(4) Denote $\psi_{m}=m\frac{2\pi}{M(K+1)}$. According to (\ref{eq:eigen_of_G}),
we have 

\begin{eqnarray*}
\lambda^{(l)}(\phi_{m}) & = & \lambda^{(1)}(\phi_{m}-\frac{(l-1)2\pi}{K+1})\\
 & = & \lambda^{(1)}(\frac{(K+1)m}{(K+1)M}2\pi-\frac{(l-1)M}{(K+1)M}2\pi)\\
 & = & \lambda^{(1)}(\psi_{(K+1)m}-\frac{(l-1)M}{(K+1)M}2\pi).
\end{eqnarray*}
Thus, according to periodic property of eigenvalue function $\lambda^{(1)}(\phi)$
with period $2\pi$, we have

\begin{eqnarray}
 &  & \left\{ \left.\lambda^{(1)}(\psi_{(K+1)m}-\frac{(l-1)M}{(K+1)M}2\pi)\right|m=0,1,\cdots,M-1\right\} \label{eq:relative}\\
 & = & \left\{ \left.\lambda^{(1)}(\psi_{(K+1)m+(K+1-mod((l-1)M),K+1)})\right|m=0,1,\cdots,M-1\right\} .\nonumber 
\end{eqnarray}
It can be easily checked that

\begin{eqnarray}
\left\{ \left.mod((l-1)M,K+1)\right|l=1,2,\cdots,K+1\right\}  & =\begin{cases}
\left\{ 0\right\} , & if\,\,\, mod(M,K+1)=0,\\
\left\{ 0,1,2\cdots,K\right\} , & else.
\end{cases}\label{eq:mod}
\end{eqnarray}
Thus, taking (\ref{eq:relative})(\ref{eq:mod}) and noting that
\begin{eqnarray*}
\left\{ \left.\lambda^{(1)}(\psi_{(K+1)m+k})\right|m=0,1,\cdots,M-1,k=0,1,2\cdots,K\right\}  & = & \left\{ \left.\lambda^{(1)}(\psi_{j})\right|j=0,1,2\cdots,M(K+1)\right\} ,
\end{eqnarray*}
we have

\begin{eqnarray*}
 &  & \left\{ \left.\lambda^{(1)}(\psi_{(K+1)m}-\frac{(l-1)M}{(K+1)M}2\pi)\right|m=0,1,\cdots,M-1;l=1,2,\cdots,K+1\right\} \\
 & = & \begin{cases}
\left\{ \left.\lambda^{(1)}(\psi_{(K+1)m})\right|m=0,1,\cdots,M-1\right\} ,\,\, & if\,\,\, mod(M,K+1)=0,\\
\left\{ \left.\lambda^{(1)}(\psi_{j})\right|j=0,1,2\cdots,M(K+1)\right\}  & else.
\end{cases}
\end{eqnarray*}
Therefore,
\begin{eqnarray*}
SG & = & \left\{ \left.\lambda^{(1)}(\phi_{m}-\frac{(l-1)}{K+1}2\pi)\right|m=0,1,\cdots,M-1;l=1,2,\cdots,K+1\right\} \\
 & = & \begin{cases}
\left\{ \left.\lambda^{(1)}(\psi_{(K+1)m})\right|m=0,1,\cdots,M-1\right\} ,\,\, & if\,\,\, mod(M,K+1)=0,\\
\left\{ \left.\lambda^{(1)}(\psi_{j})\right|j=0,1,2\cdots,M(K+1)\right\}  & else.
\end{cases}
\end{eqnarray*}

\section*{Appendix C. An example to explain properties of the eigenvalues in
Theorem 2.1}

We take a third-order WCNS for example to show the properties of the
eigenvalues of local discrete matrices and the unique spectral curve. 

A third-order WCNS reads

\begin{eqnarray*}
\frac{\partial u_{i}}{\partial t} & = & -\frac{1}{\Delta x}\left[\frac{4}{3}\left(u_{i+\frac{1}{2}}-u_{i-\frac{1}{2}}\right)-\frac{1}{6}\left(u_{i+1}-u_{i-1}\right)\right]\\
 & = & -\frac{1}{\Delta x}\left[\frac{4}{3}\left(u_{i+\frac{1}{2}}^{L}-u_{i-\frac{1}{2}}^{L}\right)-\frac{1}{6}\left(u_{i+1}-u_{i-1}\right)\right]\\
 & = & -\frac{1}{\Delta x}\left[\frac{4}{3}\left(\frac{6u_{i}+3u_{i+1}-u_{i-1}}{8}-\frac{6u_{i-1}+3u_{i}-u_{i-2}}{8}\right)-\frac{1}{6}\left(u_{i+1}-u_{i-1}\right)\right]\\
 & = & -\frac{1}{3\Delta x}\left(\frac{1}{2}u_{i-2}-3u_{j-1}+\frac{3}{2}u_{i}+u_{i+1}\right),
\end{eqnarray*}
where $\Delta x=x_{i+1/2}-x_{i-1/2}$.

Thus, the third-order WCNS can be written as the first form in (\ref{eq:Form1})
with 
\begin{eqnarray*}
E & = & -\frac{1}{K+1}\left[\begin{array}{cccccccccc}
\frac{3}{2} & 1 &  &  &  &  &  &  & \frac{1}{2} & -3\\
-3 & \frac{3}{2} & 1 &  &  &  &  &  &  & \frac{1}{2}\\
\frac{1}{2} & -3 & \frac{3}{2} & 1\\
 & \frac{1}{2} & -3 & \frac{3}{2} & 1\\
 &  & \frac{1}{2} & -3 & \frac{3}{2} & 1\\
 &  &  & \ddots & \ddots & \ddots & \ddots\\
 &  &  &  & \frac{1}{2} & -3 & \frac{3}{2} & 1\\
 &  &  &  &  & \frac{1}{2} & -3 & \frac{3}{2} & 1\\
 &  &  &  &  &  & \frac{1}{2} & -3 & \frac{3}{2} & 1\\
1 &  &  &  &  &  &  & \frac{1}{2} & -3 & \frac{3}{2}
\end{array}\right]_{N\times N}\\
 & = & Circ(c_{0},c_{1},\cdots,c_{N-1})=Circ(-\frac{1}{2},-\frac{1}{3},0,\cdots,0,-\frac{1}{6},1).
\end{eqnarray*}
Here $K=2$ and $N=(K+1)M$. In this case, $E$ is not only a block
circulant matrix but also a circulant matrix. According to special
property of circulant matrix\citet{Davis1979,Liqiao1988}, we can
directly obtain the spectrum of $E$ 
\begin{eqnarray}
\lambda_{j}(E) & = & f_{C}(\zeta^{j})=c_{0}+c_{1}\zeta^{j}+\cdots+c_{N-1}\left(\zeta^{j}\right)^{N-1}\nonumber \\
 & = & \left(-\frac{1}{6}\frac{1}{e^{2i\psi_{j}}}+\frac{1}{e^{i\psi_{j}}}-\frac{1}{2}-\frac{1}{3}e^{i\psi_{j}}\right),\label{eq:eig1-E1}\\
 & = & \left(-\frac{1}{2}+\frac{2}{3}cos(\psi_{j})-\frac{1}{6}cos(2\psi_{j})\right)+i\left(\frac{1}{6}sin(2\psi_{j})-\frac{4}{3}sin(\psi_{j})\right),\nonumber 
\end{eqnarray}
where $\zeta=e^{\frac{2\pi}{N}i}$, $\psi_{j}=j\frac{2\pi}{N}=j\frac{2\pi}{(K+1)M},\,\, j=0,1,2,\cdots,N-1$.

On the other side, the third-order WCNS can also be written as the
second form in (\ref{eq:Form2}) by putting three solution points
equally at each cell, 
\begin{eqnarray*}
\frac{\partial}{\partial t}\left[\begin{array}{c}
u_{j,1}\\
u_{j,2}\\
u_{j,3}
\end{array}\right] & = & \frac{1}{\Delta x}\cdot\frac{1}{(K+1)}\left\{ A\left[\begin{array}{c}
u_{j-1,1}\\
u_{j-2,2}\\
u_{j-3,3}
\end{array}\right]+B\left[\begin{array}{c}
u_{j,1}\\
u_{j,2}\\
u_{j,3}
\end{array}\right]+C\left[\begin{array}{c}
u_{j+1,1}\\
u_{j+2,2}\\
u_{j+3,3}
\end{array}\right]\right\} 
\end{eqnarray*}
with
\[
A=-\left[\begin{array}{ccc}
0 & \frac{1}{2} & -3\\
0 & 0 & \frac{1}{2}\\
0 & 0 & 0
\end{array}\right],\,\, B=-\left[\begin{array}{ccc}
\frac{3}{2} & 1 & 0\\
-3 & \frac{3}{2} & 1\\
\frac{1}{2} & -3 & \frac{3}{2}
\end{array}\right],\,\, C=-\left[\begin{array}{ccc}
0 & 0 & 0\\
0 & 0 & 0\\
1 & 0 & 0
\end{array}\right].
\]

Then, the matrix $G_{m}$ is

\begin{eqnarray*}
G_{m} & = & \left(Ae^{-i\phi_{m}(K+1)}+B+Ce^{i\phi_{m}(K+1)}\right)/(K+1)\\
 & =- & \frac{1}{(K+1)}\left[\begin{array}{ccc}
\frac{3}{2} & \frac{1}{2}e^{-i\phi_{m}(K+1)}+1 & -3e^{-i\phi_{m}(K+1)}\\
-3 & \frac{3}{2} & 1+\frac{1}{2}e^{-i\phi_{m}(K+1)}\\
\frac{1}{2}+e^{i\phi_{m}(K+1)} & -3 & \frac{3}{2}
\end{array}\right],
\end{eqnarray*}
where $\phi_{m}=m\frac{2\pi}{M},\,\, m=0,1,2,\cdots,M-1$. $G_{m}$
has $K+1=3$ eigenvalues, which are 
\begin{eqnarray}
z_{1}(\phi_{m})=\lambda^{(1)}(G_{m}) & = & -\left(\frac{1}{6}\frac{1}{e^{2i\phi_{m}}}-\frac{1}{e^{i\phi_{m}}}+\frac{1}{2}+\frac{e^{i\phi_{m}}}{3}\right)=f(\phi_{m}),\label{eq:eig1-G1}
\end{eqnarray}
\begin{eqnarray*}
z_{2}(\phi_{m})=\lambda^{(2)}(G_{m}) & =-\left(\frac{1}{6}\frac{1}{e^{2i(\phi_{m}+\frac{2\pi}{(K+1)})}}-\frac{1}{e^{i(\phi_{m}+\frac{2\pi}{(K+1)})}}+\frac{1}{2}+\frac{e^{i(\phi_{m}+\frac{2\pi}{(K+1)})}}{3}\right) & =f(\phi_{m}+\frac{2\pi}{(K+1)}),
\end{eqnarray*}
\begin{eqnarray*}
z_{3}(\phi_{m})=\lambda^{(3)}(G_{m}) & =-\left(\frac{1}{6}\frac{1}{e^{2i(\phi_{m}+\frac{2\cdot2\pi}{(K+1)})}}-\frac{1}{e^{i(\phi_{m}+\frac{2\cdot2\pi}{(K+1)})}}+\frac{1}{2}+\frac{e^{i(\phi_{m}+\frac{2\cdot2\pi}{(K+1)})}}{3}\right) & =f(\phi_{m}+2\times\frac{2\pi}{(K+1)}).
\end{eqnarray*}
Thus, $\left\{ G_{0},G_{1},\cdots,G_{M-1}\right\} $ has $(K+1)M$
eigenvalues, which can be clarified as $(K+1)$th groups, 

\begin{eqnarray*}
Group^{(1)} & =\left\{ z_{1}(\phi_{m}),m=0,1,2,\cdots,M-1\right\} = & \left\{ \lambda^{(1)}(G_{m}),\,\, m=0,1,2,\cdots,M-1\right\} ,\\
Group^{(2)} & =\left\{ z_{2}(\phi_{m}),m=0,1,2,\cdots,M-1\right\} = & \left\{ \lambda^{(2)}(G_{m}),\,\, m=0,1,2,\cdots,M-1\right\} ,\\
Group^{(3)} & =\left\{ z_{3}(\phi_{m}),m=0,1,2,\cdots,M-1\right\} = & \left\{ \lambda^{(3)}(G_{m}),\,\, m=0,1,2,\cdots,M-1\right\} .
\end{eqnarray*}

Compare (\ref{eq:eig1-G1}) with (\ref{eq:eig1-E1}), we can find
that $Group^{(1)}$ has the same eigenvalue functions as $Spec\left(E\right)$
obtained by circulant matrix or Fourier analysis. In addition, if
$M=(K+1)L$, then $Group^{(1)}=Group^{(2)}=Group^{(3)}$. Curves of
$Group^{(2)}$ and $Group^{(3)}$ can be obtained by taking translation
transformation of the eigenvalue curve of $Group^{(1)}$, as shown
in Fig. \ref{fig:Spectral-WCNS}(a).

Here we also draw imaginary part of eigenvalues from $\{H_{0},H_{1},H_{2}\}$,
as shown in Fig. \ref{fig:Spectral-WCNS}(b). We can see that three
eigenvalue curves obtained from $\{H_{0},H_{1},H_{2}\}$ correspond
to the first, second and third part of the spectrum $Spec(E)$ obtained
directly from property of circulant matrix or Fourier analysis.

\begin{center}
\begin{figure}
\begin{centering}
\subfloat[Imaginary part of all eigenvalues computed from matrices $\left\{ G_{0},G_{1},\cdots G_{M-1}\right\} $]{\includegraphics[width=0.48\textwidth]{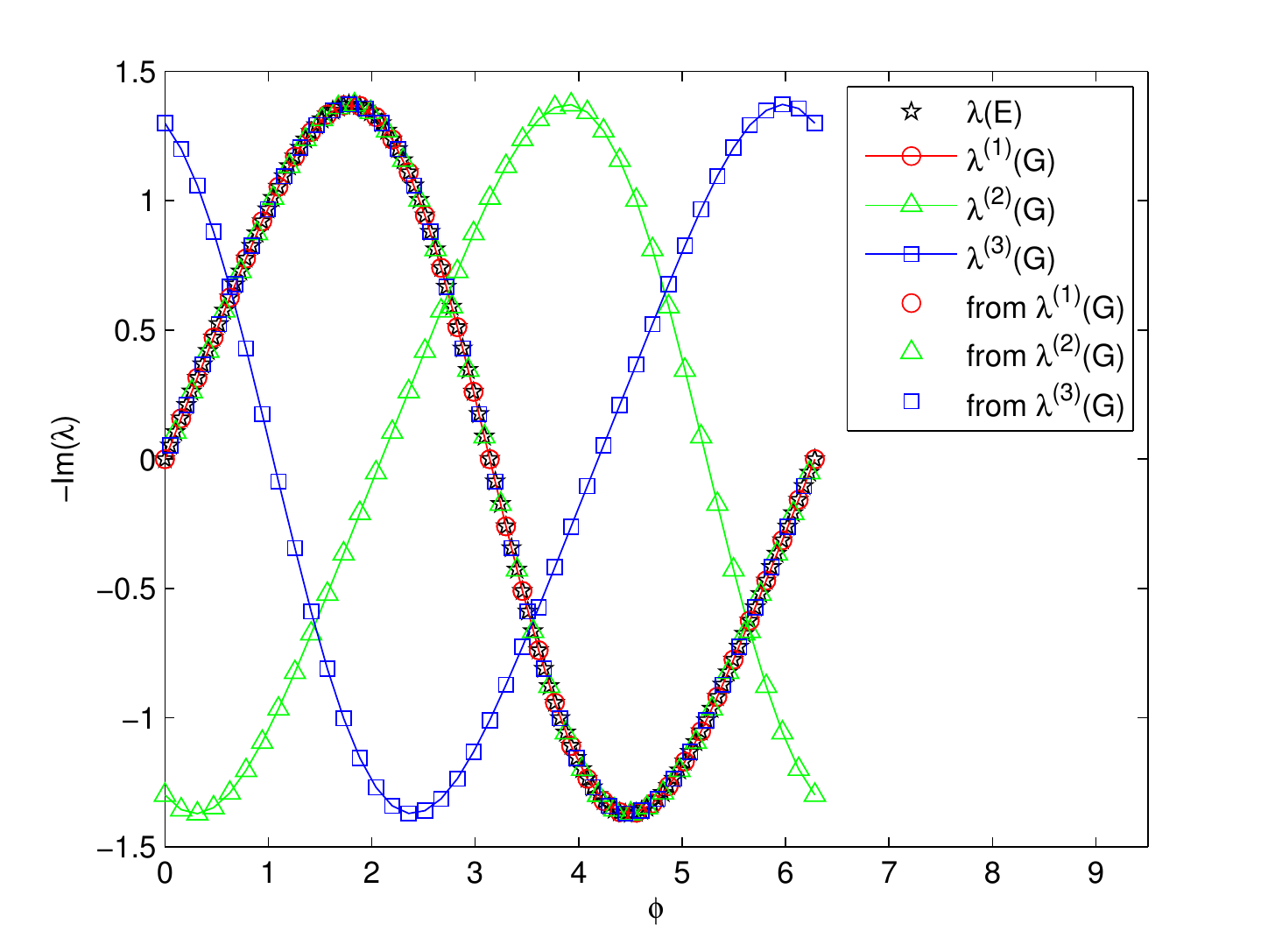}}\subfloat[Imaginary part of all eigenvalues computed from matrices $\left\{ H_{0},H_{1},\cdots H_{M-1}\right\} $]{\includegraphics[width=0.48\textwidth]{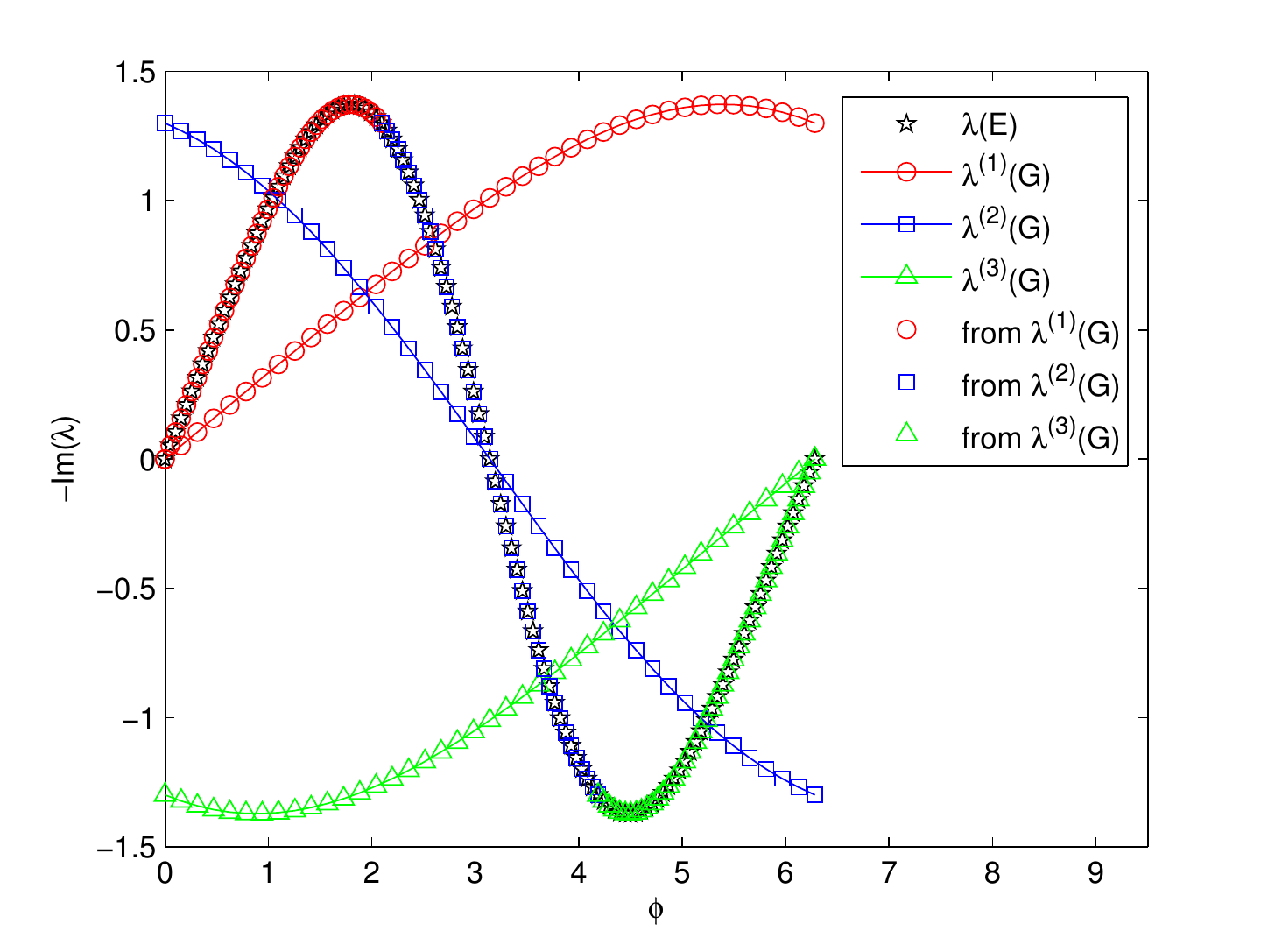}}
\par\end{centering}

\caption{Imaginary part of eigenvalues computed from two kinds of local matrices
for the third-order WCNS, where $z_{1}$, $z_{2}$ and $z_{3}$ are
the three eigenvalues of local matrices \label{fig:Spectral-WCNS}}
\end{figure}

\par\end{center}

\section*{Appendix D. Comparisons on spectrum of different high-order schemes}

The spectrum of CNNW is compared with WCNS and CPR by computing all
eigenvalues from the matrix $G$ with $M=40$.

Fig. \ref{fig:comparison-of-eigenvalues} shows eigenvalues in the
complex plane. we can see that all eigenvalues of each scheme have
negative real part, which illustrate that CNNW, CPR and WCNS are stable.
In addition, the three groups of eigenvalues (noted by $z1$, $z2$
and $z3$ in the Fig. \ref{fig:comparison-of-eigenvalues}) for third-order
schemes are different since $mod(M,3)\neq0$ and the five groups of
eigenvalues for fifth-order schemes are the same since $mod(M,5)=0$.
These results agree with Theorem \ref{thm:eigenvalues}. Dispersion
and dissipation relations in one period are shown for third-order
schemes in Fig. \ref{fig:comparison-of-dispersion-third-order} and
for fifth-order schemes in Fig. \ref{fig:comparison-of-dispersion-fifth-order}.

\begin{center}
\begin{figure}
\subfloat[C3NNW3]{\includegraphics[width=0.48\textwidth]{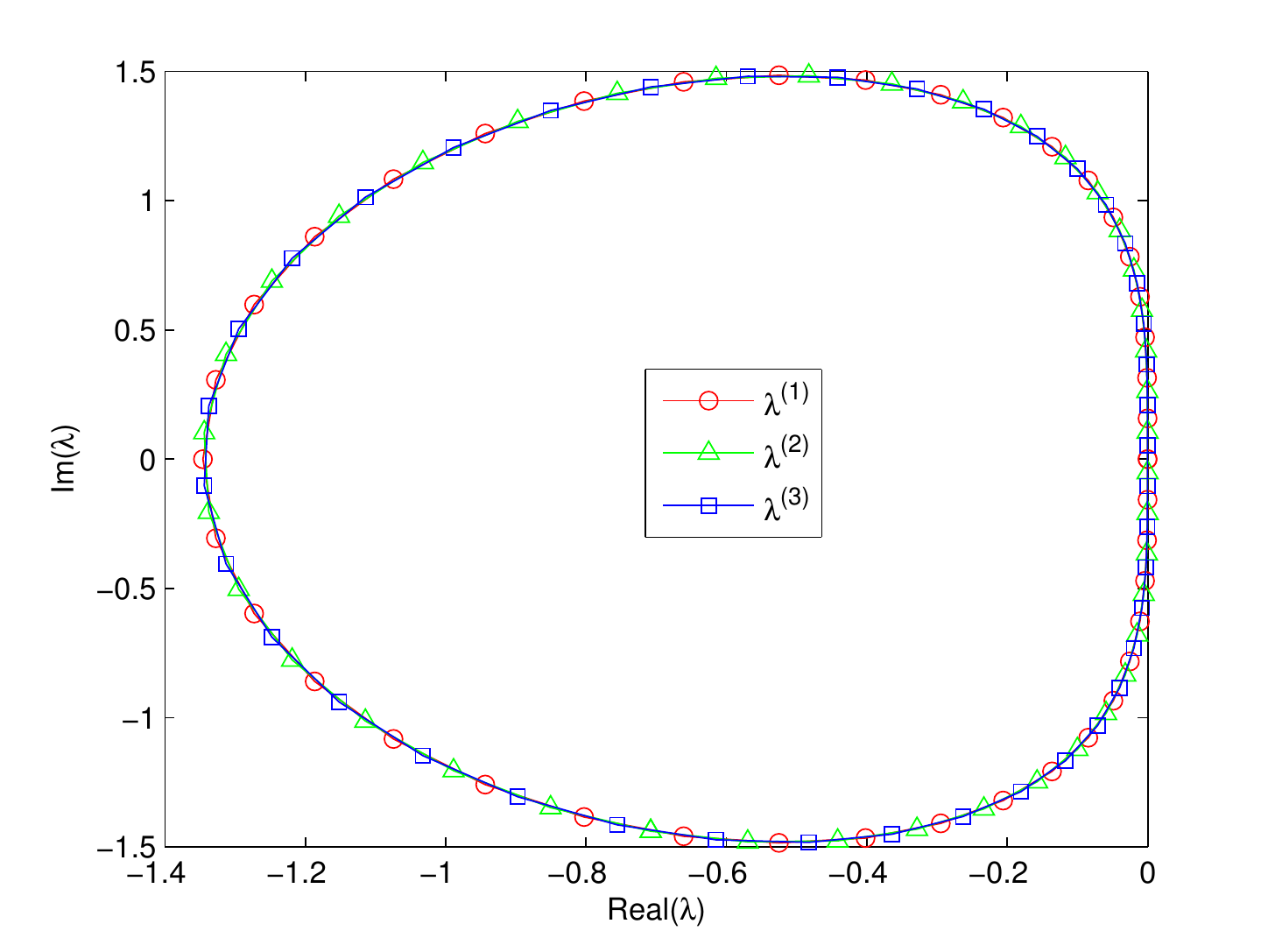}}\subfloat[C5NNW5]{\includegraphics[width=0.48\textwidth]{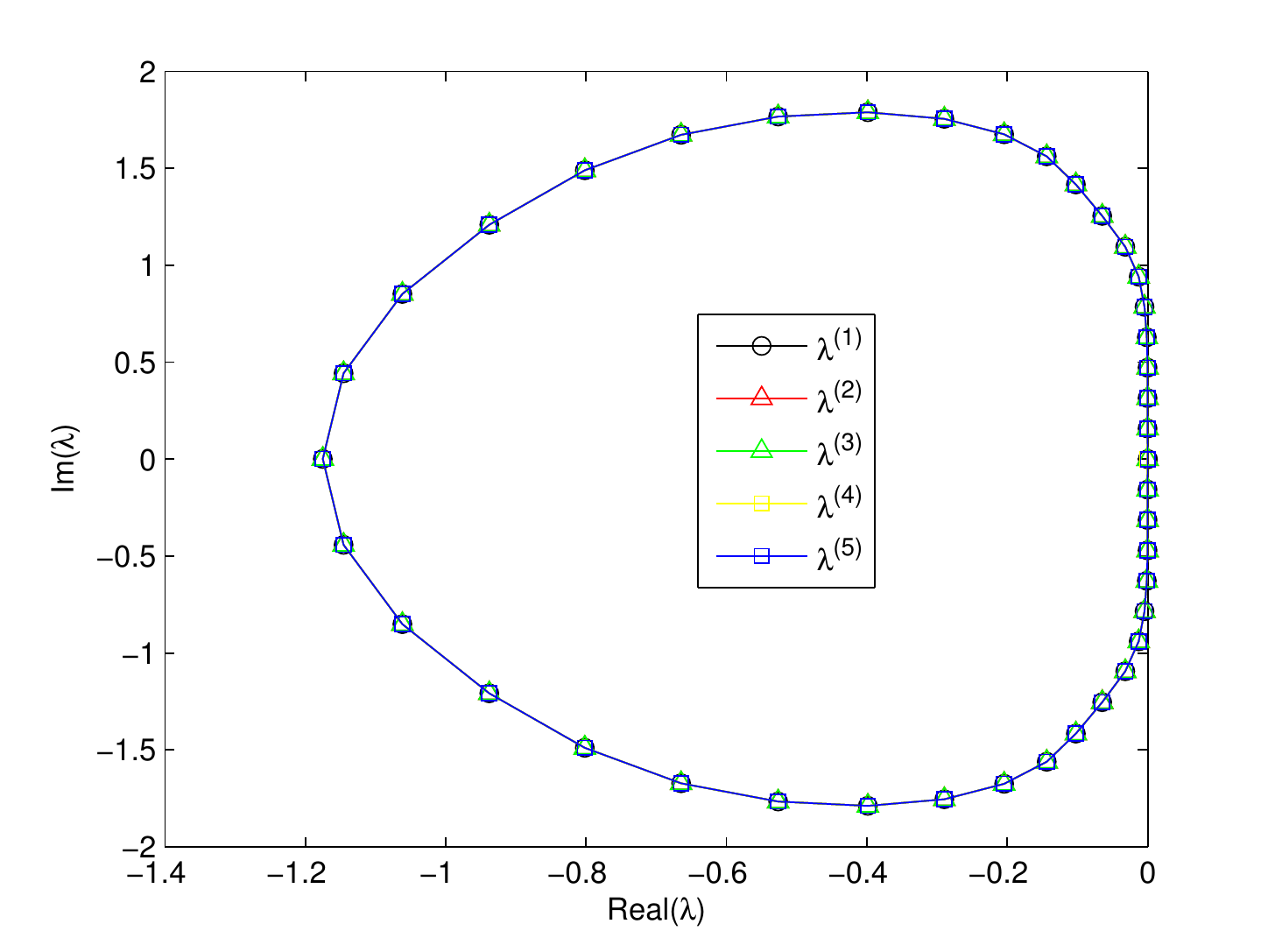}}

\subfloat[WCNS3]{\includegraphics[width=0.48\textwidth]{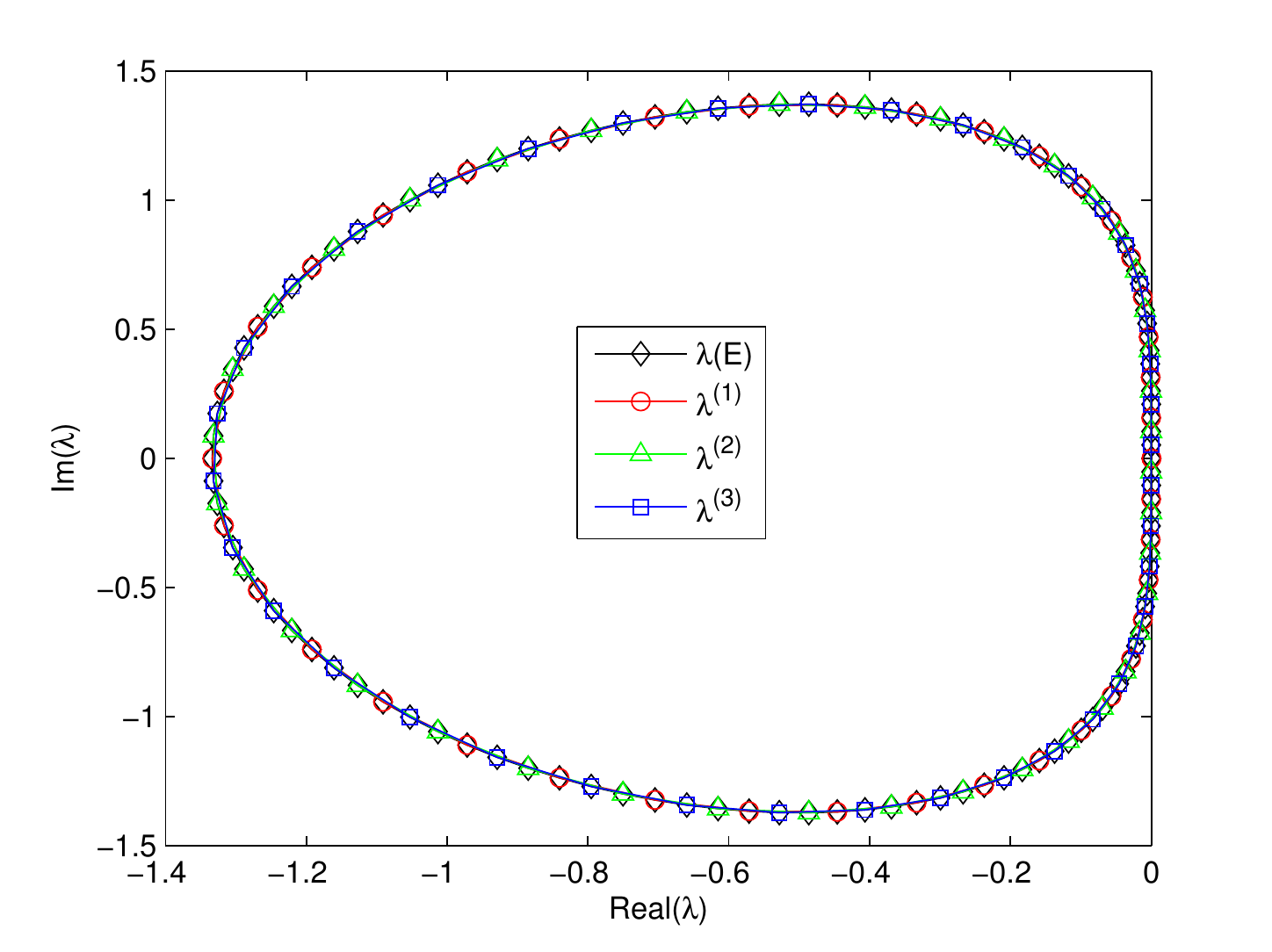}}\subfloat[WCNS5]{\includegraphics[width=0.48\textwidth]{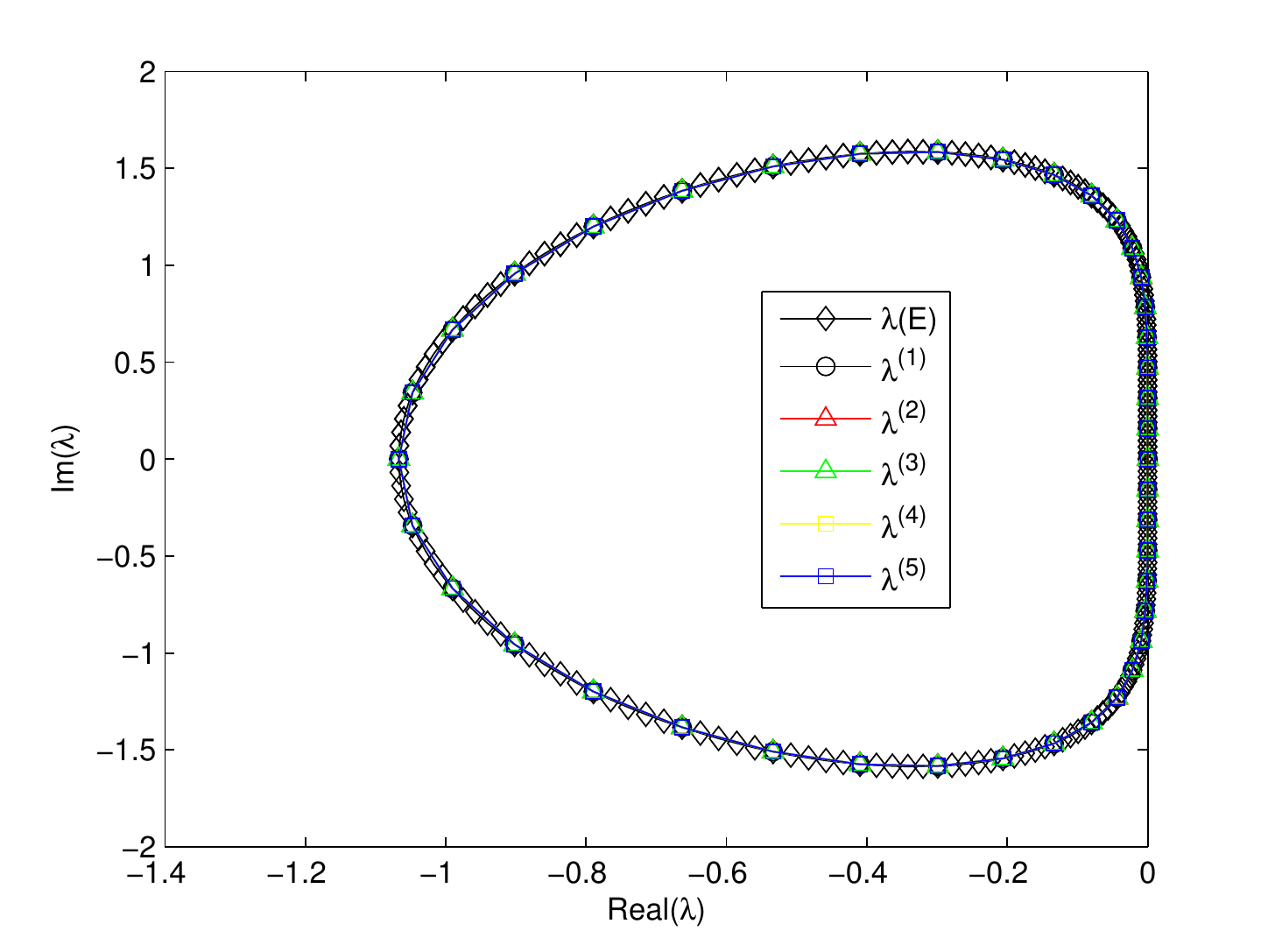}}

\subfloat[CPR-DG3]{\includegraphics[width=0.48\textwidth]{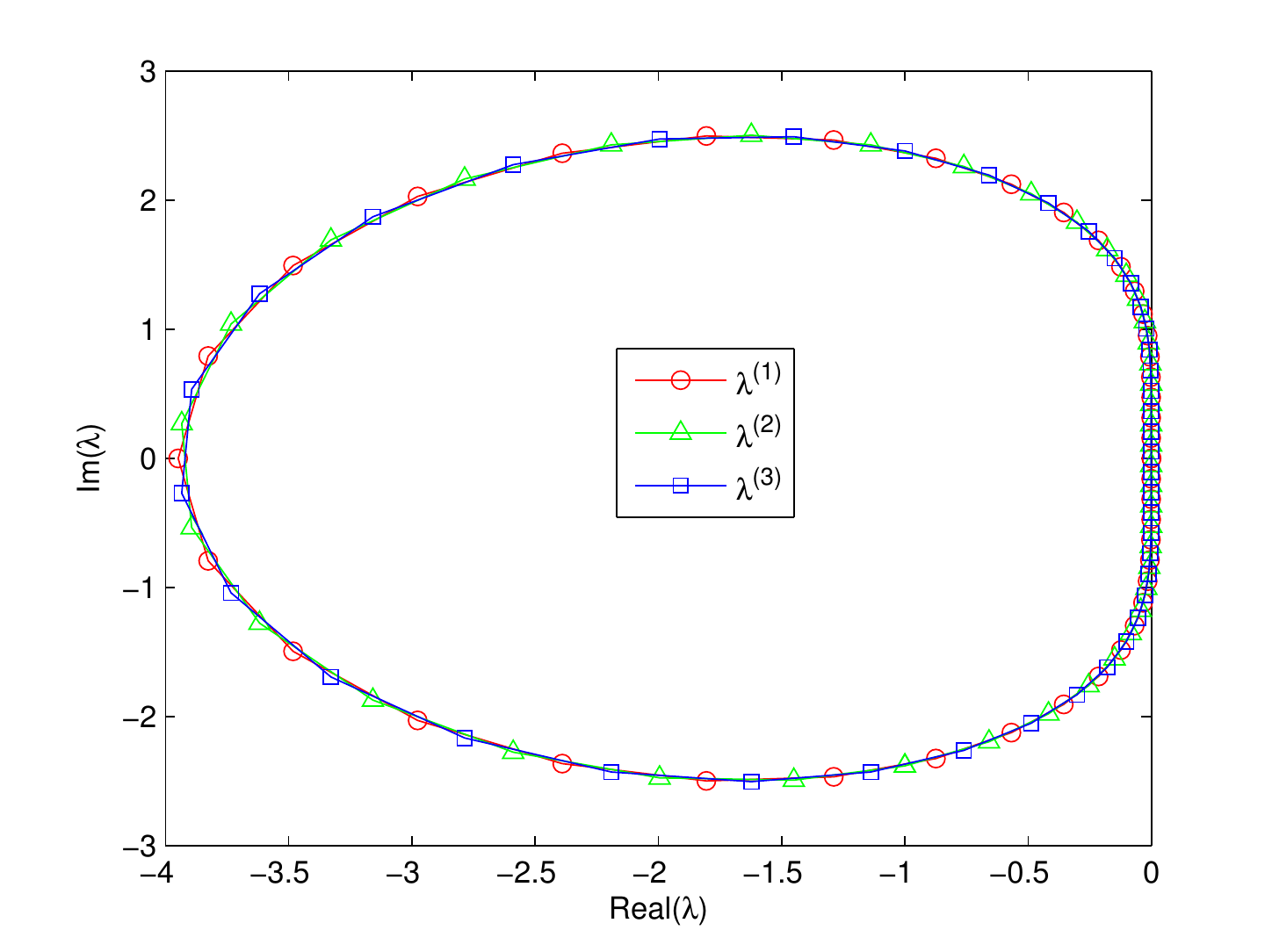}}\subfloat[CPR-DG5]{\includegraphics[width=0.48\textwidth]{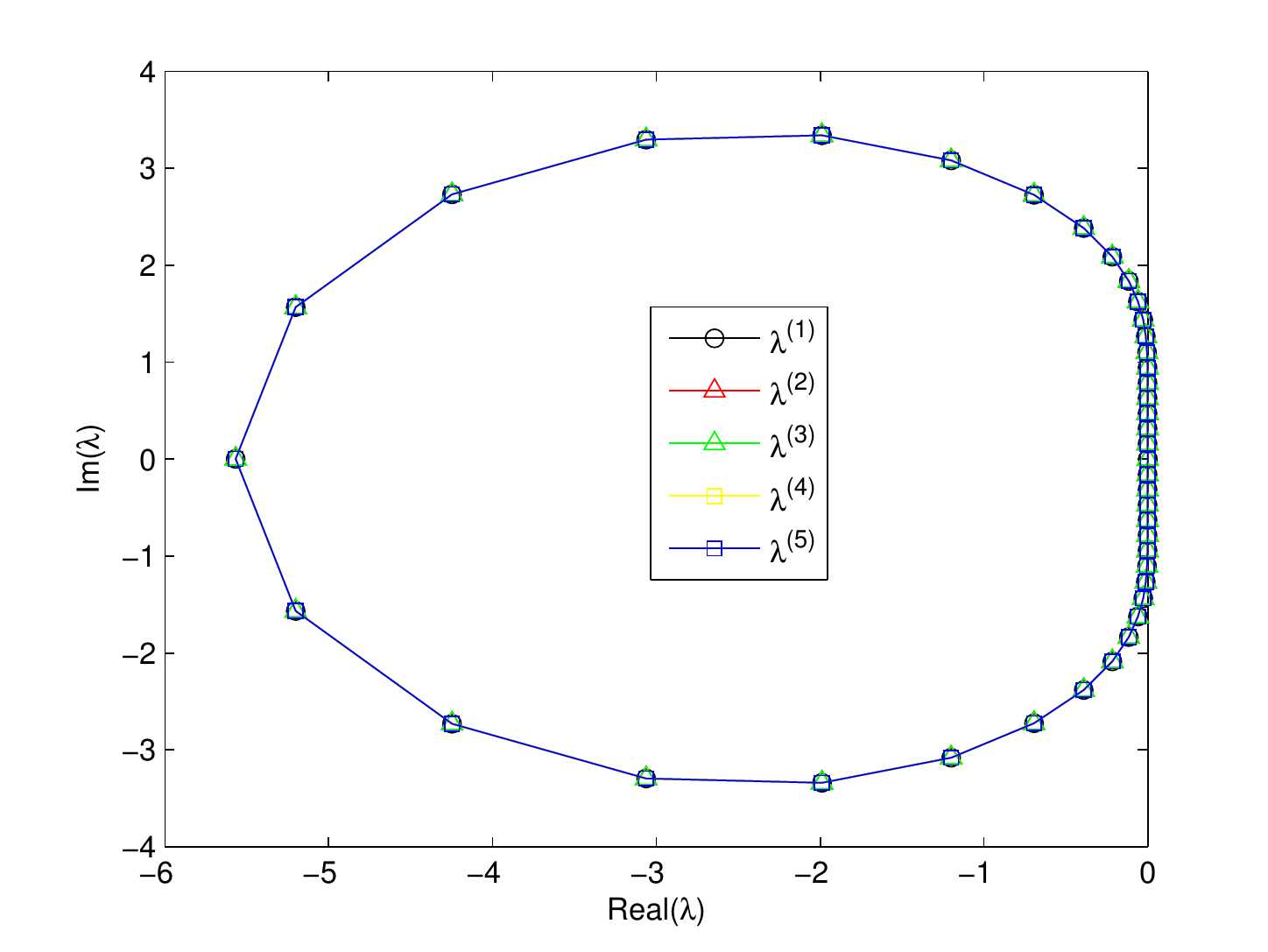}}

\caption{Comparison of the eigenvalue spectrum $\left\{ \lambda_{m}^{(l)}\left|\lambda_{m}^{(l)}=\lambda^{(l)}(G(\phi_{m})),l=1,2,\cdots K+1;m=0,1,2,\cdots,M\right.\right\} $
where $\phi_{m}=m\frac{2\pi}{M}$ for different high-order schemes
($K=2$ for third-order schemes and $K=4$ for fifth-order schemes)
and $M=40$.\label{fig:comparison-of-eigenvalues}}
\end{figure}
\begin{figure}
\subfloat[C3NNW3, dispersion]{\includegraphics[width=0.48\textwidth]{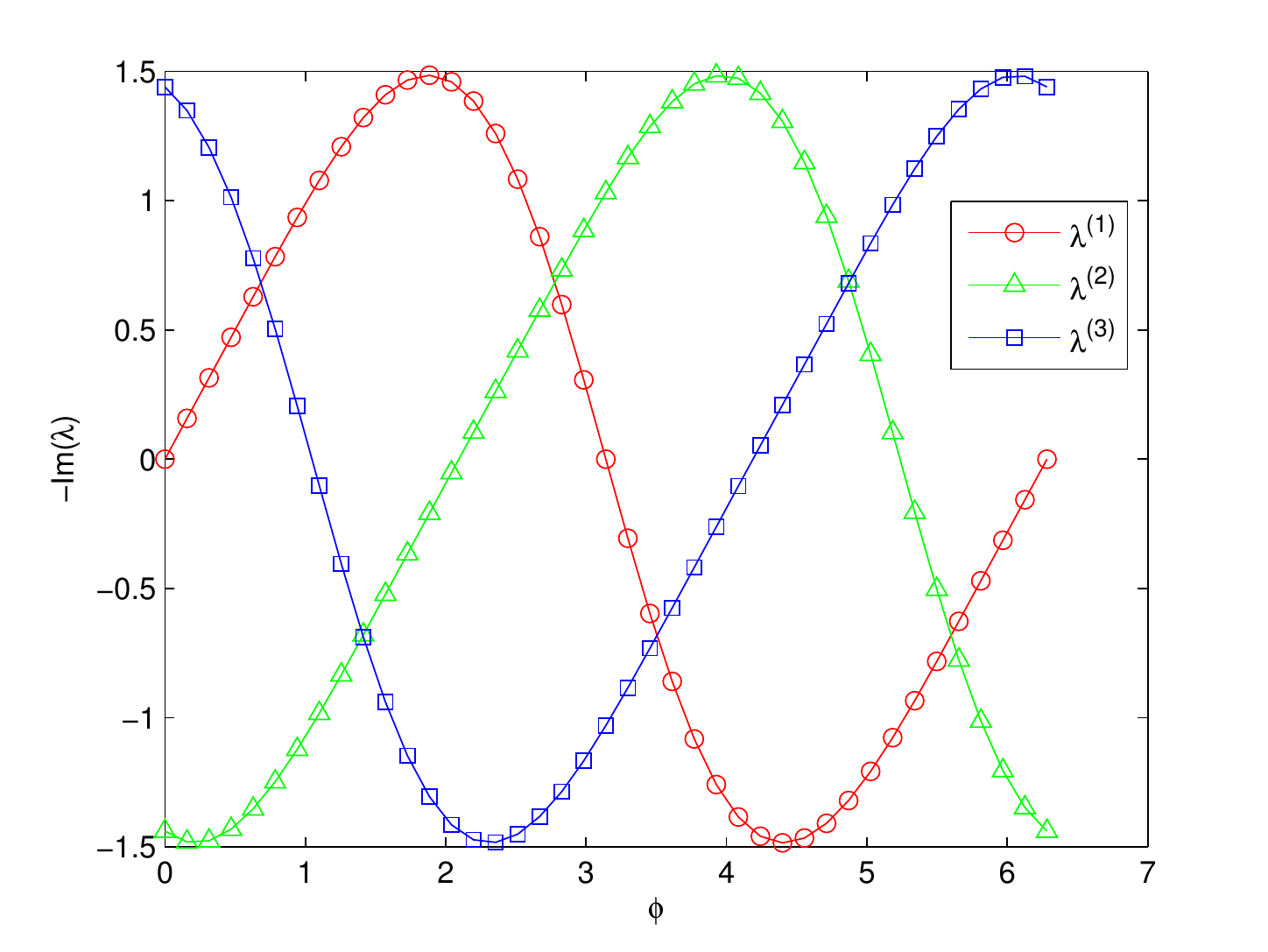}}\subfloat[C3NNW3, dissipation]{\includegraphics[width=0.48\textwidth]{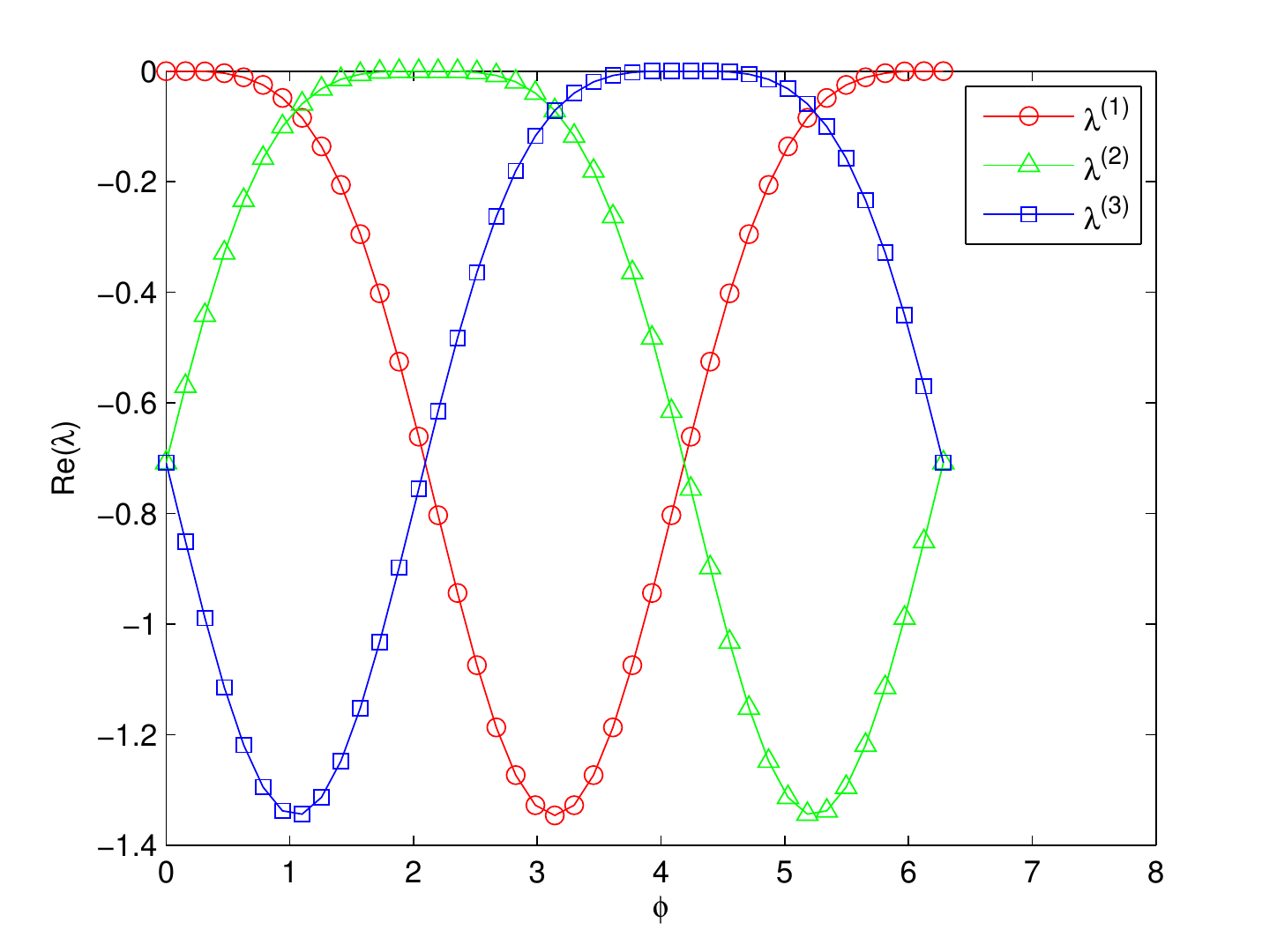}}

\subfloat[WCNS3, dispersion]{\includegraphics[width=0.48\textwidth]{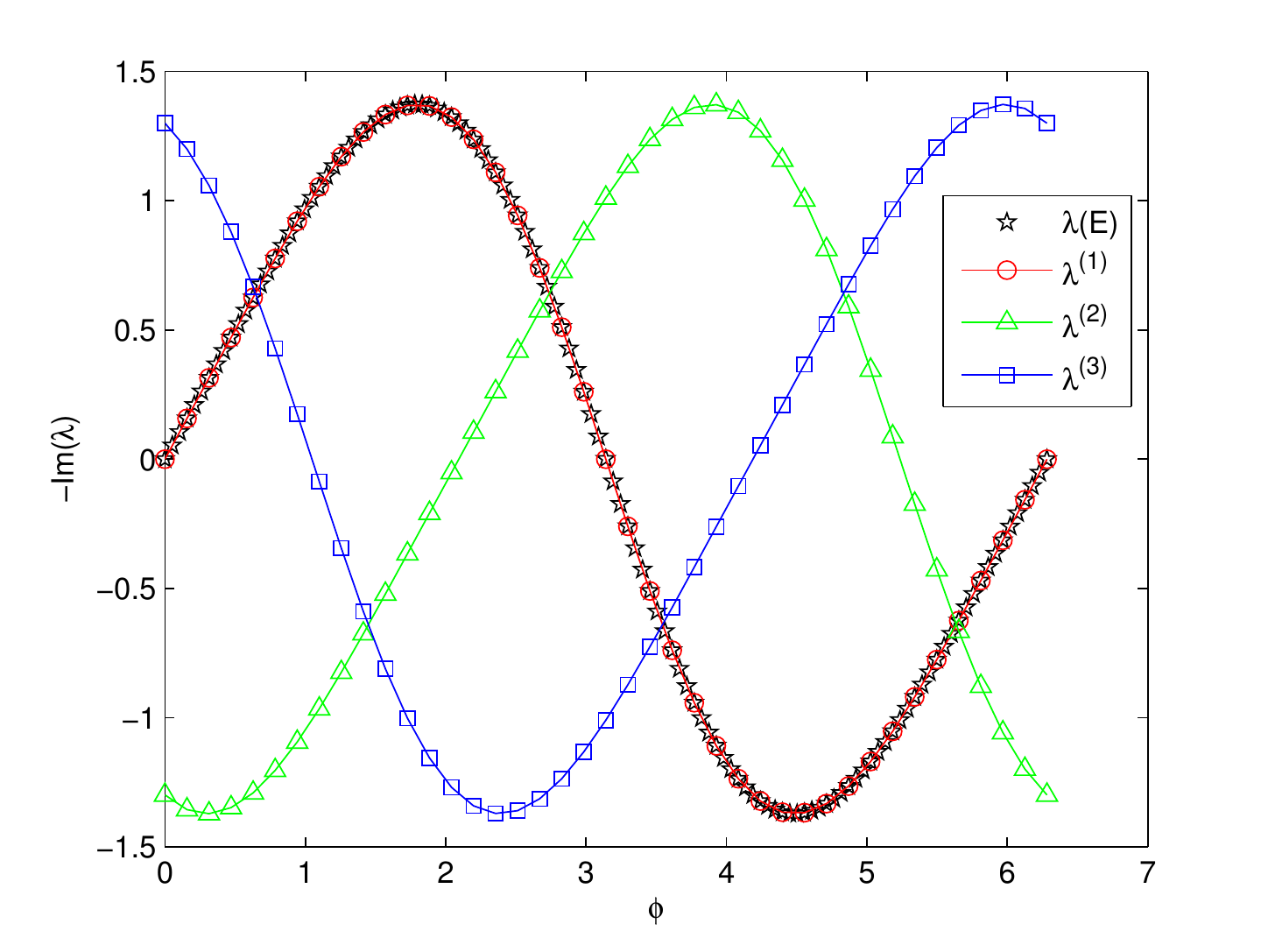}}\subfloat[WCNS3, dissipation]{\includegraphics[width=0.48\textwidth]{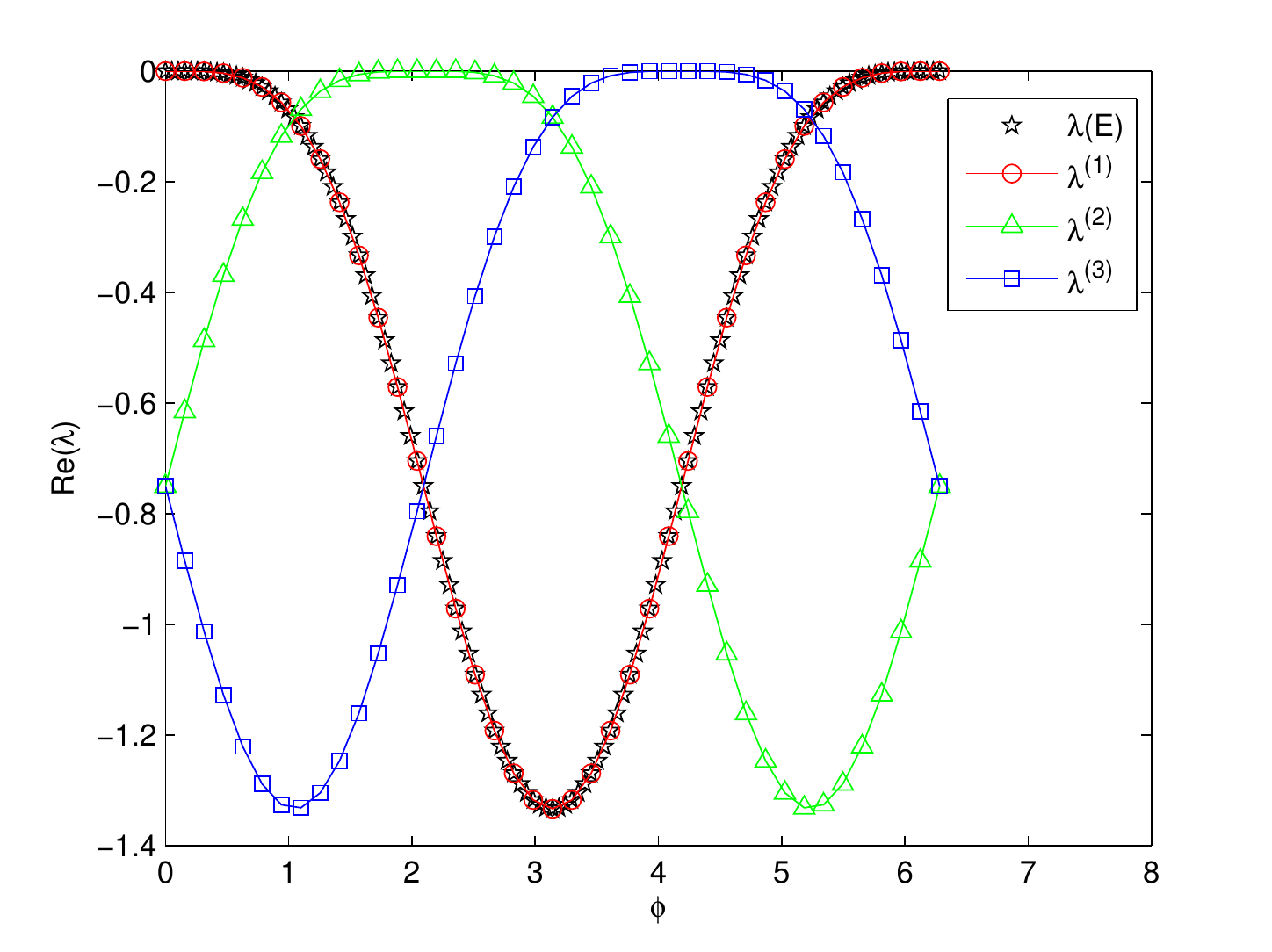}}

\subfloat[CPR3, dispersion]{\includegraphics[width=0.48\textwidth]{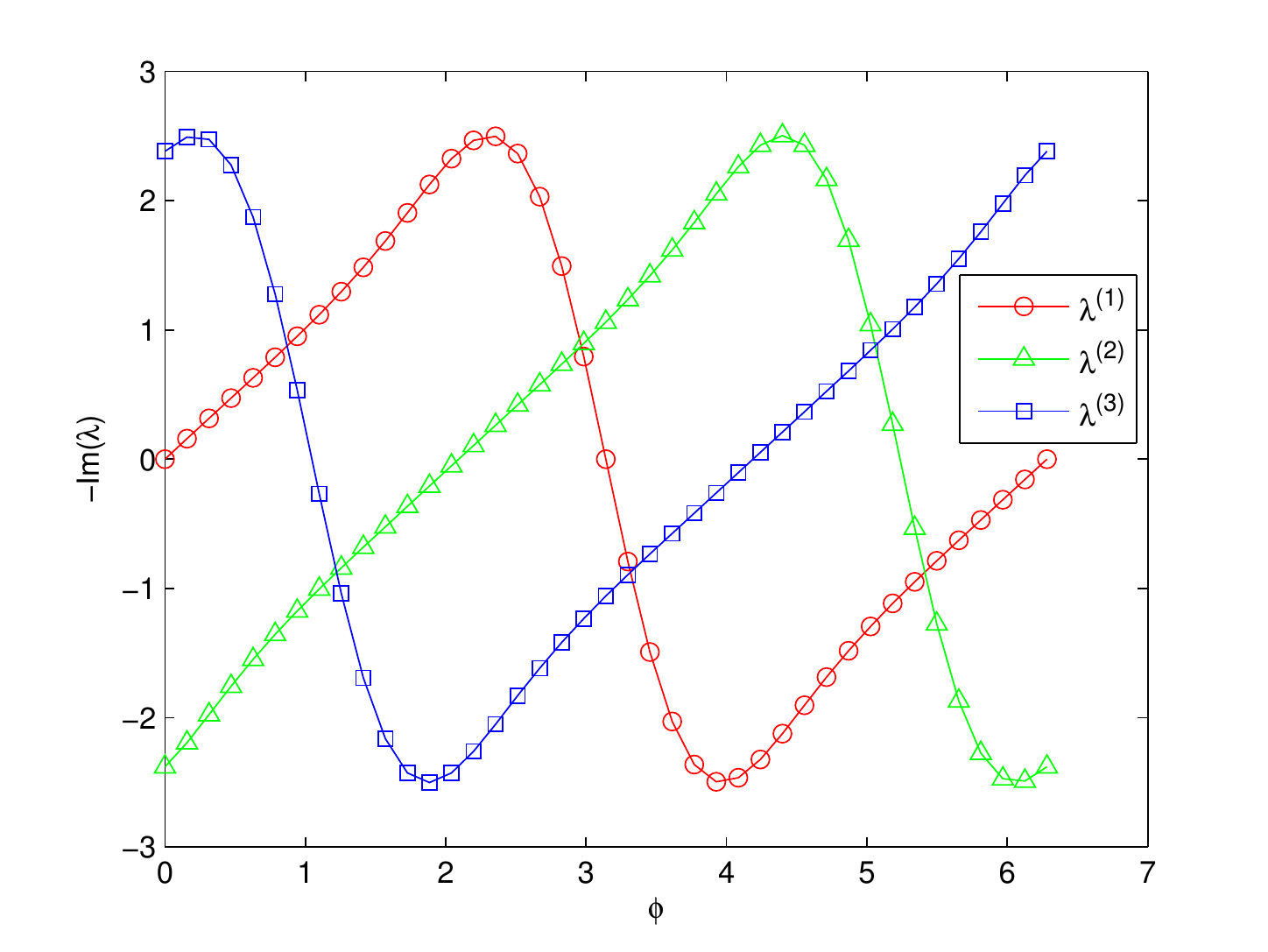}}\subfloat[CPR3, dissipation]{\includegraphics[width=0.48\textwidth]{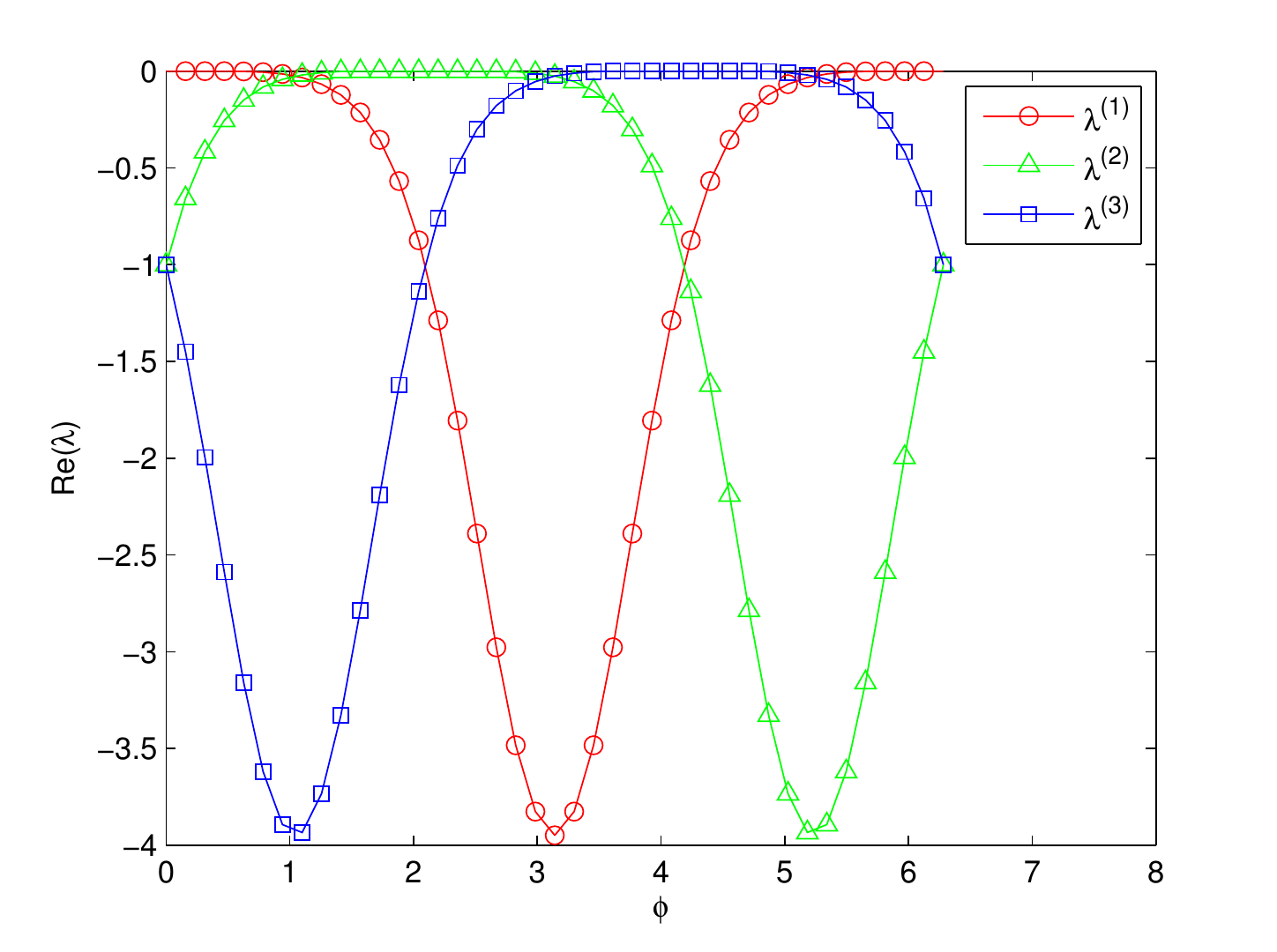}}

\caption{Comparison of dispersion (left) and dissipation (right) for third-order
schemes, where $\left\{ \lambda_{m}^{(l)}\left|\lambda_{m}^{(l)}=\lambda^{(l)}(G(\phi_{m})),l=1,2,\cdots K+1;m=0,1,2,\cdots,M\right.\right\} $,
$\phi_{m}=m\frac{2\pi}{M}$, $K=2$ and $M=40$.\label{fig:comparison-of-dispersion-third-order}}
\end{figure}

\par\end{center}

\section*{Appendix E. NNW2 interpolation}

Now we propose second-order nonuniform nonlinear weighted (NNW2) interpolation.
Consider the stencil with three nonuniformly spaced solution points
$\left\{ u_{1},u_{2},u_{3}\right\} $. Take the second solution point
for example and set $u_{1}=u_{i,sp_{1}}$, $u_{2}=u_{i,sp_{2}}$,
$u_{3}=u_{i,sp_{3}}$, as shown in in Fig. \ref{fig:C2NNW2}. The
values at flux points $u_{i,fp_{2}}^{R}$ and $u_{i,fp_{3}}^{L}$,
denoted by $u_{A}^{R}$ and $u_{B}^{L}$, can be interpolated from
$\left\{ u_{1},u_{2},u_{3}\right\} $ by following procedure. 

(1) Get $u_{A}^{(1)}$ and $u_{B}^{(1)}$ by inverse distance weighted
interpolation,

\[
u_{A}^{(1)}=\omega_{1}u_{1}+\omega_{2}u_{2},\quad\omega_{1}=\frac{\left(1/\Delta\xi_{1}\right)}{\left(1/\Delta\xi_{1}\right)+\left(1/\Delta\xi_{2}\right)},\quad\omega_{2}=\frac{\left(1/\Delta\xi_{2}\right)}{\left(1/\Delta\xi_{1}\right)+\left(1/\Delta\xi_{2}\right)};
\]
\[
u_{B}^{(1)}=\omega_{3}u_{2}+\omega_{4}u_{3},\quad\omega_{3}=\frac{\left(1/\Delta\xi_{3}\right)}{\left(1/\Delta\xi_{3}\right)+\left(1/\Delta\xi_{4}\right)},\quad\omega_{4}=\frac{\left(1/\Delta\xi_{4}\right)}{\left(1/\Delta\xi_{3}\right)+\left(1/\Delta\xi_{4}\right)};
\]
where $\Delta\xi_{1}=\xi_{fp2}-\xi_{sp1}$, $\Delta\xi_{2}=\xi_{sp2}-\xi_{fp2}$,
$\Delta\xi_{3}=\xi_{fp3}-\xi_{sp2}$, $\Delta\xi_{4}=\xi_{sp3}-\xi_{fp3}$,
as shown in Fig. \ref{fig:C2NNW2}.

(2) Calculate the gradient of $\frac{\partial u}{\partial\xi}$ with
values $\left\{ u_{A}^{(1)},u_{2},u_{B}^{(1)}\right\} $ based on
the distances between solution points and flux points, we have 

\begin{equation}
\frac{\partial u}{\partial\xi}=\omega_{5}\left(\frac{\partial u}{\partial\xi}\right)^{(1)}+\omega_{6}\left(\frac{\partial u}{\partial\xi}\right)^{(2)},\label{eq:weighted gradient}
\end{equation}
 where

\[
\omega_{5}=\frac{\left(1/\Delta\xi_{2}\right)}{\left(1/\Delta\xi_{2}\right)+\left(1/\Delta\xi_{3}\right)},\quad\omega_{6}=\frac{\left(1/\Delta\xi_{3}\right)}{\left(1/\Delta\xi_{2}\right)+\left(1/\Delta\xi_{3}\right)},\quad\left(\frac{\partial u}{\partial\xi}\right)^{(1)}=\frac{u_{2}-u_{A}^{(1)}}{\Delta\xi_{2}},\quad\left(\frac{\partial u}{\partial\xi}\right)^{(2)}=\frac{u_{B}^{(1)}-u_{2}}{\Delta\xi_{3}}.
\]
(3) Compute $u_{A}^{(2)}$ and $u_{B}^{(2)}$ based on $u_{2}$ and
the gradient $\frac{\partial u}{\partial\xi}$, we obtain

\[
u_{A}^{(2)}=u_{2}-\frac{\partial u}{\partial\xi}\Delta\xi_{2},\quad u_{B}^{(2)}=u_{2}+\frac{\partial u}{\partial\xi}\Delta\xi_{3}.
\]

(4) Add limiter to control numerical oscillation. $u_{A}^{R}$ and
$u_{B}^{L}$ are obtained by linear reconstruction with a limiter,

\begin{eqnarray}
u_{A}^{R} & = & u_{2}-\phi\frac{\partial u}{\partial\xi}\Delta\xi_{2},\quad u_{B}^{L}=u_{2}+\phi\frac{\partial u}{\partial\xi}\Delta\xi_{3}.\label{eq:NNW2 interpolation}
\end{eqnarray}
Here we take the following Birth limiter \citet{Birth1989}

\[
\phi=min\{lim(u_{A}^{(2)}),lim(u_{B}^{(2)})\},
\]
where

\begin{eqnarray*}
lim(u)= & \begin{cases}
min\{1,\frac{M-u_{2}}{u-u_{2}}\}, & if\,\, u>u_{2},\\
min\{1,\frac{m-u_{2}}{u-u_{2}}\}, & if\,\, u<u_{2},\\
1, & if\,\, u=u_{2},
\end{cases}
\end{eqnarray*}
with $m=min\{u_{1},u_{2},u_{3}\}$ and $M=max\{u_{1},u_{2},u_{3}\}$.

\bibliographystyle{unsrtnat}
\bibliography{Ref11}

\end{document}